\documentclass[a4paper,11pt]{amsart}

\usepackage[utf8]{inputenc}
\usepackage[T1]{fontenc}
\usepackage{lmodern}
\usepackage[english]{babel}
\usepackage{microtype}
\usepackage{amsmath,amssymb,amsfonts,amsthm}
\usepackage{cases}
\usepackage{mathtools,accents}
\usepackage{mathrsfs}
\usepackage{aliascnt}
\usepackage{bm}
\usepackage{bbm}
\usepackage{tikz}
\usepackage[citecolor=blue,colorlinks]{hyperref}
\usepackage{verbatim}
\usepackage{enumerate}
\usepackage{tikz-cd}
\usepackage{fullpage}
\usepackage{comment}
\usepackage{booktabs}
\usepackage{float}
\usetikzlibrary{calc}
\usetikzlibrary{arrows.meta}
\makeatletter
\g@addto@macro\@floatboxreset\centering
\makeatother

\makeatletter
\def\newaliasedtheorem#1[#2]#3{
\newaliascnt{#1\@alt}{#2}
\newtheorem{#1}[#1\@alt]{#3}
\expandafter\newcommand\csname #1\@altname\endcsname{#3}
}
\makeatother

\newtheorem{theorem}{Theorem}[section]
\newtheorem{lemma}[theorem]{Lemma}
\newtheorem{corollary}[theorem]{Corollary}
\newtheorem{proposition}[theorem]{Proposition}

\newtheorem{remark}[theorem]{Remark}

\theoremstyle{definition}
\newtheorem{definition}[theorem]{Definition}

\newtheorem{question}[theorem]{Question}

\numberwithin{equation}{section}

\newcommand{\p}{\partial}

\newcommand{\bl}{\bm{\ell}}

\newcommand{\C}{\mathbb{C}}

\newcommand{\cK}{\mathcal{K}}

\newcommand{\cL}{\mathcal{L}}

\newcommand{\R}{\mathbb{R}}

\newcommand{\CP}{\mathbb{CP}}

\newcommand{\mv}{\mathfrak v}
\newcommand{\mw}{\mathfrak w}

\newcommand{\mI}{\mathfrak I}

\newcommand{\mR}{\mathfrak R}

\newcommand{\bZ}{\mathbb Z}

\DeclareMathOperator{\diam}{diam}
\DeclareMathOperator{\dvol}{dvol}

\DeclareMathOperator{\Id}{Id}

\DeclareMathOperator{\Ric}{Ric}
\DeclareMathOperator{\Rm}{Rm}

\DeclareMathOperator{\SO}{SO}

\DeclareMathOperator{\SL}{SL}

\DeclareMathOperator{\Vol}{Vol}

\DeclareMathOperator{\Ae}{AE}

\DeclareMathOperator{\deP}{dP}
\DeclareMathOperator{\ALF}{ALF}

\DeclareMathOperator{\EH}{EH}

\DeclareMathOperator{\ALE}{ALE}

\DeclareMathOperator{\AF}{AF}

\newcommand{\bg}{\underline g}
\newcommand{\dd}{\mathrm d}
\DeclareMathOperator{\Scal}{Scal}
\DeclareMathOperator{\Area}{Area}
\DeclareMathOperator{\Hess}{Hess}

\DeclareMathOperator{\tr}{tr}
\DeclareMathOperator{\diver}{div}

\newcommand{\T}{\mathbb T}
\newcommand{\Z}{\mathbb Z}

\theoremstyle{plain}

\title{Compact toric Einstein four-manifolds}
\author{
Mingyang Li\textsuperscript{*}}
\thanks{\textsuperscript{*} Simons Center for Geometry and Physics, Stony Brook University; \url{mingyang.li@scgp.stonybrook.edu}}
\author{Song Sun\textsuperscript{\dag}}
\thanks{\dag\  Institute for Advanced Study in Mathematics, Zhejiang University; \url{songsun@zju.edu.cn}}

\begin{document}
\begin{abstract}
This is a continuation of our previous paper on toric gravitational instantons. We develop a systematic theory of compact simply-connected Einstein four-manifolds with toric symmetry. Applications include
\begin{itemize}
\item The construction of simply-connected positive Einstein four-manifolds with arbitrarily large $b_2$. Our method is based on local gluing combined with a global virtual counting argument. The key idea involves a particular design of rod structures to enable a gluing construction and restrict possible degenerations. 

\item A  topological characterization of algebraically special toric Einstein metrics in terms of $b_2^+\leq 1$. This uses special curvature identities for $W^+$ in the toric setting. 
\item  Uniqueness,  up to scaling and isometry, of toric Einstein metrics on four-manifolds admitting an algebraically special Einstein metric. 
\item Diffeomorphism classification of compact toric Einstein four-manifolds in terms of an improved Hitchin-Thorpe inequality $\chi\geq3|\tau|$.
\item Smooth compactness of simply connected toric Einstein four-manifolds with bounded Euler number. 
\end{itemize}
In particular, our construction recovers all the classical toric Einstein metrics, including the Page and Chen-LeBrun-Weber metrics, starting from the round sphere through moduli spaces of conical Einstein metrics.
\end{abstract}
\maketitle
\setcounter{tocdepth}{1}
\tableofcontents

\section{Introduction}

In this paper, we construct new four-dimensional simply-connected compact Einstein manifolds with positive Ricci curvature. Previously, the known examples are given in the following list
\begin{itemize}
\item Round $\mathbb S^4$.
\item K\"ahler-Einstein metrics on del Pezzo surfaces. The underlying manifolds are diffeomorphic to $\CP^2$, $\CP^1\times \CP^1$, or $\CP^2\# k\overline{\CP^2}$ for $3\leq k\leq 8$.
\item The Page metric on $\CP^2\# \overline{\CP^2}$.
\item The Chen-LeBrun-Weber metric on $\CP^2\# 2\overline{\CP^2}$.
\end{itemize}
All these examples have $b_2\leq9$ and occur on only finitely many diffeomorphism types. They are also \emph{algebraically special}: for a suitable choice of orientation, the self-dual Weyl curvature $W^+$, viewed as an endomorphism of $\Lambda^+$, has at most two distinct eigenvalues at every point. For a compact, simply connected Einstein four-manifold, Derdzi\'nski's theorem implies that this condition forces the metric to be either anti-self-dual or Hermitian with respect to an orientation-compatible complex structure \cite{Derdzinski}. Combining this with the classification results of Hitchin and LeBrun \cite{Hitchin,LeBrun}, the examples listed above exhaust, up to scaling and isometry, the compact, simply connected, algebraically special Einstein four-manifolds with positive Ricci curvature. Complex geometry, through twistor theory or the theory of complex surfaces, plays a crucial role in both the construction and the classification of these metrics.

Among the examples above, the toric Einstein metrics, that is, those admitting an effective isometric $\mathbb T^2$-action, are the following six, up to scaling and isometry:
\begin{itemize}
\item The round metric on $\mathbb S^4$.
\item The Fubini-Study metric on $\CP^2$.
\item The round product metric on $\mathbb S^2\times\mathbb S^2$.
\item The Page metric on $\CP^2\#\overline{\CP^2}$.
\item The Chen-LeBrun-Weber metric on $\CP^2\#2\overline{\CP^2}$.
\item The $\deP_3$ metric: the K\"ahler-Einstein metric on $\CP^2\#3\overline{\CP^2}$.
\end{itemize}

These metrics were constructed by different methods. The first three are explicit and homogeneous. The Page metric is also explicit; its cohomogeneity-one symmetry reduces the Einstein equation to a system of ordinary differential equations~\cite{Page}. The last two have cohomogeneity two, and no closed formulas are known. The Chen-LeBrun-Weber metric was obtained by constructing an extremal K\"ahler metric and making a conformal change~\cite{CLW}, whereas the existence of the $\deP_3$ metric follows from the theory of the complex Monge-Amp\`ere equation~\cite{Siu,Tian-Yau}. For numerical approximations and computations of geometric quantities, see~\cite{DHHKW,BunchDonaldson,HallMurphyNumerics}. See also~\cite{BHHPW} for a computer-assisted approximation of the $\deP_3$ metric with rigorous error estimates.

The following are the existence results we prove in this paper

\begin{theorem}[First Existence Theorem]\label{thm:main}
The following manifolds admit toric Einstein metrics:

\begin{itemize}
\item $\#r(\mathbb S^2\times \mathbb S^2)\ (r\geq 0)$ (when $r=0$, this is understood as $\mathbb S^4$);
\item $\bigl(\CP^2\#k\overline{\CP^2}\bigr) \#r(\mathbb S^2\times \mathbb S^2)$, for $0\leq k\leq3$ and $r\geq0$.

\item $M\#N$, for either of the following choices:
$$
    \begin{aligned}
        M&\in\left\{ \mathbb S^4,\ \mathbb S^2\times \mathbb S^2,\ \CP^2\#k\overline{\CP^2}\ (1\leq k\leq3) \right\},\\
        N&\in\left\{ \CP^2,\ \CP^2\#\overline{\CP^2},\ 2\CP^2\#\overline{\CP^2},\ 3\CP^2\#\overline{\CP^2} \right\};
    \end{aligned}
$$
or
$$
    \begin{aligned}
        M&\in\left\{ \mathbb S^4,\ \mathbb S^2\times \mathbb S^2,\ \CP^2\#k\overline{\CP^2}\ (0\leq k\leq2) \right\},\\
        N&\in\left\{ \overline{\CP^2},\ \CP^2\#\overline{\CP^2},\ \CP^2\#2\overline{\CP^2},\ \CP^2\#3\overline{\CP^2} \right\}.
    \end{aligned}
$$

\end{itemize}
\end{theorem}

A mild extension of our strategy, developed by AI, gives
the following more general existence theorem and a complete classification of diffeomorphism types for simply connected compact toric Einstein four manifolds. 

\begin{theorem}[=Theorem \ref{thm:improved HT}+Theorem \ref{thm:toric-existence-criterion}, General Existence]\label{thm:general}
Let $X$ be a closed simply connected smooth four-manifold admitting an effective $\mathbb T^2$-action. Then $X$ admits a toric Einstein metric (possibly for a different $\mathbb T^2$-action) if and only if
\begin{equation}\label{eq:toric exisetnce}
    \chi(X)\geq3|\tau(X)|.
\end{equation}
\end{theorem}

\begin{remark}
Independently of the present work, Ozuch~\cite{OzuchPositive} constructs positive toric Einstein metrics on $\#k(\mathbb S^2\times\mathbb S^2)$ and $k\CP^2\#(2k+1)\overline{\CP^2}$ for $k\geq2$. His construction varies a cone angle along a generic torus orbit and obtains smooth metrics by changing the cone angle along the torus and removes the conical singularity at selected angles by changing the angular period. Our construction uses gluing at torus fixed points followed by a degree argument for cone angles along rod spheres. \footnote{After the first version of this paper was posted on the arXiv, two more constructions \cite{Lafuente-Broder, ChuStern} of toric Einstein metrics on $\#k(\mathbb S^2\times\mathbb S^2)$ appeared. Both constructions share features of the geometric setup in \cite{OzuchPositive}, using a rod structure with discrete symmetry. \cite{Lafuente-Broder} also varies the cone angle along a generic torus orbit, using a degree argument, while \cite{ChuStern} uses a variational approach. }
\end{remark}
Some manifolds occur more than once in the list in Theorem \ref{thm:main}, up to orientation reversal, because they arise from different constructions. The following uniqueness theorem shows that the resulting metrics agree up to scaling and isometry in the cases below. It is a key ingredient in proving Theorem \ref{thm:main}. We also characterize algebraic speciality topologically.
For a fixed orientation, a toric Einstein metric is
anti-self-dual precisely when $b_2^+=0$, and conformally
K\"ahler precisely when $b_2^+=1$;
see Theorem~\ref{thm:toric-weyl-rigidity}.

\begin{theorem}[=Corollary \ref{cor:toric-small-betti-classification}, Uniqueness for small $b_2^+$] \label{thm:uniqueness}
On each of the manifolds
$$
    \mathbb S^4,\qquad \mathbb S^2\times \mathbb S^2,\qquad \CP^2\#k\overline{\CP^2}\quad(0\leq k\leq3),
$$
there is exactly one toric Einstein metric, up to scaling and isometry, which is the known metric listed above.
\end{theorem}
\begin{remark}
    In general, on a give smooth compact four-manifold with $b_2^+\geq 2$, there may be more than one toric Einstein metric up to scaling and isometry. For example, on $2\CP^2\# 2\overline{\CP^2}$, there are at least two non-equivalent toric Einstein metrics, which correspond to two non-equivalent torus actions. See Remark \ref{rmk:more than one}.
\end{remark}
\begin{remark}
    For a fixed torus action on a closed simply connected four-manifold, Liu proved uniqueness of an invariant Einstein metric inducing a prescribed quotient metric \cite[Theorem~1.2]{Liu}. In July 2026, Ramiro Lafuente informed the second author that he and Kyle Broder had proved uniqueness of toric Einstein metrics on $\mathbb S^4$.
\end{remark}

In particular, starting from the round metric on $\mathbb S^4$, our constructions recover the other five algebraically special toric Einstein metrics through moduli spaces of conical Einstein metrics. This provides a uniform construction of these examples.

\begin{figure}[!htb]
\centering
\small
\definecolor{RFblue}{RGB}{0,114,178}
\definecolor{RForange}{RGB}{213,94,0}
\tikzset{
  rfold/.style={draw=RFblue,text=RFblue},
  rfnew/.style={draw=RForange,text=RForange},
  rfdots/.style={draw=RForange,text=RForange,densely dotted}
}

\newcommand{\RFpicture}[2]{\begin{tikzpicture}[x=#1cm,y=#1cm,>=stealth]
    \edef\RFdata{#2}\foreach[count=\RFi from 0] \RFx/\RFy/\RFstyle/\RFlabel/\RFlx/\RFly/\RFmx/\RFmy in \RFdata {
      \coordinate (rf\RFi) at (\RFx,\RFy);
      \xdef\RFlast{\RFi}}
    \fill[black!4] (rf0) \foreach \RFi in {1,...,\RFlast} {-- (rf\RFi)} -- cycle;
    \foreach[count=\RFi from 0] \RFx/\RFy/\RFstyle/\RFlabel/\RFlx/\RFly/\RFmx/\RFmy in \RFdata {
      \pgfmathtruncatemacro{\RFj}{mod(\RFi+1,\RFlast+1)}
      \draw[\RFstyle,line width=.8pt,line cap=round,line join=round] (rf\RFi)--(rf\RFj);
      \def\RFskip{X}\ifx\RFlabel\RFskip\else
        \coordinate (rfmid) at (\RFmx,\RFmy);
        \coordinate (rflabel) at (\RFlx,\RFly);
        \draw[\RFstyle,->,line width=.35pt] (rfmid)--($(rfmid)!.11cm!(rflabel)$);
        \node[\RFstyle,draw=none,font=\footnotesize,inner sep=0pt] at (rflabel) {$\RFlabel$};
        \fill[black!65] (rf\RFi) circle[radius=.65pt];
      \fi
    }
    \path[use as bounding box] (-.48,-.48) rectangle (3.48,3.48);
  \end{tikzpicture}}
\newcommand{\RFpanelwidth}{.26\textwidth}
\newcommand{\RFpanel}[3]{\begin{minipage}[c]{\RFpanelwidth}\centering
    \RFpicture{#1}{#2}\par\vspace{.3mm}
    {\small #3\par}
  \end{minipage}}

\newcommand{\RFsphereOne}{0/0/rfold/B/1.5/-0.32/1.5/0,
3/0/rfnew/A/3.32/1.5/3/1.5,
3/3/rfnew/B/1.5/3.32/1.5/3,
0/3/rfold/A/-0.32/1.5/0/1.5}

\newcommand{\RFsphereTwo}{0/0/rfold/B/1.5/-0.32/1.5/0,
3/0/rfnew/A/3.32/0.75/3/0.75,
3/1.5/rfnew/B/2.25/1.82/2.25/1.5,
1.5/1.5/rfnew/A/1.82/2.25/1.5/2.25,
1.5/3/rfnew/B/0.75/3.32/0.75/3,
0/3/rfold/A/-0.32/1.5/0/1.5}

\newcommand{\RFsphereMany}{0/0/rfold/B/1.5/-0.32/1.5/0,
3/0/rfnew/A/3.32/0.3/3/0.3,
3/0.6/rfnew/B/2.7/0.92/2.7/0.6,
2.4/0.6/rfdots/X/2.72/0.9/2.4/0.9,
2.4/1.2/rfdots/X/2.1/1.52/2.1/1.2,
1.8/1.2/rfdots/X/2.12/1.5/1.8/1.5,
1.8/1.8/rfdots/X/1.5/2.12/1.5/1.8,
1.2/1.8/rfdots/X/1.52/2.1/1.2/2.1,
1.2/2.4/rfdots/X/0.9/2.72/0.9/2.4,
0.6/2.4/rfnew/A/0.92/2.7/0.6/2.7,
0.6/3/rfnew/B/0.3/3.32/0.3/3,
0/3/rfold/A/-0.32/1.5/0/1.5}

\newcommand{\RFprojectiveOne}{3/0/rfold/B/2.075/-0.32/2.075/0,
1.15/0/rfnew/A/0.92/0.2875/1.15/0.2875,
1.15/1.15/rfnew/B/0.2875/0.92/0.2875/1.15,
0/1.15/rfold/A/-0.32/2.075/0/2.075,
0/3/rfold/C/1.72627/1.72627/1.5/1.5}

\newcommand{\RFprojectiveTwo}{3/0/rfold/B/2.325/-0.32/2.325/0,
1.65/0/rfnew/A/1.42/0.20625/1.65/0.20625,
1.65/0.825/rfnew/B/1.03125/0.595/1.03125/0.825,
0.825/0.825/rfnew/A/0.595/1.03125/0.825/1.03125,
0.825/1.65/rfnew/B/0.20625/1.42/0.20625/1.65,
0/1.65/rfold/A/-0.32/2.325/0/2.325,
0/3/rfold/C/1.72627/1.72627/1.5/1.5}

\newcommand{\RFprojectiveMany}{3/0/rfold/B/2.55/-0.32/2.55/0,
2.1/0/rfnew/A/1.87/0.175/2.1/0.175,
2.1/0.7/rfnew/B/1.575/0.47/1.575/0.7,
1.4/0.7/rfdots/X/1.08/1.05/1.4/1.05,
1.4/1.4/rfdots/X/1.05/1.08/1.05/1.4,
0.7/1.4/rfnew/A/0.47/1.575/0.7/1.575,
0.7/2.1/rfnew/B/0.175/1.87/0.175/2.1,
0/2.1/rfold/A/-0.32/2.55/0/2.55,
0/3/rfold/C/1.72627/1.72627/1.5/1.5}

\setlength{\tabcolsep}{1mm}
\renewcommand{\arraystretch}{1}
\begin{tabular}{@{}c@{\hspace{2mm}}ccc@{}}
 & $r=1$ & $r=2$ & $r\geq3$ \\[1mm]
\shortstack{From\\[.5mm]$\mathbb S^4$} &
\RFpanel{.70}{\RFsphereOne}{$\mathbb S^2\times \mathbb S^2$} &
\RFpanel{.70}{\RFsphereTwo}{$\#2(\mathbb S^2\times \mathbb S^2)$} &
\RFpanel{.70}{\RFsphereMany}{$\#r(\mathbb S^2\times \mathbb S^2)$} \\[1mm]
\shortstack{From\\[.5mm]$\mathbb{CP}^2$} &
\RFpanel{.70}{\RFprojectiveOne}{$\mathbb{CP}^2\#(\mathbb S^2\times \mathbb S^2)$} &
\RFpanel{.70}{\RFprojectiveTwo}{$\mathbb{CP}^2\#2(\mathbb S^2\times \mathbb S^2)$} &
\RFpanel{.70}{\RFprojectiveMany}{$\mathbb{CP}^2\#r(\mathbb S^2\times \mathbb S^2)$}
\end{tabular}
\caption{Rod structure for two infinite families in Theorem~\ref{thm:main}.}
\label{fig:rod-infinite}
\end{figure}
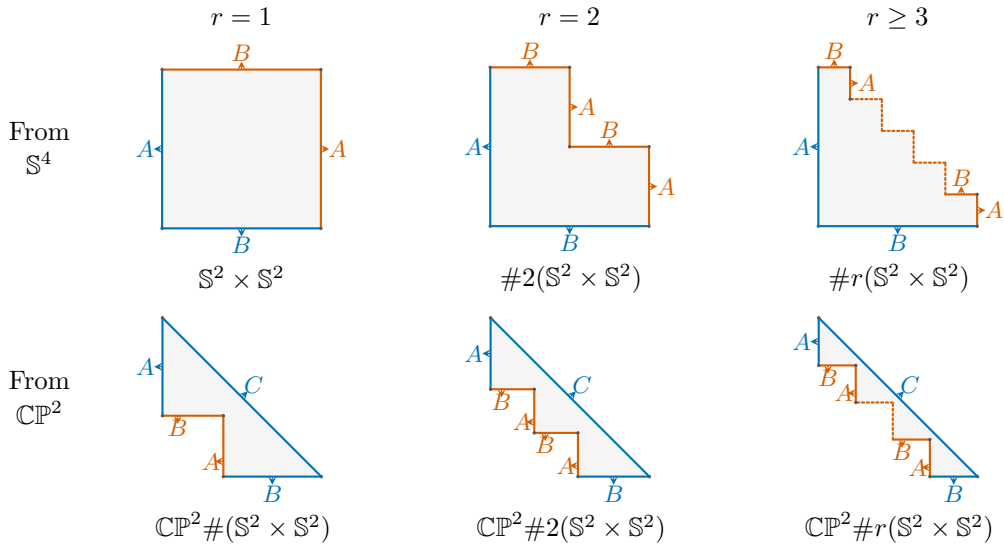

Our construction uses toric metrics.  We obtain the examples above by choosing suitable torus actions. Their orbit spaces are topological disks: a circle subgroup collapses along each boundary edge, and the vertices correspond to torus fixed points. We use only the three circle subgroups with primitive generators
$$
    A=2\pi e_1,\qquad B=2\pi e_2,\qquad C=A+B,
$$
where $e_1,e_2$ are the standard basis vectors of $\mathbb R^2$ and the torus lattice is $2\pi\mathbb Z^2$. The cyclic list of boundary labels describes the action; its precise formulation as a \emph{rod structure} is given in Section~\ref{sec:rod structure}. Any two distinct labels form a basis of the torus lattice, so distinct labels on adjacent edges give smooth local models at the fixed points.

The first construction replaces the two consecutive labels at a fixed point by a longer alternating sequence:
\begin{equation}\label{eq:pointed-insertion}
    B,A\quad\longmapsto\quad B,\underbrace{A,B,\ldots,A,B}_{2n\ \mathrm{new\ edges}},A, \qquad n\geq1.
\end{equation}
Topologically, this replaces $M$ by $M\#n(\mathbb S^2\times \mathbb S^2)$. Applying the alternating construction to the six classical toric Einstein metrics gives the first two families in Theorem~\ref{thm:main}. Figure~\ref{fig:rod-infinite} illustrates the cases starting from $\mathbb S^4$ and $\CP^2$.

For the finite constructions, we use the four sequences
\begin{equation}\label{eq:hermitian-words}
    \mathcal W=\{BCA,\ BACA,\ BCACA,\ BCABCA\}.
\end{equation}
Each replaces the pair $B,A$, with its first and last labels identified with the original edges. The intermediate edges correspond to one, two, three, or four invariant spheres. These four replacements give the connected-sum summands
in the third part of Theorem~\ref{thm:main}, with the orientation determined by the lattice identification.

For either construction, the end circles are identified with the two circles at the chosen fixed point by a torus automorphism. We call the replacement \emph{admissible} if the resulting cyclic list still uses only $\{A,B,C\}$ after one global change of lattice basis.

The existence result Theorem \ref{thm:main} follows from a more general dichotomy theorem below.

\begin{definition}[Toric nondegeneracy]
\label{def:toric-nondegeneracy}
A smooth toric Einstein metric $g$ is \emph{torically nondegenerate} if its space of $\mathbb T^2$-invariant transverse-traceless infinitesimal Einstein deformations is trivial:
\begin{equation}\label{eq:toric-kernel}
    K_g^{\T^2} =\left\{ k\in C^\infty(S^2T^*M)^{\T^2}: \delta_gk=0,\quad \tr_gk=0,\quad (\nabla_g^*\nabla_g-2\mathring R_g)k=0 \right\} =\{0\}.
\end{equation}
Here the transverse-traceless conditions are $\delta_gk=\tr_gk=0$, $\delta_g=-\diver_g$, and $\mathring R_g$ is the curvature action on symmetric two-tensors.
\end{definition}

\begin{theorem}[Dichotomy Theorem]
\label{thm:dichotomy}
Let $(M,g_0)$ be a compact toric Einstein four-manifold. Suppose that $g_0$ is torically nondegenerate and its rod directions belong to $\{A,B,C\}$ after one global change of lattice basis. Fix a torus fixed point $p_0$. Then at least one of the following holds:
\begin{enumerate}
\item Every toric manifold obtained by an admissible replacement at $p_0$ of the form \eqref{eq:pointed-insertion} or \eqref{eq:hermitian-words} admits a toric Einstein metric.

\item $M$ admits a toric Einstein metric which is not isometric to any positive constant multiple of $g_0$.
\end{enumerate}
\end{theorem}

The six toric Einstein metrics listed above are torically nondegenerate by Theorem~\ref{thm:toric-small-betti-nondegeneracy}, and their standard torus actions have rod directions in $\{A,B,C\}$ after a change of lattice basis. Theorem~\ref{thm:uniqueness} excludes the second alternative for each of these metrics. Applying the first alternative to the admissible replacements described above yields Theorem~\ref{thm:main}.

Our proof is motivated by the construction of toric gravitational instantons by our previous work  \cite{LiSun}. In the Hermitian setting, there is a striking parallel between compact Einstein four-manifolds and AF or ALF gravitational instantons \cite{BG,LiHermitian}. It is natural to ask whether such a relationship extends beyond the Hermitian setting.

For toric gravitational instantons, the Ricci-flat equation is an integrable system and reduces to an axisymmetric harmonic map equation. General harmonic map theory provides solutions with prescribed rod data, but the resulting metrics usually have conical singularities. In \cite{LiSun} we choose suitable rod structures and vary the rod lengths, using harmonic map bubbling and degree theory to make all cone angles equal to $2\pi$. In the compact case, the quotient metric and its conformal structure are also unknown. Although the Einstein-Hilbert functional admits a useful variational formulation, this does not directly give an existence theory for conical Einstein metrics.

The use of varying cone angles also recalls the solution of the K\"ahler-Einstein problem by Chen-Donaldson-Sun \cite{CDSIII}. In the K\"ahler setting, the moment map interpretation and geodesic convexity provide a clearer variational picture, accompanied by uniqueness modulo automorphisms and local continuation in the cone angle. A central difficulty there is to control degenerations without toric symmetry and use stability to identify the limits. For the general Einstein equation, these structures are unavailable even under toric symmetry, so we work directly with the moduli space of conical Einstein metrics.

We prescribe cone angle $2\pi$ along some invariant spheres and a common variable angle $2\pi\beta$ along the others. After normalization, this moduli space has expected dimension one. We construct a regular end by gluing a conical asymptotically Euclidean Ricci-flat metric to $g_0$ at a fixed point, with $\beta$ tending to a  value $\beta_*>1$. The essential property we need about the gravitational instantons is that, at fixed cone angles, its decaying $\mathbb T^2$-invariant infinitesimal deformations consist only of scaling. Prescribing total area removes this freedom, while varying the common cone angle removes the corresponding obstruction. Together with toric nondegeneracy of $g_0$, this gives the required gluing family.

Our choice of the sequences $W$ is motivated by explicit descriptions of the corresponding toric conical $\Ae$ gravitational instantons. For the alternating sequence, the two rod directions force the Gram matrix to be diagonal, reducing the harmonic-map equation to a scalar axisymmetric Laplace equation. For the four sequences in \eqref{eq:hermitian-words}, the gravitational instantons are Hermitian for a suitable orientation and admit a description by the Tod ansatz. These descriptions allow us to analyze the cone-angle map and prove the uniqueness and deformation properties needed for gluing; see Appendix~\ref{app:ae-models}.

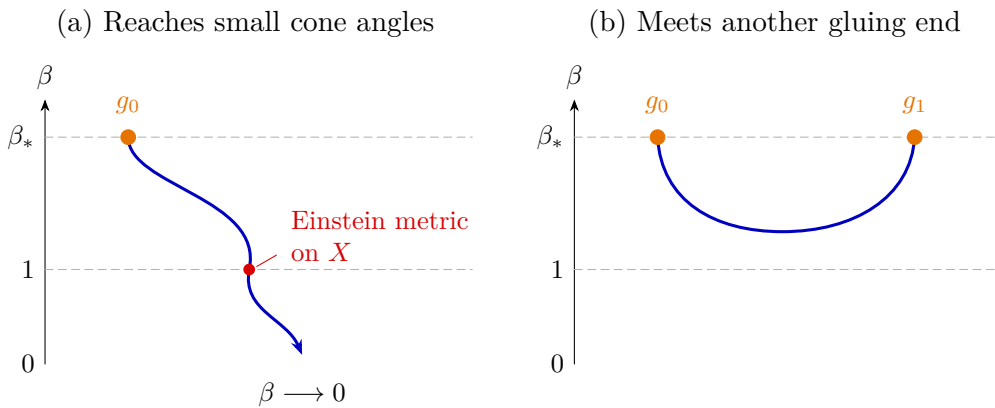
\begin{figure}[H]
\centering
\begin{tikzpicture}[
    x=1cm,y=1cm,font=\small,
    line cap=round,line join=round,
    axis/.style={thin,-{Stealth[length=1.6mm]}},
    level/.style={gray!60,thin,densely dashed},
    curve/.style={blue!75!black,line width=1.15pt},
    endpoint/.style={
        circle,draw=orange!90!black,fill=orange!90!black,
        inner sep=0pt,minimum size=5.5pt
    }]

\foreach \x in {0,7} {
  \begin{scope}[shift={(\x,0)}]
    \draw[axis] (0,0) -- (0,3.5) node[above] {$\beta$};
    \draw[level] (0,3) -- (5.65,3);
    \draw[level] (0,1.25) -- (5.65,1.25);
    \node[left] at (0,3) {$\beta_*$};
    \node[left] at (0,1.25) {$1$};
    \node[left] at (0,0) {$0$};
  \end{scope}
}

\begin{scope}
  \node[font=\normalsize] at (2.65,4.5)
    {(a) Reaches small cone angles};

  \draw[curve,-{Stealth[length=2mm]}]
    (1.1,3)
    .. controls (1.1,2.35) and (2.9,2.2) .. (2.7,1.25)
    .. controls (2.6,.7) and (3.2,.6) .. (3.4,.12);

  \node[endpoint] at (1.1,3) {};
  \node[above=5pt,text=orange!90!black]
    at (1.1,3) {$g_0$};

  \fill[red!85!black] (2.7,1.25) circle (2.2pt);
  \node[anchor=west,align=left,text=red!85!black]
    at (3.12,1.7) {Einstein metric\\on $X$};
  \draw[thin,red!85!black] (3.06,1.44) -- (2.8,1.29);

  \node[below=3pt] at (3.4,0) {$\beta\longrightarrow0$};
\end{scope}

\begin{scope}[xshift=7cm]
  \node[font=\normalsize] at (2.65,4.5)
    {(b) Meets another gluing end};

  \draw[curve]
    (1.1,3)
    .. controls (1.15,2.15) and (1.85,1.75) .. (2.75,1.75)
    .. controls (3.7,1.75) and (4.45,2.2) .. (4.5,3);

  \node[endpoint] at (1.1,3) {};
  \node[endpoint] at (4.5,3) {};
  \node[above=5pt,text=orange!90!black]
    at (1.1,3) {$g_0$};
  \node[above=5pt,text=orange!90!black]
    at (4.5,3) {$g_1$};
\end{scope}

\end{tikzpicture}\caption{Two possible scenario for the moduli space}
\label{fig:cone-angle-continuation}
\end{figure}

The global argument requires diameter bounds, noncollapsing estimates, and a precise analysis of curvature concentration when $\beta$ stays in a compact subset of $(0,\infty)$. For the chosen rod structures, the relevant degenerations are described by reverse gluing near a nondegenerate compact limit. We also prove an upper bound for $\beta$, using Ricci-flat limits and curvature-energy identities. Thus the moduli space cannot escape through unbounded cone angles.

If the relevant component were a regular curve, the picture would be that it either reaches arbitrarily small positive angles, crossing $\beta=1$, or meets another gluing end; see Figure~\ref{fig:cone-angle-continuation}. The actual moduli space can be singular and have several ends. We therefore construct an oriented Kuranishi structure and compare the boundary signs at the gluing ends. All ends arising from a compact metric isometric to $g_0$ have the same sign, so they cannot constitute the entire boundary. Consequently, unless $M$ admits another Einstein metric, the moduli space contains a metric for every $\beta\in(0,1]$. At $\beta=1$ all conical singularities disappear. The uniqueness results for algebraically special toric Einstein metrics then give the existence results in Theorem \ref{thm:main}.

\

The proof of Theorem~\ref{thm:uniqueness}, given in Section~\ref{sec:weyl-rigidity}, combines curvature identities with a variational argument. Maximum principles and an integrated identity for the Weyl curvature force anti-self-duality when $b_2^+=0$ and conformal K\"ahler geometry when $b_2^+=1$. In the spherical case, applying the curvature identities in both orientations yields $W^+=W^-=0$, so the metric has constant positive sectional curvature. In the Hermitian cases, the integrated boundary fluxes, which are proportional to the areas of the invariant spheres, determine the moment polygon and the affine conformal factor. After fixing the total flux to remove scaling, the square root of the volume is a strictly concave function of these fluxes. This follows by expressing the volume as a minimum of quadratic forms of Lorentzian signature. The Einstein equation makes this function stationary, so the fluxes are uniquely determined. With the polygon and conformal factor fixed, an integration-by-parts identity then identifies the metrics, proving uniqueness up to scaling and isometry.

\

The proof of Theorem~\ref{thm:general} is a mild extension of the strategy used in the proof of Theorem~\ref{thm:main}. The main additional technical step is to carry out gluing over a finite-dimensional stabilization of the compact moduli space, removing the need to assume toric nondegeneracy of the intermediate compact metrics. More precisely, we use a count modulo 2 of normalized toric Einstein metrics with prescribed rod data. When the rod directions belong to $\{A,B,C\}$, our compactness results and Ricci-flat rigidity already give a compact moduli space of smooth toric Einstein metrics with a fixed rod structure. Finite-dimensional stabilization then defines the count even in the presence of degenerate metrics. We choose the perturbing tensors to be divergence-free and supported away from the chosen fixed point, allowing gluing over the resulting smooth solution space. Gluing identifies the small-area boundary count with the count on the contracted manifold, while compactness and the cone-angle bound give a cobordism to the smooth boundary. Thus the count is preserved under the prescribed corner replacements. Every cyclic rod structure with directions in $\{A,B,C\}$ reduces to a bigon or triangle, whose count is one. Counting all equivalence classes also gives the advantage of avoiding any additional discrete symmetry assumption on the intermediate compact metrics, which is needed if one wants to prove Theorem \ref{thm:main} using mod 2 count. Finally, the topological classification realizes realizes every closed simply connected smooth four-manifold admitting an effective torus action and satisfying $\chi\geq3|\tau|$ by such a rod structure, while the improved Hitchin-Thorpe inequality gives necessity.

\

The paper is organized as follows. Section~\ref{sec:rod structure} introduces compact rod structures, toric conical metrics, and the reduced Einstein equations. Section~\ref{sec:a priori estimates} establishes the geometric estimates and in particular, the a priori diameter bound and the angle control implies non-collapsing theorem. Section~\ref{sec:moduli} develops the local deformation theory, the orientation of the moduli space, and the finite-dimensional perturbation argument. Section~\ref{sec:end structure} studies degeneration and bubbling, proves compactness, and establishes the gluing construction and its converse. Section~\ref{sec:angle-bounds} gives an upper bound for the common cone angle. Section~\ref{sec:weyl-rigidity} proves the classification, uniqueness, and toric nondegeneracy results for algebraically special metrics, together with the improved Hitchin-Thorpe inequality. Section~\ref{sec:degree} determines the gluing signs and completes the proofs of the existence and dichotomy theorems. Section~\ref{sec:toric-existence-criterion} extends the counting argument to allow degenerate compact metrics and proves the general existence criterion. Section~\ref{sec:further-questions} discusses open problems.  Appendices~\ref{app:asymptotics of toric gravitational instantons} and~\ref{app:ae-models} establish the asymptotic classification and the uniqueness and deformation results for the gravitational instantons used in the argument. Appendix~\ref{app:fixed-rod-smooth-compactness} proves smooth compactness for simply connected compact toric Einstein four-manifolds with bounded Euler characteristic, normalized by $\Ric_g=3g$. Together with the construction in Section~\ref{sec:toric-existence-criterion}, this defines a virtual count modulo two for every compact rod structure.

\

\noindent\textbf{Acknowledgements.}

\

We are grateful to  Dr. Tianyue Liu for explaining his interesting paper on toric Einstein four-manifolds to us in the summer of 2025, and to Christoph B\"ohm and  Jeff Viaclovsky for pushing us to extend our work on toric gravitational instantons to the compact setting. The second author thanks Ramiro Lafuente for related discussions on toric Einstein metrics during the AMSI Winter School in Brisbane in July 2026. In particular, Lafuente informed him that he and K.~Broder had a proof of the uniqueness of toric Einstein metrics on $\mathbb S^4$. The first author thanks the IASM for its hospitality during the summer of 2026. Parts of this work were carried out during these visits.

\

\noindent\textbf{AI Usage Statement.}

\

This project began with the attempt to extend our previous work \cite{LiSun} to the compact setting. In particular, we were interested in studying the variational problem, before clever (advanced) AI tools became available. We could not construct conical Einstein metrics directly through this approach, by prescribing the rod fluxes, and regarded the problem as a project for the years ahead. When AI tools became available, we tried the problem as an experiment in mathematical collaboration with AI, but these initial attempts did not produce an existence argument. They did, however, reveal difficulties with prescribing the rod fluxes that led us to set this approach aside.

The authors then had the idea of studying the moduli space of conical Einstein metrics directly, without any extrinsic parametrization. This is inspired by gauge-theoretic moduli-space methods, particularly the use of compactness and boundary counts in Donaldson theory. Our initial idea was to glue a double Eguchi–Hanson gravitational instanton to the Chen–LeBrun–Weber metric and vary a single cone angle to obtain a new Einstein metric on $2\CP^2\#2\overline{\CP^2}$. We identified three essential ingredients: gluing, weak compactness and noncollapsing for conical toric Einstein metrics, and control of bubbling, so that the only possible degenerations in the range of cone angles under consideration would reverse the gluing construction. This led to the continuation picture in Figure~\ref{fig:cone-angle-continuation}. If successful, the argument would give either an Einstein metric on $2\CP^2\#2\overline{\CP^2}$ or another toric Einstein metric on $\CP^2\#2\overline{\CP^2}$, not isometric to a positive multiple of the Chen–LeBrun–Weber metric. Numerical experiments carried out with AI suggested that the first possibility should hold.

At this stage, we expected the gluing analysis should follow from  the work of Biquard and Ozuch, and the singular moduli space to be treated using standard Kuranishi methods, as in the moduli-space arguments underlying Donaldson’s famous diagonalization theorem. However, weak compactness and noncollapsing required substantially more work than we initially anticipated. Being aware of the increasing power of AI, we tried to seek further assistance from AI. Although the individual calculations were not especially difficult, finding the differential inequalities in Section~\ref{subsec:diff inequality} involved considerable trial and error. We regard AI’s contribution to finding these inequalities as its most significant contribution to this project.

AI also contributed substantially to the uniqueness argument. The authors were less familiar with calculations involving $W^+$. After repeated corrections and challenges from the authors, AI found the identities used in the uniqueness proof. Uniqueness excludes the second alternative above and gives the desired existence result. From the curvature identities developed in this argument, the authors subsequently observed the improved Hitchin–Thorpe inequality and formulated a conjecture describing the diffeomorphism types of compact four-manifolds admitting toric Einstein metrics.

The authors also designed the rod structures for the more general existence results in Theorem~\ref{thm:main}. We chose structures with at most three rod directions to control the possible degenerations and to use explicit gravitational instantons in the gluing construction. The particular choice of the replacement sequence $W$ is guided by our own understanding of toric gravitational instantons. In particular, there are two classes of explicit families of toric gravitational instantons suitable for gluing, one from the case with two rod directions (this gives the infinite families) and one from Hermitian AE gravitational instantons (this gives the finite families). The proofs of rigidity and infinitesimal rigidity for these gravitational instantons, though fairly elementary computations, were developed by AI.

Finally, the authors asked AI to extend this strategy to prove the conjectured classification of diffeomorphism types. The proof of this final extension and the compactness result in Appendix \ref{app:fixed-rod-smooth-compactness} was developed entirely by AI, using the results and strategy established above. 

Overall, this paper was developed with substantial AI assistance and represents our first active experiment in working with AI. Beyond language editing, mathematical writing, and figure preparation, AI acted as a ``super assistant'', responding continually to our ideas and helping us test and develop them. We believe the project could in principle have been completed without AI, but would have required years of sustained effort and experimentation.

\section{Compact rod structures and toric conical Einstein metrics}\label{sec:rod structure}
\subsection{Compact rod structures and topology}
We will always fix the standard lattice $\Lambda=(2\pi\mathbb Z)^2\subset\R^2$ and denote $\T^2=\R^2/\Lambda$. There are standard angular coordinates $\varphi_1,\varphi_2$ on $\mathbb T^2$ with primitive period $2\pi$. In particular,
\begin{equation}\label{eq:torusnormalization}
    \int_{\T^2}\dd\varphi_1\dd\varphi_2=(2\pi)^2.
\end{equation}

\begin{definition}\label{def:rod}
A compact rod structure consists $\mathfrak R$ of the following data
\begin{itemize}
\item an oriented closed disk $B$;
\item $m\geq 2$ distinct points $\mathfrak p_1,\ldots,\mathfrak p_m\in\partial B$ in positive cyclic order, with $\mathfrak p_{m+1}:=\mathfrak p_1$;
\item a primitive vector $\mv_i\in\Lambda/\{\pm1\}$ on each positively oriented boundary arc $\mI_i=[\mathfrak p_i,\mathfrak p_{i+1}]$,
\end{itemize}
such that for $i=1, \cdots, m$, $\mv_i$ and $\mv_{i+1}$ form a basis of $\Lambda$.

Each $\mathfrak p_i$ is called a \emph{turning point} (sometimes also referred as a \emph{vertex}), each $\mI_i$ is called a \emph{rod} and each $\mv_i$ is called a \emph{rod vector}.

Two compact rod structures are equivalent if, after a cyclic relabelling, there are an orientation-preserving homeomorphism $f:B\to B'$ and an orientation-preserving lattice automorphism $\mathfrak L:\Lambda\to\Lambda$ such that
$$
    f(\mI_i)=\mI_i', \qquad \mathfrak L(\mv_i)=\mv_i' \quad\text{in }\Lambda/\{\pm1\}.
$$
\end{definition}
\begin{remark}
This definition is the counterpart of the rod structure defined in \cite{LiSun} for studying gravitational instantons. For our geometric purpose in this paper, we choose a simplified version than \cite{LiSun}. In particular, we fix the lattice $\Lambda$ and avoid discussing enhancements etc.
\end{remark}

Given a compact rod structure $\mR$, we set
$$
    X_{\mathfrak R}=(B\times\T^2)/\sim,
$$
where over the interior of $\mI_i$ the circle generated by $\mv_i$ is collapsed and over every $\mathfrak p_i$ the full torus is collapsed. The local standard models are $\R^2\times\T^2$ over the interior, $\C\times S^1\times\R$ over a rod, and $\C^2$ over a turning point. Then $X_{\mR}$ is a smooth four-manifold endowed with the obvious $\mathbb T^2$-action, and the obvious quotient map $\pi: X_{\mR}\rightarrow B$. We denote
$$
    X_{\mR}^\circ:=\pi^{-1}(B^\circ), \qquad \Sigma_i:=\pi^{-1}(\mI_i), \qquad p_i:=\pi^{-1}(\mathfrak p_i).
$$
Then
$$
    X_{\mR} =X_{\mR}^{\circ}\sqcup \left(\bigcup_{i=1}^{m}\Sigma_i\right).
$$
and each $\Sigma_i=\pi^{-1}(\mI_i)$ is an embedded two-sphere. For simplicity sometimes we will refer to such a two sphere as a \emph{rod sphere}. Each turning point $\mathfrak p_i\in\partial B$ corresponds to a unique torus fixed point $p_i\in X$, with $\pi(p_i)=\mathfrak p_i$.

We orient $X_{\mR}$ by requiring that $\pi^*\nu\wedge \dd\varphi_1\wedge \dd\varphi_2$ to be positive on $X_{\mR}^\circ$ for every positive local volume form $\nu$ on $B$. Notice that if $\mR'$ is equivalent to $\mR$, then there is an orientation-preserving diffeomorphism between $X_\mR$ and $X_{\mR'}$, which intertwines the torus actions through a lattice automorphism.

We now describe the local models for the $\T^2$-action on $X_{\mR}$. Over a point of $B^\circ$, a $\T^2$-invariant neighborhood has coordinates
$$
    U\simeq D^2_b\times\T^2_{\psi_1,\psi_2},\qquad (\varphi_1,\varphi_2)\cdot(b,\psi_1,\psi_2) =(b,\psi_1+\varphi_1,\psi_2+\varphi_2).
$$
Over an interior point of $\mI_i$, choose a signed representative $\mv_i$ and complete it to a lattice basis $(\mv_i,\mw_i)$. Let $\theta_1,\theta_2$ be the corresponding angular coordinates of period $2\pi$, so that
$$
    \partial_{\theta_1} =\sum_{\alpha=1}^2 \frac{(\mv_i)_\alpha}{2\pi}\partial_{\varphi_\alpha}, \qquad \partial_{\theta_2} =\sum_{\alpha=1}^2 \frac{(\mw_i)_\alpha}{2\pi}\partial_{\varphi_\alpha}.
$$
A $\mathbb T^2$-invariant neighborhood of the corresponding circle orbit has coordinates
$$
    U\simeq D^2_z\times(-\epsilon,\epsilon)_s\times S^1_\psi, \qquad (\theta_1,\theta_2)\cdot(z,s,\psi) =(e^{\sqrt{-1}\theta_1}z,s,\psi+\theta_2),
$$
with $\Sigma_i\cap U=\{z=0\}$.

At $p_i$, choose signed representatives $(\mv_{i-1},\mv_i)$ forming a lattice basis, and use its angular coordinates $\theta_1,\theta_2$ of period $2\pi$. An invariant neighborhood has coordinates
$$
    U\simeq D^2_{z_1}\times D^2_{z_2}, \qquad (\theta_1,\theta_2)\cdot(z_1,z_2) =(e^{\sqrt{-1}\theta_1}z_1,e^{\sqrt{-1}\theta_2}z_2),
$$
with
$$
    \Sigma_{i-1}\cap U=\{z_1=0\}, \qquad \Sigma_i\cap U=\{z_2=0\}.
$$

Conversely, every closed, simply connected, oriented smooth four-manifold with an effective smooth $\T^2$-action is equivariantly diffeomorphic to $X_{\mR}$ for some compact rod structure $\mR$, unique up to equivalence, with the orientation of $B$ chosen to agree with the given orientation of the manifold. The rod data determines the action up to equivariant diffeomorphism \cite{OrlikRaymond}. Moreover, the underlying oriented smooth manifold is either $\mathbb S^4$ or a connected sum of copies of $\mathbb S^2\times \mathbb S^2$, $\mathbb{CP}^2$, and $\overline{\mathbb{CP}}^{ 2}$.

We collect some easy topological facts about $X_\mR$. Recall that we fix the lattice to be $\Lambda=(2\pi\mathbb{Z})^2$.

\begin{itemize}
\item $\pi_1(X_\mR)=\{1\}$;
\item $\chi(X_\mR)=m$;
\item For each $i$, the self-intersection $e_i$ of $\Sigma_i$ is given by
$$
    (2\pi)^{-6}\det(\mv_{i-1}, \mv_{i})\det(\mv_i, \mv_{i+1})\det(\mv_{i-1}, \mv_{i+1}).
$$
\item $3\tau(X_\mR)=\sum_i e_i$.
\end{itemize}
The last formula is a consequence of the signature theorem. When $m=2$ it holds trivially. For $m\geq3$, the stable tangent-bundle splitting \cite[Theorem~6.6]{DavisJanuszkiewicz} gives
$$
    p_1(TX)=\sum_{i=1}^m(\operatorname{PD}[\Sigma_i])^2.
$$

\subsection{Compact toric conical metrics}
Roughly speaking, Riemannian metrics with cone singularities along a smooth codimension-two submanifold $S$ are locally modelled, to leading order, on the product of a Euclidean space and a two-dimensional flat cone. The regularity required beyond this model depends on the particular setting. In K\"ahler geometry, a by-now standard analytic formulation for cone singularities along a smooth divisor is given by Donaldson \cite{DonaldsonCone}. Conical K\"ahler-Einstein metrics played a central role in the solution of the K\"ahler-Einstein problem for Fano manifolds \cite{CDS1,CDS2,CDSIII}. For the toric metrics considered here, we use the formulation of conical regularity in terms of the rod structure developed in \cite{LiSun}.

\begin{definition}\label{def:conical}
Given a compact rod structure $\mR$, a \emph{toric conical metric} on $X_{\mR}$ is a $\mathbb T^2$-invariant Riemannian metric $g$ on $X_{\mR}^{\circ}$ whose metric completion is naturally identified with $X_{\mR}$ and which satisfies the following conditions for constants $\beta_i>0$.

In the adapted local coordinates, lift each collapsing angular coordinate to $\mathbb R$. We require the following
\begin{itemize}
\item Near an interior point of $\mI_i$, the lifted metric is the pullback of a smooth Riemannian metric across $\{z=0\}$ under
 $$
    (r,\theta,s,\psi) \longmapsto (re^{\sqrt{-1}\beta_i\theta},s,\psi).
$$

 \item Near $p_i$, the lifted metric is the pullback of a smooth Riemannian metric across both coordinate axes under
 $$
    (r_1,\theta_1,r_2,\theta_2) \longmapsto \bigl(r_1e^{\sqrt{-1}\beta_{i-1}\theta_1}, r_2e^{\sqrt{-1}\beta_i\theta_2}\bigr).
$$
\end{itemize}
These smooth extensions are called the \emph{local developments} of $g$. The cone angle along $\Sigma_i$ is $2\pi\beta_i$.
\end{definition}

The cone angle along the rod $\mI_i$ is $2\pi\beta_i$. When $\beta_i=1$ for all $i$, then the metric is an ordinary smooth metric. The conical singularities we defined here are very mild-they arise from an identification with an incorrect period, similar to orbifold singularities. In particular, the curvature is bounded across the conical singularities.

A toric conical Einstein metric is by definition a toric conical metric which satisfies the Einstein equation
$$
    \Ric=\lambda g
$$
on $X_{\mR}^\circ$.

The following is well-known \cite{Dammerman, Liu}. For simplicity we give simple explanation following \cite{LiSun}

\begin{lemma}\label{lem:polar}
A toric conical Einstein metric on $X_\mR$ can be globally written over $X_{\mR}^\circ$ as
\begin{equation}\label{eq:block}
    g=\bg+G_{\alpha\beta}\dd\varphi_\alpha\dd\varphi_\beta, \qquad
\end{equation}
where $G=\rho\Phi$, with $\rho>0$ on $B^\circ$ and
$$
    \Phi:B^\circ\longrightarrow \mathcal H=SL(2,\mathbb R)/SO(2).
$$
We identify $\mathcal H$ with the space of positive definite symmetric $2\times2$ matrices of determinant one.
\end{lemma}
\begin{proof}
We follow the derivation in \cite{LiSun}. Let $X_j=\partial_{\varphi_j}$, $j=1,2$, be the Killing fields of period $2\pi$, and set $\eta_j=g(X_j,\cdot)$.  The Killing identities and the Einstein equation imply that $*(\eta_j\wedge\dd\eta_j)$ is closed. Thus locally there are twist potentials $\mathcal E_j$ satisfying
$$
    \dd \mathcal E_j=*(\eta_j\wedge\dd\eta_j).
$$
Set
$$
    \Theta_j=*(\eta_1\wedge\eta_2\wedge\dd\eta_j).
$$
As in \cite{LiSun}, $\Theta_1=\mathcal L_{X_2}\mathcal E_1$
$$
    \dd\Theta_1=\dd\mathcal L_{X_2}\mathcal{E}_1=\mathcal L_{X_2} (*(\eta_1\wedge \dd\eta_1))=0,
$$
and similarly, $\dd\Theta_2=0$.

Now the $\mathbb T^2$-action has at least one fixed point, at which both $\Theta_1$ and $\Theta_2$ tend to zero. It follows that $\Theta_1\equiv\Theta_2\equiv 0$. This implies that the horizontal distribution is integrable. Since $B^\circ$ is simply connected, its flat connection has a global horizontal section, giving the trivialization in \eqref{eq:block}. Finally, set $\rho=\sqrt{\det G}$ and $\Phi=\rho^{-1}G$.
\end{proof}

In this paper we will endow $\mathcal H$ with the standard Riemannian metric with constant negative curvature $-2$:
$$
    g_{\mathcal H}=\frac{1}{4}\tr(\Phi^{-1}\dd\Phi)^2.
$$

\subsection{Reduction of the Einstein equation}

For compactly supported variations over $B^\circ$, consider the Einstein-Hilbert functional
\begin{equation}\label{eqn:EH}
    \EH_\lambda(g)=\int_X(\Scal_g-2\lambda)\dvol_g.
\end{equation}
It is straightforward to compute the scalar curvature
\begin{equation}
    \Scal_g=\Scal_{\bg}-2\rho^{-1}\Delta_{\bg}\rho +\tfrac12\rho^{-2}|\dd\rho|_{\bg}^2 -\tfrac14\tr(\Phi^{-1}\dd\Phi)^2_{\bg}.\label{eq:scalarreduction}
\end{equation}

\begin{proposition}\label{prop:reduced system}
For a metric of the form \eqref{eq:block}, the Einstein equation $\Ric(g)=\lambda g$ is equivalent on $B^\circ$ to the following system:
\begin{align}
    &\Delta_{\bg}\rho=-2\lambda\rho, &&\text{(fiber trace)},\label{eq:rho}\\
    &\dd_{\bg}^{*}\bigl(\rho\Phi^{-1}\dd\Phi\bigr)=0, &&\text{(fiber trace-free)},\label{eq:map}\\
    &K_{\bg} =-\frac12|\dd\log\rho|_{\bg}^{2} +\frac18\tr\bigl((G^{-1}\dd G)^2\bigr)_{\bg}, &&\text{(base trace)},\label{eq:nu}\\
    &\left( \Hess_{\bg}\log\rho +\frac14\tr\bigl(G^{-1}\dd G\otimes G^{-1}\dd G\bigr) \right)^{\circ}=0, &&\text{(base trace-free)}.\label{eq:stress}
\end{align}
Here $K_{\bg}$ denotes the Gaussian curvature of $\bg$.
\end{proposition}
\begin{proof}

Variation of $\bg$ in the conformal direction gives~\eqref{eq:rho}. Variation of $\Phi$ gives~\eqref{eq:map}. For a compactly supported variation of $\rho$, its coefficient is
$$
    2K_{\bg}-\rho^{-1}\Delta_{\bg}\rho +\frac{|\dd\rho|_{\bg}^2}{2\rho^2} -\frac14\tr\bigl((\Phi^{-1}\dd\Phi)^2\bigr)_{\bg} -2\lambda.
$$
Using~\eqref{eq:rho} and the fact
$$
    \tr\bigl((G^{-1}\dd G)^2\bigr)_{\bg} =2|\dd\log\rho|_{\bg}^2 +\tr\bigl((\Phi^{-1}\dd\Phi)^2\bigr)_{\bg},
$$
we get~\eqref{eq:nu}.

Finally, a compactly supported trace-free variation of $\bg$ gives
$$
    \left( \Hess_{\bg}\rho -\frac{\dd\rho\otimes\dd\rho}{2\rho} +\frac{\rho}{4} \tr(\Phi^{-1}\dd\Phi\otimes\Phi^{-1}\dd\Phi) \right)^\circ=0.
$$
Here the Hessian term comes from integrating by parts in the variation of $\int_B\rho R_{\bg} \dd A_{\bg}$. Dividing by $\rho$, and using
$$
    \rho^{-1}\Hess_{\bg}\rho =\Hess_{\bg}\log\rho +\dd\log\rho\otimes\dd\log\rho,
$$
gives \eqref{eq:stress}. The mixed horizontal-vertical Ricci components vanish for \eqref{eq:block}, so the four equations together are equivalent to the Einstein equation.
\end{proof}

\begin{remark}
In \cite{LiSun}, for the Ricci-flat case (i.e., $\lambda=0$), \eqref{eq:rho} implies that $\rho$ is a harmonic function on the base, so $\rho+\sqrt{-1}z$, for a harmonic conjugate $z$, gives a natural local conformal coordinate on the base. This fact allows us to write the whole system in terms of an axisymmetric harmonic map from a domain in $\R^3$ to $\mathcal H$.

When $\lambda\neq 0$, $\Phi$ is still a weighted harmonic map, and $\rho$ is an eigenfunction for $\Delta_{\underline g}$. So we no longer have natural local conformal coordinates. This is a main difficulty to directly solve the above PDE system.
\end{remark}

\subsection{Rod asymptotics}
We now derive the rod asymptotics for toric conical Einstein metrics. We first specify the coordinates used near rods and turning points.

Whenever the quotient metric $\bg$ is smooth up to the interior of a rod $\mI_i$, choose an arclength coordinate $s$ along the rod and extend it by keeping it constant along the inward unit-speed geodesics perpendicular to the rod. For sufficiently small $r\geq0$, where $r$ is the distance to the rod, this gives local Fermi coordinates $(r,s)$ satisfying
$$
    \bg=\dd r^2+h(r,s)^2\dd s^2,\qquad h(0,s)=1.
$$

At a turning point $\mathfrak p_i$, choose signed representatives $(\mv_{i-1},\mv_i)$ forming a basis of $\Lambda$. The local model for the torus action gives smooth equivariant coordinates $(z_1,z_2)$ near the corresponding fixed point $p_i$, in which the two circle factors act by separate rotations and $\Sigma_{i-1}=\{z_1=0\}$, $\Sigma_i=\{z_2=0\}$. Writing
$$
    z_1=xe^{\sqrt{-1}\theta_1},\qquad z_2=ye^{\sqrt{-1}\theta_2},\qquad x,y\geq0,
$$
gives local quotient coordinates $(x,y)$, with the two rods corresponding to the coordinate axes.

In the following local expansions, we express the Gram matrix in the adapted angular coordinates:
$$
    G_{ab}=g(\partial_{\theta_a},\partial_{\theta_b}).
$$
The function $\rho=\sqrt{\det G}$ is unchanged by these integral changes of angular basis.

\begin{lemma}\label{lem:rodasymptotics}
Given a compact rod structure $\mR$, and data $(\underline g, \rho, \Phi)$ which solves \eqref{eq:rho}-\eqref{eq:stress}, then it determines a toric conical Einstein metric with cone angle $2\pi\beta_i$ along $\Sigma_i$, if and only if the pair $(\underline g,\rho\Phi)$ satisfies the following boundary conditions.
\begin{enumerate}
\item[(1)] The metric $\underline g$ extends smoothly across the interior of each $\mI_i$.
\item[(2)] At an interior point of $\mI_i$, choose the basis $(\mv_i,\mw_i)$ of $\Lambda$ and let $(r,s)$ be local Fermi coordinates on $B$. Then we have
$$
    \bg=\dd r^2+f(r^2,s)^2\dd s^2
$$
    for a smooth nonzero function $f$ with $f(0, s)=1$. Moreover, choose an equivariant chart adapted to these coordinates,
$$
    (z,s,\theta_2),\qquad z=re^{\sqrt{-1}\theta_1},
$$
    such that $\Sigma_i=\{z=0\}$, then we require that the Gram matrix in these angular coordinates satisfies
\begin{align*}
    G_{11} &=r^2\bigl(\beta_i^2+r^2a(r^2,s)\bigr),\\
    G_{12} &=r^2b(r^2,s), \\
    G_{22} &=c(r^2,s),
\end{align*}
    Here $a, b, c$ are all local smooth functions and $c(0, s)>0$.
\item[(3)] At a fixed point $p_i$, choose the basis $(\mv_{i-1},\mv_{i})$ of $\Lambda$ and local equivariant coordinates
$$
    z_1=xe^{\sqrt{-1}\theta_1},\qquad z_2=ye^{\sqrt{-1}\theta_2}, \qquad x,y\geq0,
$$
    so that $\Sigma_{i-1}=\{z_1=0\}$ and $\Sigma_{i}=\{z_2=0\}$. Then we require $\underline g$ and $G$ to take the form
\begin{align}
    & \bg      =a(x^2,y^2)\dd x^2 +2xyb(x^2,y^2)\dd x\dd y +c(x^2,y^2)\dd y^2,\nonumber\\
    & G        =\begin{pmatrix} x^2p(x^2,y^2)    &x^2y^2q(x^2,y^2)\\
    x^2y^2q(x^2,y^2) &y^2r(x^2,y^2)
\end{pmatrix},\label{eq:cornerexpansion}\\
& p(0,y^2)=\beta_{i-1}^2a(0,y^2), \qquad r(x^2,0)=\beta_i^2c(x^2,0),\nonumber
\end{align}
    where $a, b, c, p, q, r$ are all smooth local functions, and $a(0, 0)>0, c(0, 0)>0$.
\end{enumerate}
\end{lemma}

\begin{proof}
For a point in the interior of a rod $\mI_i$, and let $(r, s)$ be the Fermi coordinates on $B$. Consider an adapted equivariant chart $(r, \theta_1, s, \theta_2)$, if we develop the collapsing $S^1$ factor by setting
$$
    \vartheta=\beta_i\theta_1,
$$
and denote
$$
    X=r\cos\vartheta,\quad Y=r\sin\vartheta,\qquad \alpha=X\dd Y-Y\dd X=r^2\dd\vartheta.
$$
Then the smoothness and rotational invariance of the new metric precisely mean that it can be written as
$$
    \widetilde g =\dd X^2+\dd Y^2+q(r^2,s)^2\dd s^2 +A(r^2,s)\alpha^2 +2B(r^2,s)\alpha \dd\theta_2 +C(r^2,s)\dd\theta_2^2,
$$
where $q,A,B,C$ are smooth, $q(0,s)=1$, and $C(0,s)>0$. Substituting $\alpha=\beta_i r^2\dd\theta_1$ gives
$$
    \bg=\dd r^2+q(r^2,s)^2\dd s^2,\quad G_{11}=r^2\bigl(\beta_i^2+ \beta_i^2r^2A(r^2,s)\bigr),\quad G_{12}=\beta_i r^2B(r^2,s),\quad G_{22}=C(r^2,s).
$$
The discussion at a fixed point is similar and we omit it.
\end{proof}

As a consequence, for a toric conical Einstein metric, the quotient metric $\underline g$ extends smoothly to $B$ as a Riemannian metric on a surface with corners. The interior of each rod is a geodesic boundary arc, and adjacent arcs meet at right angles at each turning point $\mathfrak p_j$. Moreover, we have the boundary asymptotics for $\rho$ as a function on $B$

\begin{equation}\label{eq:rhoends}
    \begin{cases}
        \rho=r\kappa(s)+O(r^3),\quad \kappa(s)=\beta_i\sqrt{c(0,s)}>0, &\text{along an open rod},\\
        \rho=xy\tau(x^2,y^2),\quad \tau(0,0)>0, &\text{at a turning point}.
    \end{cases}
\end{equation}
where $\kappa$ and $\tau$ are smooth local functions.

The induced metric on $\Sigma_i$ is given by $\dd s^2+\beta_i^{-2}\kappa^2\dd\theta_2^2$, with Gaussian curvature
$$
    K_{\Sigma_i}=-\frac{\kappa''}{\kappa}.
$$
In the local smooth redeveloping model, $\Sigma_i$ is a fixed-point surface and is therefore totally geodesic. By Lemma \ref{lem:polar} the horizontal section is also totally geodesic. Thus the sectional curvature of a plane spanned by $\partial_s$ and a normal direction equals $K_{\bg}$. Taking the Ricci curvature in the unit $\partial_s$ direction yields
\begin{equation}\label{eqn:relation ambient and sigma curvature}
    \lambda=K_{\Sigma_i}+2K_{\bg} =-\frac{\kappa''}{\kappa}+2K_{\bg}.
\end{equation}

The area of each $\Sigma_i$ can be expressed as the rod flux of $\rho$:
\begin{equation}\label{eq:rod-flux}
    \Area_g(\Sigma_i)=-\frac{2\pi}{\beta_i}\int_{\mI_i}\partial_n\rho \dd s_{\bg}
\end{equation}
where $n$ is the outward unit normal to the rod in $\bg$.. Indeed, in Fermi coordinates $n=-\partial_r$ and $-\partial_n\rho=\beta_i\sqrt{c(0,s)}$, whereas
$$
    \Area_g(\Sigma_i) =2\pi\int_{\mI_i}\sqrt{c(0,s)} \dd s_{\bg}.
$$
The total volume of $(X, g)$ is given by
\begin{equation}\label{eq:volumeangle}
    \Vol_g(X)=\lambda^{-1}\pi\sum_i\beta_i\Area_g(\Sigma_i).
\end{equation}
This follows from the divergence theorem and \eqref{eq:rho}:
$$
    \lambda \Vol_g(X)=(2\pi)^2\lambda\int_B\rho \dd\mu_{\bg}= -2\pi^2\int_B\Delta_{\bg}\rho \dd\mu_{\bg} =-2\pi^2\sum_i\int_{\mI_i}\partial_n\rho \dd s_{\bg} =\pi\sum_i\beta_i\Area_g(\Sigma_i).
$$

We also have a version of Bishop-Gromov volume comparison
\begin{equation}\label{geom:volumeupper}
    \operatorname{Vol}_g(X) \leq \frac{24\pi^2}{\lambda^2} \min_{1\leq i\leq m}\beta_i\beta_{i+1}, \qquad \beta_{m+1}=\beta_1.
\end{equation}
Indeed, by scaling the lattice, we can arrange that all the cone angles are below $2\pi$, so the convexity argument of \cite[Lemma~4.43]{LiSun} permits the usual Bonnet-Myers and Bishop-Gromov comparison. At the torus fixed point $\Sigma_i\cap\Sigma_{i+1}$, the  volume density is $\beta_i\beta_{i+1}$, while its total volume is $\operatorname{Vol}_g(X)$.

\section{A priori estimates and geometric bounds}\label{sec:a priori estimates}

Throughout this section we will work with a compact toric conical Einstein metric $(X_\mR, g)$ with Einstein constant $\lambda>0$. We will derive some a priori estimates and geometric bounds. All the quantities on $B$ are calculated using the quotient metric $\bg$. For simplicity of discussion we introduce some notations
$$
    u=\log\rho,\qquad U=\frac{1}{2}|\dd u|^2,\qquad K=K_{\bg},\qquad T_{ij}=\frac14\tr(\Phi^{-1}\Phi_i\Phi^{-1}\Phi_j).
$$
So $T$ is non-negative. Denote
$$
    e=\frac{1}{2}\tr T=\frac{1}{8}\tr(\Phi^{-1}\dd\Phi)^2\geq 0
$$
and
$$
    \det T:=\frac{\det(T_{ij})}{\det(\underline g_{ij})}\geq 0.
$$
Notice that
\begin{equation}\label{eqn:det trace}
    2\det T+|T|^2=(\tr T)^2.
\end{equation}
Denote the weighted Laplacian operator
$$
    \cL=\Delta_{\bg}+\langle \dd u,\dd \cdot \rangle.
$$
The reduced system in Proposition \ref{prop:reduced system} implies the following equations
\begin{align}
    &\Delta u=-2\lambda-2U \label{geom:eigenfunction},\\
    &2e=2K+U,\label{geom:scalartrace}\\
    &\Hess u=(K-\lambda)\bg-\tfrac12\dd u\otimes \dd u-T, \label{geom:hessu} \\
    & \tau(\Phi)+\dd\Phi(\nabla u)=0. \label{geom:harmonic map}
\end{align}
Here $\tau(\Phi):=\tr(\nabla \dd\Phi)$ is the tension field of the map $\Phi$.

\subsection{A differential inequality}\label{subsec:diff inequality}
We now prove a key differential inequality. One motivation is Lu's interior curvature estimate for special K\"ahler metrics \cite{Lu}. Locally, these correspond to the above system with constant $\rho$ and $\lambda=0$. Then in this situation \eqref{geom:hessu} gives $T=K\bg$, thus $\Phi$ is conformal and harmonic. On $\{K>0\}$, the metric $K\bg=\Phi^*g_{\mathcal H}$ has curvature $-2$, and hence
$$
    \Delta_{\bg}\log K=6K.
$$
Such an equation gives a priori interior estimates on $K$, without extra assumptions. In our Einstein setting, we have a more complicated differential inequality.

\begin{lemma}\label{geom:scalar}
For a toric conical Einstein metric, at a point in $B^\circ$ where $\dd\Phi\neq 0$, we have
\begin{equation}\label{geom:scalarineq}
    \Delta K+\frac{3U+2\lambda}{2e}\langle \dd u,\dd K\rangle \geq \frac{|\dd K|^2}{e} +\frac{K(3K-\lambda)(2U+K+\lambda)}{e} +3\det T.
\end{equation}
\end{lemma}

\begin{proof}
First, the weighted Bochner formula for a function gives
$$
    \cL U=|\Hess u|^2+\langle \nabla u, \nabla \cL u\rangle+(\Ric-\Hess u)(\nabla u, \nabla u).
$$
Since $\cL u=-2\lambda$, and the weighted Ricci curvature is
\begin{equation}\label{eqn:weighted Ricci}
    \Ric_{\bg}-\Hess u=\lambda\bg+\frac12 \dd u\otimes \dd u+T,
\end{equation}
we obtain
\begin{equation}\label{geom:bochneru}
    \cL U=|\Hess u|^2+2\lambda U+2U^2+T(\nabla u,\nabla u).
\end{equation}
Similarly, the weighted Bochner formula for the density of $\Phi$ gives, using \eqref{geom:harmonic map} and $g_{\mathcal H}$ has curvature $-2$,
\begin{align}
    \mathcal L e ={} &|\nabla\dd\Phi|^2 +\langle\Ric_{\bg}-\Hess_{\bg}u,T\rangle_{\bg}+2(\tr T)^2-2|T|^2\\
    =                &|\nabla \dd\Phi|^2+\lambda\tr T+\tfrac12T(\nabla u,\nabla u) +2(\tr T)^2-|T|^2. \label{geom:bochnermap-general}
\end{align}
Using \eqref{geom:scalartrace}, \eqref{geom:bochneru}, we obtain
\begin{equation}\label{geom:LK}
    \cL K=|\nabla \dd\Phi|^2-\tfrac12|\Hess u|^2+2\lambda K +2(\tr T)^2-|T|^2-U^2.
\end{equation}
By \eqref{geom:hessu} we have
$$
    \dd U=(K-\lambda-U)\dd u-T(\nabla u,\cdot).
$$
Together with \eqref{geom:scalartrace} and \eqref{geom:harmonic map}, this gives
$$
    \dd K+\tfrac12(K-\lambda-U)\dd u =\sum_{j=1}^2 \left\langle(\nabla\dd\Phi)^\circ(\cdot,e_j), \dd\Phi(e_j)\right\rangle,
$$
where $e_1,e_2$ is any local orthonormal frame and $(\nabla \dd\Phi)^\circ$ denotes the trace-free part of $\nabla \dd\Phi$.

Now we recall an elementary fact. For a symmetric trace-free $(2\times 2)$ matrix $A$, we have
$$
    A^2=\frac{1}{2}\|A\|^2\Id
$$
so
$$
    |Av|^2=\frac{1}{2}\|A\|^2|v|^2.
$$
Apply this  to each component of $(\nabla\dd\Phi)^\circ$ in an orthonormal frame of $\Phi^*T\mathcal H$, and then use Cauchy-Schwarz in the target index, we obtain
$$
    \begin{aligned}
        \left|\dd K+\tfrac12(K-\lambda-U)\dd u\right|^2 &= \left|\sum_{j=1}^2 \left\langle (\nabla\dd\Phi)^\circ(\cdot,e_j),\dd\Phi(e_j) \right\rangle\right|^2\\
        &\leq \tfrac12|(\nabla\dd\Phi)^\circ|^2|\dd\Phi|^2 =e|(\nabla\dd\Phi)^\circ|^2.
    \end{aligned}
$$
Finally,
$$
    |\nabla\dd\Phi|^2 =|(\nabla\dd\Phi)^\circ|^2 +\tfrac12|\tau(\Phi)|^2 =|(\nabla\dd\Phi)^\circ|^2 +\tfrac12T(\nabla u,\nabla u).
$$
Consequently, on $\{e>0\}$, we have
\begin{equation}\label{geom:refinedbochner}
    |\nabla\dd\Phi|^2\geq \frac1e\left|\nabla K+\tfrac12(K-\lambda-U)\nabla u\right|^2 +\tfrac12T(\nabla u,\nabla u).
\end{equation}
On the other hand, by \eqref{geom:hessu}
$$
    |\Hess u|^2=2(K-\lambda)^2-2(K-\lambda)U-2(K-\lambda)\tr T +U^2+|T|^2+T(\nabla u,\nabla u).
$$
Inserting these into \eqref{geom:LK}, by \eqref{eqn:det trace} and \eqref{geom:scalartrace} and a straightforward manipulation we obtain \eqref{geom:scalarineq}.
\end{proof}

A first consequence of the differential inequality is a curvature upper bound on $B$.

\begin{proposition}\label{geom:Kupper}
We have
\begin{equation}\label{geom:Kconcavity}
    K\leq \frac{\lambda}{3}.
\end{equation}
\end{proposition}

\begin{proof}
Let $K_0$ be the maximum of $K$ on $B$. Suppose that $K_0>\lambda/3$. Then $e=K+U/2>0$ at every maximum. An interior maximum contradicts \eqref{geom:scalarineq}, since $\dd K=0$, $\Delta_{\bg}K\leq0$, and its right-hand side is strictly positive. It remains to bound a boundary maximum. Since $e=K+U/2$,
$$
    \frac{3U+2\lambda}{2e} =3+\frac{\lambda-3K}{e}.
$$
At an interior point of a rod, by \eqref{eq:rhoends} we have $u_r=r^{-1}+O(r)$, $u_s=O(1)$, and $K=O(1)$. Thus $U=(2r^2)^{-1}+O(1)$, and
$$
    e=K+\frac{U}{2}=(4r^2)^{-1}+O(1).
$$
Since $K$ is even in $r$, at a maximum $(0,s_0)$ we have
$$
    \lim_{r\to0} \left(\Delta_{\bg}K+ \frac{3U+2\lambda}{2e}\langle \dd u,\dd K\rangle_{\bg}\right)(r,s_0) =4K_{rr}(0,s_0)+K_{ss}(0,s_0)\leq0.
$$
On the other hand, the second term on the right-hand side of \eqref{geom:scalarineq} tends to $4K_0(3K_0-\lambda)>0$, while the other terms are nonnegative. This is a contradiction.

At a vertex, choose the coordinates $(x,y)$ as in Lemma \ref{lem:rodasymptotics}, with
$$
    \rho=xy a(x^2,y^2),\qquad \bg(0)=dx^2+dy^2.
$$
Along $x=y=t\to 0$, we have
$$
    u_x=t^{-1}+O(t),\quad u_y=t^{-1}+O(t),\quad U=t^{-2}+O(1),\quad e=(2t^2)^{-1}+O(1).
$$
Since $K$ is even in both variables,
$$
    K_x=tK_{xx}(0)+O(t^3),\qquad K_y=tK_{yy}(0)+O(t^3).
$$
Consequently,
$$
    \Delta_{\bg}K+ \frac{3U+2\lambda}{2e}\langle \dd u,\dd K\rangle_{\bg} \to 4\bigl(K_{xx}(0)+K_{yy}(0)\bigr)\leq0.
$$
The second term again tends to $4K_0(3K_0-\lambda)>0$, and the remaining terms are nonnegative. This contradiction proves $K\leq\lambda/3$ on $B$.
\end{proof}

By \eqref{eqn:relation ambient and sigma curvature}, it follows that
\begin{corollary}
Along the interior of a rod, we have
\begin{equation}\label{geom:rodidentity}
    \kappa''\leq-\frac\lambda3\kappa<0.
\end{equation}
\end{corollary}

A consequence of the weighted Bochner formula gives
\begin{proposition}
For $\mathfrak p\in B^\circ$, we have
\begin{equation}\label{geom:principalgradient}
    U(\mathfrak p)\leq\frac{12}{d(\mathfrak p, \p B)^2}.
\end{equation}
\end{proposition}

\begin{proof}

Fix $\mathfrak p\in B^\circ$ and $R<d_{\bg}(\mathfrak p,\partial B)$, and let $r=d_{\bg}(\mathfrak p,\cdot)$. At a smooth point of $r$, the Bochner formula applying to the equation $\mathcal L(|\nabla r|^2)=0$ and \eqref{eqn:weighted Ricci} give
$$
    \begin{aligned}
        \partial_r(\mathcal Lr) &=-|\Hess_{\bg}r|^2 -(\Ric_{\bg}-\Hess_{\bg}u)(\partial_r,\partial_r)\\
        &\leq-(\Delta_{\bg}r)^2-\frac12(u_r)^2 \leq-\frac13(\mathcal Lr)^2.
    \end{aligned}
$$
Since $\mathcal Lr=r^{-1}+O(1)$ near $\mathfrak p$, a scalar Riccati comparison gives $\mathcal Lr\leq3r^{-1}$. So
$
    \mathcal L(r^2)\leq8.
$
Set
$$
    F=(R^2-r^2)^2U.
$$
The weighted Bochner formula \eqref{geom:bochneru} implies that
$
    \mathcal LU\geq 2U^2.
$
If $F$ has a positive maximum at a smooth point of $r^2$, then we have
$$
    \begin{aligned}
        0\geq\mathcal LF \geq2(R^2-r^2)^2U^2-16(R^2-r^2)U-24r^2U.
    \end{aligned}
$$
So
$$
    F\leq8(R^2-r^2)+12r^2\leq12R^2.
$$
If the maximum lies in the cut locus, the same estimate holds by a standard barrier argument.

Now evaluate at $\mathfrak p$ and letting $R$ increase to $d_{\bg}(\mathfrak p,\partial B)$ we get the conclusion
\end{proof}
By \eqref{geom:scalartrace} we get an interior lower bound on $K$ and and upper bound on the map density $e$.
\begin{corollary}
We have
\begin{equation}\label{eqn:Klower}
    K\geq - \frac{6}{d(\mathfrak p, \p B)^2}, \qquad \  e\leq \frac{\lambda}{3}+\frac{6}{d(\mathfrak p, \p B)^2}.
\end{equation}
\end{corollary}

Together with Proposition \ref{geom:Kupper}, we obtain uniform interior bound on $K$. Using standard elliptic theory for the reduced system we also obtain a uniform bound on the curvature of $g$ away from the rods. Later we will prove a stronger bound including the interior of the rods (see Proposition \ref{geom:vertexcurvature}).

\subsection{Diameter and injectivity radius estimate}
\begin{proposition}\label{geom:diameter}
For a toric conical Einstein metric $(X_{\mR}, g)$ with $\lambda>0$ we have the uniform estimate
\begin{equation}\label{geom:diameterbounds}
    \frac{\pi}{2\sqrt{2\lambda}}\leq\operatorname{diam}(X_{\mR},g) \leq2\pi\sqrt{\frac 3{\lambda}}.
\end{equation}
\end{proposition}
\begin{remark}
The upper bound is expected from Myers's theorem in the smooth setting. The proof below is a more direct argument without dealing the conical singularities, but gives slightly weaker bound which suffices for our purpose.
\end{remark}
Denote
$$
    M:=\max_B\rho.
$$

\begin{lemma}
We have a uniform bound
\begin{equation}\label{geom:rhogradient}
    |\dd\rho|^2+2\lambda\rho^2\leq2\lambda M^2.
\end{equation}
\end{lemma}

\begin{proof}
Denote $P=|\dd\rho|^2+2\lambda\rho^2$. By the Bochner formula, we have
$$
    \tfrac12\Delta P=|\Hess\rho|^2+K|\dd\rho|^2-4\lambda^2\rho^2.
$$
Suppose otherwise, then at a maximum $\mathfrak q\in B$ we have $P>2\lambda M^2$. At a corner point $\mathfrak p_i$ we have $\rho=xya(x^2,y^2)$, so $P=0$. In particular, $\mathfrak q$ can not be a corner point.

If $\mathfrak q$ is in $B^\circ$, then it follows that $\dd\rho\neq 0$ at $\mathfrak q$. So we can use an orthonormal frame $(e_1, e_2)$ at $\mathfrak q$, with $e_1=\nabla\rho/|\dd\rho|$. At $\mathfrak q$ we have $\dd P=0$. For any tangent vector $X$, differentiating $P$ gives
$$
    \dd P(X)=2\Hess \rho(\nabla\rho,X)+4\lambda\rho\dd\rho(X).
$$ 
Then the fact that $\Delta\rho=-2\lambda\rho$ gives
$$
    \Hess\rho(e_1,e_1)=-2\lambda\rho,\qquad \Hess\rho(e_1,e_2)=\Hess\rho(e_2,e_2)=0.
$$
In particular, $\Hess u(e_2, e_2)=0$ and it follows from \eqref{geom:hessu} that $K\geq\lambda>0$ at $\mathfrak q$. It follows that $\Delta P>0$ at $\mathfrak q$. Contradiction.

Suppose now that $\mathfrak q$ is in the interior of a rod. We work with a local Fermi coordinates. Along the rod we have $P(0,s)=\kappa(s)^2$. So $\kappa'=0$, $\kappa''\leq0$ at $\mathfrak q$. The equation $\Delta\rho=-2\lambda\rho$ gives
$$
    \rho_{rrr}-K\kappa+\kappa''=-2\lambda\kappa.
$$
Differentiating $P$ and then using \eqref{eqn:relation ambient and sigma curvature} gives that at $\mathfrak q$
$$
    P_r=0,\qquad P_{rr}=2\kappa\rho_{rrr}+2(\kappa')^2+4\lambda\kappa^2 =\lambda\kappa^2-\kappa\kappa''>0.
$$
Contradiction.
\end{proof}

\begin{proof}[Proof of Proposition \ref{geom:diameter}]
For $\epsilon>0$, \eqref{geom:rhogradient} implies
$$
    |\dd\arcsin(\frac{\rho}{M+\epsilon})|\leq\sqrt{2\lambda}.
$$
Integrate from a maximum point $\mathfrak q$ of $\rho$ to the boundary and let $\epsilon\rightarrow0$. This gives
$$
    d(\mathfrak q,\partial B)\geq\frac{\pi}{2\sqrt{2\lambda}}
$$
and hence the lower diameter bound.

For the upper bound, we first notice that
\begin{equation}\label{eqn:hessian bound sqrt rho}
    \Hess\sqrt\rho=\frac{\sqrt\rho}{2}\big[(K-\lambda)\bg-T\big] \leq-\frac\lambda3\sqrt\rho \bg.
\end{equation}
Along a minimizing geodesic in $B^\circ$ from a maximum point of $\sqrt\rho$,  from the Hessian bound we see that
$$
    \sqrt{\rho(t)}\leq \sqrt{\rho(0)}\cos(\sqrt{{\lambda}/{3}} t).
$$
This implies that
$$
    \operatorname{diam}B\leq\pi\sqrt{{3}/{\lambda}}.
$$
On $X_\mR$, the distance to a torus fixed point equals distance to its image in the quotient. Joining two points through any fixed point therefore gives $\operatorname{diam}(X_\mR, g)\leq2\operatorname{diam}B$.
\end{proof}

\begin{lemma}\label{geom:basecomparison}
We have
$$
    \operatorname{inj}(\mathfrak p)\geq d(\mathfrak p,\partial B), \quad \sup_{\mathfrak p\in B^\circ}\operatorname{inj}(\mathfrak p)\geq \frac{\pi}{2\sqrt{2\lambda}}.
$$
Moreover,
\begin{equation}\label{eqn:diameter lower bound}
    \operatorname{diam}(\partial B) \geq \frac{\pi}{\sqrt{2\lambda}}.
\end{equation}
\end{lemma}

\begin{proof}
Arguing as in the proof of Proposition \ref{geom:diameter} we know that every geodesic starting in $B^\circ$ reaches $\partial B$ before time $\pi\sqrt{\frac{3}{\lambda}}$. Fix $\mathfrak{p}\in B^\circ$, and let $D_{\mathfrak{p}}\subset T_{\mathfrak{p}}B$ be the star-shaped domain of vectors whose radial geodesics remain in $B^\circ$. Since $K\leq \frac{\lambda}{3}$, Jacobi field comparison shows that these geodesics have no conjugate points before reaching the boundary. Hence $\exp_{\mathfrak{p}}:D_{\mathfrak{p}}\to B^\circ$ is a local diffeomorphism. Since the rods are totally geodesic and the vertices have angle $\pi/2$, a geodesic starting in $B^\circ$ cannot touch $\partial B$ and then continue within $B$. By continuous dependence on initial data, a sequence approaching $\partial D_{\mathfrak p}$ therefore has its endpoints approaching $\partial B$. Since $D_{\mathfrak p}$ is bounded, $\exp_{\mathfrak p}:D_{\mathfrak p}\to B^\circ$ is proper. It is also onto, since a minimizing curve between two interior points stays in $B^\circ$.  As $D_{\mathfrak{p}}$ is connected and $B^\circ$ is simply connected, $\exp_{\mathfrak{p}}$ is a diffeomorphism. In particular, every interior geodesic segment minimizes, and
$$
    \operatorname{inj}(\mathfrak{p})\geq d(\mathfrak{p},\partial B).
$$
Again as in the proof of Proposition \ref{geom:diameter} we have for a maximum point $\mathfrak q$ of $\rho$,
$$
    d(\mathfrak q,\partial B)\geq\frac{\pi}{2\sqrt{2\lambda}}.
$$
Extend a geodesic through $\mathfrak q$ in both directions until it reaches boundary points $\mathfrak q_1$ and $\mathfrak q_2$. Its interior segments minimize, so, by continuity of distance, the entire segment from $\mathfrak q_1$ to $\mathfrak q_2$ minimizes. Both halves have length at least $d(\mathfrak q,\partial B)$. Therefore
$$
    \operatorname{diam}(\partial B) \geq d(\mathfrak q_1,\mathfrak q_2) \geq 2d(\mathfrak q,\partial B) \geq\frac{\pi}{\sqrt{2\lambda}}.
$$
\end{proof}

The proof of Lemma~\ref{geom:basecomparison} shows that $B^\circ$ is geodesically convex and has no cut or conjugate points. Since $K\leq\lambda/3$, the usual Rauch comparison gives $\operatorname{CAT}(\lambda/3)$ triangle comparison on $B^\circ$. Passing to its metric completion $B$ we get

\begin{corollary}
\label{geom:cat}
The quotient $(B,d_{\bg})$ is a $\operatorname{CAT}(\lambda/3)$ space.
\end{corollary}

\subsection{Curvature estimate away from fixed points}
Let $\mathcal V$ denote the set of torus fixed points in $M$. The following provides a uniform curvature bound for the curvature of conical Einstein metrics, away from $\mathcal V$.

\begin{proposition}\label{geom:vertexcurvature}
There is a universal constant $C$ such that for a compact toric conical Einstein manifold $(M,g)$ and $x\notin \mathcal V$ we have
\begin{equation}\label{geom:vertexcurvatureeq}
    |\operatorname{Rm}_g|(x) \leq C\max\left\{\frac{\lambda}{3}, d_g(x,\mathcal V)^{-2}\right\}.
\end{equation}
\end{proposition}
Here the curvature at a point in $\Sigma_i$ is to be understood in terms of local smooth development, as in the definition of toric conical metrics Definition \ref{def:conical}.

We first prove a lemma. Recall that the interior of each rod is a totally geodesic boundary arc of $(B,\bg)$. For $\mathfrak q\in I_i^\circ$, let $\nu_{\mathfrak q}$ denote the inward unit normal and define
$$
    D_r^+(\mathfrak q) := \left\{ v\in T_{\mathfrak q}B: |v|_{\bg}<r,\quad \bg(v,\nu_{\mathfrak q})\geq 0 \right\}.
$$
We define the \emph{injectivity radius at $\mathfrak q$}, denoted by $\operatorname{inj}_{B,\bg}(\mathfrak q)$, to be the supremum of all $r>0$ such that the exponential map
$$
    \exp_{\mathfrak q}: D_r^+(\mathfrak q) \longrightarrow B_{\bg}(\mathfrak q,r)
$$
is well-defined and is a diffeomorphism of manifolds with boundary. Notice that the ball $B_{\bg}(\mathfrak q, r)$ does not intersect any other rod for $r<\operatorname{inj}_{B,\bg}(\mathfrak q)$. We will also use the double of $(B,\bg)$, obtained by gluing two copies of $B$ by the identity along $\partial B$. By Lemma \ref{lem:rodasymptotics}, near the interior of each rod the quotient metric has the form
$$
    \bg=\dd r^2+f(r^2,s)^2\dd s^2.
$$
Allowing $r$ to take both signs therefore extends the metric smoothly across the rod. At each turning point, the corner expansion extends smoothly under both coordinate reflections, and the doubled neighborhood is the quotient of this smooth extension by $(x,y)\mapsto(-x,-y)$. Thus the double is a compact surface without boundary, smooth away from the turning points, with cone angle $\pi$ at each turning point. We call this the even double of $(B,\bg)$.

\begin{lemma}\label{geom:rodhalfdisk}
Suppose a point $\mathfrak q$ lies in the interior of a rod $\mI_i$, and suppose that $\overline B_{\bg}(\mathfrak q,2R)$ contains no turning points and that
$$
    |\operatorname{Rm}_g|\leq c \quad\text{on}\quad \pi^{-1}\bigl(B_{\bg}(\mathfrak q,2R)\cap B^\circ\bigr),
$$
Then we have
$$
    \operatorname{inj}_{B,\bg}({\mathfrak q})\geq \min\left\{R,\frac{\pi}{2\sqrt c}\right\},
$$
\end{lemma}

\begin{proof}
Denote $r= \min\left\{R,\frac{\pi}{2\sqrt c}\right\}$. Suppose $\operatorname{inj}_{B,\bg}({\mathfrak q})<r$, then one of the following must occur
\begin{itemize}
\item There is a unit-speed geodesic segment $\gamma:[0, \ell]\rightarrow B$, for $\ell<r$, going from $\mathfrak q$ and reaches the interior of another rod orthogonally.
\item There is a unit-speed geodesic segment $\gamma:[0, \ell]\rightarrow B$, for $\ell<2r$, going from $\mathfrak q$ back to $\mI_i$.
\end{itemize}
In both cases, we have $\gamma((0,\ell))\subset B^\circ$.

To see this, if another rod meets $B_{\bg}(\mathfrak q,r)$, a shortest segment to the other rods gives the first case. Otherwise, we pass to the even base double, where the curvature bound excludes conjugate points before $r$. A closest cut point at distance $s<r$ produces a geodesic loop of length $2s$, contained in the doubled radius-$s$ ball. The loop cannot remain on the rod, which is an interval; its first excursion into one copy of $B$ therefore returns to $\mI_i$, giving the second case.

Let $\widetilde\gamma:[0,\ell]\to X_\mR$ be a horizontal lift of $\gamma$, so that $\pi\circ\widetilde\gamma=\gamma$. Since the horizontal sections are totally geodesic, $\widetilde\gamma$ is a unit-speed geodesic. Let $J$ be the Killing Jacobi field along $\widetilde\gamma$ generated by the circle vanishing on $\mI_i$. Then $J(0)=0$, and $J$ is nonzero in the interior. At $t=0$, this is understood in the smooth local development. We claim that at $t=\ell$, $\langle D_tJ,J\rangle=0.$ Again this is understood in the smooth local development. In the second case, the claim follows from the fact that $J(\ell)=0$. In the first case, notice that at the other end point, the Killing field is tangent to the corresponding totally geodesic fixed surface, whereas the lifted geodesic is normal to it. The claim then follows from $\langle D_tJ,J\rangle =-\langle\nabla_JJ,\dot{\widetilde\gamma}\rangle=0$.

Integration of the Jacobi equation gives
$$
    \int_0^\ell |D_tJ|^2 dt \leq c\int_0^\ell |J|^2 dt.
$$
Standard Jacobi field comparison gives
$$
    \ell\geq \frac{\pi}{2\sqrt c}
$$
in the first case, and the second case, the stronger endpoint condition $J(0)=J(\ell)=0$ gives
$$
    \ell\geq \frac{\pi}{\sqrt c}.
$$
In both cases we are lead to a contradiction.
\end{proof}

\begin{proof}[Proof of Proposition \ref{geom:vertexcurvature}]

Suppose that \eqref{geom:vertexcurvatureeq} fails. Then there is a sequence of toric conical Einstein metrics $(X_j, g_j)$, with $\operatorname{Ric}_{g_j}=\lambda_j g_j$, such that
$$
    M_j:= \max_{x\in X_j} a_j(x)^2|\Rm_{g_j}|(x) \longrightarrow\infty, \qquad a_j(x):= \min\left\{ d_{g_j}(x,V_j),\sqrt{{3}/{\lambda_j}} \right\}.
$$
Choose a maximizing point $x_j$, and set
$$
    Q_j:=|\operatorname{Rm}_{g_j}|(x_j), \qquad \widehat g_j:=Q_jg_j, \qquad R_j:=\frac12\sqrt{Q_j} a_j(x_j) =\frac12\sqrt{M_j}\rightarrow\infty.
$$
Since $a_j$ is $1$-Lipschitz, we have
$$
    |\operatorname{Rm}_{\widehat g_j}|(x_j)=1, \qquad \sup_{B_{\widehat g_j}(x_j,R_j)} |\operatorname{Rm}_{\widehat g_j}|\leq4.
$$
Moreover,
$$
    d_{\widehat g_j}(x_j,\mathcal V_j)\geq2R_j\rightarrow\infty, \qquad \operatorname{Ric}_{\widehat g_j} =\widehat\lambda_j\widehat g_j, \qquad \widehat\lambda_j:=\frac{\lambda_j}{Q_j} \leq\frac{3}{M_j}\longrightarrow0.
$$

Now we divide into two cases:

\

\noindent {\bf Case (1).} All the fixed surfaces escape to infinity.

In this case the limiting procedure does not see conical singularities. For sufficiently small $\delta>0$, independent of $j$, the balls $B_{\widehat g_j}(x_j,2\delta)$ are smooth. The local covering consequence of \cite[Theorem~1.3 and conditions~(1.1.4)-(1.1.5)]{CFG} gives a uniform volume lower bound on the universal covers:
$$
    \Vol_{\widetilde g_j} B_{\widetilde g_j}(\widetilde x_j,\delta) \geq c\delta^4.
$$
Here $\widetilde g_j$ is the lift of $\widehat g_j$. The curvature bound and interior Einstein regularity give smooth convergence (after passing to a subsequence) on smaller lifted balls to a Ricci-flat metric $\widetilde g_\infty$, with
$$
    |\Rm_{\widetilde g_\infty}|(\widetilde x_\infty)=1.
$$
We now show that $\widetilde g_\infty$ is flat. Here $U_j,K_j,e_j,T_j$ are associated with the rescaled metric $\widehat g_j$. For every fixed $A>0$, the distance from $B_{\widehat g_j}(x_j,A)$ to the union of the fixed surfaces tends to infinity. Thus \eqref{geom:principalgradient} and \eqref{eqn:Klower}, together with $\widehat\lambda_j\to0$, give
$$
    U_j\longrightarrow0,\qquad e_j\longrightarrow0
$$
uniformly on this ball. By \eqref{geom:scalartrace}, $K_j=e_j-U_j/2\to0$ uniformly as well. Finally, since $T_j$ is nonnegative,
$$
    |T_j|_{\underline{\widehat g}_j} \leq\operatorname{tr}_{\underline{\widehat g}_j}T_j =2e_j\longrightarrow0.
$$
Normalize the constant real angular basis so that $G_j(x_j)=I$. By integrating along horizontal directions it follows that $G_j\to I$ uniformly on every fixed smaller neighborhood. Thus, in the smooth local equivariant limit, the metric becomes
$$
    \widetilde g_\infty=\bg_\infty+(d\theta^1)^2+(d\theta^2)^2.
$$
Since $K_j\to0$, the base metric $\bg_\infty$ is flat. So $\widetilde g$ is also flat, contradiction.

\

\noindent {\bf Case (2).} There exists a fixed surface $\Sigma_j$ within bounded distance to $x_j$.

In this case away from the fixed surfaces we  still have the limit of local universal cover around as above, which is Ricci-flat. But we can not deduce a contradiction directly. Instead, the flatness comes from studying the local convergence near a point in $\Sigma_j$.

We choose a point $\mathfrak q_j\in \partial B_j$ closest to $\pi_j(x_j)$, and let $\mI_j$ be the rod containing $\mathfrak q_j$. After passing to a subsequence, we have
$$
    d_{\widehat{\bg}_j}(\pi_j(x_j),\mathfrak q_j)\leq D
$$
for some fixed $D<\infty$. Let $q_j$ be the endpoint of a horizontal lift, starting at $x_j$, of a shortest segment from $\pi_j(x_j)$ to $\mathfrak q_j$.

We parametrize $\mI_j$ by arclength $s$, with $s(\mathfrak q_j)=0$. Since all vertices escape to infinity, there are intervals $[-t_j,t_j]\subset\mI_j$, with $t_j\to\infty$ and $t_j+D<R_j$, such that on $\pi_j^{-1}([-t_j, t_j])$ the curvature of $\widehat g_j$ is bounded by $4$. Notice that by \eqref{geom:rodidentity}, $\kappa_j$ is positive and concave. So we have
$$
    1-\frac{|s|}{t_j}\leq \frac{\kappa_j(s)}{\kappa_j(0)} \leq1+\frac{|s|}{t_j}, \qquad \frac{|\kappa_j'(0)|}{\kappa_j(0)}\leq\frac1{t_j}.
$$
Thus $\kappa_j(s)/\kappa_j(0)\to1$ uniformly on each fixed arclength interval.

By Lemma~\ref{geom:rodhalfdisk}, a quotient half-disk of radius $\delta>0$ about $\mathfrak q_j$, with $\delta$ independent of $j$, meets no rod other than $\mI_j$. Let $g_j^*$ be the smooth angular development across this rod. We normalize the vanishing Killing field to have unit normal rotation and a complementary Killing field to have unit length at $q_j$. Subtract a constant multiple of the vanishing field from the complementary field so that its derivative has zero normal rotational component at $q_j$. This does not change its value there. Its remaining derivative is controlled by $|\kappa_j'(0)|/\kappa_j(0)$.

The developed metrics $g_j^*$ have uniformly bounded curvature. Applying the local covering argument from Case~(1), we obtain smooth convergence on smaller lifted balls to a Ricci-flat metric $g_\infty$. The normalized values and derivatives of the Killing fields, together with $\nabla^2V=\Rm*V$, give smooth convergence of both fields. Thus the limit satisfies the reduced Einstein equations with $\lambda=0$. In this angular normalization, the boundary coefficient of $\rho$ is $\kappa_j(s)/\kappa_j(0)$, so $\kappa=1$ along the limiting rod. Moreover, the curvature upper bound and \eqref{eqn:relation ambient and sigma curvature} give
$$
    K\leq0, \qquad 0=2K-\frac{\kappa''}{\kappa} \quad\text{on the rod}.
$$
So we must have $K=0$ along the rod.

We claim that $K$ vanishes in a collar of the rod. Notice that we have \eqref{geom:scalarineq} for $\lambda=0$:
\begin{equation}\label{eqn:diff ineq Ricci flat}
    \Delta_{\bg_\infty}K+\frac{3U}{2e}\langle \dd u,\dd K\rangle_{\bg_\infty} \geq \frac{|dK|^2}{e} +\frac{3K^2(2U+K)}{e} +3\det T.
\end{equation}
Since $e\geq0$ and $2e=2K+U$, we have $2U+K\geq0$. By the rod expansion,
$$
    U=\frac1{2r^2}+O(1), \qquad e=K+\frac U2=\frac1{4r^2}+O(1),
$$
so $e>0$ in a sufficiently small punctured collar.

To apply the strong maximum principle across the rod, we use Fermi coordinates $(r,s)$, where $r\geq0$ is distance to the rod. The smooth rod expansions give
$$
    \bg_\infty=\dd r^2+f(r^2,s)^2 \dd s^2,\qquad \rho=r+O(r^3),\qquad U=\frac1{2r^2}+O(1),\qquad e=\frac1{4r^2}+O(1).
$$
Consequently, the operator on the left of \eqref{eqn:diff ineq Ricci flat} can be written as
$$
    \mathcal A =\partial_{rr}+\frac3r\partial_r +a(r^2,s)r\partial_r +b(r^2,s)\partial_{ss} +c(r^2,s)\partial_s,
$$
where $a,b,c$ are smooth and $b>0$. Regard $r$ as the radius $|x|$ in $\mathbb R^4$. Since $K$ is smooth and even in $r$, the function $\widehat K(x,s):=K(|x|,s)$ is smooth, and $\mathcal A$ becomes the smooth uniformly elliptic operator
$$
    \widehat{\mathcal A} =\Delta_{\mathbb R^4} +a(|x|^2,s)x\cdot\nabla_x +b(|x|^2,s)\partial_{ss} +c(|x|^2,s)\partial_s.
$$
The inequality $\widehat{\mathcal A}(\widehat K)\geq0$ extends across $x=0$ by continuity. Since $\widehat K\leq0$ and $\widehat K(0,s)=0$, the strong maximum principle implies $K\equiv0$ in a collar of the rod.

We now show directly that the four-dimensional limit metric $g_\infty$ is flat there. Since the base is flat and the rod is totally geodesic, its Fermi coordinates satisfy
$$
    \bg_\infty=\dd r^2+\dd s^2.
$$
The function $\rho$ satisfies
$$
    \Delta_{\bg_\infty}\rho=0,\qquad \rho(0,s)=0,\qquad \partial_r\rho(0,s)=1.
$$
Uniqueness for these harmonic Cauchy data gives $\rho=r$, hence $u=\log r$.

The reduced base equation is now
$$
    \operatorname{Hess}_{\bg_\infty}u =-\frac12 \dd u\otimes \dd u-T.
$$
Evaluating on $\partial_s$ gives $T_{ss}=0$. Since $T=\Phi^*g_{\mathcal H}$, this implies $\partial_s\Phi=0$. Together with $\partial_s\rho=0$ and $G=\rho\Phi$, we obtain $\partial_sG=0$. This means that we can write
$$
    g_\infty=\dd s^2+h,\qquad h=\dd r^2+G_{ab}(r) \dd\theta^a \dd\theta^b.
$$
Since $g_\infty$ is Ricci-flat, the three-dimensional metric $h$ is Ricci-flat and therefore flat. Thus $g_\infty$ is flat in the collar.

Let $\gamma_j$ be the shortest quotient geodesic from $\mathfrak q_j$ to $\pi_j(x_j)$, parametrized by arclength from $\mathfrak q_j$. Since $\mathfrak q_j$ is a closest boundary point,
$$
    d_{\underline{\widehat g}_j} (\gamma_j(t),\partial B_j)=t.
$$
Its horizontal lift has length at most $D$. Near $q_j$ use the smooth angular development above. Outside a fixed small initial segment, the displayed identity gives a uniform positive distance from the fixed surfaces, so the local covering argument applies on balls of uniform radius.

Cover the lifted segment by a uniformly bounded number of overlapping neighborhoods. Choose lifts successively along the segment. On smaller neighborhoods of the overlap points, the lifted metrics are related by local isometries. The curvature bounds, local noncollapse, and Einstein regularity allow both the metrics and these local isometries to converge after passing to a subsequence.

The first limiting neighborhood contains the flat collar and is therefore flat by analyticity. Each successive Ricci-flat limiting neighborhood contains an open subset locally isometric to the preceding one. Analyticity therefore propagates flatness along the finite chain, including to the limit of $x_j$. This contradicts $|\Rm_{\widehat g_j}|(x_j)=1$.
\end{proof}

\subsection{Noncollapsing under angle bound}

In this subsection, we prove a volume noncollapsing theorem.
\begin{theorem}\label{geom:noncollapse}
Fix a compact rod structure $\mR$ and fix $0<\beta_-\leq\beta_+<\infty$. There exists $\sigma>0$ such that any compact toric conical Einstein metric $g$ with all cone angles contained in $[\beta_-,\beta_+]$ satisfies
$$
    \Vol(X_{\mR},g)\geq \sigma\lambda^{-2}.
$$
\end{theorem}

We will need the following elementary fact.
\begin{lemma}
\label{geom:rodmatrixcone} Denote
$$
    J=\begin{pmatrix} 0 &-1\\
    1 &0
\end{pmatrix}, \qquad A_i=(J\mv_i)(J\mv_i)^\top.
$$
Then at every point $\mathfrak{p}\in B^\circ$, we can write
\begin{equation}\label{geom:rodmatrices}
    G(\mathfrak{p})=\sum_i c_i(\mathfrak{p})A_i,\qquad c_i(\mathfrak{p})\geq0.
\end{equation}
Moreover,
\begin{equation}\label{geom:detmatrices}
    \det G=\sum_{i<j}c_ic_j\det(\mv_i,\mv_j)^2, \qquad G(v,v)=\sum_i c_i\det(v,\mv_i)^2.
\end{equation}
\end{lemma}

\begin{proof}
With our definition, $A_i\mathfrak v_i=0$. Notice that $\Phi$ is a weighted harmonic map into $\mathcal H$:
$$
    \dd^*(\rho\Phi^{-1}\dd\Phi)=0.
$$
Since $\det\Phi=1$, the matrix-valued one-form $\Phi^{-1}\dd\Phi$ is traceless. We apply $\mathcal L$ to matrix-valued functions entrywise in the fixed torus basis. Thus
$$
    (\mathcal L\Phi)_{\alpha\beta}=\Delta_{\bg}\Phi_{\alpha\beta}+\langle \dd u,\dd\Phi_{\alpha\beta}\rangle_{\bg}.
$$
Choose a local $\bg$-orthonormal frame $(e_1,e_2)$ and write $\Phi_i=\dd\Phi(e_i)$. Expanding the above equation gives
$$
    \Phi^{-1}\cL\Phi=\sum_i(\Phi^{-1}\Phi_i)^2.
$$
By the Cayley-Hamilton identity for traceless $2\times2$ matrices,
$$
    A^2=\frac12\tr(A^2)\Id,
$$
and hence
$$
    \cL\Phi=\frac12\tr\bigl((\Phi^{-1}\dd\Phi)^2\bigr)\Phi=4e\Phi.
$$
Let $\ell$ be a linear functional on the space of symmetric matrices such that $\ell(A_i)\geq0$ for every $i$. Then
$$
    (\cL-4e)\ell(\Phi)=0.
$$
At the interior of a rod $I_i$, the rod asymptotics in Lemma \ref{lem:rodasymptotics} gives
$$
    G=c(s)A_i+O(r^2),\qquad \rho\sim r,
$$
so $\ell(\Phi)$ has nonnegative lower limit. At a vertex between $I_i$ and $I_j$, the corner expansion gives
$$
    G=r_1^2aA_j+r_2^2bA_i+O(r_1^2r_2^2), \qquad a,b>0,
$$
and hence $\ell(\Phi)\geq-Cr_1r_2$. Thus a negative minimum of $\ell(\Phi)$ would be attained in $B^\circ$, contradicting the maximum principle. Therefore $\Phi$, and hence $G$, lies in the closed cone generated by the $A_i$ at each point, proving \eqref{geom:rodmatrices}.

Now we have
$$
    G(v, v)=\sum_i c_i (v^\top J\mv_i)^2=\sum_i c_i \det(v, \mv_i)^2.
$$
Write $J\mv_i=(a_i, b_i)^\top$. Then
$$
    G=\sum_i c_i(J\mv_i)(J\mv_i)^\top = \begin{pmatrix} \sum_i c_i a_i^2  & \sum_i c_i a_ib_i\\
    \sum_i c_i a_ib_i & \sum_i c_i b_i^2
\end{pmatrix}.
$$
Hence
$$
    \begin{aligned}
        \det G &= \left(\sum_i c_i a_i^2\right) \left(\sum_j c_j b_j^2\right) - \left(\sum_i c_i a_ib_i\right)^2\\
        &= \sum_{i<j}c_ic_j \left(a_i^2b_j^2+a_j^2b_i^2-2a_ib_i a_jb_j\right)\\
        &= \sum_{i<j}c_ic_j(a_ib_j-a_jb_i)^2\\
        &=\sum_{i<j}c_ic_j\det(\mv_i,\mv_j)^2.
    \end{aligned}
$$
\end{proof}

\begin{proof}[Proof of Theorem \ref{geom:noncollapse}]
The proof is by contradiction. The rod directions identify a fixed collapsing circle. After unwrapping this circle, its normalized length converges to a positive function satisfying the equation \eqref{geom:positiveeigenfunction}. The rod expansions and the volume estimate near the limiting vertices justify an integration by parts that contradicts this equation.

After scaling, we may assume $\lambda=3$. Choose a fixed integer $N$ sufficiently large and quotient by the subgroup $\T^2[N]\cong(\mathbb Z/N\mathbb Z)^2$. This divides all cone angles by $N$ and the volume by $N^2$, while leaving the quotient metric on the base unchanged. Relabeling the lattice, it suffices to prove the theorem when all cone angles are strictly less than $2\pi$, with a uniform positive lower bound. The local convexity argument in \cite[Lemma 4.43]{LiSun} then applies. Together with Bishop-Gromov comparison and the diameter bound \eqref{geom:diameterbounds}, this gives uniform covering bounds for the spaces $(X_{\mR},g)$ and their torus quotients.

Suppose the conclusion fails. There is then a sequence $g_j$ such that $\Vol(X_{\mR},g_j)\to0$. After passing to a subsequence, the bases converge in the Gromov-Hausdorff sense,
$$
    B_j\longrightarrow B_\infty,
$$
and the cone angles on each rod converge to positive limits. Let $\mathcal S$ be the limit of the vertex sets $\pi_j(\mathcal V_j)$, and let $\mathcal F$ be the limit of the union of rods. The estimates in Lemma \ref{geom:basecomparison}, Proposition \ref{geom:vertexcurvature}, and Lemma \ref{geom:rodhalfdisk} imply that $B_\infty\setminus\mathcal S$ is a smooth surface with geodesic boundary $\mathcal F\setminus\mathcal S$, and the convergence is locally smooth up to this boundary. The boundary diameter lower bound \eqref{eqn:diameter lower bound} implies that $\mathcal F\setminus\mathcal S\neq\varnothing$.

Choose an interior point $\mathfrak s_\infty$ of a surviving rod, with rod vector $\mv_1$, and let $Q$ be the connected component of $B_\infty\setminus\mathcal S$ containing $\mathfrak s_\infty$. Write $Q^\circ=Q\setminus\mathcal F$. Moving a fixed small distance normally inward from $\mathfrak s_\infty$, choose $\mathfrak q_\infty\in Q^\circ$ and corresponding points $\mathfrak q_j\in B_j$. The Killing field generated by $\mv_1$ vanishes on the rod, and its normal derivative there is controlled above and below by the cone angle bounds. The curvature bound and the Jacobi equation therefore give
$$
    0<c\leq G_j(\mv_1,\mv_1)(\mathfrak q_j)\leq C.
$$
Thus the $\mv_1$ direction does not collapse at $\mathfrak q_j$.

On the other hand, $\rho_j=\sqrt{\det G_j}$ satisfies $-\Delta_{\bg_j}\rho_j=6\rho_j$. On a fixed small neighborhood $U_j$ of $\mathfrak q_j$, the base metrics converge smoothly and the Harnack inequality gives
$$
    \rho_j(\mathfrak q_j)\leq C\int_{U_j}\rho_j\dd A_{\bg_j}\leq C\Vol(X_{\mR},g_j)\longrightarrow0.
$$
We now use Lemma \ref{geom:rodmatrixcone} to identify a fixed collapsing direction. Write
$$
    G_j=\sum_i c_{ij}A_i,\qquad c_{ij}\geq0,
$$
combining the terms corresponding to parallel rod vectors. At $\mathfrak q_j$, Equation \eqref{geom:detmatrices} gives
$$
    c_{ij}\det(\mv_1,\mv_i)^2\leq G_j(\mv_1,\mv_1)\leq C
$$
for every $\mv_i$ not parallel to $\mv_1$. These coefficients are therefore uniformly bounded. For the remaining coefficient, denoted by $c_{1j}$, the same equation gives
$$
    \det G_j\geq c_{1j}G_j(\mv_1,\mv_1),
$$
so $c_{1j}(\mathfrak q_j)\to0$. Consequently, after passing to a subsequence, $G_j(\mathfrak q_j)$ converges to a nonzero rank one matrix $G_*$ in the cone generated by the $A_i$. The determinant identity shows that two generators with nonparallel rod vectors cannot both have positive coefficients in $G_*$. Hence
$$
    G_*=aA_{i_0},\qquad a>0.
$$
Set $\mathfrak w=\mv_{i_0}$. Since $G_*(\mv_1,\mv_1)>0$, the vectors $\mathfrak w$ and $\mv_1$ are not parallel.

Complete $\mathfrak w$ to a fixed integral basis $(\mathfrak w,\mathfrak u)$ of $\Lambda$. Using the period-one angular parameters generated by this basis, write
\begin{equation}\label{geom:circledecomp}
    g_j=\bg_j+u_j^2(\dd\theta+a_j\dd\psi)^2+f_j^2\dd\psi^2,
\end{equation}
where
$$
    u_j^2=G_j(\mathfrak w,\mathfrak w),\qquad a_j=\frac{G_j(\mathfrak w,\mathfrak u)}{G_j(\mathfrak w,\mathfrak w)},\qquad f_j^2=G_j(\mathfrak u,\mathfrak u)-\frac{G_j(\mathfrak w,\mathfrak u)^2}{G_j(\mathfrak w,\mathfrak w)}.
$$
In particular, if $D=|\det(\mathfrak w,\mathfrak u)|$, then $u_jf_j=D\rho_j$. Here $D$ is fixed and only records the angular normalization.

Write $\mv_i=p_i\mathfrak w+q_i\mathfrak u$. The decomposition into the $A_i$ gives
$$
    a_j=-\frac{\sum_{q_i\neq0}c_{ij}p_iq_i}{\sum_{q_i\neq0}c_{ij}q_i^2},\qquad |a_j|\leq\max_{q_i\neq0}\left|\frac{p_i}{q_i}\right|.
$$
Set
$$
    \varepsilon_j=u_j(\mathfrak q_j)\longrightarrow0.
$$
Since $\mv_1=p_1\mathfrak w+q_1\mathfrak u$ with $q_1\neq0$, we have
$$
    G_j(\mv_1,\mv_1)(\mathfrak q_j)=\varepsilon_j^2(p_1+q_1a_j(\mathfrak q_j))^2+q_1^2f_j(\mathfrak q_j)^2.
$$
The first term tends to zero, so
$$
    C^{-1}\leq f_j(\mathfrak q_j)\leq C.
$$

We next propagate these estimates through $Q^\circ$. The interior estimates leading to \eqref{eqn:Klower}, together with the equation for $\rho_j$, give uniform bounds for the logarithmic variation of the quadratic forms $G_j$ on compact subsets of $Q^\circ$. Integrating along fixed paths from $\mathfrak q_\infty$, and using
$$
    u_j^2=G_j(\mathfrak w,\mathfrak w),\qquad f_j^2=\inf_{t\in\mathbb R}G_j(t\mathfrak w+\mathfrak u,t\mathfrak w+\mathfrak u),
$$
we obtain, for every compact subset $K\subset Q^\circ$,
$$
    C_K^{-1}\leq\varepsilon_j^{-1}u_j\leq C_K,\qquad C_K^{-1}\leq f_j\leq C_K.
$$
Thus the same fixed direction $\mathfrak w$ collapses throughout the interior, whereas the quotient direction remains noncollapsed.

Lift the $\theta$ coordinate locally to $\mathbb R$ and set $\theta=\varepsilon_j^{-1}\tau$. With $V_j=\varepsilon_j^{-1}u_j$, the developed metric is
$$
    \bg_j+V_j^2(\dd\tau+\varepsilon_ja_j\dd\psi)^2+f_j^2\dd\psi^2.
$$
The preceding bounds give local noncollapse on every compact subset of $Q^\circ$. The curvature estimates and Einstein regularity therefore give smooth local convergence after passing to a subsequence. Since $a_j$ is uniformly bounded, the limit takes the diagonal form
$$
    g_\infty=\bg+V^2\dd\tau^2+f^2\dd\psi^2,
$$
where $V,f>0$ are defined throughout $Q^\circ$. These developments are local, but the functions $V$ and $f$ agree on overlaps because they are limits of the globally defined functions $V_j$ and $f_j$.

The developed metrics still satisfy $\operatorname{Ric}=3g$. Equations \eqref{eq:rho} and \eqref{eq:map} therefore give
$$
    \Delta_{\bg}(Vf)=-6Vf,\qquad \operatorname{div}_{\bg}(V\nabla f-f\nabla V)=0.
$$
Combining these equations, we obtain
\begin{equation}\label{geom:positiveeigenfunction}
    -\operatorname{div}_{\bg}(f\nabla V)=3fV.
\end{equation}
We will integrate this equation against a cutoff times $V^{-1}$. For this purpose, we first show that $V$ remains positive along every surviving rod in $Q$, and then control the contribution near $\mathcal S$.

Fix an arbitrary compact segment of a surviving rod in $\partial Q$, with rod vector
$$
    \mv=p\mathfrak w+q\mathfrak u.
$$
We claim that $q\neq0$. Otherwise, since both vectors are primitive, $\mv=\pm\mathfrak w$. Choose points a fixed small normal distance $r_0>0$ from this rod. The Jacobi equation and the positive lower bound for the cone angles would give
$$
    u_j\geq cr_0
$$
at these points. But the interior estimates already proved give $u_j\leq C\varepsilon_j\to0$ there. This is a contradiction. Thus $\mathfrak w$ is not parallel to any surviving rod vector in $Q$.

\

We now prove that near this rod,
\begin{equation}\label{geom:limitingrod}
    V=V_0(s)+O(r^2),\qquad V_0(s)>0,\qquad f=\frac{2\pi\beta}{|q|}r+O(r^3),
\end{equation}
where $r$ is the normal distance to the rod and $2\pi\beta$ is its limiting cone angle.

Fix a compact segment of a surviving rod, with rod vector $\mv=p\mathfrak w+q\mathfrak u$, where $q\neq0$. Normalize its normal cone angle using $\beta_j^{-1}\mv$ and unwrap the $\mathfrak w$ direction as above, with $\theta=\varepsilon_j^{-1}\tau$. The curvature bound and Lemma \ref{geom:rodhalfdisk} give a developed collar of uniform size. The interior noncollapse and volume comparison give noncollapse up to the rod, so Einstein regularity gives smooth convergence across the fixed surface.

Let $K_{\mathfrak w,j}$ and $K_{\mv,j}$ denote the corresponding Killing fields. Their interior convergence and the Killing equation give smooth convergence across the fixed surface,
$$
    \varepsilon_j^{-1}K_{\mathfrak w,j}\longrightarrow K_{\mathfrak w,\infty},\qquad K_{\mv,j}\longrightarrow K_{\mv,\infty}.
$$
On the interior, $K_{\mathfrak w,\infty}=\partial_\tau$ and $K_{\mv,\infty}=q\partial_\psi$. These fields commute, and the diagonal limit metric makes them orthogonal.

We claim that $K_{\mathfrak w,\infty}$ does not vanish on the fixed surface. Otherwise, at a zero, commutation with the rod rotation implies that its derivative preserves the normal plane and acts there as a rotation. Orthogonality to $K_{\mv,\infty}$ forces this normal rotation to vanish. Its flow would then fix the inward normal geodesic pointwise, contradicting $|K_{\mathfrak w,\infty}|=V>0$ in the interior.

Thus $V=|K_{\mathfrak w,\infty}|$ is smooth, positive, and rotationally invariant near the fixed surface. Consequently,
\begin{equation}\label{Vasymptotics}
    V=V_0(s)+O(r^2),\qquad V_0(s)>0,\qquad \partial_rV=O(r).
\end{equation}
The rod rotation satisfies $|K_{\mv,\infty}|=2\pi\beta r+O(r^3)$, where $2\pi\beta$ is the limiting cone angle. Since $|K_{\mv,\infty}|=|q|f$, we also obtain
\begin{equation}\label{fasymptotics}
    f=\frac{2\pi\beta}{|q|}r+O(r^3).
\end{equation}
These expansions hold on every compact segment of every surviving rod in $Q$.

It remains to control the neighborhoods of $\mathcal S$. On $B_j$, set
$$
    \nu_j=f_j\dd A_{\bg_j},\qquad \mathfrak r_j=d_{\bg_j}(\cdot,\pi_j(\mathcal V_j)).
$$
We claim that
\begin{equation}\label{geom:capvolume}
    \nu_j\{\mathfrak r_j<r\}\leq Cr^3,\qquad \nu_j(B_j)\leq C.
\end{equation}
Consider the metric quotient $Y_j=X_j/S^1_{\mathfrak w}$. On its regular part, its metric is
$$
    \widehat g_j=\bg_j+f_j^2\dd\psi^2.
$$
Since $\psi$ has period one, the volume of the inverse image of a measurable subset of $B_j$ equals its $\nu_j$ measure. Let $P_j\subset Y_j$ be the images of the torus fixed points. The orbit over each vertex is a single point, so the distance to $P_j$ in $Y_j$ equals the distance to $\pi_j(\mathcal V_j)$ in $B_j$.

The spaces $Y_j$ inherit uniform metric doubling bounds from $X_j$. Since the number of vertices is fixed, each annulus
$$
    \{y\in Y_j:r/2\leq d(y,P_j)<r\}
$$
can therefore be covered by a uniformly bounded number of balls of radius $r/16$. On slightly enlarged balls, the curvature estimate \eqref{geom:vertexcurvatureeq} and O'Neill's formula give a sectional curvature lower bound $-Cr^{-2}$, with the uniform diameter bound absorbed into $C$. Away from $P_j$, a rod parallel to $\mathfrak w$ gives a smooth even boundary of $Y_j$, while a rod $\mv=p\mathfrak w+q\mathfrak u$ with $q\neq0$ gives a conical curve of angle $2\pi\beta_j/|q|<2\pi$. Reflection across the smooth boundary and local volume comparison across these conical curves show that each covering ball has volume at most $Cr^3$. Consequently, the annulus has volume at most $Cr^3$. Summing over dyadic annuli proves the first estimate in \eqref{geom:capvolume}. Taking $r$ larger than the uniformly bounded diameter gives the second.

Let $\mathcal N_r$ be the $r$-neighborhood of $\mathcal S$ intersected with $Q$. Smooth convergence on compact subsets of $Q^\circ$, followed by exhaustion, now gives
$$
    \nu=f\dd A_{\bg},\qquad 0<\nu(Q)<\infty,\qquad \nu(\mathcal N_r)\leq Cr^3.
$$
Here positivity follows from $f>0$ on the nonempty open set $Q^\circ$.

We can now integrate \eqref{geom:positiveeigenfunction}. Choose a Lipschitz cutoff $\chi_r$ which vanishes on $\mathcal N_r$, equals one on $Q\setminus\mathcal N_{2r}$, and satisfies $|\dd\chi_r|\leq Cr^{-1}$. Along each surviving rod, \eqref{Vasymptotics} and \eqref{fasymptotics} give
$$
    \frac{f}{V}\partial_rV=O(r^2).
$$
Thus there is no boundary contribution from the rods when we multiply \eqref{geom:positiveeigenfunction} by $\chi_r^2/V$ and integrate by parts. The cutoff removes the limiting vertices. We obtain
\begin{equation}\label{integration by parts}
    \begin{aligned}
        3\int_Q f\chi_r^2\dd A_{\bg} &=\int_Q f\left\langle\dd V,\dd\left(\frac{\chi_r^2}{V}\right)\right\rangle\dd A_{\bg}\\
        &=\int_Q f\left(|\dd\chi_r|^2-\left|\dd\chi_r-\frac{\chi_r}{V}\dd V\right|^2\right)\dd A_{\bg}\\
        &\leq\int_Q f|\dd\chi_r|^2\dd A_{\bg}\leq Cr^{-2}\nu(\mathcal N_{2r})\leq Cr.
    \end{aligned}
\end{equation}
As $r\to0$, the right-hand side tends to zero, whereas the left-hand side tends to $3\nu(Q)>0$. This contradiction proves the volume lower bound for $\lambda=3$. Rescaling gives the stated bound.
\end{proof}

\subsection{Hitchin-Thorpe inequality}
For compact Einstein four-manifolds, a well-known topological restriction is given by the Hitchin-Thorpe inequality. This was generalized in \cite[Theorem~A]{AtiyahLeBrun}. to Einstein four-manifolds with a cone singularity along a smooth surface in a suitable sense. In our toric conical setting, we will prove a version below.

Let $(X_{\mR},g)$ be a compact toric conical Einstein four-manifold with $\Ric_g=\lambda g$. Let $\Sigma_1,\ldots,\Sigma_m$ be associated 2-spheres along which $g$ has cone angle $2\pi\beta_i$, and set $e_i=[\Sigma_i]^2$.

Define the modified Euler characteristic and signature by
\begin{align}
    \chi_c(g) &=\chi(X)-2\sum_{i=1}^m(1-\beta_i) +\sum_{i=1}^m(1-\beta_i)(1-\beta_{i+1}) =\sum_{i=1}^m\beta_i\beta_{i+1}, \label{eq:chic}\\
    \tau_c(g) &=\tau(X)-\frac13\sum_{i=1}^m(1-\beta_i^2)e_i =\frac13\sum_{i=1}^m e_i\beta_i^2. \label{eq:tauc}
\end{align}
Set
$$
    I_\pm(g)=\int_{X^\circ}|W_g^\pm|^2 \mathrm \dvol_g.
$$

\begin{proposition}
\label{geom:hitchinthorpe}
We have
\begin{equation}\label{eq:ht}
    I_\pm(g)=2\pi^2\bigl(2\chi_c(g)\pm3\tau_c(g)\bigr)-\frac{\lambda^2}{3}\operatorname{Vol}_g(X).
\end{equation}
In particular, if $\lambda\neq0$, then
\begin{equation}\label{geom:angleinequality}
    2\sum_{i=1}^m\beta_i\beta_{i+1}>\left|\sum_{i=1}^m e_i\beta_i^2\right|.
\end{equation}
\end{proposition}

\begin{proof}
Let $\mathcal E(g)$ and $p_1(g)$ denote the normalized Euler and first Pontryagin forms, and set $\delta_i=1-\beta_i$. We compute their integrals by smoothing the conical singularities and applying the ordinary Chern-Gauss-Bonnet and signature formulas.

We first arrange that the metric is exactly a product of flat cones near each torus fixed point. In equivariant normal coordinates in a local development, the metric differs from its Euclidean tangent metric by $O(r^2)$, with the corresponding derivative bounds. A cutoff at scale $\eta$ therefore gives a conical metric $h_\eta$ which agrees with the flat product model within distance $\eta$ of each torus fixed point and agrees with $g$ outside distance $2\eta$. These metrics have the same cone angles and uniformly bounded curvature on the regular locus. Since the modified regions have volume $O(\eta^4)$, the regular Euler and Pontryagin integrals of $h_\eta$ converge to those of $g$ as $\eta\to0$.

Fix $\eta$ and write $h=h_\eta$. Smooth the normal cones at scale $\varepsilon$, using product smoothings near the torus fixed points, to obtain smooth metrics $h_\varepsilon$ on $X$. These agree with $h$ outside shrinking neighborhoods of the singular surfaces. The ordinary characteristic number formulas give
$$
    \int_X\mathcal E(h_\varepsilon)=\chi(X),\qquad \int_Xp_1(h_\varepsilon)=3\tau(X).
$$
As $\varepsilon\to0$, the contributions away from the singular surfaces converge to the corresponding integrals over $X^\circ$. We now compute the additional contributions concentrated along the surface interiors and at their intersections.

Consider first the interior of $\Sigma_i$. Smoothing a normal cone of angle $2\pi\beta_i$ gives a smooth disk $N_\varepsilon$ which agrees with the cone near its boundary. The boundary has total geodesic curvature $2\pi\beta_i$, so the two-dimensional Gauss-Bonnet formula gives
$$
    \frac1{2\pi}\int_{N_\varepsilon}K_N\dd A_N =1-\frac1{2\pi}\int_{\partial N_\varepsilon}k_N\dd s =1-\beta_i=\delta_i.
$$
For a product of this disk with a surface patch $U\subset\Sigma_i^\circ$, the Euler form is the product of the two normalized curvature forms. Hence its integral is
$$
    \left(\frac1{2\pi}\int_{N_\varepsilon}K_N\dd A_N\right) \left(\frac1{2\pi}\int_UK_{\Sigma_i,h}\dd A_h\right) =\frac{\delta_i}{2\pi}\int_UK_{\Sigma_i,h}\dd A_h.
$$
The local edge-cone calculation in \cite[Section~2, proof of Theorem~2.1]{AtiyahLeBrun} gives the same limiting contribution for $h$. Here we first work away from the torus fixed points and then exhaust $\Sigma_i^\circ$.

The induced metric on $\Sigma_i$ has intrinsic cone angles $2\pi\beta_{i-1}$ and $2\pi\beta_{i+1}$ at its two poles, since these are the intersections with the neighboring singular surfaces. Thus
$$
    \frac1{2\pi}\int_{\Sigma_i^\circ}K_{\Sigma_i,h}\dd A_h =2-\delta_{i-1}-\delta_{i+1}.
$$
The Euler contribution concentrated along the interior of $\Sigma_i$ is therefore
$$
    \delta_i(2-\delta_{i-1}-\delta_{i+1}).
$$

At the torus fixed point $\Sigma_i\cap\Sigma_{i+1}$, the metric $h$ is an exact product of two flat cones. The product smoothing has normalized curvature integrals $\delta_i$ and $\delta_{i+1}$ in the two factors. Its Euler integral is consequently $\delta_i\delta_{i+1}$. Combining the regular integral with the surface and fixed-point contributions gives
$$
    \chi(X)=\int_{X^\circ}\mathcal E(h)+\sum_i\delta_i(2-\delta_{i-1}-\delta_{i+1})+\sum_i\delta_i\delta_{i+1}.
$$
Letting $\eta\to0$ and using $\chi(X)=m$, we obtain
$$
    \begin{aligned}
        \int_{X^\circ}\mathcal E(g) &=\chi(X)-\sum_i\delta_i(2-\delta_{i-1}-\delta_{i+1})-\sum_i\delta_i\delta_{i+1}\\
        &=\chi(X)-2\sum_i\delta_i+\sum_i\delta_i\delta_{i+1}\\
        &=\sum_i\beta_i\beta_{i+1} =\chi_c(g).
    \end{aligned}
$$

For the first Pontryagin form, the local normal-bundle calculation in \cite[Theorem~2.2 and its proof]{AtiyahLeBrun} gives the contribution $(1-\beta_i^2)e_i$ along $\Sigma_i$. Here $e_i=[\Sigma_i]^2$ is the Euler number of its oriented normal bundle. There is no additional fixed-point contribution: for the product smoothing, the first Pontryagin form is the sum of the squares of the two normalized curvature forms, and each square vanishes because it is pulled back from a two-dimensional factor. Taking the smoothing limit and then letting $\eta\to0$ therefore gives
$$
    3\tau(X)=\int_{X^\circ}p_1(g)+\sum_i(1-\beta_i^2)e_i.
$$
Using the topological identity $3\tau(X)=\sum_i e_i$, we conclude that
$$
    \int_{X^\circ}p_1(g)=3\tau(X)-\sum_i(1-\beta_i^2)e_i=\sum_i e_i\beta_i^2=3\tau_c(g).
$$

Finally, since $\Ric_g=\lambda g$, the scalar curvature is $4\lambda$. From  the Chern-Weil formulas we see
$$
    \begin{aligned}
        8\pi^2\chi_c(g)&=I_+(g)+I_-(g)+\frac{2\lambda^2}{3}\operatorname{Vol}_g(X),\\
        12\pi^2\tau_c(g)&=I_+(g)-I_-(g).
    \end{aligned}
$$
Taking half the sum and difference yields
$$
    I_\pm(g)=2\pi^2\bigl(2\chi_c(g)\pm3\tau_c(g)\bigr)-\frac{\lambda^2}{3}\operatorname{Vol}_g(X),
$$
which proves \eqref{eq:ht}. Since $I_\pm(g)\geq0$, we have
$$
    2\chi_c(g)\pm3\tau_c(g)\geq\frac{\lambda^2}{6\pi^2}\operatorname{Vol}_g(X).
$$
For $\lambda\neq0$, the right-hand side is positive. Substituting the expressions for $\chi_c(g)$ and $\tau_c(g)$ gives \eqref{geom:angleinequality}.
\end{proof}

\section{Local structure of the moduli space }
\label{sec:moduli}
From now on in this paper we work with the setting of Theorem \ref{thm:dichotomy}. Let $\mathfrak R_0$ be a compact rod structure with $m_0$ rods and a chosen turning point $\mathfrak p$. After a global integral change of lattice basis, possibly reversing its orientation, assume that the incident rod directions are $B,A$ in boundary order and all rod directions belong to $\{A,B,C\}$, with $A,B,C$ as defined in the introduction.

Replace the corner $(B,A)$ by one of the following sequences of rod directions:
\begin{equation}\label{eq:alternating-compact-word}
    W\in\left\{
    \begin{gathered}
        (B,(A,B)^n,A)\quad(n\geq1),\qquad (B,C,A),\\
        (B,A,C,A),\qquad (B,C,B,A),\\
        (B,C,A,C,A),\qquad (B,C,B,C,A),\\
        (B,C,A,B,C,A)
    \end{gathered}
    \right\}.
\end{equation}
Set
$$
    \mathfrak R=\mathfrak R_0\#_{\mathfrak p}W,\qquad M=X_{\mathfrak R_0},\qquad X=X_{\mathfrak R}.
$$
Number the rods of $\mathfrak R$ in positive cyclic order so that $\mathfrak I_1,\ldots,\mathfrak I_q$  correspond to the intermediate directions of $W$. There are $m=m_0+q$ rods altogether. As in Section \ref{sec:rod structure}, write
$$
    \pi:X\longrightarrow B,\qquad \Sigma_i=\pi^{-1}(\mathfrak I_i),\qquad e_i=[\Sigma_i]^2,
$$
with $\mathfrak p_i\in B$ and $p_i=\pi^{-1}(\mathfrak p_i)$ denoting the turning points and the corresponding torus fixed points, respectively. Set
$$
    \mathcal S=\bigcup_{\alpha=1}^q\Sigma_\alpha, \qquad \sigma=\sum_{\alpha=1}^q e_\alpha.
$$
The orientations of $M$ and $X$ are determined by the oriented quotient and the normalized angular coordinates, as in Section~\ref{sec:rod structure}. All connected-sum identifications and self-intersection numbers below use these orientations. If the lattice normalization reverses the original lattice orientation, then these four-manifold orientations are opposite to the original ones. To express the connected-sum identifications in the original orientation of $M$, reverse the orientation of each listed summand. The self-intersection numbers and signatures then change sign.

\subsection{Topology of the toric manifolds}
It is easy to check that with the above choice of rod structures, we get exactly the smooth toric four-manifolds listed in Theorem \ref{thm:main}.

Indeed, the connected-sum identifications follow directly from the neighborhoods of the corresponding spheres. For the $(AB)$ alternating sequence, two consecutive spheres have self-intersection zero and meet transversely once. Their regular neighborhood is $(\mathbb S^2\times \mathbb S^2)\setminus B^4$, as seen from the two coordinate spheres in $\mathbb S^2\times \mathbb S^2$. Replacing this neighborhood by a four-ball removes one pair $(A,B)$ from the sequence. Repeating this operation gives
$$
    X\cong M\#n(\mathbb S^2\times \mathbb S^2).
$$
For each remaining sequence, the self-intersection formula gives successive deletions of rods with self-intersection $\pm1$, reducing the sequence to $(B,A)$. Each $+1$ deletion contributes a $\CP^2$ summand, and each $-1$ deletion contributes an $\overline{\CP^2}$ summand. Counting these deletions gives the connected-sum identifications below.

$$
    \begin{array}{c|c}
        W&\text{connected sum by}\\
        \hline
        (B,C,A)&\overline{\CP^2}\\
        (B,A,C,A),\ (B,C,B,A)&\CP^2\#\overline{\CP^2}\\
        (B,C,A,C,A),\ (B,C,B,C,A)&\CP^2\#2\overline{\CP^2}\\
        (B,C,A,B,C,A)&\CP^2\#3\overline{\CP^2}
    \end{array}
$$
The two sequences grouped in each row give the same oriented connected sum summand, but their insertions into a fixed rod structure may yield inequivalent $\T^2$-actions. For the four-manifolds $\#r(\mathbb S^2\times \mathbb S^2)$ and $(\CP^2\# k\overline{\CP^2})\#r(\mathbb S^2\times \mathbb S^2)$ with $0\leq k\leq 3$ in Theorem \ref{thm:main}, we can obtain them by starting with $\mathbb S^4$ or $\CP^2\# k\overline{\CP^2}$ with $0\leq k\leq 3$, then insert with $W=(B,(A,B)^r,A)$. For the third situation in Theorem \ref{thm:main}, we may begin with $\mathbb S^4$, $\mathbb S^2\times \mathbb S^2$, or $\CP^2\#k\overline{\CP^2}$ with either of its orientaion, where $0\leq k\leq 3$, then insert the rest of sequences.

For all the sequences we get
\begin{equation}\label{eq:insertion-self-intersections}
    \begin{array}{c|c|c|c}
        W&q&(e_1,\ldots,e_q)&\sigma\\
        \hline
        (B,(A,B)^n,A)&2n&(0,\ldots,0)&0\\
        (B,C,A)&1&(-1)&-1\\
        (B,A,C,A)&2&(1,0)&1\\
        (B,C,B,A)&2&(0,1)&1\\
        (B,C,A,C,A)&3&(-1,0,0)&-1\\
        (B,C,B,C,A)&3&(0,0,-1)&-1\\
        (B,C,A,B,C,A)&4&(-1,-1,-1,-1)&-4
    \end{array}
\end{equation}
In particular, for $[(\beta, g)]\in \mathcal M$, we have
\begin{equation}\label{eq:alternating-characteristic}
    \chi_c(g)=(q-1)\beta^2+2\beta+m_0-1,\qquad \tau_c(g)=\tau(X)+\frac{\beta^2-1}{3}\sigma.
\end{equation}

\begin{remark}\label{rmk:q=1}
For $W=(B,C,A)$, we have $q=1$ and $\sigma=-1$, so the characteristic formulas and the Hitchin-Thorpe inequality give an upper bound on $\beta$. For the remaining sequences, these inequalities alone do not exclude $\beta\to\infty$. The upper bound is more difficult to obtain and will be proved in Section~\ref{sec:angle-bounds}.
\end{remark}

\subsection{Definition of the moduli space}
 We consider toric conical Einstein metrics with
$$
    \Ric(g)=3g,\qquad \beta_i= \begin{cases} \beta, &1\le i\le q,\\
    1,     &q<i\le m.
\end{cases}
$$

Let $\mathcal Z$ be the space of pairs $(\beta,g)$, where $\beta>0$ and $g$ is such a metric satisfying $\Ric(g)=3g$. We equip $\mathcal Z$ with the topology of convergence of cone angles and $C^\infty$ convergence in local developments. Recording $\beta$ and the developed tensor coefficients in a fixed finite collection of charts realizes $\mathcal Z$ as a subspace of the product of $(0,\infty)$ with finitely many separable $C^\infty$ spaces. In particular, $\mathcal Z$ is metrizable and second countable.

We call an orientation-preserving diffeomorphism $F:X\to X$ \emph{admissible} if it commutes with the fixed $\mathbb T^2$-action and preserves each rod sphere individually. We say two pairs $(\beta,g)$ and $(\beta',g')$ are \emph{equivalent} if $\beta=\beta'$ and $g'=F^*g$ for an admissible diffeomorphism $F$. We denote the resulting quotient of $\mathcal Z$ by $\mathcal M$, equip it with the quotient topology, and write $\varpi:\mathcal Z\to\mathcal M$ for the quotient projection.

The stabilizer of every pair in $\mathcal Z$ is precisely the group of torus translations. Indeed, an equivariant orientation-preserving isometry induces an orientation-preserving isometry of the orbit space preserving each inserted rod, so it fixes the two turning points bounding such a rod individually. Its restriction to that rod is therefore the identity. Preservation of inward normal geodesics then implies that it is the identity near the interior of the rod, and hence on the whole orbit space. On the principal orbit set, an equivariant map covering the identity has the form $x\mapsto a(\pi(x))\cdot x$, where $a$ maps the interior of the orbit space to $\mathbb T^2$. Since an isometry preserves the orthogonal complements of the torus orbits, $\dd a=0$. Thus $a$ is constant, and the isometry is a torus translation. In particular, the stabilizer acts trivially on invariant tensors.

We next show that the equivalence relation is closed. Suppose that
$$
    (\beta_j,g_j)\longrightarrow(\beta,g), \qquad (\beta_j,F_j^*g_j)\longrightarrow(\beta,h),
$$
where each $F_j$ is an admissible diffeomorphism. Smooth convergence of the metrics in local developments gives uniform bounds for the first derivatives of $F_j$ and $F_j^{-1}$, and the isometry equation gives bounds for all higher derivatives. Compactness of $X$ therefore gives a subsequence for which the maps and their inverses converge smoothly in local developments to a diffeomorphism $F$ and its inverse. The limit preserves orientation, commutes with the torus action, and preserves each $\Sigma_\alpha$, so it is an admissible diffeomorphism. Passing to the limit gives $h=F^*g$, proving closedness.

The projection $\varpi$ is open. Indeed, for every open set $U\subset\mathcal Z$,
$$
    \varpi^{-1}(\varpi(U))=\bigcup_F F^*U,
$$
and the right-hand side is open because each pullback map is a homeomorphism. Closedness of the equivalence relation now implies that $\mathcal M$ is Hausdorff: inequivalent pairs have neighborhoods containing no equivalent pair, whose images under $\varpi$ are disjoint open neighborhoods. Finally, since $\mathcal Z$ is second countable and $\varpi$ is open, the images of a countable basis for $\mathcal Z$ form a countable basis for $\mathcal M$.

\subsection{The linearized equation}
We first specify the analytic conventions used below. For fixed cone angles, an invariant tensor is called \emph{smooth} if it is smooth on the regular locus and extends smoothly in each local angular development, including the simultaneous developments at the torus fixed points. All tensor spaces below consist of $\T^2$-invariant tensors with this regularity. These conventions apply to arbitrary compact toric conical metrics, including the smooth case.

In this paper we will use the following divergence gauge for the Einstein
equation. For related formulations, see
\cite[Remark~4.6, equation~(4.25)]{AndersonSymmetries}
and \cite[Section~2.10]{BiquardDGGA}.

\begin{definition}[Divergence gauge]
\label{def:divergence-gauge}
We fix a reference metric $g$. We say that a nearby metric $h$ is in \emph{divergence gauge relative to $g$} if
$$
 \delta_hg=0.
$$
\end{definition}

Denote the Einstein tensor by
$$
    \mathscr E_\lambda(g) =\Ric(g)-\tfrac12\Scal_g g+\lambda g,
$$
and define its linearization in divergence gauge by
$$
    \mathscr P_g =D\mathscr E_{\lambda,g}+\delta_g^*\delta_g,
$$
where $\delta_g=-\diver_g$ and $\delta_g^*$ is its formal $L^2$ adjoint. In dimension four, $\mathscr E_\lambda(g)=0$ is equivalent to $\Ric(g)=\lambda g$.
For the metrics in $\mathcal M$, we take $\lambda=3$.
The operator $\mathscr P_g$ is elliptic, and its coefficients are smooth.

Fix $(\beta_0,g)$ representing a point of $\mathcal M$. For each integer $s\geq0$, let $L^2_s$ be the completion of smooth $\T^2$-invariant symmetric two-tensors in the norm
$$
    \|k\|_{L^2_s}^2 =\sum_{j=0}^s\|\nabla_g^j k\|_{L^2}^2, \qquad L^2=L^2_0.
$$
Unless otherwise stated, all tensor spaces and orthogonal complements below are defined using this fixed metric $g$. Since $g$ is Einstein, we have
$$
    2\mathscr P_gk =(\nabla_g^*\nabla_g-2\mathring R_g)k -\nabla_g^2\tr_gk -\bigl(\nabla_g^*\nabla_g\tr_gk+\delta_g\delta_gk-3\tr_gk\bigr)g.
$$
On transverse-traceless tensors,  we have
$$
    \mathscr P_g=\tfrac12 (\nabla_g^*\nabla_g-2\mathring R_g),
$$

\begin{lemma}\label{lem:conical-linear}
The operator
$$
    \mathscr P_g:L^2_2\subset L^2\longrightarrow L^2
$$
is self-adjoint and Fredholm of index zero. Its kernel is
\begin{equation}\label{eq:conical-einstein-kernel}
    \cK_g= \left\{ k\in L^2_2: \delta_gk=0,\quad \operatorname{tr}_gk=0,\quad (\nabla_g^*\nabla_g-2\mathring R_g)k=0 \right\},
\end{equation}
and every element of $\cK_g$ is smooth. For every integer $s\geq0$, the map $\mathscr P_g:L^2_{s+2}\to L^2_s$ is Fredholm of index zero with the same kernel.
\end{lemma}

\begin{proof}
At an Einstein metric, $D\mathscr E_{3,g}$ is formally self-adjoint, and hence so is $\mathscr P_g=D\mathscr E_{3,g}+\delta_g^*\delta_g$. Let $k,h$ be smooth invariant symmetric two-tensors. We integrate by parts on the complement of shrinking neighborhoods of the fixed surfaces and torus fixed points. Smoothness in local developments implies that $k,h$ and their first covariant derivatives remain bounded, so the boundary terms tend to zero as these neighborhoods shrink. Therefore
$$
    \int_X\langle\mathscr P_gk,h\rangle_g\dvol_g=\int_X\langle k,\mathscr P_gh\rangle_g\dvol_g.
$$
By density, this identity also holds for $k,h\in L^2_2$, proving that $\mathscr P_g$ is symmetric on this domain.

To prove self-adjointness, we must identify the domain of the adjoint $\mathscr P_g^*$ of the unbounded operator $\mathscr P_g:L^2_2\subset L^2\to L^2$. Let $k\in\operatorname{Dom}(\mathscr P_g^*)$. By definition, there exists $f=\mathscr P_g^*k\in L^2$ such that
$$
    \langle k,\mathscr P_gh\rangle_{L^2}=\langle f,h\rangle_{L^2} \qquad\text{for every }h\in L^2_2.
$$
At this stage, $k$ is only known to belong to $L^2$. We will use this identity to prove that $k\in L^2_2$. Lift $k$ and $f$ to an angularly invariant local development. Given an arbitrary compactly supported smooth test tensor $\varphi$ there, its angular average is smooth and invariant, descends to the corresponding conical neighborhood, and extends by zero to a global smooth invariant tensor. It can therefore be used as $h$ in the adjoint identity. Since $\mathscr P_g$ commutes with angular averaging and both $k$ and $f$ are invariant, the same identity holds for the arbitrary test tensor $\varphi$. Formal self-adjointness of the differential expression then means precisely that the lifted equation $\mathscr P_gk=f$ holds in the distributional sense throughout the development, including across the axes. Since $f\in L^2$, interior elliptic regularity gives $k\in L^2_{2,\mathrm{loc}}$ in each development. A finite covering and invariant smooth approximation yield $k\in L^2_2$. Thus $\operatorname{Dom}(\mathscr P_g^*)\subset L^2_2$, while the reverse inclusion follows from symmetry proved above. Hence $\operatorname{Dom}(\mathscr P_g^*)=L^2_2$ and $\mathscr P_g^*=\mathscr P_g$ on this domain.

The local elliptic estimates give
$$
    \|k\|_{L^2_2} \leq C_g\bigl(\|\mathscr P_gk\|_{L^2} +\|k\|_{L^2}\bigr).
$$
Together with the compact inclusion $L^2_2\hookrightarrow L^2$, this proves the Fredholm property. Self-adjointness gives index zero, and elliptic regularity gives the corresponding statements at higher Sobolev orders.

The linearized contracted Bianchi identity gives
\begin{equation}\label{eq:linearized-divergence}
    \delta_g\mathscr P_gk=\delta_g\delta_g^*\delta_gk.
\end{equation}
For $k\in\ker\mathscr P_g$, integration by parts first gives $\delta_g^*\delta_gk=0$ and then $\delta_gk=0$. For a divergence-free tensor, the trace equation is
\begin{equation}\label{eq:linearized-trace}
    \operatorname{tr}_g(\mathscr P_gk) =-(\nabla_g^*\nabla_g-3)\operatorname{tr}_gk.
\end{equation}
The integrated Bochner identity is valid for invariant functions, again by integration by parts in local developments. Since $\Ric(g)=3g$ and $\dim X=4$, the standard Lichnerowicz eigenvalue estimate gives $\lambda_1\geq4$, where $\lambda_1$ denotes the first positive eigenvalue of the scalar Laplacian $\nabla_g^*\nabla_g$ on invariant functions. Consequently $\nabla_g^*\nabla_g-3$ is invertible on invariant functions, and $\operatorname{tr}_gk=0$.
\end{proof}

The ellipticity statement also applies to non-Einstein reference metrics. On a compact manifold, $\mathscr P_g$ differs from its formal adjoint by an operator of order zero. Its symmetrization therefore has the same elliptic principal symbol and is self-adjoint on $L^2_2$ domain. Since the difference is compact as a map $L^2_2\to L^2$, the Fredholm index of $\mathscr P_g$ is also zero.

\subsection{The local structure of the moduli space}

We first identify the space of tensors when $\beta$ varies. Choose $\mathbb T^2$-invariant rod and corner charts and a subordinate partition of unity $\{\chi_a\}$, smooth in local developments and independent of $\beta$. Recall that $\mathcal S=\bigcup_{\alpha=1}^q\Sigma_\alpha$ is the union of the inserted rod spheres. In each normal plane to a component of $\mathcal S$, define the operator $J_{\beta,a}$ by keeping the tensor coefficients and replacing
$$
    \beta_0 r \dd\theta\longmapsto\beta r \dd\theta
$$
in both covariant factors, while keeping the tangential part unchanged. At a fixed point, make the replacement in each normal plane corresponding to a component of $\mathcal S$. Thus both angles vary when both incident spheres belong to $\mathcal S$, and only one varies otherwise. On charts away from $\mathcal S$ we do not make any change. Now set
$$
    J_\beta=\sum_a\chi_aJ_{\beta,a}, \qquad b_\beta=J_\beta g.
$$
For $\beta$ sufficiently close to $\beta_0$, $J_\beta$ is an isomorphism between the corresponding tensor spaces and $b_\beta$ is a conical metric with cone angle $\beta$ along $\mathcal S$. Both depend real-analytically on $\beta$, with $J_{\beta_0}=\mathrm{Id}$ and $b_{\beta_0}=g$. They may be chosen to agree with the identity and $g$, respectively, outside a fixed neighborhood of $\mathcal S$.

Throughout this subsection, all Sobolev spaces $L^2_s$ are defined using the fixed reference metric $g$, with cone angle $\beta_0$ along $\mathcal S$. In particular, the domain space $L^2_8$ and the target space $L^2_6$ below are independent of $\beta$; the maps $J_\beta$ and $J_\beta^{-1}$ transport tensors between the reference spaces and the corresponding spaces for cone angle $\beta$. Write a nearby metric with cone angle $\beta$ as
$$
    h=b_\beta+J_\beta k, \qquad k\in L^2_8.
$$

We define the map with fixed domain and target by
\begin{equation}\label{eq:fixedtarget}
 \Psi(\beta,k)
 =J_\beta^{-1}
 \{\mathscr E_3(h)-\delta_h^*\delta_hb_\beta\}.
\end{equation}
It is real-analytic from a neighborhood of $(\beta_0,0)$ in $\mathbb R\times L^2_8$ to $L^2_6$, and
$$
    D_k\Psi(\beta_0,0)=\mathscr P_g.
$$

\begin{proposition}\label{prop:localmoduli}
A neighborhood of $[(\beta_0,g)]$ in $\mathcal M$ is homeomorphic to the zero set of a real-analytic map
$$
    \kappa_g:U_g\subset\cK_g\times\mathbb R \longrightarrow\cK_g,
$$
where $U_g$ is a neighborhood of $(0,0)$. Under this homeomorphism, a point $(\xi,t)\in\kappa_g^{-1}(0)$ represents a metric with common cone angle $2\pi(\beta_0+t)$ along $\mathcal S$.
\end{proposition}

\begin{proof}
At a zero of $\Psi$, set $\omega=\delta_hb_\beta$. The contracted Bianchi identity gives $\delta_h\delta_h^*\omega=0$. Integration by parts yields $\delta_h^*\omega=0$ and hence
$$
    \|\omega\|_{L^2(h)}^2 =\langle b_\beta,\delta_h^*\omega\rangle_{L^2(h)}=0.
$$
Thus $\Psi(\beta,k)=0$ is equivalent to $\Ric(h)=3h$ together with $\delta_hb_\beta=0$.

To set a nearby invariant metric $h$ in this gauge, solve
$$
    \delta_h(f^*b_\beta)=0
$$
for an equivariant diffeomorphism $f$ near the identity. At $(\beta,h,f)=(\beta_0,g,\mathrm{Id})$, the derivative in $f$, after identifying vector fields and one-forms using $g$, is $2\delta_g\delta_g^*$. Its kernel consists of the invariant Killing fields, namely the infinitesimal torus translations. The equation is automatically orthogonal to the torus Killing fields for $h$. Restricting to their orthogonal complement and applying the implicit function theorem therefore gives $f$, locally unique modulo torus translations. The metric $(f^{-1})^*h$ satisfies the required gauge condition.

Let $\Pi_g$ denote the $L^2$-orthogonal projection onto $\cK_g$. By Lemma~\ref{lem:conical-linear},
$$
    \mathscr P_g:L^2_8\cap\cK_g^\perp \longrightarrow L^2_6\cap\cK_g^\perp
$$
is an isomorphism. Writing $k=\xi+v$, with $\xi\in\cK_g$ and $v\perp\cK_g$, the implicit function theorem solves
$$
    (1-\Pi_g)\Psi(\beta_0+t,\xi+v)=0
$$
uniquely as $v=v(\xi,t)$. Set
$$
    \kappa_g(\xi,t) =\Pi_g\Psi(\beta_0+t,\xi+v(\xi,t)).
$$
The remaining equation is precisely $\kappa_g=0$. As observed above, the stabilizer consists of torus translations, which act trivially on invariant tensors. Thus the local moduli space is the zero set of $\kappa_g$ itself, with no further quotient. Elliptic regularity upgrades the solutions and their convergence to $C^\infty$.
\end{proof}

Proposition~\ref{prop:localmoduli} identifies neighborhoods in $\mathcal M$ with finite-dimensional real analytic zero sets. Such sets are locally connected, so every connected component of $\mathcal M$ is open. Extending Definition \ref{def:toric-nondegeneracy} to conical metrics, we call $g$ \emph{torically nondegenerate} if $\cK_g=0$. At such a point, the local model in Proposition~\ref{prop:localmoduli} is an interval, and the moduli space is a one-dimensional manifold near $[(\beta,g)]$, with local coordinate $\beta$.

\subsection{Orientation of the moduli space}

Let $\mathcal D_g$ be the full linearization of \eqref{eq:fixedtarget} at $(\beta_0,0)$:
$$
    \mathcal D_g:\mathbb R\oplus L^2_8\longrightarrow L^2_6, \qquad \mathcal D_g(t,k)=tB_g+\mathscr P_gk, \qquad B_g=\partial_\beta\Psi(\beta_0,0).
$$
It is Fredholm of index one. Its cokernel is represented, using $L^2$ pairings, by
\begin{equation}\label{eq:moduli-cokernel}
    \left\{\tau\in\cK_g: \langle B_g,\tau\rangle_{L^2}=0\right\}.
\end{equation}
In particular, every cokernel element is represented by a smooth invariant transverse-traceless tensor. We call a point \emph{regular} if $\mathcal D_g$ is surjective. Every nondegenerate point is regular, but regularity alone does not assert that $\beta$ is a local coordinate.

The determinant lines
$$
    \det\mathcal D_g =\Lambda^{\max}\ker\mathcal D_g\otimes \bigl(\Lambda^{\max}\operatorname{coker}\mathcal D_g\bigr)^*
$$
form a line bundle over $\mathcal M$. To describe its local topology, choose a finite-dimensional subspace $E\subset L^2_6$ such that the stabilization
$$
    \mathcal D_{h,E}(v,e)=\mathcal D_hv+e
$$
is surjective for $h$ near $g$ in a fixed chart. Its kernels form a vector bundle, and stabilization gives
$$
    \det\mathcal D_h \cong\Lambda^{\max}\ker\mathcal D_{h,E} \otimes(\Lambda^{\max}E)^*.
$$
The standard transition isomorphisms make these descriptions independent of the stabilization and the divergence gauge; see \cite{Zinger}. They descend to the moduli space because torus translations act trivially on invariant tensors. At a regular point, the kernel of $\mathcal D_g$ is the tangent space of the local one-dimensional moduli space, and hence
$$
    \det\mathcal D_g\cong T_{[(\beta_0,g)]}\mathcal M.
$$

Let $L_g$ be the restriction of $\mathscr P_g$ to invariant transverse-traceless tensors. Thus
$$
    L_g=\tfrac12(\nabla_g^*\nabla_g-2\mathring R_g).
$$
This operator is self-adjoint and bounded below, with compact resolvent. Let $\mu(g)$ be the number of its negative eigenvalues, counted with multiplicity.

\begin{lemma}\label{lem:degree-orientation}
The determinant line bundle $\det\mathcal D\to{\mathcal M}$ has a natural orientation, determined by the convention that increasing $\beta$ orients the angle parameter. At a nondegenerate metric $g$, the induced orientation is $(-1)^{\mu(g)}$ times the orientation defined by increasing $\beta$.
\end{lemma}

\begin{proof}

Suppose $\mathscr P_gk=\gamma k$. By \eqref{eq:linearized-divergence},
$$
    \delta_g\delta_g^*(\delta_gk)=\gamma\delta_gk.
$$
The form $\delta_gk$ is orthogonal to the $\mathbb T^2$-invariant Killing fields. On this orthogonal complement, $\delta_g\delta_g^*$ has a positive spectral gap. Hence $\delta_gk=0$ when $|\gamma|$ is sufficiently small. The trace equation \eqref{eq:linearized-trace} then gives
$$
    (\nabla_g^*\nabla_g-3+\gamma)\tr_gk=0.
$$
The spectrum of $\nabla^*_g\nabla_g$ on $\mathbb T^2$-invariant functions is contained in $\{0\}\cup[4,\infty)$, so $\tr_gk=0$ whenever $|\gamma|<1$. Thus every eigentensor of $\mathscr P_g$ with eigenvalue sufficiently close to zero is transverse-traceless. These spectral gaps persist locally. Choose $\varepsilon>0$ sufficiently small that every eigentensor of $\mathscr P_g$ with eigenvalue in $(-\varepsilon,\varepsilon)$ is transverse-traceless and neither endpoint lies in the spectrum. After shrinking the neighborhood of $g$, the same properties hold for nearby Einstein metrics $h$ in the chosen chart. Let $V_h$ be the sum of the eigenspaces in this interval. These spaces are the same for $\mathscr P_h$ and $L_h$. Their spectral projections vary continuously, and hence the spaces $V_h$ form a finite-dimensional vector bundle of constant rank. Eigenvalues may cross zero, changing the kernel, but no eigenvalue crosses the endpoints of the chosen interval. The operators $\mathscr P_h$ and $L_h$ have the same restriction $A_h$ to $V_h$, and both are invertible on their respective complementary spectral subspaces. Thus their determinant lines both reduce to $\det A_h$. The exact sequence
$$
    0\longrightarrow\ker A_h\longrightarrow V_h\xrightarrow{A_h}V_h\longrightarrow\operatorname{coker}A_h\longrightarrow0
$$
gives the natural identification
$$
    \det A_h\cong\Lambda^{\max}V_h\otimes(\Lambda^{\max}V_h)^*.
$$
These identifications are continuous in $h$, including where the kernel dimension changes, and therefore identify the determinant line bundles
$$
    \det\mathscr P\cong\det L
$$
over this neighborhood. The comparison is independent of the chosen spectral window, since changing the window adds the same invertible spectral blocks to both reductions.

Choose $\Lambda>0$ outside the spectrum of $L_g$, and let $E_g$ be the sum of its eigenspaces with eigenvalues less than $\Lambda$. This space is finite-dimensional, and $L_g$ is positive and invertible on $E_g^\perp$. Finite-dimensional reduction gives
\begin{equation}\label{eq:orientation with LAmbda}
    \Phi_{g,\Lambda}:\det L_g\cong\Lambda^{\max}E_g\otimes(\Lambda^{\max}E_g)^*.
\end{equation}
In the following, $\displaystyle\det_{\mathrm{num}}$ denotes the ordinary numerical determinant of a finite-dimensional endomorphism. The isomorphism $\Phi_{g,\Lambda}$ uses the restriction $A_g=L_g|_{E_g}$. When $A_g$ is invertible, so is $L_g$, and its zero kernel and cokernel give $\det L_g=\Lambda^0\{0\}\otimes(\Lambda^0\{0\})^*$ with canonical generator $\mathbf 1_g=1\otimes1^*$. The isomorphism \eqref{eq:orientation with LAmbda} sends this generator to
$$
    \Phi_{g,\Lambda}(\mathbf 1_g) =u\otimes((\Lambda^{\max}A_g)u)^* =(\det_{\mathrm{num}}A_g)^{-1}(u\otimes u^*),
$$
where $u$ is any nonzero top wedge of $E_g$ and $w^*$ denotes the dual element taking value $1$ on $w$. Orient the right-hand side of \eqref{eq:orientation with LAmbda} by its canonical pairing. Increasing $\Lambda$ adds only positive eigenspaces, so this orientation is independent of the cutoff. At a torically nondegenerate metric, its sign relative to the canonical generator of $\det L_g$ is therefore
$$
    \operatorname{sgn}\det_{\mathrm{num}}(L_g|_{E_g})=(-1)^{\mu(g)}.
$$

Finally, the homotopy
$$
    \mathcal D_{g,s}(t,k)=\mathscr P_gk+s tB_g, \qquad 0\leq s\leq1,
$$
transports the product orientation of $\det\mathscr P_g\otimes\Lambda^1\mathbb R$ to $\det\mathcal D_g$, where the last factor $\Lambda^1\R$ is oriented by increasing $\beta$. Equivalently, one may use the determinant isomorphism for adding the angle parameter; see \cite{Zinger}. The spectral pairing is preserved by equivariant pullback, and changes of reference family give homotopic identifications. Thus these orientations descend and agree on overlapping charts. At a nondegenerate point, projection of $\ker\mathcal D_g$ onto the angle parameter is an isomorphism, giving the asserted sign.
\end{proof}

\subsection{The oriented boundary count}
\label{subsec:boundary-count}

We now equip $\mathcal M$ with the determinant-line orientation of Lemma~\ref{lem:degree-orientation}. Define the area function on $\mathcal M$
$$
    \mathcal A([(\beta,g)]) =\sum_{\alpha=1}^q\Area_g(\Sigma_\alpha).
$$
For $\eta>0$, set
$$
    \mathcal M_\eta=\mathcal A^{-1}([\eta,\infty)), \qquad \partial\mathcal M_\eta:=\mathcal A^{-1}(\eta).
$$

\begin{proposition}\label{prop:parity}
Let $\eta>0$ be such that every point of $\partial\mathcal M_\eta$ is regular and $\dd\mathcal A\neq0$ there. If $\mathcal M_\eta$ is compact, then $\partial\mathcal M_\eta$ is finite and
$$
    \sum_{x\in\partial\mathcal M_\eta} \operatorname{sign}_{\partial\mathcal M_\eta}(x)=0.
$$
The signs are the boundary signs induced by the determinant-line orientation, with outward direction $-\partial_a$ for the local coordinate $a=\mathcal A$. In particular, $\partial\mathcal M_\eta$ has even cardinality.
\end{proposition}

\begin{proof}
The set $\partial\mathcal M_\eta$ is closed in the compact space $\mathcal M_\eta$, hence compact. At each of its points, regularity makes the moduli space locally one-dimensional, and $\dd\mathcal A\neq0$ makes $\mathcal A$ a local coordinate. Thus every point of $\partial\mathcal M_\eta$ is isolated, so this set is finite.

Choose neighborhoods of $\partial\mathcal M_\eta$ on which the linearization is surjective. The set of points in $\mathcal M_\eta$ where it is not surjective is compact and disjoint from these neighborhoods. By \eqref{eq:moduli-cokernel}, its cokernels are represented by smooth invariant transverse-traceless tensors. Compactness and openness of surjectivity allow us to choose finitely many such tensors whose extensions stabilize the equation near this set.

We construct the extensions equivariantly. Let $\tau_\nu$ be a smooth invariant tensor representing a selected element of $\operatorname{coker}\mathcal D_{g_\nu}$ at $(\beta_\nu,g_\nu)$. Write $b_{\nu,\beta}$ and $J_{\nu,\beta}$ for the reference family and tensor identifications based at this pair, with $b_{\nu,\beta_\nu}=g_\nu$ and $J_{\nu,\beta_\nu}=\mathrm{Id}$. Applying the earlier gauge-fixing argument with $\beta$ as a parameter gives a local slice $S_\nu$ of pairs $(\beta,h)$ near $(\beta_\nu,g_\nu)$, where
$$
    h=b_{\nu,\beta}+J_{\nu,\beta}k, \qquad \delta_hb_{\nu,\beta}=0.
$$
The gauge-fixing diffeomorphisms preserve $\beta$, while the reference metric $b_{\nu,\beta}$ varies with it. For $(\beta,h)\in S_\nu$, set
$$
    \widetilde\sigma_\nu(h)=J_{\nu,\beta}\tau_\nu.
$$
After shrinking $S_\nu$, compactness of isometries together with local gauge uniqueness implies that any admissible diffeomorphism relating two points of $S_\nu$ is a torus translation. Since these translations act trivially on invariant tensors, the rule
$$
    \widetilde\sigma_\nu(F^*h)=F^*\widetilde\sigma_\nu(h)
$$
defines a well-defined extension to the saturation of $S_\nu$. Choose a cutoff $\chi_\nu$ on $S_\nu$, equal to one near $(\beta_\nu,g_\nu)$ and supported in a smaller slice neighborhood, and extend it invariantly by $\chi_\nu(F^*h)=\chi_\nu(h)$. Choose the corresponding supports in the quotient disjoint from the boundary neighborhoods.

The transported tensors $\widetilde\sigma_\nu(h)$ need not be divergence-free with respect to $h$. To ensure this property, for an invariant symmetric tensor $\tau$ define
$$
    \mathsf Q_h\tau =\tau-\delta_h^*(\delta_h\delta_h^*)_\perp^{-1}\delta_h\tau.
$$
Here the inverse is taken on invariant one-forms orthogonal to the torus Killing fields, with orthogonality defined using $h$. Since $\delta_h\tau$ is automatically orthogonal to these fields, the inverse is well defined. The operator $\mathsf Q_h$ is the $L^2(h)$-orthogonal projection onto $\ker\delta_h$; in particular, it fixes every divergence-free tensor. Set
$$
    \sigma_\nu(h)=\chi_\nu(h)\mathsf Q_h\widetilde\sigma_\nu(h),
$$
and extend this section by zero outside its support. For each fixed $h$, the cutoff $\chi_\nu(h)$ is a scalar constant on $X$. Thus multiplication by the cutoff preserves divergence-freeness. Naturality of $\mathsf Q_h$ under pullback, together with the equivariance of $\widetilde\sigma_\nu$ and invariance of $\chi_\nu$, also preserves equivariance. Consequently,
$$
    \delta_h\sigma_\nu(h)=0, \qquad \sigma_\nu(F^*h)=F^*\sigma_\nu(h).
$$
These identities hold throughout the domain, including where the cutoff varies.

For the required differentiability, take the divergence gauges and cutoffs on $L^2_4$ neighborhoods. On $L^2_8$ metrics the gauge-fixing diffeomorphisms depend smoothly with values in $L^2_9$, and changes of representative are $C^2$ with values in $L^2_6$. The cutoffs are consequently $C^2$ in every divergence gauge. Pullback of the smooth selected tensors and the elliptic projection above give $C^2$ sections $\sigma_\nu:L^2_8\to L^2_6$.

We use these sections to perturb the Einstein equation. To verify local surjectivity near $(\beta_\nu,g_\nu)$, first consider the perturbation directions associated with this reference pair. For $h=b_{\nu,\beta}+J_{\nu,\beta}k$, the perturbed equation in the divergence gauge is
$$
    J_{\nu,\beta}^{-1}\left\{ \mathscr E_3(h)-\delta_h^*\delta_hb_{\nu,\beta} -\sum_{j=1}^{d_\nu}u_j\sigma_{\nu,j}(h) \right\}=0, \qquad u\in\mathbb R^{d_\nu}.
$$
Write $\mathcal D_h$ for the full linearization of the unperturbed gauged equation in this reference chart. At $u=0$, the linearization of the displayed perturbed equation, with variations $(t,\ell)$ in $(\beta,k)$ and $v$ in $u$, is
$$
    (t,\ell,v)\longmapsto \mathcal D_h(t,\ell) -\sum_{j=1}^{d_\nu}v_jJ_{\nu,\beta}^{-1}\sigma_{\nu,j}(h).
$$
At $(\beta_\nu,g_\nu)$, we have $J_{\nu,\beta_\nu}=\mathrm{Id}$ and $\sigma_{\nu,j}(g_\nu)=\tau_{\nu,j}$, so this becomes
$$
    (t,\ell,v)\longmapsto \mathcal D_{g_\nu}(t,\ell) -\sum_{j=1}^{d_\nu}v_j\tau_{\nu,j}.
$$
This operator is surjective because the classes of the $\tau_{\nu,j}$ span $\operatorname{coker}\mathcal D_{g_\nu}$. By openness of surjectivity, the augmented linearization remains surjective for nearby $(\beta,h)$. Equivalently,
$$
    L^2_6 =\operatorname{im}\mathcal D_h +\operatorname{span}\left\{ J_{\nu,\beta}^{-1}\sigma_{\nu,j}(h): 1\leq j\leq d_\nu \right\}.
$$
Thus the classes of these tensors span $\operatorname{coker}\mathcal D_h$ throughout a sufficiently small neighborhood of $(\beta_\nu,g_\nu)$.

Combining the perturbations from all reference neighborhoods gives the equivariant equation
$$
    \mathscr E_3(h) =\sum_\nu\sum_{j=1}^{d_\nu}u_{\nu,j}\sigma_{\nu,j}(h), \qquad u\in\mathbb R^N, \qquad N=\sum_\nu d_\nu,
$$
with the common cone angle allowed to vary. In any divergence gauge, writing $b_\beta$ and $J_\beta$ for its reference data, the gauged equation is
$$
    J_\beta^{-1}\left\{ \mathscr E_3(h)-\delta_h^*\delta_hb_\beta -\sum_\nu\sum_{j=1}^{d_\nu}u_{\nu,j}\sigma_{\nu,j}(h) \right\}=0.
$$
Since the perturbations are divergence-free, the Bianchi identity removes the gauge term exactly as in Proposition~\ref{prop:localmoduli}. Equivariance therefore ensures that solutions in different divergence gauges represent the same geometric equation.

At every unperturbed solution $([(\beta,g)],0)$ with $[(\beta,g)]\in\mathcal M_\eta$, the augmented linearization is surjective of index $N+1$. By openness of surjectivity, it remains surjective in sufficiently small neighborhoods of these points. In each fixed divergence gauge, the implicit function theorem therefore gives a local $C^2$ solution manifold of dimension $N+1$.

We next check compatibility under changes of divergence gauge. On an overlap, let $T_{\mu\nu}$ denote the change of representative from gauge $\nu$ to gauge $\mu$, preserving the angle and perturbation parameters. By the previously established estimates, this map is $C^2$ with tensor coordinates in $L^2_6$. On each local solution manifold, choose coordinates from linear functions of the angle and perturbation parameters and $L^2$ pairings with smooth invariant tensors, with differentials independent on the augmented kernel. These coordinate functions extend to continuous affine linear functionals on the $L^2_6$ coordinate space, so $T_{\mu\nu}$ induces $C^2$ coordinate transitions. Slice uniqueness modulo torus translations makes these transitions well defined and consistent on triple overlaps. Thus the local solution spaces glue to a $C^2$ manifold $W$ near $\mathcal M_\eta\times\{0\}$, with the $C^2$ parameter projection
$$
    \pi:W\longrightarrow\mathbb R^N, \qquad \pi([(\beta,h)],u)=u.
$$

Near each point of $\partial\mathcal M_\eta$, all perturbations vanish, so the solution space is a product of the original regular solution interval and a parameter ball. On this interval, $a=\mathcal A$ is a local coordinate. Restricting to $a\geq\eta$ and shrinking the neighborhood and parameter ball gives a $C^2$ manifold with boundary, still denoted by $W$, with the natural identifications
$$
    \pi^{-1}(0)=\mathcal M_\eta\times\{0\}, \qquad \partial W=\partial\mathcal M_\eta\times B_\rho(0).
$$
Under the second identification, $\pi|_{\partial W}$ is projection onto the parameter ball.

Along the zero fiber, stabilization gives
$$
    \Lambda^{N+1}T_{(h,0)}W \cong\det\mathcal D_h\otimes\Lambda^N\mathbb R^N.
$$
The determinant-line orientation and the standard orientation of $\mathbb R^N$ orient $W$ along the zero fiber. Since an orientation of a real vector bundle along a closed subset of a manifold extends to a neighborhood of that subset, this orients $W$ after shrinking around the zero fiber and induces orientations on the regular fibres of $\pi$.

Choose a relatively open neighborhood $V$ of the zero fiber with compact closure in $W$. The relative boundary $\overline V\setminus V$ is compact and disjoint from the zero fiber. By compactness, there is $0<\varepsilon<\rho$ such that
$$
    |\pi|>\varepsilon \quad\text{on }\overline V\setminus V, \qquad \partial\mathcal M_\eta\times B_\varepsilon(0)\subset V.
$$
Choose a regular value $u\in B_\varepsilon(0)$ by Sard's theorem. Then
$$
    \pi^{-1}(u)\cap\overline V
$$
is a compact oriented one-manifold whose boundary is naturally identified with $\partial\mathcal M_\eta$. The product structure near the cuts preserves their boundary signs. Since the signed boundary count of a compact oriented one-manifold is zero, the proposition follows.
\end{proof}
\begin{remark}\label{rem:independent-cone-angles}
The local moduli and orientation arguments also apply when the cone angles vary independently. Fix an ordering of all invariant spheres $\Sigma_1,\ldots,\Sigma_m$, preserve this order in the equivalence relation, and write $\boldsymbol\beta=(\beta_1,\ldots,\beta_m)\in(0,\infty)^m$. The fixed-angle analysis is unchanged, and the identifications $J_{\boldsymbol\beta}$ are constructed by making the angular replacement independently in each normal plane, including at corners. The local Kuranishi model is then
$$
    \kappa_g:U_g\subset\cK_g\times\mathbb R^m \longrightarrow\cK_g.
$$
The full linearization $\mathcal D_g$ has Fredholm index $m$, and the determinant-line construction gives
$$
    \det\mathcal D_g \cong\det\mathscr P_g\otimes\Lambda^m\mathbb R^m.
$$
Orienting the parameter space by $\dd\beta_1\wedge\cdots\wedge \dd\beta_m$ therefore orients the determinant line bundle. At a nondegenerate metric, the cone angles are local coordinates on the moduli space, and the induced orientation is $(-1)^{\mu(g)}$ times their standard orientation.
\end{remark}

\section{End structure of the moduli space}
\label{sec:end structure}

In this section the goal is to understand the end structure of the moduli space $\mathcal M$. We will consider a sequence $[(\beta_j, g_j)]\in \mathcal M$ satisfying the bound
$$
    1\leq\beta_j\leq\beta_+<\infty.
$$
So passing to a subsequence, we may assume $\beta_j\to\beta_\infty$.

We study the limiting behavior of $g_j$ under the above angle bound. To describe convergence near conical singularities, we use the angular developments of Definition~\ref{def:conical}. We say that $g_j$ converges \emph{smoothly in local developments} if the cone angles converge and the developed metrics converge in $C^\infty$ on fixed smaller coordinate neighborhoods, including where the radial coordinates vanish.

The main compactness Theorem \ref{thm:compactness} says that in our setting, there are only two alternatives passing to a subsequence

\begin{itemize}
\item $g_j$ converges smoothly to a limit toric conical Einstein metric $g_\infty$ on $X_{\mR}$;
\item $g_j$ converges smoothly, away from a point, to a smooth toric Einstein metric on $X_{\mR_0}$, and the only non-flat bubble is a $\Ae$ toric conical gravitational instanton, whose rod structure is determined by $W$.
\end{itemize}
We also use the gluing construction to prove a converse in the second alternative. Furthermore, we also establish quantitative convergence estimates and prove the degeneration and the gluing construction exactly match.

\subsection{Preliminary results}

Write $\mathcal V=\{p_1,\ldots,p_m\}$ for the torus fixed points. By Proposition~\ref{geom:diameter}, Theorem~\ref{geom:noncollapse} and Proposition~\ref{geom:vertexcurvature}, we have the following geometric bounds
\begin{equation}\label{scd:bounds}
    \diam(X_{\mR},g_j)\leq D,\qquad \Vol_{g_j}(X_{\mR})\geq v>0,\qquad |\Rm_{g_j}|(x)\leq C\max\{1,d_{g_j}(x,\mathcal V)^{-2}\}.
\end{equation}
The computation in the proof of Proposition~\ref{geom:hitchinthorpe} also provides a uniform bound on the curvature energy
\begin{equation}\label{deg:energy}
    \int_{X_{\mR}}|\Rm_{g_j}|^2 \dvol_{g_j} \leq C\sum_{\alpha=1}^m \beta_{\alpha,j}\beta_{\alpha+1,j} \leq C m\beta_+^2.
\end{equation}
As in the proof of the angle bound in Section~\ref{sec:a priori estimates}, we can pass to a fixed finite quotient so that the cone angles are all smaller than $2\pi$. Then we have volume comparison and the Sobolev inequality on the quotient space. From these we get a uniform volume ratio bound
\begin{equation}\label{scd:volume}
    cr^4\leq\Vol_{g_j}B_{g_j}(x,r)\leq Cr^4 \quad(0<r\leq1)
\end{equation}
and a uniform Sobolev inequality
\begin{equation}\label{scd:Sobolev}
    \|f\|_{L^4}^2\leq C\bigl(\|\dd f\|_{L^2}^2+\|f\|_{L^2}^2\bigr).
\end{equation}

\begin{proposition}
\label{prop:bounded-curvature-compactness} If
$$
    \sup_j\max_{X_{\mR}}|\Rm_{g_j}|<\infty,
$$
then, after passing to a subsequence and pulling back by $\T^2$-equivariant diffeomorphisms preserving each $\Sigma_\alpha$ and $p_\alpha$, the metrics $g_j$ converge smoothly in local developments to a toric conical Einstein metric $g_\infty$ with the same rod structure $\mR$ and cone angle $2\pi\beta_\infty$ along $\mathcal S$.
\end{proposition}
\begin{proof}
By Lemma~\ref{geom:rodhalfdisk}, the rods have neighborhoods of uniform size away from the turning points.
The same argument bounds below the distances between distinct turning points and between a turning point and a nonincident rod.
Together with Lemma~\ref{geom:basecomparison} and the angle bounds, these estimates provide uniformly noncollapsed local developments of uniform size.

Then the conclusion follows from standard compactness for Einstein metrics. 
\end{proof}

The following generalizes the standard $\epsilon$-regularity theorem to our toric conical setting.
\begin{lemma}\label{scd:epsilon}
There are $\varepsilon_0,C>0$ such that, for $0<r\leq\tfrac12$,
$$
    \int_{B_{g_j}(x,2r)}|\Rm_{g_j}|^2\dvol_{g_j}<\varepsilon_0 \Longrightarrow\sup_{B_{g_j}(x,r)}|\Rm_{g_j}| \leq Cr^{-2} \left(\int_{B_{g_j}(x,2r)}|\Rm_{g_j}|^2\dvol_{g_j}\right)^{1/2}.
$$
\end{lemma}

\begin{proof}
The Einstein curvature equation gives away from the conical singularities the differential inequality
$$
    \Delta_{g_j}|\Rm_{g_j}|\geq-C|\Rm_{g_j}|^2.
$$
We claim that this inequality holds weakly globally on $X$. This follows from the fact that $|\Rm_{g_j}|$ is bounded on $X$ for each fixed $j$ and the conical locus $\mathcal S$ has zero $W^{1,2}$-capacity.

Then the uniform bound on $|\Rm_{g_j}|$ follows from the Sobolev inequality and the usual Moser iteration.
\end{proof}

\subsection{Bubble classification}
\label{sec:degeneration-bubble} Now suppose that
$$
    Q_j:=\max_{X}|\Rm_{g_j}|\longrightarrow\infty, \qquad \ell_j:=Q_j^{-1/2}.
$$
Choose a maximum point $x_j$. By \eqref{scd:bounds}, after passing to a subsequence there is a fixed point $p_i$, independent of $j$, such that
$$
    d_{g_j}(x_j,p_i)\leq C\ell_j.
$$
Set $\widehat g_j=\ell_j^{-2}g_j$ and use $p_i$ as base point. Under this rescaling, the reduced data $\underline g_j$, $G_j$, and $\rho_j$ are multiplied by $\ell_j^{-2}$, while $\Phi_j$ is unchanged.

By \eqref{scd:volume} and scaling,
$$
    cR^4\leq\Vol_{\widehat g_j}B_{\widehat g_j}(x,R)\leq CR^4,\qquad 0<R\leq\ell_j^{-1}.
$$
Thus the rescaled sequence is uniformly noncollapsed and satisfies $|\Rm_{\widehat g_j}|\leq1$. After passing to a subsequence, $(X_{\mR},\widehat g_j,p_i)$ converges in the pointed Gromov-Hausdorff sense to $(Y,d_Y,p_\infty)$. The local equivariant compactness argument of Proposition~\ref{prop:bounded-curvature-compactness} gives smooth convergence in local developments on bounded subsets, retaining the fixed torus action and its primitive isotropy weights. Since the rescaled Einstein constants are $3\ell_j^2\to0$, the distance $d_Y$ is induced by a complete toric conical Ricci-flat metric $h$.

The points $x_j$ remain at bounded rescaled distance from $p_i$, so their subsequential limit satisfies $|\Rm_h|=1$. Thus $\max_Y|\Rm_h|=1$, and $h$ is not flat. Moreover, it has  finite curvature energy and Euclidean volume growth.

We write $(N,\underline h)=Y/\T^2$.
The equivariant convergence is locally smooth on $N^\circ$ and preserves the minimizing property of interior geodesic segments.
Since the rods are totally geodesic, the corners have angle $\pi/2$, and $K_{\underline h}\leq0$ by Proposition~\ref{geom:Kupper}, the argument of Lemma~\ref{geom:basecomparison} implies that the exponential map identifies $N^\circ$ diffeomorphically with a star-shaped domain.
Consequently, $N$ is contractible.

\begin{lemma}\label{scd:proper}
The limit $(Y,h)$ is proper in the sense of Definition~\ref{def: app gravitational instanton proper}.
\end{lemma}

\begin{proof}
The limit $(Y,h)$ is complete and integrable, while $N^\circ$ is simply connected, $\partial N\neq\varnothing$, and $K_{\underline h}\leq0$.
The finite curvature energy and Lemma~\ref{scd:epsilon}, after rescaling and passing to the limit, imply quadratic curvature decay.
So Lemma~\ref{geom:proper-ricciflat} implies that $(Y,h)$ is proper.
\end{proof}

\begin{proposition}\label{deg:end}
The limit $(Y,h)$ is a proper toric conical gravitational instanton with distinct infinite rod directions. Its end is $\Ae$.
\end{proposition}

\begin{proof}
By Lemma~\ref{scd:proper}, $(Y,h)$ is a proper toric conical gravitational instanton. The Euclidean volume growth case of Appendix~\ref{app:asymptotics of toric gravitational instantons} gives distinct infinite rod directions $\mv_-,\mv_+$ and $\ALE$ asymptotic. In our case, their primitive vectors form a basis of $\Lambda$, so the end is actually $\Ae$.
\end{proof}

\begin{proposition}\label{prop:bubble}
The limit $(Y,h)$ is $\Ae$ and $\beta_\infty>1$. Its finite rods are exactly the rods $\mI_1,\ldots,\mI_q$ appeared in $W$.
\end{proposition}

\begin{proof}
Let $s,t$ be the original indices of the two infinite rods, and write $\beta_{\alpha,\infty}=\lim_j\beta_{\alpha,j}$. Then $\beta_{\alpha, \infty}$ is either $1$ or $\beta_\infty$.
By Proposition~\ref{deg:end}, the volume density at infinity of $Y$ is $\beta_{s,\infty}\beta_{t,\infty}$. By Bishop-Gromov comparison (as in \cite{LiSun}) we know that at a fixed point $p_\alpha$ in $Y$,
\begin{equation}\label{eq:density}
    \beta_{\alpha-1,\infty}\beta_{\alpha,\infty} > \beta_{s,\infty}\beta_{t,\infty}.
\end{equation}
The inequality is strict since $Y$ is non-flat. In particular, $\beta_\infty>1$.

If both infinite rods had cone angle $2\pi\beta_\infty$, the right-hand side of \eqref{eq:density} would be $\beta_\infty^2$, but the left-hand side is at most $\beta_\infty^2$. Contradiction. If exactly one had angle $2\pi\beta_\infty$, the fixed point incident to the smooth infinite rod would have left-hand side at most $\beta_\infty$, again contradicting \eqref{eq:density}. Thus both infinite rods are smooth. Moreover, \eqref{eq:density} also excludes the possibility of two adjacent smooth rods.

The rods of $Y$ occur consecutively in $\mR$. The finite rod adjacent to either infinite rod must therefore be inserted. Since the inserted rods form a single consecutive block, the chain inherited on $Y$ must be of the form
$$
    \mI_m,\mI_1,\ldots,\mI_q,\mI_{q+1}.
$$
\end{proof}

\subsection{The compact limit}
\label{sec:degeneration-compact-limit}
We return to the original sequence $g_j$. The additional rescalings considered below either send the conical locus to infinity or collapse it to a point. We will show that the resulting limits are complete toric Ricci-flat spaces with smooth rod interiors and at most isolated orbifold singularities. Their finite curvature energy and Euclidean volume growth give $\ALE$ asymptotics. The restriction of the rod directions to $\{A,B,C\}$ makes the group at infinity trivial, since any two distinct directions among $A,B,C$ are represented by primitive vectors forming a basis of $\Lambda$. Thus these limits are $\Ae$. The following simple rigidity statement then excludes any additional non-flat limit.

\begin{lemma}\label{deg:smooth-ae-rigidity}
A complete toric Ricci-flat $\ALE$ four-orbifold with isolated singularities and rod directions in $\{A,B,C\}$ is Euclidean.
\end{lemma}

\begin{proof}
If $\mv_-,\mv_+$ are the primitive infinite rod vectors, the group at infinity is
$$
    \Lambda/(\mathbb Z\mv_-+\mathbb Z\mv_+).
$$
Its finiteness implies that the two directions are distinct. The constraint on the rod vectors therefore makes this group trivial, so the end is AE. Bishop-Gromov comparison then implies that the orbifold is Euclidean.
\end{proof}

Recall that $\mathcal S=\bigcup_{\alpha=1}^q\Sigma_\alpha\subset X_{\mR}$ is the union of the inserted totally geodesic two-spheres, and that $p\in M$ is the original torus fixed point replaced by this chain.

\begin{theorem}\label{thm:compactness}
After passing to a subsequence, either $g_j$ converges smoothly in local developments on $X_{\mR}$, or $\mathcal S$ shrinks to $p$ and the compact Gromov-Hausdorff limit is a smooth toric Einstein metric on $M$, with smooth convergence away from $p$. In the latter case, the limit $(Y,h)$ in Proposition~\ref{prop:bubble} is the only non-flat bubble.
\end{theorem}

\begin{proof}
The bounded curvature case follows from Proposition~\ref{prop:bounded-curvature-compactness}. Otherwise, Proposition~\ref{prop:bubble} gives the AE limit $(Y,h)$ at scale $\ell_j$, containing all the inserted spheres. Write $q\in\mathcal S\subset X_{\mR}$ for the torus fixed point used as the base point for this rescaled convergence. In particular,
$$
    \diam_{g_j}\mathcal S=O(\ell_j).
$$
We first claim that
\begin{equation}\label{deg:exterior-curvature}
    d_{g_j}(x,\mathcal S)^2|\Rm_{g_j}|(x)\leq C.
\end{equation}
Otherwise, by a point selection and rescaling by curvature we would get a complete smooth nonflat Ricci-flat limit: the distance to $\mathcal S$ tends to infinity in this rescaling. Proposition~\ref{geom:vertexcurvature} supplies a torus fixed point at bounded rescaled distance. The limit has finite curvature energy and Euclidean volume growth, so it is $\ALE$ by \cite{BKN}. Its rod vectors belong to $\{A,B,C\}$, contradicting Lemma~\ref{deg:smooth-ae-rigidity}.

Using the geometric bounds we obtain a compact Gromov-Hausdorff limit
$$
    (X_{\mR},g_j,q)\longrightarrow(X_\infty,p).
$$
Since $\diam_{g_j}\mathcal S\to0$, the whole chain $\mathcal S$ converges to $p$. By estimate \eqref{deg:exterior-curvature} we have smooth convergence to an Einstein metric $g_\infty$ away from $p$. Finiteness of the curvature energy and Lemma~\ref{scd:epsilon} imply
$$
    r^2\sup_{B(p,2r)\setminus B(p,r)}|\Rm_{g_\infty}|\longrightarrow0 \qquad\text{as }r\to0.
$$
By volume comparison on the fixed finite quotient we see that tangent cones at $p$ are flat. The argument of \cite[proof of Theorem~5.7]{DonaldsonSun} implies that $g_\infty$ extends across $p$ as an Einstein orbifold metric.

We next show that no fixed point outside $\mathcal S$ converges to $p$. Otherwise, after passing to a subsequence, some fixed point $p_\alpha\notin\mathcal S$ satisfies
$$
    s_j=d_{g_j}(q,p_\alpha)\longrightarrow0.
$$
Proposition~\ref{prop:bubble} implies $s_j/\ell_j\to\infty$. Rescaling by $s_j^{-2}$, \eqref{deg:exterior-curvature} and taking limit we then get a complete Ricci-flat $\ALE$ orbifold. Again by Lemma~\ref{deg:smooth-ae-rigidity} this limit has to be Euclidean. However, its torus action has two distinct fixed points, given by the limits of $q$ and $p_\alpha$, which are at distance one. This is impossible.

So every remaining rod and turning point persists in the limit. The two rods incident to $p$ have primitive vectors which form a basis of $\Lambda$, so the local orbifold group at $p$ is trivial. The compact limit is therefore smooth, and its rod structure is $\mR_0$. The construction in Section~\ref{sec:rod structure} identifies $(X_\infty,p)$ with $(M,p)$.

Finally, we exclude further nonflat bubbles. For fixed $R$ large and $\delta>0$ small, using the equivariant convergence we can choose two invariant three-spheres enclosing $\mathcal S$.
The inner sphere and its collar, after rescaling by $\ell_j^{-2}$, converge smoothly to a radius-$R$ coordinate sphere and its collar in the AE end of $Y$. The outer sphere and its collar converge smoothly to $\partial B_{g_\infty}(p,\delta)$ and its collar. We denote by $A_j(R,\delta)$ the region between these spheres. For sufficiently large $j$, this region contains neither conical points nor torus fixed points, and its two boundary rods have directions $B,A$. Thus the rod construction identifies its closure equivariantly with $S^3\times[0,1]$.

Since $g_j$ is Einstein and $\chi(A_j(R,\delta))=0$, the Chern-Gauss-Bonnet formula expresses its curvature energy in terms of the two boundary contributions.
By smooth convergence of the boundary collars, followed by $R\to\infty$ and $\delta\to0$, these contributions approach the Euclidean values with opposite signs.
Consequently,
\begin{equation}\label{deg:neck-energy}
    \lim_{\delta\to0}\lim_{R\to\infty}\limsup_{j\to\infty} \int_{A_j(R,\delta)}|\Rm_{g_j}|^2 \dvol_{g_j}=0.
\end{equation}
Using the $\varepsilon$-regularity we see no nonflat bubbles can arise in this region.
This proves that $(Y,h)$ is the only nonflat bubble.
\end{proof}

In particular, if $\beta_\infty=1$, then Proposition~\ref{prop:bubble} excludes bubbling, so $g_j$ converges smoothly in local developments to a smooth toric Einstein metric on $X_{\mR}$.

\subsection{A gluing construction}
\label{sec:gluing}

Let $g_0$ be a smooth compact toric Einstein metric on $M$ satisfying $\Ric_{g_0}=3g_0$.
Let $(Y,h)$ be the $\Ae$ toric conical Ricci-flat gravitational instanton associated with $W$, as constructed in Appendix~\ref{app:ae-models}. Its two infinite rods are smooth, and its finite rods have a common cone angle $2\pi\beta_*$, where $\beta_*>1$. In equivariant coordinates at infinity,
\begin{equation}\label{glue:model-expansion}
    h=g_{\mathrm{Euc}}+O_\infty(r^{-2}).
\end{equation}
We also denote by $\Sigma_\alpha$ the invariant spheres corresponding to the finite rods of $Y$.
For a metric $g$ on either $X_{\mR}$ or $Y$, we write
$$
    A(g):=\sum_{\alpha=1}^q\Area_g(\Sigma_\alpha).
$$
We fix a scale for $h$ and set $A_*=A(h)$.

\begin{theorem}\label{thm:glue}
If $g_0$ is torically nondegenerate, then there is a smooth family of toric conical Einstein metrics $g_a$ on $X$, defined for sufficiently small $a>0$, such that
\begin{align*}
    \Ric(g_a)&=3g_a,\\
    A(g_a)&=a.
\end{align*}
Each $g_a$ is smooth away from $\mathcal S$ and has common cone angle $2\pi\beta(a)$ along $\mathcal S$, where $\beta(a)\to\beta_*$ as $a\to0$.
Furthermore, $g_a$ converges smoothly to $g_0$ on compact subsets of $M\setminus\{p\}$, while $a^{-1}A_*g_a$ converges to $h$ smoothly in local developments on compact subsets of $Y$.
\end{theorem}

We choose equivariant normal coordinates at $p$ so that
\begin{equation}\label{glue:compact-expansion}
    g_0=g_{\mathrm{Euc}}+O_\infty(r^2).
\end{equation}
The toric nondegeneracy of $g_0$ and Lemma~\ref{lem:conical-linear} imply that $\mathscr P_{g_0}$ has trivial $\T^2$-invariant kernel.

On $(Y,h)$, Proposition~\ref{prop:ae-deformations} identifies decaying invariant infinitesimal Ricci-flat deformations with variations of the finite Weyl rod lengths, modulo equivariant gauge transformations.
At fixed cone angles, the kernel is generated by the gauged scaling deformation.
We allow the common cone angle to vary and prescribe the total area; then Proposition~\ref{prop:ae-augmented} implies that the augmented operator is invertible in the weighted spaces specified there.

We set $a=\varepsilon^2A_*$ and $r_*=\varepsilon^{2/3}$.
We identify the normal coordinates $x$ near $p$ with the $\Ae$ coordinates $y$ by $x=\varepsilon y$.
Using a smooth radial cutoff across $r_*\leq|x|\leq2r_*$, we interpolate between $\varepsilon^2h$ and $g_0$ to obtain a metric $b_\varepsilon$ that equals $\varepsilon^2h$ on the inner region and $g_0$ on the outer region.
Since $\mathcal S$ lies entirely in the inner region, we have
$$
    A(b_\varepsilon)=\varepsilon^2A_*.
$$

We construct the tensor identifications of Section~\ref{sec:moduli} from cone angle $2\pi\beta_*$, using charts and a partition of unity on a fixed compact neighborhood $U$ of $\mathcal S\subset Y$, transferred to its rescaled image $U_\varepsilon\subset X$.
We denote these identifications by $J_\beta$, suppressing their dependence on $\varepsilon$, and define
$$
    b_{\varepsilon,\beta}=J_\beta b_\varepsilon.
$$
Then $J_{\beta_*}=\mathrm{Id}$, and $J_\beta$ is the identity outside $U_\varepsilon$.
Thus $b_{\varepsilon,\beta}$ depends real-analytically on $\beta$, has cone angle $2\pi\beta$ along $\mathcal S$, and agrees with $b_\varepsilon$ outside $U_\varepsilon$.

We set
\begin{align*}
    \sigma&=\varepsilon(\beta-\beta_*),\\
    g&=b_{\varepsilon,\beta}+J_\beta k,
\end{align*}
where $k$ is a $\T^2$-invariant symmetric two-tensor expressed in the local developments of $b_\varepsilon$.
We solve the gauged Einstein equation with prescribed area:
\begin{equation}\label{eq:gaugedglue}
    \mathcal F_\varepsilon(k,\sigma)=\left(J_\beta^{-1}\bigl\{\mathscr E_3(g)-\delta_g^*\delta_g b_{\varepsilon,\beta}\bigr\},\ \varepsilon^{-1}\bigl(A(g)-\varepsilon^2A_*\bigr)\right)=0.
\end{equation}

We choose a smooth positive weight $w_\varepsilon$ satisfying, uniformly in $\varepsilon$,
\begin{align*}
    C^{-1}(\varepsilon+r)&\leq w_\varepsilon\leq C(\varepsilon+r) &&\text{near }p,\\
    C^{-1}&\leq w_\varepsilon\leq C &&\text{away from }p.
\end{align*}
Here $r$ is the common radial coordinate on the overlap, extended smoothly over the rescaled compact part of $Y$ with $0\leq r\leq C\varepsilon$.

We fix $0<\alpha<1$ and denote by
\begin{align*}
    X_\varepsilon=C^{2,\alpha}_{-1},\qquad
    Y_\varepsilon=C^{0,\alpha}_{-3}
\end{align*}
the weighted H\"older spaces of $\T^2$-invariant symmetric two-tensors.
Their norms consist of
$$
    \sum_{j=0}^2\sup_X w_\varepsilon^{j+1}|\nabla^j k|
    \qquad\text{and}\qquad
    \sup_X w_\varepsilon^3|f|,
$$
respectively, together with the corresponding scaled H\"older seminorms.
We use $b_\varepsilon$ to define the norms and covariant derivatives, and as before measure regularity in its local developments.
Then $\mathcal F_\varepsilon$ is defined on a neighborhood of $(0,0)$ in $X_\varepsilon\times\mathbb R$ and takes values in $Y_\varepsilon\times\mathbb R$.
We denote the derivative at $(0,0)$ by
\begin{equation}\label{eq:augmented}
    \mathcal L_\varepsilon(k,\sigma):=D\mathcal F_\varepsilon(0,0)(k,\sigma)=\left(\mathscr P_{b_\varepsilon}k+\sigma B_\varepsilon,\ a_\varepsilon(k)+c_\varepsilon\sigma\right),
\end{equation}
where
\begin{align*}
    &B_\varepsilon=\varepsilon^{-1}\left.\partial_\beta\bigl(J_\beta^{-1}\mathscr E_3(b_{\varepsilon,\beta})\bigr)\right|_{\beta=\beta_*},\\&
    a_\varepsilon(k)=\varepsilon^{-1}\dd A_{b_\varepsilon}(k),\\&
    c_\varepsilon=\varepsilon^{-2}\left.\partial_\beta A(b_{\varepsilon,\beta})\right|_{\beta=\beta_*}.
\end{align*}
Indeed, since $\delta_g g=0$, we have
$$
    -\delta_g^*\delta_g b_{\varepsilon,\beta}=\delta_g^*\delta_g(J_\beta k).
$$
Thus the gauge term vanishes for $k=0$ at every $\beta$, while its derivative in $k$ at $(0,0)$ is $\delta_{b_\varepsilon}^*\delta_{b_\varepsilon}k$.

\begin{proposition}\label{prop:glueinverse}
For $\varepsilon>0$ sufficiently small, $\mathcal L_\varepsilon$ is invertible and with a uniformly bounded inverse.
\end{proposition}

\begin{proof}
By rescaling the local Schauder estimates, we obtain
\begin{equation}\label{eq:glue-schauder}
    \|k\|_{X_\varepsilon}\leq C\left(\|\mathscr P_{b_\varepsilon}k\|_{Y_\varepsilon}+\sup_X w_\varepsilon|k|_{b_\varepsilon}\right),
\end{equation}
with $C$ independent of $\varepsilon$.
Indeed, the rescaled metrics have uniformly controlled coefficients on the compact part of $Y$, on the overlap annuli, and on the remaining compact part of $M$.
Over rods and turning points, we apply the estimates in local developments; see also \cite[p.~26, equation~(54)]{BiquardI}.
The same rescaling on the support of the angle variation implies that $\|B_\varepsilon\|_{Y_\varepsilon}\leq C$.

We claim that
$$
    \|k\|_{X_\varepsilon}+|\sigma|\leq C\|\mathcal L_\varepsilon(k,\sigma)\|_{Y_\varepsilon\times\mathbb R}.
$$
Otherwise, we can choose $\varepsilon_i\to0$ and pairs $(k_i,\sigma_i)$ such that
\begin{align*}
    &\|k_i\|_{X_{\varepsilon_i}}+|\sigma_i|=1,\\
    &\mathcal L_{\varepsilon_i}(k_i,\sigma_i)=(f_i,t_i),\\
    &\|f_i\|_{Y_{\varepsilon_i}}+|t_i|\longrightarrow0.
\end{align*}
After passing to a subsequence, we may assume that $\sigma_i\to\sigma_\infty$. We first consider the limit on $Y$. We denote by $J_\beta^Y$ the tensor identifications on $Y$ used to construct $J_\beta$, and set
\begin{align*}
    &b_\beta^Y=J_\beta^Yh,\\
    &\dot b=\left.\partial_\beta b_\beta^Y\right|_{\beta=\beta_*},\\
    &c_0=\left.\partial_\beta A(b_\beta^Y)\right|_{\beta=\beta_*}.
\end{align*}
The variation $\dot b$ is supported in a fixed neighborhood of $\mathcal S$ and varies the common cone angle. By construction, the transported tensors $(J_\beta^Y)^{-1}\mathscr E_0(b_\beta^Y)$ are smooth in the fixed-angle developments. Since $\mathscr E_0(h)=0$, we have 
$$
    D\mathscr E_{0,h}(\dot b)=\left.\partial_\beta\bigl((J_\beta^Y)^{-1}\mathscr E_0(b_\beta^Y)\bigr)\right|_{\beta=\beta_*}.
$$
We denote by $\iota_i$ the gluing embeddings from expanding subsets of $Y$ into $X$, and set
$$
    \ell_i=\varepsilon_i^{-1}\iota_i^*k_i.
$$
The weighted bounds imply that, after passing to a subsequence, $\ell_i$ converges in $C^2$ on compact subsets to a tensor $\ell\in C^{2,\alpha}_{-1}(Y)$. On each fixed compact subset, we have $\iota_i^*b_{\varepsilon_i,\beta}=\varepsilon_i^2b_\beta^Y$ for all sufficiently large $i$, so the rescaled equations are
\begin{align*}
    &\mathscr P_h\ell_i+3\varepsilon_i^2\ell_i+\sigma_iD\mathscr E_{0,h}(\dot b)=\varepsilon_i\iota_i^*f_i,\\
    &\dd A_h(\ell_i)+c_0\sigma_i=t_i.
\end{align*}
Here the angle derivative of $(J_\beta^Y)^{-1}(3\varepsilon_i^2b_\beta^Y)$ vanishes because this expression equals $3\varepsilon_i^2h$. The weighted bound on $f_i$ implies that $\varepsilon_i\iota_i^*f_i\to0$ locally. Therefore, passing to the limit, we obtain
\begin{align*}
    &\mathscr P_h\ell+\sigma_\infty D\mathscr E_{0,h}(\dot b)=0,\\
    &\dd A_h(\ell)+c_0\sigma_\infty=0.
\end{align*}
Proposition~\ref{prop:ae-augmented} then implies that $\ell=0$ and $\sigma_\infty=0$.

Consequently, $\|\mathscr P_{b_{\varepsilon_i}}k_i\|_{Y_{\varepsilon_i}}\to0$. Since $\|k_i\|_{X_{\varepsilon_i}}\to1$, estimate~\eqref{eq:glue-schauder} implies that there are points $x_i\in X$ and a constant $c>0$ such that
$$
    w_{\varepsilon_i}(x_i)|k_i(x_i)|_{b_{\varepsilon_i}}\geq c.
$$
The convergence $\ell_i\to0$ excludes points lying in the rescaled image of any fixed compact subset of $Y$.

If $x_i$ stays away from $p$, local compactness implies that a subsequence of $k_i$ converges on compact subsets of $M\setminus\{p\}$ to a nonzero invariant tensor $k_0$ satisfying
\begin{align*}
    &\mathscr P_{g_0}k_0=0,\\
    &|k_0|\leq Cr^{-1}.
\end{align*}
Local elliptic estimates also imply that $|\nabla k_0|\leq Cr^{-2}$. Thus the boundary terms in Green's formula on a sphere of radius $r$ are $O(r)$, so the equation extends distributionally across $p$. Elliptic regularity then implies that $k_0$ is smooth on $M$. Toric nondegeneracy and Lemma~\ref{lem:conical-linear} imply that $k_0=0$, a contradiction.

In the remaining case, $\rho_i=r(x_i)\to0$ and $\varepsilon_i/\rho_i\to0$. We use the coordinate dilations $d_i(z)=\rho_i z$. After passing to a subsequence, on compact subsets of $\mathbb R^4\setminus\{0\}$ we have
\begin{align*}
    &\rho_i^{-2}d_i^*b_{\varepsilon_i}\longrightarrow g_{\mathrm{Euc}},\\
    &\rho_i^{-1}d_i^*k_i\longrightarrow u,
\end{align*}
where $u$ is nonzero at some point of the unit sphere and satisfies
\begin{align*}
    &\mathscr P_0u=0,\\
    &|u(z)|\leq C|z|^{-1}.
\end{align*}
Here $\mathscr P_0$ is the operator at $g_{\mathrm{Euc}}$ with Einstein constant zero. By similar arguments as above we know that $u$ extends smoothly across the origin. Hence $u$ is bounded on $\mathbb R^4$, and interior derivative estimates on balls of arbitrarily large radius imply that $u$ is constant. Its decay at infinity then forces $u=0$, again a contradiction. This proves the uniform estimate.

For each fixed $\varepsilon>0$, the elliptic theory of Section~\ref{sec:moduli} implies that $\mathcal L_\varepsilon$ is Fredholm of index zero. Hence the uniform estimate implies that $\mathcal L_\varepsilon$ is invertible with a uniformly bounded inverse.
\end{proof}

\begin{proof}[Proof of Theorem~\ref{thm:glue}]
Since $A(b_\varepsilon)=\varepsilon^2A_*$, the area component of $\mathcal F_\varepsilon(0,0)$ vanishes. Using the expansions \eqref{glue:compact-expansion} and \eqref{glue:model-expansion} and the cutoff estimates, we obtain
$$
    \|\mathcal F_\varepsilon(0,0)\|_{Y_\varepsilon\times\R}\leq C\bigl(\varepsilon^3+\varepsilon^2r_*^{-1}+r_*^3\bigr)\leq C\varepsilon^{4/3},
$$
where $r_*=\varepsilon^{2/3}$. For $z=(k,\sigma)$, we set
\begin{align*}
    &N_\varepsilon(z)=\|k\|_{X_\varepsilon}+|\sigma|,\\
    &Q_\varepsilon(z)=\mathcal F_\varepsilon(z)-\mathcal F_\varepsilon(0)-\mathcal L_\varepsilon z.
\end{align*}
There are constants $c,C>0$, independent of $\varepsilon$, such that
$$
    \|Q_\varepsilon(z)-Q_\varepsilon(z')\|_{Y_\varepsilon\times\R}\leq C\varepsilon^{-1}\bigl(N_\varepsilon(z)+N_\varepsilon(z')\bigr)N_\varepsilon(z-z')
$$
whenever $N_\varepsilon(z),N_\varepsilon(z')\leq c\varepsilon$. For the metric part of the equation, we obtain this estimate by rescaling the local developments to unit size and using $w_\varepsilon\geq C^{-1}\varepsilon$ and $\beta-\beta_*=\varepsilon^{-1}\sigma$. For the normalized area part of the equation, we rescale to $Y$ and use the fact that area scales by $\varepsilon^2$.

By Proposition~\ref{prop:glueinverse}, for a fixed $R$ large and all $\varepsilon>0$ small, the map
$$
    z\longmapsto-\mathcal L_\varepsilon^{-1}\bigl(\mathcal F_\varepsilon(0)+Q_\varepsilon(z)\bigr)
$$
preserves the ball $N_\varepsilon(z)\leq R\varepsilon^{4/3}$ and has Lipschitz constant at most $CR\varepsilon^{1/3}<1$. Thus we obtain a unique solution $z_\varepsilon=(k_\varepsilon,\sigma_\varepsilon)$ in this ball. We set $\beta_\varepsilon=\beta_*+\varepsilon^{-1}\sigma_\varepsilon$, so
\begin{align*}
    &\|k_\varepsilon\|_{X_\varepsilon}+|\sigma_\varepsilon|\leq C\varepsilon^{4/3},\\
    &|\beta_\varepsilon-\beta_*|\leq C\varepsilon^{1/3}.
\end{align*}
We now set
\begin{align*}
    b=b_{\varepsilon,\beta_\varepsilon},  \qquad g_\varepsilon=b+J_{\beta_\varepsilon}k_\varepsilon.
\end{align*}
By elliptic regularity for the gauged equation we see that $g_\varepsilon$ is smooth in local developments. To remove the gauge term, we set $\omega=\delta_{g_\varepsilon}b$ and write the metric part of the equation as
$$
    \mathscr E_3(g_\varepsilon)-\delta_{g_\varepsilon}^*\omega=0.
$$
By the contracted Bianchi identity, we have $\delta_{g_\varepsilon}\delta_{g_\varepsilon}^*\omega=0$. Integrating by parts, we therefore obtain
\begin{align*}
    &\|\delta_{g_\varepsilon}^*\omega\|_{L^2(g_\varepsilon)}^2=\langle\delta_{g_\varepsilon}\delta_{g_\varepsilon}^*\omega,\omega\rangle_{L^2(g_\varepsilon)}=0,\\
    &\|\omega\|_{L^2(g_\varepsilon)}^2=\langle b,\delta_{g_\varepsilon}^*\omega\rangle_{L^2(g_\varepsilon)}=0.
\end{align*}
Hence we have $\Ric(g_\varepsilon)=3g_\varepsilon$ and $\delta_{g_\varepsilon}b=0$.

On compact subsets of $M\setminus\{p\}$, the reference metric equals $g_0$ for sufficiently small $\varepsilon$, so the weighted estimate implies that $g_\varepsilon\to g_0$ in $C^{2,\alpha}$. We denote by $\iota_\varepsilon$ the gluing embeddings from expanding subsets of $Y$.
For every compact $K\subset Y$, we have
$$
    \|\varepsilon^{-2}\iota_\varepsilon^*k_\varepsilon\|_{C^{2,\alpha}(K,h)}\leq C_K\varepsilon^{-1}\|k_\varepsilon\|_{X_\varepsilon}\leq C_K\varepsilon^{1/3}.
$$
Together with $\beta_\varepsilon\to\beta_*$ and the construction of $b_{\varepsilon,\beta}$, we obtain $\varepsilon^{-2}\iota_\varepsilon^*g_\varepsilon\to h$ in $C^{2,\alpha}$ in local developments.
Elliptic estimates then upgrade both convergences to smooth convergence.

Finally, the nonlinear estimate implies that
$$
    \bigl\|\mathcal L_\varepsilon^{-1}\bigl(D\mathcal F_\varepsilon(z_\varepsilon)-\mathcal L_\varepsilon\bigr)\bigr\|_{\mathrm{op}}\leq C\varepsilon^{1/3}.
$$
Thus $D\mathcal F_\varepsilon(z_\varepsilon)$ is invertible, and the implicit-function theorem implies smooth dependence on $\varepsilon>0$ and local uniqueness in the chosen gauge.
Since $A(g_\varepsilon)=\varepsilon^2A_*$, we obtain the asserted family $g_a$ by setting $a=\varepsilon^2A_*$.
\end{proof}
\begin{remark}
Gluing Einstein metrics in dimension four generally involves finite-dimensional obstructions. In our setting, we remove these obstructions using toric symmetry and variation of the common cone angle.

In Biquard's desingularization of $A_1$ Einstein orbifold singularities using Eguchi-Hanson metrics, the decaying model kernel is three-dimensional, and the leading obstruction requires $\det R^+_{g_0}(p)=0$; see \cite{BiquardI}. Under the standard torus action, only the scaling direction is invariant.

For our $\Ae$ models, Proposition~\ref{prop:ae-deformations} implies that every decaying invariant infinitesimal Ricci-flat deformation at fixed cone angles is a scale deformation modulo equivariant coordinate changes. By the nonzero angle-scale pairing in Lemma~\ref{lem:pair}, varying the common cone angle removes the cokernel obstruction, while prescribing the total area removes the scale kernel. Thus we obtain the augmented inverse in Proposition~\ref{prop:ae-augmented}. Hence, when $g_0$ is torically nondegenerate, the conical gluing construction requires no additional pointwise curvature condition at $p$, and the cone angle is determined as part of the solution.
\end{remark}
\subsection{Degeneration versus gluing}

Consider a degenerating sequence $g_j$ as in Theorem~\ref{thm:compactness}, with compact limit $(M,g_0)$. Throughout this subsection, we assume $g_0$ is torically nondegenerate. We use the notation of Section~\ref{sec:gluing} and set
\begin{equation}\label{scd:area-scale}
    \varepsilon_j^2=\frac{A(g_j)}{A_*}.
\end{equation}
Since the rescaled convergence includes every $\Sigma_\alpha$, this scale is comparable to the curvature scale. The uniqueness  in Appendix~\ref{app:ae-models} implies that
$$
    \varepsilon_j^{-2}g_j\longrightarrow h \quad\text{on compact subsets of }Y, \qquad \beta_j\longrightarrow\beta_*,
$$
with convergence smooth in local developments.

The proof of the following proposition is a mild adaptation of Ozuch's arguments \cite{OzuchI,OzuchII} to our toric conical setting.

\begin{proposition}\label{prop:converse}
For all $j$ large, the metric $g_j$ is equivariantly isometric to $g_{A(g_j)}$ obtained in Theorem~\ref{thm:glue}, by a diffeomorphism preserving each $\Sigma_\alpha$.
\end{proposition}

We denote $|u|_{m,\alpha;w_\varepsilon}$ for the local $C^{m,\alpha}$ norm at scale $w_\varepsilon$, with derivatives of order $\ell$ multiplied by $w_\varepsilon^\ell$, and set
$$
    H_\varepsilon(k)=\sup_X|k|_{2,\alpha;w_\varepsilon}.
$$
Again, we use $b_\varepsilon$ and its local developments to measure these norms. The maps $J_\beta$ identify tensors at the varying cone angles, as in Section~\ref{sec:gluing}.

\begin{lemma}\label{scd:gauge}
After equivariant pullbacks preserving each $\Sigma_\alpha$, we can write
\begin{align*}
    &g_j=b_{\varepsilon_j,\beta_j}+J_{\beta_j}k_j,\\
    &\delta_{g_j}b_{\varepsilon_j,\beta_j}=0,\\
    &H_{\varepsilon_j}(k_j)\longrightarrow0.
\end{align*}
\end{lemma}

\begin{proof}
We denote by $A_j(R,s)$ the annulus between the images of the radius-$R$ sphere in $Y$ and the radius-$s$ sphere about $p$. By Theorem~\ref{thm:compactness}, this annulus is equivariantly diffeomorphic to $S^3\times(0,1)$, with flat limits at every intermediate scale. The Chern-Gauss-Bonnet argument, together with \eqref{glue:compact-expansion} and \eqref{glue:model-expansion}, implies that
\begin{equation}\label{scd:neck-energy}
    \limsup_{j\to\infty}\int_{A_j(R,s)}|\Rm_{g_j}|^2\dvol_{g_j}\leq C(s^4+R^{-4}).
\end{equation}
We apply Ozuch's annular coordinate estimates on this smooth region. The simultaneous half-turn in the two circle directions $A,B$ acts freely, since the only isotropy circles have directions $A$ and $B$. Its quotient therefore has flat model $\mathbb R^4/\{\pm\mathrm{Id}\}$, so we can apply \cite[Propositions~7-8 and Remark~16]{OzuchI} on the quotient and lift the resulting coordinates.

We choose $L_j\to\infty$ sufficiently slowly that $\varepsilon_jL_j^3\to0$ and the convergence on the two limiting regions, through six derivatives and including the boundary terms in \eqref{scd:neck-energy}, has error at most $L_j^{-M}$ for a fixed sufficiently large $M$. Taking $R=L_j$ and $s=L_j^{-1}$, we obtain
\begin{equation}\label{scd:cmc}
    |g_j-g_{\mathrm{Euc}}|_{4,\alpha;r}\leq CL_j^{-2}\left\{\left(\frac{\varepsilon_jL_j}{r}\right)^\gamma+(rL_j)^\gamma\right\}
\end{equation}
on a slightly smaller annulus, for some $\gamma>0$. Its endpoint radii are comparable to $\varepsilon_jL_j$ and $L_j^{-1}$. The prefactor follows from \eqref{scd:neck-energy}, while the decay follows from the cited estimates.

We choose these coordinates equivariantly. We use the foliation by constant-mean-curvature hypersurfaces on the half-turn quotient constructed in \cite[Proposition~7]{OzuchI}. Each leaf is locally unique among nearby normal graphs with the same mean curvature. Since the torus acts by isometries and is connected, we see that it preserves every leaf. Lifting this foliation to $A_j(R,s)$, we obtain a foliation by invariant three-spheres. On the lift of the central leaf used in Ozuch's construction, the spectral subspace near the first nonzero eigenvalue of the round sphere provides an equivariant identification with $S^3$, which we extend along the normal flow. The rod directions determine the two rotation planes. After adjusting by torus rotations, which extend to isometries of $h$ and $g_0$, the transition maps to the converging end coordinates are $O(L_j^{-2})$ close to the identity in scaled norms. We interpolate these transition maps on the end collars, so the coordinates agree there without changing the estimate.

We fix $0<\nu<\min\{\gamma,\tfrac12\}$ and set
$$
    \omega_\varepsilon=w_\varepsilon^\nu+(\varepsilon/w_\varepsilon)^\nu.
$$
We use the weighted norms of \cite[Definitions~2.5-2.6]{OzuchII}, restricted to invariant tensors and vector fields. In our notation, the tensor norm is
$$
    \|u\|_{m,*}=\inf_T\left\{|T|+\sup_X\omega_\varepsilon^{-1}|u-\chi T|_{m,\alpha;w_\varepsilon}\right\},
$$
where $T$ ranges over invariant constant symmetric two-tensors in the annular coordinates. Here $\chi$ equals one on the middle annulus and has transitions at a fixed large radius in $Y$ and a fixed small radius in $M$. We include the trace direction, since no volume constraint is imposed.

Dividing \eqref{scd:cmc} by $\omega_{\varepsilon_j}$, we obtain a bound of $CL_j^{\nu-2}$.
For the reference metric, we have
$$
    \frac{r^2+(\varepsilon_j/r)^2}{\omega_{\varepsilon_j}}\leq C\left\{r^{2-\nu}+(\varepsilon_j/r)^{2-\nu}\right\}.
$$
The choice $\varepsilon_jL_j^3\to0$ places the two end collars in the regions where the reference metric equals $\varepsilon_j^2h$ and $g_0$, respectively. Together with convergence on these regions and the collar identifications, we therefore obtain
\begin{equation}\label{scd:proximity}
    \left\|J_{\beta_j}^{-1}(g_j-b_{\varepsilon_j,\beta_j})\right\|_{4,*}\longrightarrow0.
\end{equation}

We next verify the model estimates needed for the gauge argument. On $(M,g_0)$ and each reference metric, an invariant Killing field vanishes at a torus fixed point, and its derivative is skew-adjoint and commutes with the two normal rotations. Subtracting a torus field, we obtain a Killing field with zero first jet, which must vanish. Thus the invariant Killing fields are precisely the torus fields.

On $(Y,h)$, we claim that there is no nonzero invariant solution of
\begin{align*}
    &\delta_h\delta_h^*v=0,\\
    &|v|=O(r^{1-\nu}).
\end{align*}
Indeed, the order-two asymptotic expansion and scaled elliptic estimates imply that
$$
    \delta_{g_{\mathrm{Euc}}}\delta_{g_{\mathrm{Euc}}}^*v=O(r^{-3-\nu}).
$$
Applying the Euclidean fundamental solution to a cutoff of $v$, we obtain
$$
    v=v_0+O(r^{-1-\nu}),
$$
with corresponding derivative estimates, where $v_0$ is a constant vector. Torus invariance forces $v_0=0$. The boundary term in Green's formula is therefore $O(r^{-2\nu})$ at infinity, while the terms along the conical strata vanish by regularity in local developments. Hence we have $\delta_h^*v=0$, and the Ricci-flat Bochner formula and decay imply that $v=0$.

For vector fields, we use the weight $w_\varepsilon\omega_\varepsilon$, with cutoffs of invariant linear fields included separately, and the image norm uses $\omega_\varepsilon/w_\varepsilon$. The flat annular estimate in \cite[Proposition~4]{OzuchII} applies after the same half-turn quotient. We can therefore follow the proof of \cite[Lemma~3.4]{OzuchII} in local developments. The compact estimate controls the linear coefficients on the annulus; after subtracting these terms, we apply the preceding estimate on $Y$ to the remaining sublinear field. The annular estimate controls the remaining error.
Since the errors from the reference metric are bounded by
$$
    r^{2-\nu}+(\varepsilon/r)^{2-\nu},
$$
we absorb them by choosing the inner model radius large and the outer compact radius small.  Using the Fredholm theory and integration by parts in local developments, we get a uniform inverse for $\delta_b\delta_b^*$ on the complement of the torus fields, where $b=b_{\varepsilon,\beta_*}$.  Changing $\beta$ changes the rescaled coefficients by $O(|\beta-\beta_*|)$, so the estimate also holds for $b=b_{\varepsilon,\beta}$ when $\beta$ is close to $\beta_*$.

We now set $b=b_{\varepsilon_j,\beta_j}$ and $g=g_j$, and apply the fixed-point argument in \cite[proof of Proposition~10]{OzuchII} to
$$
    v\longmapsto\delta_g\bigl((\exp_v)^*b\bigr),
$$
where the exponential map is taken with respect to $b$. At $g=b$, the derivative at zero is $2\delta_b\delta_b^*$, with vector fields identified with one-forms using $b$. By \eqref{scd:proximity}, the derivative remains invertible between the corresponding torus complements. We take the target complement with respect to the $L^2(g)$ pairing. Indeed, for every torus field $K$ and every symmetric tensor $u$ smooth in local developments, we have
$$
    \int_X(\delta_g u)(K)\dvol_g=\int_X\langle u,\delta_g^*K\rangle_g\dvol_g=0.
$$
Thus the image already lies in this complement, and the projected equation is the full equation.

By rescaling the local developments and treating the constant tensor terms separately, we obtain the required nonlinear estimates. Together with the uniform inverse, these allow us to apply the fixed-point argument and find an equivariant diffeomorphism $\psi_j$ such that
\begin{align*}
    &\delta_{g_j}(\psi_j^*b_{\varepsilon_j,\beta_j})=0,\\
    &\left\|J_{\beta_j}^{-1}\bigl(\psi_j^*b_{\varepsilon_j,\beta_j}-b_{\varepsilon_j,\beta_j}\bigr)\right\|_{2,*}\longrightarrow0.
\end{align*}
Using \eqref{scd:proximity} and the estimates for $\psi_j$, after pulling back by $\psi_j^{-1}$ we obtain
\begin{align*}
    &g_j=b_{\varepsilon_j,\beta_j}+J_{\beta_j}k_j,\\
    &\delta_{g_j}b_{\varepsilon_j,\beta_j}=0,\\
    &\|k_j\|_{2,*}\longrightarrow0.
\end{align*}
Since $\omega_\varepsilon$ is uniformly bounded, we have $H_{\varepsilon_j}(k_j)\to0$. Finally, each $\psi_j$ is equivariantly isotopic to the identity and therefore preserves every $\Sigma_\alpha$.
\end{proof}

\begin{proof}[Proof of Proposition~\ref{prop:converse}]
We use the representatives from Lemma~\ref{scd:gauge} and set
\begin{align*}
    &\sigma_j=\varepsilon_j(\beta_j-\beta_*),\\
    &z_j=(k_j,\sigma_j).
\end{align*}
The pullbacks preserve each $\Sigma_\alpha$ and hence the total area. Thus the gauge condition and \eqref{scd:area-scale} imply that $\mathcal F_{\varepsilon_j}(z_j)=0$. The local estimates in the proof of Theorem~\ref{thm:glue} imply that
\begin{equation}\label{scd:tame}
    \bigl\|(D\mathcal F_\varepsilon(z)-\mathcal L_\varepsilon)\dot z\bigr\|_{Y_\varepsilon\times\mathbb R}\leq C\bigl(H_\varepsilon(k)+|\beta-\beta_*|\bigr)\bigl(\|\dot k\|_{X_\varepsilon}+|\dot\sigma|\bigr),
\end{equation}
provided $H_\varepsilon(k)+|\beta-\beta_*|$ is sufficiently small, where $z=(k,\varepsilon(\beta-\beta_*))$ and $\dot z=(\dot k,\dot\sigma)$. Indeed, this quantity controls the difference of the metric coefficients and their first two derivatives after rescaling to scale $w_\varepsilon$. The angle variation is supported where $w_\varepsilon$ is comparable to $\varepsilon$, so the normalization of $\sigma$ implies the same estimate for the angle derivative. For the area component, we obtain the estimate by rescaling to $Y$.

Now we can write
$$
    \mathcal L_{\varepsilon_j}z_j=-\mathcal F_{\varepsilon_j}(0)-\int_0^1\bigl(D\mathcal F_{\varepsilon_j}(tz_j)-\mathcal L_{\varepsilon_j}\bigr)z_j\dd t.
$$
By Lemma~\ref{scd:gauge} and since $\beta_j\to\beta_*$, we know the coefficient in \eqref{scd:tame} tends to zero uniformly along this segment. Using Proposition~\ref{prop:glueinverse} and the reference error estimate from Theorem~\ref{thm:glue}, we therefore obtain
$$
    \|k_j\|_{X_{\varepsilon_j}}+|\sigma_j|\leq C\varepsilon_j^{4/3}+o(1)\bigl(\|k_j\|_{X_{\varepsilon_j}}+|\sigma_j|\bigr).
$$
After absorbing the last term, we have
$$
    \|k_j\|_{X_{\varepsilon_j}}+\varepsilon_j|\beta_j-\beta_*|\leq C\varepsilon_j^{4/3}.
$$
The solution constructed in Theorem~\ref{thm:glue} with total area $\varepsilon_j^2A_*$ satisfies the same bound. Since $H_\varepsilon(k)\leq C\varepsilon^{-1}\|k\|_{X_\varepsilon}$, the coefficient in \eqref{scd:tame} is $O(\varepsilon_j^{1/3})$ along the segment joining the two solutions. Applying  the uniform inverse estimate once more, we see that their difference must vanish. Hence $g_j$ is equivariantly isometric to $g_{A(g_j)}$, as asserted.
\end{proof}
\section{An upper bound on the cone angle}\label{sec:angle-bounds}
We continue to adopt the setting in Section \ref{sec:moduli}. In this section, the main goal is

\begin{theorem}\label{thm:upperangle}
There is a $\beta_{\max}>0$ depending only on the rod structures $\mR$ such that for any $[(\beta, g)]\in \mathcal M$, we have
$$
    \beta\le\beta_{\max}
$$
\end{theorem}

Before giving the proof, we outline the heuristic argument. The angle bound comes from comparing the curvature energy required by the Ricci-flat bubbles with the energy available on the compact manifold. If $\beta\to\infty$, rescaling the angular periods by $\beta^{-1}$ makes the designated cone angles equal to $2\pi$, while the cone angles along the original rods and the total volume tend to zero. The resulting metrics therefore collapse. The internal torus fixed points are smooth, and rescaling around them produces nonflat $\AF$ or $\ALF$ gravitational instantons. In the Gauss-Bonnet formula, the conical corrections account for the disappearing contribution of the original rods. We then see that the Ricci-flat bubbles exhaust the available Weyl curvature energy. The signature formula imposes a further constraint: the twisting of the circle fibration at each $\ALF$ end contributes a sign, and these contributions must balance the signatures of the limits. Our choice of rod structures makes this balance impossible. For the alternating sequence, the odd number of internal fixed points forces an odd number of $\ALF$ ends. Their signs cannot sum to zero, as required by the vanishing limiting signature integral. For the other chosen sequences, the intersection forms lead to analogous contradictions. Thus splitting into several limits cannot remove the obstruction. In the following we make this argument precise.

\

Given $[(\beta, g)]\in \mathcal M$. To prove the above theorem we always assume $\beta> 1$. As in Remark \ref{rmk:q=1}, the upper bound is automatic if $q=1$, so we will assume $q\geq 2$ throughout this section. To apply Riemannian convergence theory, it is helpful to change the angular periods by defining
$$
    \widehat g=\underline g+\beta^{-2}G_{ab} \mathrm d\varphi_a\mathrm d\varphi_b.
$$
Then $\widehat g$ is again a toric conical Einstein metric, and each cone angle is multiplied by $\beta^{-1}$. Thus the rods $\mI_1, \cdots, \mI_q$ become smooth, while the remaining rods have cone angle $2\pi\beta^{-1}<2\pi$. Moreover we have
\begin{equation}\label{eq:infinite-integral-scaling}
    \operatorname{Vol}_{\widehat g}(X) =\beta^{-2}\operatorname{Vol}_g(X), \qquad I_\pm(\widehat g)=\beta^{-2}I_\pm(g).
\end{equation}
Suppose $\beta_j\to\infty$. Then by \eqref{geom:volumeupper} and \eqref{eq:alternating-characteristic} we have
\begin{align}
    \operatorname{Vol}_{\widehat g_j}(X) &\leq \frac{8\pi^2}{3\beta_j^2}\rightarrow0, \label{eq:infinite-volume}\\
    \chi_c(\widehat g_j)                 &=q-1+\frac{2}{\beta_j}+\frac{m_0-1}{\beta_j^2}, \nonumber\\
    \tau_c(\widehat g_j)                 &=\frac{\sigma}{3} +\frac{1}{\beta_j^2} \left(\tau(X)-\frac{\sigma}{3}\right), \label{eq:infinite-characteristic}\\
    I_\pm(\widehat g_j)                  &\longrightarrow 2\pi^2\bigl(2(q-1)\pm\sigma\bigr). \label{eq:infinite-budgets}
\end{align}
In particular, the $L^2$ curvature energies of $\widehat g_j$ are uniformly bounded.

Let $\mathcal S_0$ be the union of the spheres corresponding to the original rods in $\mathfrak R_0$, and denote the internal fixed points by $p_j^1,\ldots,p_j^{q-1}$. These points are smooth for $\widehat g_j$. Write $\omega_4$ for the volume of the Euclidean unit ball and, for each $1\leq a\leq q-1$, choose $r_j^a>0$ such that
\begin{equation}\label{eq:alternating-volume-radius}
    \operatorname{Vol}_{\widehat g_j} B_{\widehat g_j}(p_j^a,r_j^a) =\frac34\omega_4(r_j^a)^4.
\end{equation}
Such a radius exists because Bishop-Gromov volume comparison holds for $\widehat g_j$ (since the cone angles are all at most $2\pi)$ and the total volume goes to zero. Moreover, \eqref{eq:infinite-volume} implies $r_j^a\leq C\beta_j^{-1/2}\to0$. Set
$$
    h_j^a=(r_j^a)^{-2}\widehat g_j.
$$

\begin{lemma}\label{lem:alternating-smooth-limits}
After passing to a subsequence, each $(X,h_j^a,p_j^a)$ converges smoothly and equivariantly to a complete noncompact nonflat toric Ricci-flat four-manifold with finite curvature energy.
\end{lemma}

\begin{proof}
We fix $a$ and write $p_j=p_j^a$, $r_j=r_j^a$, and $h_j=h_j^a$. At $z\in\mathcal S_0$, the volume density is at most $\beta_j^{-1}$. By the normalization \eqref{eq:alternating-volume-radius} and Bishop-Gromov comparison, we have
$$
    \frac34\omega_4r_j^4
    \leq \operatorname{Vol}_{\widehat g_j}B_{\widehat g_j}(z,d_{\widehat g_j}(p_j,z)+r_j)
    \leq \frac{\omega_4}{\beta_j}\bigl(d_{\widehat g_j}(p_j,z)+r_j\bigr)^4.
$$
Hence
\begin{equation}\label{eq:alternating-old-escape}
    d_{h_j}(p_j,\mathcal S_0)
    \geq \left(\frac{3\beta_j}{4}\right)^{1/4}-1
    \longrightarrow\infty.
\end{equation}
Thus every fixed-radius ball about $p_j$ is smooth for large $j$. Since $B_{h_j}(p_j,1)\subset B_{h_j}(x,R+1)$ for $x\in B_{h_j}(p_j,R)$, another application of volume comparison implies
\begin{equation}\label{eq:alternating-local-noncollapse}
    \operatorname{Vol}_{h_j}B_{h_j}(x,t)
    \geq \frac{3\omega_4}{4(R+1)^4}t^4,
    \qquad x\in B_{h_j}(p_j,R),\quad 0<t\leq R+1.
\end{equation}

Suppose that curvature is unbounded on some fixed-radius ball. As in the proof of Theorem~\ref{thm:compactness}, point selection and rescaling by curvature produce a complete smooth nonflat Ricci-flat limit. The bound \eqref{eq:alternating-local-noncollapse} implies Euclidean volume growth in this further limit, and the uniform curvature energy bound implies finite curvature energy. It follows from Proposition~\ref{geom:vertexcurvature} that there is a smooth torus fixed point at bounded distance. We therefore obtain an effective limiting torus action with rod directions in $\{A,B,C\}$. By \cite[Theorem~1.5]{BKN}, the limit is $\ALE$, contradicting Lemma~\ref{deg:smooth-ae-rigidity}.

We conclude that curvature is uniformly bounded on every fixed-radius ball. Together with \eqref{eq:alternating-local-noncollapse}, smooth equivariant compactness yields a complete toric limit $(Y,h,p_\infty)$. Since $\Ric_{h_j}=3r_j^2h_j$ and $r_j\to0$, the limit is Ricci-flat, and its curvature energy is finite by lower semicontinuity. By \eqref{eq:alternating-old-escape},  $Y$ is noncompact.

The limiting torus action fixes $p_\infty$. If $h$ were flat, this torus action would force it to be Euclidean, which is a contradiction to the fact that 

$$
    \operatorname{Vol}_hB_h(p_\infty,1)=\frac34\omega_4.
$$
Finally, since there are only finitely many indices $a$, we may choose a common subsequence.
\end{proof}

Equivariant convergence preserves stabilizers and identifies their fixed sets on compact subsets. Hence the torus fixed points in the limit are precisely the limits of internal fixed points that remain at bounded rescaled distance. Using the local uniqueness of the surface fixed by each stabilizer circle and the equivariant charts at torus fixed points, we identify each finite limiting rod with a single designated rod, with matching endpoints. Thus every finite chain in the limit corresponds to a consecutive segment of the designated chain.

Let $(Y,h)$ be a limit from Lemma~\ref{lem:alternating-smooth-limits}, and write
$$
    (B_\infty,\underline h)=Y/\mathbb T^2, \qquad \mathcal V_\infty=Y^{\mathbb T^2}.
$$
The set $\mathcal V_\infty$ is finite. Moreover, the rods are smooth and the rod vectors belong to $\{A,B,C\}$.

\begin{proposition}
\label{asym:prop:bubbleproper}
The limit $(Y,h)$ is a proper toric gravitational instanton of $\AF$ or $\ALF$ type. The rods of $Y/\mathbb T^2$, in boundary order, correspond to a consecutive part of the designated chain.
\end{proposition}
\begin{proof}
As in Corollary~\ref{geom:cat}, the quotient $B_\infty=Y/\mathbb T^2$ is complete and has simply connected interior with no cut locus. Since $K_{\underline h}\leq0$, the same comparison argument implies that $B_\infty$ is $\operatorname{CAT}(0)$. Its boundary is nonempty because the basepoint is fixed by the torus. By Lemma~\ref{lem:polar}, the horizontal distribution is integrable. Moreover, the torus fixed-point set is finite, so \eqref{geom:vertexcurvatureeq} implies quadratic curvature decay at infinity. Lemma~\ref{geom:proper-ricciflat} therefore implies that $(Y,h)$ is proper and identifies $B_\infty$ with the closed half-plane.

By \eqref{geom:rodidentity}, each rod in $(X/\mathbb T^2,\underline h_j^a)$ has length at most $\pi/r_j^a$. Since this quotient is $\operatorname{CAT}((r_j^a)^2)$ by Corollary~\ref{geom:cat}, every rod is minimizing. Thus a rod cannot escape to infinity between two points at bounded rescaled distance, so each contributes at most one connected piece to $\partial B_\infty$. Together with the correspondence of rods and endpoints established above, we conclude that the limiting rod sequence is a consecutive segment of the designated chain. Theorem~\ref{thm:appendix classification of asymptotics} implies that $(Y,h)$ has $\ALE$, $\ALF$ or $\AF$ asymptotics. In the $\ALE$ case, Lemma~\ref{deg:smooth-ae-rigidity} would force $(Y,h)$ to be Euclidean, contrary to its nonflatness. Hence $(Y,h)$ is of $\AF$ or $\ALF$ type.
\end{proof}

We denote
$$
    I_\pm(Y)=\int_Y|W_h^\pm|^2 \dvol_h.
$$
Let $v$ be the number of torus fixed points in $Y$. Set $\varepsilon(Y)=0$ for an AF end; for an ALF end, let $\varepsilon(Y)=\pm1$ be the oriented Euler number of the circle bundle at infinity.

\begin{lemma}\label{lem:alternating-bubble-energies}
We have $\chi(Y)=v$, and $\tau(Y)$ is the signature of the intersection matrix of the finite rod spheres. Moreover,
\begin{equation}\label{eq:limiting-bubble-energies}
    I_\pm(Y)=4\pi^2 \left(v\pm\left(\frac32\tau(Y)+\varepsilon(Y)\right)\right).
\end{equation}
If only two rod directions occur, then $\tau(Y)=0$ and
\begin{equation}\label{eq:alternating-bubble-energies}
    I_\pm(Y)=4\pi^2\bigl(v\pm\varepsilon(Y)\bigr), \qquad \varepsilon(Y)= \begin{cases}
        0,    &v\text{ even},\\
        \pm1, &v\text{ odd}.
    \end{cases}
\end{equation}
\end{lemma}

\begin{proof}
By the rod construction, we have $\chi(Y)=v$, and the $v-1$ finite rod spheres form a basis of $H_2(Y;\mathbb Z)$; see \cite[Section~4]{Nilsson}. Hence $\tau(Y)$ is the signature of their intersection matrix.  Equal extreme directions give an AF end; distinct directions form a lattice basis and give an ALF end with Euler number $\pm1$.

The AF and ALF Chern-Weil formulas give
\begin{align}
    I_+(Y)+I_-(Y)                 &=8\pi^2v, \label{eq:af-alf-gb}\\
    \frac{I_+(Y)-I_-(Y)}{12\pi^2} &=\tau(Y)+\frac23\varepsilon(Y), \label{eq:af-alf-signature}
\end{align}
see \cite[Section~6.2, equations~(6.19)-(6.25)]{ChenLi} and \cite[Section~4.1]{Nilsson}. For a screw AF end, the signature boundary term also vanishes: equatorial reflection commutes with the screw rotation and reverses the orientation of the model cross-section. Adding and subtracting the two identities gives \eqref{eq:limiting-bubble-energies}.

With only two directions, the rods alternate and every finite sphere has self-intersection zero. Changing the signs of alternate sphere classes negates the intersection matrix, so $\tau(Y)=0$. The extreme directions agree exactly when $v$ is even, proving \eqref{eq:alternating-bubble-energies}.
\end{proof}

\begin{proof}[Proof of Theorem~\ref{thm:upperangle}]
We recall that $q\ge2$ and consider the sequence $\beta_j\to\infty$ above. After passing to a subsequence and relabelling the internal fixed points, we may assume that $r_j^1\geq\cdots\geq r_j^{q-1}$ for every $j$. We may also assume that $d_{\widehat g_j}(p_j^a,p_j^b)/r_j^c$ converges in $[0,\infty]$ for every $1\leq a,b,c\leq q-1$. At each step, we choose a remaining internal fixed point with largest radius $r_j^a$, take its pointed limit, and remove all internal fixed points at bounded distance on that scale. If $p_j^a$ is chosen before $p_j^b$, then $r_j^b\leq r_j^a$ and $d_{\widehat g_j}(p_j^a,p_j^b)/r_j^a\to\infty$. Hence, for every fixed $R>0$, the balls $B_{\widehat g_j}(p_j^a,Rr_j^a)$ centered at the chosen points are pairwise disjoint for all sufficiently large $j$. The same separation prevents a fixed point assigned to one limit from appearing in another. Since every internal fixed point is assigned at some step, the fixed-point correspondence established above shows that each occurs in exactly one selected limit.

We denote these limits by $Y_1,\ldots,Y_t$ and their numbers of torus fixed points by $v_1,\ldots,v_t$.
Thus
\begin{equation}\label{eq:vertex-partition}
    \sum_{\alpha=1}^t v_\alpha=q-1.
\end{equation}
By smooth convergence, scale invariance of the Weyl energies, and disjointness, we obtain after exhaustion
\begin{equation}\label{eq:weyl-semicontinuity}
    \sum_{\alpha=1}^t I_\pm(Y_\alpha)
    \leq \lim_{j\to\infty}I_\pm(\widehat g_j)
    =2\pi^2\bigl(2(q-1)\pm\sigma\bigr).
\end{equation}
Substituting \eqref{eq:limiting-bubble-energies} and using \eqref{eq:vertex-partition}, the $+$ and $-$ inequalities bound the following sum above and below by $\sigma$, respectively. Therefore,
\begin{equation}\label{eq:bubble-signature-balance}
    \sum_{\alpha=1}^t\bigl(3\tau(Y_\alpha)+2\varepsilon(Y_\alpha)\bigr)=\sigma.
\end{equation}

For the alternating sequence $W=(B,(A,B)^n,A)$, every limit uses only the directions $A,B$, so $\tau(Y_\alpha)=0$. Since $\sigma=0$, we obtain $\sum_\alpha\varepsilon(Y_\alpha)=0$ from \eqref{eq:bubble-signature-balance}.
On the other hand, \eqref{eq:alternating-bubble-energies} implies
$$
    \sum_\alpha\varepsilon(Y_\alpha)
    \equiv\sum_\alpha v_\alpha
    =2n-1\pmod2,
$$
which is a contradiction.

For $W=(B,A,C,A)$ or $(B,C,A,C,A)$, each limit uses only the directions $A,C$. For $W=(B,C,B,A)$ or $(B,C,B,C,A)$, every limit uses only the directions $B,C$. In all four cases, we have $\tau(Y_\alpha)=0$ and $\sigma=\pm1$. Thus \eqref{eq:bubble-signature-balance} has an even left-hand side and an odd right-hand side, again a contradiction.

Finally, we consider $W=(B,C,A,B,C,A)$. Here $q-1=3$, $\sigma=-4$, and each limiting rod sequence is a consecutive segment of $(C,A,B,C)$. The finite rod spheres inherit self-intersection $-1$ from the corresponding spheres in $X$. Hence a limit with two or three torus fixed points has intersection matrix, respectively,
$$
    (-1),\qquad
    \begin{pmatrix}
        -1&1\\
        1&-1
    \end{pmatrix},
$$
and in either case its signature is $-1$. If such a limit occurs, \eqref{eq:vertex-partition} implies that $\sum_\alpha\tau(Y_\alpha)=-1$. The left-hand side of \eqref{eq:bubble-signature-balance} is then odd, whereas $\sigma=-4$.

Otherwise, there are three limits, each with one torus fixed point. Each has signature zero and $\varepsilon(Y_\alpha)=\pm1$. Equation~\eqref{eq:bubble-signature-balance} would then require three signs to sum to $-2$, which is impossible.
\end{proof}

\section{Rigidity of toric Einstein metrics with \texorpdfstring{$b_2^+\leq1$}{b2+ <= 1}}\label{sec:weyl-rigidity}

Throughout this section, $g$ is a smooth toric Einstein metric on the compact four-manifold $X_{\mR}$, with $\Ric(g)=\lambda g$ and $\lambda>0$. We use the orientation convention of Section~\ref{sec:rod structure}. We first give a topological criterion for $g$ to be algebraically special; its proof is given in Section~\ref{subsec:toric-weyl-rigidity}.

\begin{theorem}\label{thm:toric-weyl-rigidity}
Let $g$ be a toric Einstein metric on $X_{\mR}$. Then
\begin{itemize}
\item $b_2^+(X_{\mR})=0$ if and only if $g$ is anti-self-dual, that is, $W_g^+=0$.
\item $b_2^+(X_{\mR})=1$ if and only if $g$ is conformally K\"ahler with respect to a complex structure compatible with the orientation of $X_{\mR}$.
\end{itemize}
\end{theorem}

As in the case of toric gravitational instantons studied in \cite{LiSun}, one can also formulate the criterion in terms of the \emph{degree} of the rod structure. Recall that the rod vector on $\mI_i$ is defined up to sign, so it determines a projective direction $[\mv_i]\in\mathbb{RP}^1$.

\begin{definition}\label{def:compact-rod-degree}
Let $\mv_1,\ldots,\mv_m$ be the rod vectors in positive boundary order, with $\mv_{m+1}=\mv_1$. Choose real angular lifts $\theta_1,\ldots,\theta_{m+1}\in\R$ of these projective directions successively so that
$$
    0<\theta_i-\theta_{i+1}<\pi,\qquad i=1,\ldots,m.
$$
The \emph{degree} of $\mR$ is
$$
    d(\mR)=\frac{\theta_1-\theta_{m+1}}{\pi}-1.
$$
Thus $d(\mR)+1$ counts the clockwise turns made by joining consecutive directions in $\mathbb{RP}^1$.
\end{definition}

The degree is a nonnegative integer, independent of the initial angular lift and invariant under orientation-preserving changes of lattice basis. With our orientation conventions,
\begin{equation}\label{eq:compact-rod-degree}
    d(\mR)=b_2^+(X_{\mR}).
\end{equation}
To prove this identity, choose a circle subgroup generated by a primitive lattice vector $\xi$ whose direction differs from every rod direction. Its fixed points are precisely the $m$ torus fixed points. At the vertex between $\mI_i$ and $\mI_{i+1}$, choose signed representatives so that $(\mv_i,\mv_{i+1})$ is a positive lattice basis, and write
$$
    \xi=a_i\mv_i+b_i\mv_{i+1},\qquad a_i,b_i\in\mathbb Z\setminus\{0\}.
$$
In the corresponding local complex coordinates, the circle action is
$$
    t\cdot(z_1,z_2)=(t^{a_i}z_1,t^{b_i}z_2),\qquad t\in\mathbb S^1.
$$
Thus $a_i,b_i$ are the two rotation weights. To compare orientations, write $z_k=r_ke^{\sqrt{-1}\psi_k}$. The rods $\mI_i$ and $\mI_{i+1}$ are given by $r_1=0$ and $r_2=0$, respectively, so their positive boundary order makes $\dd r_1\wedge\dd r_2$ positive on the orbit space. The positive lattice basis makes $\dd\psi_1\wedge\dd\psi_2$ positive on the torus. Our orientation of $X_{\mR}$ is therefore represented by $\dd r_1\wedge\dd r_2\wedge\dd\psi_1\wedge\dd\psi_2$, whereas the complex orientation is represented by
$$
    (r_1\dd r_1\wedge\dd\psi_1)\wedge(r_2\dd r_2\wedge\dd\psi_2)=-r_1r_2\dd r_1\wedge\dd r_2\wedge\dd\psi_1\wedge\dd\psi_2.
$$
Hence these local complex coordinates induce the orientation opposite to that of $X_{\mR}$.

We apply the isolated-fixed-point form of the Atiyah-Singer equivariant signature theorem \cite[Theorem~6.12, p.~582]{AtiyahSingerIII}. Since the circle is connected, its action on cohomology is trivial, so the equivariant signature equals the ordinary signature $\tau(X_{\mR})$. By the  fixed-point formula we get that for generic $t\in\mathbb S^1$,
$$
    \tau(X_{\mR})=-\sum_{i=1}^m\frac{1+t^{a_i}}{1-t^{a_i}}\frac{1+t^{b_i}}{1-t^{b_i}}.
$$
The minus sign in each summand accounts for the opposite local complex orientation. This is an identity of rational functions of $t$, so we may take $t\to0$. Since
$$
    \lim_{t\to0^+}\frac{1+t^a}{1-t^a}=\operatorname{sgn}(a)\qquad(a\neq0),
$$
we obtain
$$
    \tau(X_{\mR})=-\sum_{i=1}^m\operatorname{sgn}(a_ib_i).
$$
Thus the contribution of the vertex between $\mI_i$ and $\mI_{i+1}$ to this signature sum is $-\operatorname{sgn}(a_ib_i)$: it is $+1$ when the two rotation weights have opposite signs and $-1$ when they have the same sign.

This contribution is $+1$ precisely when the chosen circle direction lies in the clockwise projective arc from $\mv_i$ to $\mv_{i+1}$, and is $-1$ otherwise. There are exactly $d(\mR)+1$ such arcs containing this direction, counted with multiplicity. Therefore
$$
    \tau(X_{\mR})=2d(\mR)+2-m.
$$
Together with $\chi(X_{\mR})=m$ and $b_1(X_{\mR})=0$, this proves \eqref{eq:compact-rod-degree}. Using $3\tau(X_{\mR})=\sum_i[\Sigma_i]^2$, we also obtain
$$
    d(\mR)=\frac{m-2}{2}+\frac16\sum_{i=1}^m[\Sigma_i]^2.
$$
Thus the two cases in Theorem~\ref{thm:toric-weyl-rigidity} correspond to degrees zero and one, respectively. Reversing orientation changes the degree to $m-2-d(\mR)=b_2^-(X_{\mR})$ and gives the corresponding criteria for $W_g^-$.

Together with the classification results of Hitchin \cite{Hitchin} and LeBrun \cite{LeBrun}, Theorem~\ref{thm:toric-weyl-rigidity} gives the following classification result. In Section~\ref{subsec:hermitian-uniqueness} we will provide an alternative direct proof. 

\begin{corollary}\label{cor:toric-small-betti-classification}
If $b_2^+(X_{\mR})\leq1$, then, up to scaling and isometry, a toric Einstein metric $g$ on $X_{\mR}$ is one of the following:
\begin{itemize}
\item the round metric on $\mathbb S^4$;
\item the Fubini-Study metric on $\CP^2$;
\item the round product metric on $\mathbb S^2\times\mathbb S^2$;
\item the Page metric on $\CP^2\#\overline{\CP^2}$;
\item the Chen-LeBrun-Weber metric on $\CP^2\#2\overline{\CP^2}$;
\item the K\"ahler-Einstein metric on $\CP^2\#3\overline{\CP^2}$.
\end{itemize}
\end{corollary}

In Section~\ref{subsec:hermitian-uniqueness}, we also prove infinitesimal rigidity for algebraically special toric Einstein metrics. This essentially follows from the work of Andersson-Araneda \cite{AA}.

\begin{theorem}\label{thm:toric-small-betti-nondegeneracy}
All the metrics listed in Corollary~\ref{cor:toric-small-betti-classification} are torically nondegenerate.
\end{theorem}
\subsection{Proof of Theorem \ref{thm:toric-weyl-rigidity}}\label{subsec:toric-weyl-rigidity}

We use the notation of Sections \ref{sec:rod structure} and \ref{sec:a priori estimates}. Recall that
$$
    \begin{aligned}
        &u=\log\rho,\qquad U=\frac12|\dd u|_{\bg}^2,\qquad K=K_{\bg},\\
        &T=\Phi^*g_{\mathcal H},\qquad e=\frac12\tr_{\bg}T=K+\frac U2.
    \end{aligned}
$$
We identify the upper half-plane with the space of positive-definite symmetric matrices of determinant one by
$$
    (X,Y)\longmapsto
    \begin{pmatrix}
        X+X^{-1}Y^2&X^{-1}Y\\
        X^{-1}Y&X^{-1}
    \end{pmatrix},
    \qquad X>0.
$$
In these coordinates, the Riemannian  metric  can be written as
$$
    g_{\mathcal H}=\frac14\tr(\Phi^{-1}\dd\Phi)^2=\frac{\dd X^2+\dd Y^2}{2X^2},
$$
and $\Phi=\rho^{-1}G$ has coordinates
$$
    X=\frac{\rho}{G_{22}}, \qquad Y=\frac{G_{12}}{G_{22}}.
$$
Choose a positively oriented orthonormal frame $e_1,e_2,e_3,e_4$ locally on $X_{\mR}$, with $e_1,e_2$ horizontal and $e_3,e_4$ vertical. Set
$$
    J=R_{1234}=\langle R(e_1,e_2)e_4,e_3\rangle.
$$
\begin{lemma}\label{lem:weyl-jacobian}
The curvature component $J$ is the negative Jacobian of $\Phi$. More precisely,
\begin{equation}\label{eq:weyl-jacobian}
    J\dd A_{\bg} =\dd\left(\frac{G_{22}}{2\rho} \dd\frac{G_{12}}{G_{22}}\right) =-\frac{\dd X\wedge\dd Y}{2X^2}.
\end{equation}
\end{lemma}
\begin{proof}
We denote by $\xi_1,\xi_2$ the Killing vector fields corresponding to a positively oriented basis of the torus lattice.
We choose the orthonormal vertical frame
$$
    e_3=\frac{\sqrt{G_{22}}}{\rho}(\xi_1-Y\xi_2),\qquad e_4=\frac{\xi_2}{\sqrt{G_{22}}}.
$$
For a horizontal vector $v$, we have $2\langle\nabla_v\xi_a,\xi_b\rangle=v(G_{ab})$.
Thus the connection one-form $\alpha$ of this frame in horizontal directions satisfies
$$
    \alpha=\langle\nabla e_3,e_4\rangle
    =-\frac{\dd G_{12}-Y\dd G_{22}}{2\rho}
    =-\frac{\dd Y}{2X}.
$$
Since covariant differentiation in horizontal directions preserves the vertical distribution, we have
$$
    \nabla_v e_3=\alpha(v)e_4,\qquad \nabla_v e_4=-\alpha(v)e_3.
$$
With our convention $R_{ijkl}=\langle R(e_i,e_j)e_l,e_k\rangle$, the curvature equation therefore implies that
$$
    J\dd A_{\bg}=-\dd\alpha=-\frac{\dd X\wedge\dd Y}{2X^2}.
$$
\end{proof}
The integral of $J$ is determined by the signature of $X_{\mR}$:
\begin{lemma}\label{lem:weyl-signature}
We have
\begin{equation}\label{eq:weyl-signature}
    \int_BJ\dd A_{\bg}=-\frac\pi2\tau(X_{\mR}),
\end{equation}
and
\begin{equation}\label{eq:weyl-signature2}
    \int_B(K+J)\dd A_{\bg}=\pi(1-b_2^+(X_{\mR})).
\end{equation}
\end{lemma}
\begin{proof}
We choose a positively oriented basis of the torus lattice, with associated Killing vector fields $\xi_1,\xi_2$, such that the circle generated by $\xi_2$ differs from every rod circle.
For the connection form
$$
    \alpha=-\frac{G_{22}}{2\rho}\dd\frac{G_{12}}{G_{22}},
$$
equation~\eqref{eq:weyl-jacobian} reads $J\dd A_{\bg}=-\dd\alpha$.
Along each open rod, $G_{12}/G_{22}$ has constant boundary value with an error of order $r^2$, whereas $\rho/G_{22}$ vanishes linearly.
Thus the integrals of $\alpha$ over curves approaching the open rods tend to zero, with neighborhoods of the vertices excluded.

At a turning point, we choose radial coordinates $x,y$ with $\dd x\wedge\dd y$ positive and choose signs for the adjacent primitive rod vectors so that they form a positive lattice basis.
We express $\xi_1,\xi_2$ in this local basis as $(a,b)$ and $(c,d)$, respectively, where $ad-bc=1$ and $cd\neq0$.
The leading terms of $G_{12}$ and $G_{22}$ are $acx^2+bdy^2$ and $c^2x^2+d^2y^2$, respectively.
Consequently, on the arc $x=r\cos t$, $y=r\sin t$, we have
$$
    \left.\alpha\right|_{\mathrm{arc}}
    =\left(\frac{cd}{c^2\cos^2t+d^2\sin^2t}+O(r^2)\right)\dd t,
$$
where the error is uniform up to the endpoints by the corner expansions.
The boundary of the truncated domain traverses this arc from $\pi/2$ to $0$, so the integral of $\alpha$ over the arc tends to $-\pi\operatorname{sgn}(cd)/2$.

We apply Stokes' theorem to $J\dd A_{\bg}=-\dd\alpha$ on domains truncated near the rods and vertices.
We first let the distances from the open rods tend to zero with the corner radii fixed, and then let the corner radii tend to zero.
Each vertex therefore contributes $\pi\operatorname{sgn}(cd)/2$ to $\int_BJ\dd A_{\bg}$.
By the fixed-point expression for the signature established above, the same vertex contributes $-\operatorname{sgn}(cd)$ to $\tau(X_{\mR})$, so we obtain \eqref{eq:weyl-signature}.

Since the rods are geodesic and the $\chi(X_{\mR})$ corner angles are all $\pi/2$, the Gauss--Bonnet formula implies that
$$
    \int_BK\dd A_{\bg}=2\pi-\frac\pi2\chi(X_{\mR}).
$$
Combining this formula with \eqref{eq:weyl-signature} and the identity $\chi(X_{\mR})+\tau(X_{\mR})=2+2b_2^+(X_{\mR})$, we obtain \eqref{eq:weyl-signature2}.
\end{proof}

Lemma \ref{lem:weyl-jacobian} gives
$$
    J^2=\det_{\bg}T, \qquad |J|\leq\tfrac12\tr_{\bg}T=e.
$$
Consider the orthonormal basis of $\Lambda^+$ given by
$$
    \omega_0=\frac{e^{12}+e^{34}}{\sqrt2},\qquad \omega_1=\frac{e^{13}-e^{24}}{\sqrt2},\qquad \omega_2=\frac{e^{14}+e^{23}}{\sqrt2}.
$$
In this basis, the local reflection $\varphi\mapsto-\varphi$ acts as $\operatorname{diag}(1,-1,-1)$. Since $W^+$ commutes with this reflection, $\omega_0$ is an eigenform of $W^+$. The Einstein condition gives $R_{3434}=R_{1212}=K$, and hence
$$
    \langle\Rm_g\omega_0,\omega_0\rangle=\frac12(R_{1212}+R_{3434}+2R_{1234})=K+J.
$$
On the other hand,
$$
    \left.\Rm_g\right|_{\Lambda^+}=W^++\frac{\Scal_g}{12}\Id=W^++\frac\lambda3\Id.
$$
Thus, defining real coefficient functions $a,b$ in the chosen frame, we may write
\begin{equation}\label{eq:weyl-matrix}
    W^+=\begin{pmatrix} -w&0&0\\
        0&w/2+a&b\\
        0&b&w/2-a
    \end{pmatrix},
\end{equation}
where
$$
    w=\frac{\lambda}{3}-K-J, \qquad q=a+\sqrt{-1} b.
$$
The coefficient $q$ depends on the frame, but its squared norm
\begin{equation}\label{eq:weyl-q-norm}
    |q|^2=\frac12\left(|W^+|^2-\frac32w^2\right)
\end{equation}
is well defined.

The boundary expansions of Lemma~\ref{lem:rodasymptotics} show that $K,J,w,|q|^2$ are smooth and even in the normal variables at rods and turning points. On an open rod, choose the limiting frame so that $e_1,e_3$ span the normal plane to the fixed surface. The stabilizer circle fixes $\omega_1$ and rotates $\operatorname{span}\{\omega_0,\omega_2\}$. Since $W^+$ is self-adjoint and commutes with this action, it acts as $-w\Id$ on this two-plane. By the trace freeness we see that the eigenvalues are given by $-w,-w,2w$, so
\begin{equation}\label{eq:weyl-boundary}
    |q|^2=\frac12\left(6w^2-\frac32w^2\right) =\frac94w^2 \qquad\text{on }\partial B.
\end{equation}
The tensor $W^+$ satisfies the contracted Bianchi identity
\begin{equation}\label{eqn:contracted Bianchi}
    \delta_gW^+=0
\end{equation}
and the Weyl Weitzenb\"ock formula \cite[Section~2]{Derdzinski},
\begin{equation}\label{eqn:Weyl Weitzenbock}
    \Delta_gW^+=2\lambda W^+-6(W^+)^2+2|W^+|^2\Id.
\end{equation}
Here $\Delta_g=-\nabla_g^*\nabla_g$ on tensors. These identities give the following equations for $w$ and $q$.

\begin{lemma}\label{lem:weyl-equations}
On $B^\circ$, we have
\begin{align}
    \left|\dd w+\frac32w\dd u\right|_{\bg}^2 &=4(e+J)|q|^2,\label{eq:weyl-bianchi-norm}\\
    \Delta_{\bg}w+3\langle\dd u,\dd w\rangle_{\bg} &=4\lambda w-3w^2-4|q|^2.\label{eq:weyl-scalar}
\end{align}
\end{lemma}

\begin{proof}
Use a $\T^2$-invariant orthonormal frame and write
$$
    \bigl(\langle\nabla_{e_\alpha}e_i,e_\beta\rangle\bigr) _{\alpha,\beta=3,4} =\begin{pmatrix}u_i/2+p_i&t_i\\t_i&u_i/2-p_i\end{pmatrix}, \qquad i=1,2.
$$
Here $e_1,e_2$ are horizontal, $e_3,e_4$ are vertical, and $u_i=e_i(u)$. The functions $p_i,t_i$ describe the trace-free part of the matrix. Set
$$
    \zeta=2\bigl(p_1-t_2+\sqrt{-1}(p_2+t_1)\bigr).
$$
We first establish the identities
$$
    e=p_1^2+p_2^2+t_1^2+t_2^2,\qquad J=-2(p_1t_2-p_2t_1),\qquad |\zeta|^2=4(e+J).
$$
Take constant linear combinations of the angular Killing fields that agree with $e_3,e_4$ at the chosen point. Their Gram matrix satisfies $G=\Id$ and $\rho=1$ there, so
$$
    \begin{aligned}
        e_i(G)&=
        \begin{pmatrix}
            u_i+2p_i&2t_i\\
            2t_i&u_i-2p_i
        \end{pmatrix},\\
        e_i(\Phi)&=e_i(G)-u_i\Id =2\begin{pmatrix} p_i&t_i\\
            t_i&-p_i
        \end{pmatrix}.
    \end{aligned}
$$
Hence by the definition of the target metric we see that 
$$
    e=\frac18\sum_{i=1}^2 \tr\bigl((\Phi^{-1}e_i(\Phi))^2\bigr) =p_1^2+p_2^2+t_1^2+t_2^2.
$$
In the upper half-plane coordinates, we obtain
$$
    e_i(X)=e_i(\rho)-e_i(G_{22})=2p_i,\qquad e_i(Y)=2t_i.
$$
Since $X=1$ at the chosen point, \eqref{eq:weyl-jacobian} gives
$$
    J=-\frac12(\dd X\wedge\dd Y)(e_1,e_2) =-2(p_1t_2-p_2t_1).
$$
Finally, recalling that
$$
    \zeta=2\bigl(p_1-t_2+\sqrt{-1}(p_2+t_1)\bigr),
$$
we obtain
$$
    |\zeta|^2 =4\bigl((p_1-t_2)^2+(p_2+t_1)^2\bigr) =4(e+J).
$$

For a horizontal vector $\xi$, define
$$
    \vartheta(\xi) =2\left( \langle\nabla_\xi e_1,e_2\rangle +\langle\nabla_\xi e_3,e_4\rangle \right).
$$
The two summands are the horizontal and vertical connection forms along the base. From the definitions of $\omega_1,\omega_2$,
$$
    \nabla_\xi\omega_1=\frac{\vartheta(\xi)}2\omega_2, \qquad \nabla_\xi\omega_2=-\frac{\vartheta(\xi)}2\omega_1.
$$
The induced covariant derivative of the trace-free symmetric block represented by $q=a+\sqrt{-1} b$ is therefore
$$
    Dq=\dd q+\sqrt{-1}\vartheta q =(\dd a-\vartheta b)+\sqrt{-1}(\dd b+\vartheta a).
$$
The factor two arises because the components of a trace-free symmetric endomorphism rotate through twice the angle of the frame. With our curvature convention, the horizontal and vertical connection forms have curvatures $-K\dd A_{\bg}$ and $-J\dd A_{\bg}$, respectively. Hence
\begin{equation}\label{eq:weyl-connection}
    \dd\vartheta=-2(K+J)\dd A_{\bg}.
\end{equation}
At the chosen point, rotate the horizontal and vertical frames so that their connection forms vanish in horizontal directions. Then $D_iq=a_i+\sqrt{-1} b_i$ there, where $a_i=e_i(a)$ and $b_i=e_i(b)$. Although the coefficients of $W^+$ are $\T^2$-invariant, their covariant derivatives in vertical directions include the derivatives of the frame.

The shape coefficients give
$$
    \begin{aligned}
        \nabla_{e_3}\omega_0 &=\left(\frac{u_2}{2}+p_2+t_1\right)\omega_1 -\left(\frac{u_1}{2}+p_1-t_2\right)\omega_2,\\
        \nabla_{e_4}\omega_0 &=\left(\frac{u_1}{2}-p_1+t_2\right)\omega_1 +\left(\frac{u_2}{2}-p_2-t_1\right)\omega_2.
    \end{aligned}
$$
Differentiating $W^+\omega_0=-w\omega_0$ and pairing with $\omega_j$, for $j=1,2$, yields
$$
    \langle(\nabla_{e_\alpha}W^+)\omega_0,\omega_j\rangle =-\langle(W^++w\Id) \nabla_{e_\alpha}\omega_0,\omega_j\rangle.
$$
Now we use the contracted Bianchi identity $\delta_gW^+=0$. Its $e_2$-component, paired with $\omega_0$, reads
$$
    0=w_1 +\langle(\nabla_{e_3}W^+)\omega_0,\omega_2\rangle -\langle(\nabla_{e_4}W^+)\omega_0,\omega_1\rangle.
$$
Substituting the preceding expressions, and computing the $e_1$-component in the same way, gives
\begin{equation}\label{eq:weyl-first-bianchi}
    \begin{aligned}
        w_1+\frac32wu_1 &=2a(p_1-t_2)+2b(p_2+t_1),\\
        w_2+\frac32wu_2 &=2b(p_1-t_2)-2a(p_2+t_1).
    \end{aligned}
\end{equation}
Squaring and adding gives
$$
    \begin{aligned}
        \left|\dd w+\frac32w\dd u\right|_{\bg}^2 &=4(a^2+b^2)\bigl((p_1-t_2)^2+(p_2+t_1)^2\bigr)\\
        &=4(e+J)|q|^2.
    \end{aligned}
$$
This proves \eqref{eq:weyl-bianchi-norm}.

We next derive the second Bianchi equation in the same frame. Pairing the $e_3$-component of $\delta_gW^+=0$ with $\omega_1,\omega_2$ gives
$$
    \begin{aligned}
        0&=a_1+b_2+\frac12w_1 -\langle(\nabla_{e_4}W^+)\omega_0,\omega_1\rangle,\\
        0&=b_1-a_2+\frac12w_2 -\langle(\nabla_{e_4}W^+)\omega_0,\omega_2\rangle.
    \end{aligned}
$$
Using the preceding formulas for $\nabla_{e_4}\omega_0$ and $(\nabla_{e_4}W^+)\omega_0$, these become
$$
    \begin{aligned}
        0={}&a_1+b_2+\frac12w_1 +\left(\frac32w+a\right) \left(\frac12u_1-p_1+t_2\right) +b\left(\frac12u_2-p_2-t_1\right),\\
        0={}&b_1-a_2+\frac12w_2 +b\left(\frac12u_1-p_1+t_2\right) +\left(\frac32w-a\right) \left(\frac12u_2-p_2-t_1\right).
    \end{aligned}
$$
Substituting \eqref{eq:weyl-first-bianchi} yields
$$
    \begin{aligned}
        a_1+b_2 &=\frac32w(p_1-t_2)-\frac12(au_1+bu_2),\\
        b_1-a_2 &=\frac32w(p_2+t_1)-\frac12(bu_1-au_2).
    \end{aligned}
$$
With $\zeta=2(p_1-t_2+\sqrt{-1}(p_2+t_1))$, these combine into
\begin{equation}\label{eq:weyl-second-bianchi}
    D_1q-\sqrt{-1}D_2q =\frac34w\zeta-\frac12q(u_1-\sqrt{-1} u_2).
\end{equation}
The use of covariant derivatives makes this identity independent of the choice of oriented frames.

To prove \eqref{eq:weyl-scalar}, we use the Weyl Weitzenb\"ock formula. Differentiating $W^+\omega_0=-w\omega_0$ twice gives
\begin{equation}\label{eq:weyl-eigenvalue-laplacian}
    \begin{aligned}
        \Delta_gw &=-\langle(\Delta_gW^+)\omega_0,\omega_0\rangle+2\sum_{\alpha=1}^4 \langle(w\Id+W^+)\nabla_{e_\alpha}\omega_0, \nabla_{e_\alpha}\omega_0\rangle.
    \end{aligned}
\end{equation}
The horizontal derivatives of $\omega_0$ vanish, while
$$
    \begin{aligned}
        \nabla_{e_3}\omega_0 &=\tfrac12(u_2+\operatorname{Im}\zeta)\omega_1 -\tfrac12(u_1+\operatorname{Re}\zeta)\omega_2,\\
        \nabla_{e_4}\omega_0 &=\tfrac12(u_1-\operatorname{Re}\zeta)\omega_1 +\tfrac12(u_2-\operatorname{Im}\zeta)\omega_2.
    \end{aligned}
$$
Substituting these and \eqref{eq:weyl-first-bianchi} into the sum in \eqref{eq:weyl-eigenvalue-laplacian} gives
$$
    \begin{aligned}
        &2\sum_{\alpha=1}^4 \langle(w\Id+W^+)\nabla_{e_\alpha}\omega_0, \nabla_{e_\alpha}\omega_0\rangle\\
        &\qquad=6w(K+J)-2\langle\dd u,\dd w\rangle_{\bg}.
    \end{aligned}
$$
The first term on the right of \eqref{eq:weyl-eigenvalue-laplacian} is $2\lambda w+3w^2-4|q|^2$. Equation \eqref{eq:weyl-scalar} now follows from $K+J=\lambda/3-w$ and
$$
    \Delta_gw=\Delta_{\bg}w+\langle\dd u,\dd w\rangle_{\bg}.
$$

\end{proof}

The scalar equation implies the following positivity statement.
\begin{proposition}\label{prop:weyl-positive}
Either $W^+=0$ everywhere, or $w>0$ on $B$.
\end{proposition}

\begin{proof}
At an open rod, the operator on the left of \eqref{eq:weyl-scalar} has the form
$$
    \partial_{rr}+\frac3r\partial_r +f(r^2,t)^{-2}\partial_{tt} +r b_1(r^2,t)\partial_r+b_2(r^2,t)\partial_t.
$$
On even functions it is the restriction of a smooth uniformly elliptic operator obtained by replacing $r$ with the radial variable in $\R^4$. At a turning point it similarly has the form
$$
    a\left(\partial_{xx}+\frac3x\partial_x\right) +c\left(\partial_{yy}+\frac3y\partial_y\right) +2xyb\partial_{xy}+xd_1\partial_x+yd_2\partial_y,
$$
where the coefficients are smooth in $x^2,y^2$. Writing $x=|\mathbf{x}|$ and $y=|\mathbf{y}|$, with $\mathbf{x},\mathbf{y}\in\R^4$, gives a smooth uniformly elliptic operator; its mixed principal coefficients are $b\mathbf{x}_i\mathbf{y}_j$. The ordinary maximum principle therefore also applies at boundary extrema of $w$.

At a negative minimum, the left-hand side of \eqref{eq:weyl-scalar} is nonnegative and the right-hand side is negative. Thus $w\geq0$. Moreover,
$$
    \Delta_{\bg}w+3\langle\dd u,\dd w\rangle_{\bg} -4\lambda w=-3w^2-4|q|^2\leq0.
$$
The strong maximum principle implies that either $w>0$ everywhere or $w=0$ everywhere. In the latter case \eqref{eq:weyl-scalar} gives $q=0$, hence $W^+=0$.
\end{proof}

Suppose now that $W^+\not\equiv0$. By Proposition \ref{prop:weyl-positive}, $w>0$ on $B$. Using $\Delta_{\bg}u=-2\lambda-2U$ and $e=K+U/2$ and Equations \eqref{eq:weyl-bianchi-norm}-\eqref{eq:weyl-scalar}, we have
\begin{align}
    \Delta_{\bg}\log w &=3(K+J)+\left(\frac{4\lambda}{3}+2U\right) \left(\frac94-\frac{|q|^2}{w^2}\right), \label{eq:weyl-logw}\\
    \Delta_{\bg}\log(\rho^{3/2}w) &=3(e+J)-4\left(\frac\lambda3+\frac U2\right) \frac{|q|^2}{w^2}. \label{eq:weyl-logweight}
\end{align}

\begin{lemma}\label{lem:weyl-ratio}
On $B$, we have $|q|\leq3w/2$.
\end{lemma}

\begin{proof}
We first derive an equation for $|q|/w$ on $\{q\neq0\}$. Write $q=|q|e^{\sqrt{-1}\theta}$ locally on this set. Then
$$
    \frac{D_iq}{q} =e_i(\log|q|)+\sqrt{-1}(\dd\theta+\vartheta)(e_i).
$$
Dividing \eqref{eq:weyl-second-bianchi} by $q$ and using the conjugate of \eqref{eq:weyl-first-bianchi}, we obtain
$$
    \frac{D_1q-\sqrt{-1}D_2q}{q} =\frac{3w^2}{4|q|^2}(e_1-\sqrt{-1} e_2) \log(\rho^{3/2}w) -\frac12(u_1-\sqrt{-1} u_2).
$$
Taking real and imaginary parts we see
$$
    \begin{aligned}
        (\dd\theta+\vartheta)(e_1) &=e_2(\log|q|)+\frac12u_2 -\frac{3w^2}{4|q|^2}e_2\log(\rho^{3/2}w),\\
        (\dd\theta+\vartheta)(e_2) &=-e_1(\log|q|)-\frac12u_1 +\frac{3w^2}{4|q|^2}e_1\log(\rho^{3/2}w).
    \end{aligned}
$$
Let $\underline{\ast}$ denote the Hodge star of $\bg$. Since $\underline{\ast}e^1=e^2$ and $\underline{\ast}e^2=-e^1$, these identities are equivalent to
\begin{equation}\label{eq:weyl-phase}
    \dd\theta+\vartheta =-\underline{\ast}\dd\log|q| -\frac12\underline{\ast}\dd u +\frac{3w^2}{4|q|^2}\underline{\ast} \dd\log(\rho^{3/2}w).
\end{equation}

Taking exterior derivatives and applying \eqref{eq:weyl-connection} gives
$$
    \begin{aligned}
        &\Delta_{\bg}\log\frac{|q|}{w} -\operatorname{div}_{\bg}\left[ \left(\frac{3w^2}{4|q|^2}-\frac13\right) \nabla_{\bg}\log(\rho^{3/2}w)\right]\\
        &\qquad=2(K+J)-\frac23\Delta_{\bg}\log w\\
        &\qquad=-\frac83\left(\frac\lambda3+\frac U2\right) \left(\frac94-\frac{|q|^2}{w^2}\right).
    \end{aligned}
$$
Expanding the divergence and using \eqref{eq:weyl-logweight}, we obtain 
\begin{equation}\label{eq:weyl-ratio-equation}
    \begin{aligned}
        &\Delta_{\bg}\log\frac{|q|}{w} +\frac{3w^2}{2|q|^2} \left\langle\dd\log(\rho^{3/2}w), \dd\log\frac{|q|}{w}\right\rangle_{\bg}\\
        &\qquad=\left(\frac94-\frac{|q|^2}{w^2}\right) \left(\frac{(e+J)w^2}{|q|^2} -4\left(\frac\lambda3+\frac U2\right)\right).
    \end{aligned}
\end{equation}

By \eqref{eq:weyl-boundary}, the ratio equals $3/2$ on $\partial B$. If it had a larger maximum, that maximum would occur in $B^\circ\cap\{q\neq0\}$. The left-hand side of \eqref{eq:weyl-ratio-equation} would be nonpositive there, whereas both factors on the right would be negative. For the second factor, use $e+J\geq0$ and
$$
    \frac{(e+J)w^2}{|q|^2} -4\left(\frac\lambda3+\frac U2\right) =(e+J)\left(\frac{w^2}{|q|^2}-4\right)-4w<0.
$$
This is a contradiction.
\end{proof}

\begin{proof}[Proof of Theorem~\ref{thm:toric-weyl-rigidity}]
Suppose first that $W^+\not\equiv0$. The function $w$ is positive and smooth on the closed quotient. We claim
$$
    \int_B\Delta_{\bg}\log w\dd A_{\bg}=0.
$$
Indeed, the normal derivative vanishes on open rods, and the flux over a corner arc of radius $r$ is $O(r^2)$. Integrating \eqref{eq:weyl-logw} we get
\begin{equation}\label{eq:weyl-topological-identity}
    \int_B\left(\frac{4\lambda}{3}+2U\right) \left(\frac94-\frac{|q|^2}{w^2}\right)\dd A_{\bg} =3\pi\bigl(b_2^+(X_{\mR})-1\bigr).
\end{equation}
The integrand is nonnegative by Lemma~\ref{lem:weyl-ratio}. If $b_2^+(X_{\mR})=0$, then  \eqref{eq:weyl-topological-identity} is impossible. Hence $W^+=0$ in that case.

If $b_2^+(X_{\mR})=1$, the case $W^+=0$ is excluded by the Weitzenb\"ock formula for harmonic self-dual two-forms,
$$
    0=\nabla_g^*\nabla_g\eta+\frac{4\lambda}{3}\eta-2W^+(\eta).
$$
Thus \eqref{eq:weyl-topological-identity} applies and forces $|q|=3w/2$. By \eqref{eq:weyl-matrix},
$$
    \operatorname{spec}(W^+)=\{-w,-w,2w\},\qquad w>0.
$$
Derdzi\'nski's theorem \cite[Proposition~5(i), pp.~420-421]{Derdzinski} implies that $w^{2/3}g$ is locally K\"ahler. Since $X_{\mR}$ is simply connected, the simple eigenline of $W^+$ admits a global unit section, and the local K\"ahler forms agree after a choice of sign. Thus $w^{2/3}g$ is globally K\"ahler with respect to a complex structure compatible with the orientation.

Conversely, if $W^+=0$, the Weitzenb\"ock formula above implies that $b_2^+(X_{\mR})=0$.
If $g$ is conformally K\"ahler, the associated K\"ahler form is a nonzero self-dual harmonic two-form for $g$, so $W^+\not\equiv0$.
The K\"ahler eigenvalue structure of $W^+$ implies that $|q|/w$ takes values in $\{0,3/2\}$.
By continuity and \eqref{eq:weyl-boundary}, we have $|q|/w=3/2$ everywhere.
Equation~\eqref{eq:weyl-topological-identity} therefore implies that $b_2^+(X_{\mR})=1$.
\end{proof}

\subsection{An improved Hitchin-Thorpe inequality and topological constraint}
\label{ss:improved HT}
As an application of the above discussion, we obtain the following improvement of the Hitchin-Thorpe inequality and the restriction on self-intersection numbers for comapct toric Einstein metrics.
\begin{theorem}\label{thm:improved HT}
Let $(X_{\mR},g)$ be a compact toric Einstein four-manifold. Then
\begin{equation}\label{eq:weyl-improved-toric-hitchin-thorpe}
    \chi(X_{\mR})\geq3|\tau(X_{\mR})|.
\end{equation}
Equality holds if and only if all rod self-intersections have the same value $\varepsilon\in\{-1,1\}$. 
\end{theorem}

We rescale $g$ so that $\Ric(g)=3g$. Given a rod sphere $\Sigma_i$, write $e_i=[\Sigma_i]^2$.
\begin{proposition}\label{lem:rod-self-intersection-bound1}
We have for all $i$, 
\begin{equation}\label{eq:weyl-rod-area-euler-bound}
    |e_i|+\frac{\Area_g(\Sigma_i)}{2\pi}\leq2.
\end{equation}
In particular, $e_i\in\{-1,0,1\}$.
\end{proposition}

\begin{proof}
Let $w_+$ be the function $w$ of Section~\ref{subsec:toric-weyl-rigidity}, and let $w_-$ be its counterpart for the opposite orientation. Proposition~\ref{prop:weyl-positive}, applied in both orientations, implies that $w_\pm\geq0$. This includes the half-conformally flat cases, when the corresponding function vanishes identically.

The sphere $\Sigma_i$ is totally geodesic. At a point over the interior of its rod, let $\sigma$ and $\nu$ be the unit tangent and normal area forms, with $*\sigma=\nu$. Then $(\sigma\pm\nu)/\sqrt2$ are eigenforms of $W^\pm$ with eigenvalues $2w_\pm$, respectively. Since the Einstein curvature operator on $\Lambda^\pm$ is $W^\pm+\Id$, these identities give
$$
    \Rm_g(\sigma)=(1+w_++w_-)\sigma+(w_+-w_-)\nu.
$$
Denote by $K_i$ the intrinsic Gaussian curvature of $\Sigma_i$ and by $K_i^\perp$ the scalar curvature of its oriented normal connection. It follows from the Gauss and Ricci equations that we have
\begin{equation}\label{eq:weyl-rod-tangent-normal-curvature}
    K_i=1+w_++w_-,\qquad |K_i^\perp|=|w_+-w_-|\leq K_i-1.
\end{equation}
The normal Euler-number formula and Gauss-Bonnet therefore yield
$$
    \begin{aligned}
        2\pi|e_i|=\left|\int_{\Sigma_i}K_i^\perp\dd A_g\right|
        \leq\int_{\Sigma_i}|K_i^\perp|\dd A_g
        \leq4\pi-\Area_g(\Sigma_i).
    \end{aligned}
$$
This proves \eqref{eq:weyl-rod-area-euler-bound}. The area is positive, so $|e_i|<2$, and the integrality of $e_i$ gives the last assertion.
\end{proof}

\begin{proof}[Proof of Theorem~\ref{thm:improved HT}]
Recall that
$$
    \chi(X_{\mR})=m,\qquad 3\tau(X_{\mR})=\sum_{i=1}^m e_i.
$$
By Proposition~\ref{lem:rod-self-intersection-bound1}, each $e_i$ belongs to $\{-1,0,1\}$. Consequently,
$$
    \begin{aligned}
        3|\tau(X_{\mR})|&=\left|\sum_{i=1}^m e_i\right|\leq\sum_{i=1}^m|e_i|\leq m=\chi(X_{\mR}).
    \end{aligned}
$$
Equality in the last inequality requires $|e_i|=1$ for every $i$, while equality in the triangle inequality then requires all $e_i$ to have the same sign. Conversely, these conditions imply equality.
\end{proof}

\begin{remark}
For a prescribed rod structure,  the individual self-intersection bounds in Proposition \ref{lem:rod-self-intersection-bound1} give a stronger obstruction than the improved Hitchin-Thorpe inequality. For example, the standard toric action on the Hirzebruch surface $\mathbb F_2$ has rod self-intersections $(2,0,-2,0)$. Although $\chi=4$ and $\tau=0$ satisfy the improved Hitchin-Thorpe inequality, this action admits no invariant Einstein metric. Its underlying manifold is $\mathbb S^2\times\mathbb S^2$, which admits the round product Einstein metric with a different torus action.
\end{remark}

\subsection{Rigidity in the Hermitian case}
\label{subsec:hermitian-uniqueness}

The following uniqueness result is known by the work of Hitchin \cite{Hitchin} and LeBrun \cite{LeBrun}. We give a new proof using a concavity property of the volume functional. Let $X_{\mR}$ be one of the complex surfaces $\CP^2\# k\overline{\CP^2}$, $0\leq k\leq3$, or $\mathbb S^2\times\mathbb S^2$, with its standard complex structure $\mathcal J$ and $\T^2$-action.

\begin{proposition}\label{prop:toric-hermitian-uniqueness}
 A toric Einstein metric on $X_{\mR}$ that is Hermitian with respect to $\mathcal J$ is unique up to scaling and isometry.
\end{proposition}

By the conformally K\"ahler reduction in \cite[Proposition~1]{LeBrun}, we can write $g=f^{-2}h$, where $h$ is a $\T^2$-invariant extremal K\"ahler metric in the sense of Calabi. Its scalar curvature is affine in moment coordinates. In the non-K\"ahler case, $f$ is a positive constant multiple of $\Scal_h$, and hence is a positive affine function on the moment polygon $P$. In the K\"ahler case, we take $h=g$ and $f=1$.

For the flux and volume calculations below, we work more generally with smooth metrics $g=f^{-2}h$, where $h$ is any $\T^2$-invariant K\"ahler metric on the fixed complex surface $(X_{\mR},\mathcal J)$ and $f$ is a positive affine function on its moment polygon $P$. We can write, in the action-angle coordinates, that
$$
    \begin{aligned}
        h&=H_{ij}\dd x_i\dd x_j+H^{ij}\dd\theta_i\dd\theta_j,\\
        H_{ij}&=\partial_i\partial_j\psi,
    \end{aligned}
$$
where $(H^{ij})=H^{-1}$ and $\psi$ is a symplectic potential. Let $E_i$ be the edge of $P$ corresponding to $\mI_i$, with endpoints $z_i,z_{i+1}$ and primitive outward normal $n_i$. Write
$$
    P=\{x:\langle n_i,x\rangle\leq b_i\}, \qquad z_{i+1}-z_i=L_it_i,
$$
where $t_i$ is the primitive tangent in the positive boundary direction. We denote by $\dd\sigma$ the affine boundary measure, defined on each edge $E_i$ by $\dd\sigma=\dd s$ under the parametrization $x=z_i+s t_i$, where $0\leq s\leq L_i$.
Thus
$$
    \int_{E_i}\dd\sigma=L_i.
$$
Using the notation of Section~\ref{sec:rod structure}, we define the integrated boundary flux by
$$
    M_i=-\int_{\mI_i}\partial_n\rho\dd s_{\bg}.
$$
By \eqref{eq:rod-flux} and the conformal relation $g=f^{-2}h$, we have
$$
    M_i=\frac{\Area_g(\Sigma_i)}{2\pi}=\int_{E_i}f^{-2}\dd\sigma.
$$
Since $f$ is affine on $E_i$, we obtain
$$
    M_i=\frac{L_i}{f(z_i)f(z_{i+1})}.
$$
Denote the total boundary flux 
$$ M_{\mathrm{tot}}=\sum_iM_i, $$
and the normalized total volume 
\begin{equation}\label{eq:hermitian-flux-volume}
    \begin{aligned}
      V=\int_Pf^{-4}\dd x =\frac{\Vol_g(X_{\mR})}{(2\pi)^2},
    \end{aligned}
\end{equation}
where $\dd x$ denotes the affine area measure on $P$. 

We now show that the integrated boundary fluxes give smooth coordinates for the pairs $(P,f)$. For the parametrization below, we identify $(h,f)$ with $(s^2h,sf)$ for $s>0$, since this rescaling leaves $g=f^{-2}h$ unchanged.

\begin{lemma}\label{lem:hermitian-flux-data}
With the primitive outward normals $n_i$ fixed, the positive vectors $\mathbf M=(M_1,\ldots,M_m)$ give a smooth parametrization of the pairs $(P,f)$, up to translation and the above rescaling.
\end{lemma}
\begin{proof}
By the divergence theorem, we have 
\begin{equation}\label{eq:hermitian-flux-direction}
    \sum_iM_in_i=-2\nabla f\int_Pf^{-3}\dd x.
\end{equation}
We first assume that $f$ is constant.
After rescaling, we may assume that $f=1$, so $L_i=M_i$.
Thus the flux vector $\mathbf M$ determines the polygon $P$ uniquely up to translation.
To prove smooth dependence near this pair, we fix one labelled vertex at the origin and impose $f(0)=1$ on nearby pairs.
An infinitesimal variation of $f$ then has the form $\dot f=\xi\cdot x$, where $\xi\in\R^2$.
Differentiating the flux identities at $f=1$, we obtain
$$
    \begin{aligned}
        &\dot M_i=\dot L_i-L_i\xi\cdot(z_i+z_{i+1}),\\
        &\sum_i\dot M_in_i=-2\Area(P)\xi.
    \end{aligned}
$$
The second identity determines $\xi$ from $\dot{\mathbf M}$, and the first then determines each $\dot L_i$.
Since the allowed length variations satisfy the two closure conditions $\sum_i\dot L_i n_i=0$, the source and target of the differential both have dimension $m$.
Thus the differential is invertible, and the inverse function theorem implies smooth dependence of $(P,f)$ on $\mathbf M$ near every pair with constant $f$.

 Now we assume $f$ is not constant, then $\nabla f$ is a positive multiple of $-\sum_iM_in_i$. We can use the rescaling to fix its magnitude, and choose positively oriented orthonormal coordinates $(y_1,y_2)$ in which $f=y_2$. Denote by $t_{i, a}$ and $z_{i, a}$ the $y_a$-components of  $t_i$ and $z_i$ respectively, where $a=1,2$. Denote $r_i=1/z_{i, 2}$, $r=r_1$, and $c_i=r_i-r_1$. Then
$$
    r_{i+1}-r_i=-M_it_{i, 2},\qquad r_i=r+c_i,
$$
where the $c_i$ are determined by $\mathbf M$ and $c_1=0$. It suffices to determine $r$. Notice that, 
$$
    L_i(r)=\frac{M_i}{(r+c_i)(r+c_{i+1})}.
$$
Since $r_{m+1}=r_1$, the total displacement in the $y_2$-direction is $$\sum_iL_it_{i, 2}=\sum_i(r_{i+1}^{-1}-r_i^{-1})=0.$$ To obtain a closed polygon, we need the sum in the $y_1$-direction also vanishes, i.e., 
$$
    F(r)=\sum_i \frac{M_it_{i, 1}}{(r+c_i)(r+c_{i+1})}=0.
$$
At a zero of $F$, we have 
$$  
        F'(r)=-\sum_iL_it_{i, 1}\bigl(z_{i,2}+z_{i+1,2})\bigr)=-2\int_{\partial P}y_2\dd y_1=2\Area(P)>0.
$$
So there is at most one zero. With the direction fixed by \eqref{eq:hermitian-flux-direction},
$$
    \sum_iM_it_{i,1}>0,
$$
so $F(r)>0$ for large $r$. We now consider the limit $r\rightarrow-\min_i c_i$.
If there is a horizontal upper edge, the leading term of $F(r)$ has order $(r+\min_i c_i)^{-2}$ and a negative coefficient, since the positive boundary direction along this edge is the negative $y_1$-direction.
Otherwise, we denote by $t_-$ and $t_+$ the incoming and outgoing primitive tangents at the upper vertex.
The leading term of $F(r)$ then has order $(r+\min_i c_i)^{-1}$, with coefficient
$$
    \frac{t_{-,1}}{t_{-,2}}-\frac{t_{+,1}}{t_{+,2}}
    =\frac{\det(t_-,t_+)}{t_{-,2}t_{+,2}}<0,
$$
since $t_{-,2}>0$, $t_{+,2}<0$, and $\det(t_-,t_+)>0$.
Thus $F(r)\to-\infty$ at this endpoint.
Since $F(r)>0$ for sufficiently large $r$, continuity implies that $F$ has a zero.
Moreover, $\partial_rF>0$ at its unique zero, so the implicit function theorem implies that $r$ depends smoothly on $\mathbf M$.

\end{proof}

By Lemma~\ref{lem:hermitian-flux-data}, we can therefore regard $V$ as a smooth function $V(\mathbf M)$ for $\mathbf M\in (0,\infty)^m$. Now we fix a polygon $P$, with centroid given by the origin, and set $A=\Area(P)$.

For $\mathbf M\in\R^m$, we set
$$
    \xi=-\frac{1}{2A}\sum_iM_in_i,
    \qquad
    \eta_i=M_i+2\int_{E_i}\xi\cdot x \dd\sigma.
$$
By the divergence theorem, we have
$$
    \sum_i\eta_in_i=\sum_iM_in_i+2A\xi=0,
$$
so the vectors $\eta_i t_i$, placed successively, form a closed polygon. Its signed area is
\begin{equation}\label{eqn:area formula}
    \mathcal A(\eta) :=\frac12\sum_{i<j}\eta_i\eta_j\det(t_i,t_j).
\end{equation}

Now we write $\mathbf L=(L_1,\ldots,L_m)$, By the definition of the area form and since $\sum_i L_it_i=0$, we have $$\mathcal A(\mathbf L)=A=\frac{1}{2}\sum_i b_iL_i.$$ Now define the mixed area by
$$
    \mathcal A(P,\eta)=\frac12\bigl(\mathcal A(\mathbf L+\eta)-\mathcal A(\mathbf L)-\mathcal A(\eta)\bigr).
$$
Differentiating \eqref{eqn:area formula} we see
$$
    2\mathcal A(P,\eta)=\sum_i b_i\eta_i.
$$

\begin{lemma}\label{lem:hermitian-volume-formula}
If  $\mathbf M$ is the flux vector of a positive affine function $f$ on $P$, then
\begin{equation}\label{eq:hermitian-volume-formula}
    6V(\mathbf M)=2\mathcal A(\eta)+\frac{4}{A}\mathcal A(P,\eta)^2-12\int_P(\xi\cdot x)^2\dd x.
\end{equation}
\end{lemma}

\begin{proof}
We first assume that $f=c$ is constant.
Then $\xi=0$ and $\eta=c^{-2}\mathbf L$, so both sides of \eqref{eq:hermitian-volume-formula} equal $6c^{-4}A$.

We now assume that $f$ is nonconstant and rotate the coordinates to write $f=c+py_2$ with $p>0$.
We set
$$
    \kappa=\frac{p}{A}\int_Pf^{-3}\dd x.
$$
By \eqref{eq:hermitian-flux-direction} and the definitions of $\xi$ and $\eta$, we have
$$
    \xi=\kappa\partial_{y_2},\qquad
    \eta_i=\int_{E_i}(f^{-2}+2\kappa y_2)\dd\sigma.
$$
To compute $\mathcal A(\eta)$, we set
$$
    a_0(y_2)=f^{-2}+2\kappa y_2,\qquad
    a_1(y_2)=-\frac1{pf}+\kappa y_2^2,
$$
so that $a_1'=a_0$. Then we integrate the vector-valued one-form $a_0(y_2)(\dd y_1,\dd y_2)$ along $\partial P$.
The resulting curve is closed because $\sum_i\eta_it_i=0$.
Each edge $E_i$ is mapped into a line with displacement $\eta_it_i$, so the signed area of this curve is $\mathcal A(\eta)$.
Its second coordinate is $a_1(y_2)$ up to an additive constant, and hence
$$
    \mathcal A(\eta)=-\int_{\partial P}a_1(y_2)a_0(y_2)\dd y_1=\int_P(a_0^2+a_1a_0')\dd x.
$$
By direct computation we obtain
$$
    \mathcal A(\eta)=\int_P\left(3f^{-4}-\frac{2\kappa c^2}{p}f^{-3}+6\kappa^2y_2^2\right)\dd x
       =3V(\mathbf M)-\frac{2Ac^2\kappa^2}{p^2}+6\int_P(\xi\cdot x)^2\dd x.
$$
Using  the divergence theorem, we also obtain
$$
  2\mathcal A(P,\eta)=\int_P(2a_0+y_2a_0')\dd x=2c\int_Pf^{-3}\dd x+6\kappa\int_Py_2\dd x=\frac{2cA\kappa}{p}.
$$
Then we see \eqref{eq:hermitian-volume-formula} holds.
\end{proof}

Motivated by Lemma~\ref{lem:hermitian-volume-formula}, we define the quadratic form
\begin{equation}\label{eq:hermitian-quadratic-form}
    Q_P(\mathbf M)=2\mathcal A(\eta)+\frac{4}{A}\mathcal A(P,\eta)^2-12\int_P(\xi\cdot x)^2\dd x
\end{equation}
for  $\mathbf M\in\R^m$. If $P$ is the polygon associated with $\mathbf M$, then from the above discussion we have 
\begin{equation}\label{eqn:volume equals quadratic form}
    Q_P(\mathbf M)=6V(\mathbf M).
\end{equation}

We now restrict to the triangle, rectangle, trapezoid, pentagon, and hexagon families corresponding to the surfaces in Proposition~\ref{prop:toric-hermitian-uniqueness}.
We fix one of these families and take the minimum below over polygons $P$ with the prescribed primitive outward normals.

\begin{lemma}\label{lem:hermitian-quadratic-minimum}
Each $Q_P$ has signature $(1,m-1)$, and for every $\mathbf M\in(0,\infty)^m$ we have
\begin{equation}\label{eq:hermitian-quadratic-minimum}
    0<V(\mathbf M)=\frac16\min_P Q_P(\mathbf M),
\end{equation}
with the minimum attained at the polygon associated with $\mathbf M$ by Lemma~\ref{lem:hermitian-flux-data}.
\end{lemma}

\begin{proof}

\medskip
\noindent\textbf{Step 1. The case of hexagons.}

Using translation invariance, we may represent a hexagon by
$$
    P=\operatorname{conv}\{(c,0),(a,0),(a,b-d),(a-d,b),(0,b),(0,c)\},
$$
where $0<c<\min\{a,b\}$ and $0<d<\min\{a,b\}$.
We index its vertices and corresponding edges from $0$ to $5$ in the displayed order.
Then
$$
    \begin{aligned}
        &A=ab-\frac{c^2+d^2}{2},\\
        &\mathbf L=(a-c,b-d,d,a-d,b-c,c).
    \end{aligned}
$$

We write $x=(x_1,x_2)$ in the coordinates used to describe the vertices of $P$.
Since $x$ is affine along each edge $E_i$, we have
$$
    \int_{E_i}x\dd\sigma=\frac{L_i}{2}(z_i+z_{i+1}).
$$
We compute the integrals of $x_1,x_2,x_1^2,x_1x_2,x_2^2$ over $P$ by subtracting the corresponding integrals over the two removed corner triangles from those over $[0,a]\times[0,b]$.
After expressing these integrals in coordinates centered at the centroid of $P$ and substituting into \eqref{eq:hermitian-quadratic-form}, we obtain
\begin{equation}\label{eq:hermitian-quadratic-quotient}
    \begin{aligned}
        &Q_P(\mathbf M)=Q_0(\mathbf M)+\frac{N_{\mathbf M}(a,b,c,d)}{A},\\
        &Q_0(\mathbf M)=\sum_{i<j}M_iM_j-M_0M_3-M_1M_4-M_2M_5,
    \end{aligned}
\end{equation}
where
\begin{equation}\label{eq:hermitian-quadratic-numerator}
    \begin{aligned}
        &N_{\mathbf M}=(b-d)^2M_0M_2+b^2M_0M_3+c^2M_0M_4+d^2M_1M_3+a^2M_1M_4+(a-c)^2M_1M_5\\
        &\qquad+(a-d)^2M_2M_4+(a+b-c-d)^2M_2M_5+(b-c)^2M_3M_5.
    \end{aligned}
\end{equation}
Since $Q_0(\mathbf M)$ is independent of $P$, minimizing $Q_P(\mathbf M)$ is equivalent to minimizing $N_{\mathbf M}/A$.

We first show that $N_{\mathbf M}/A$ is stationary at the polygon associated with $\mathbf M$.
Keeping $\mathbf M$ fixed, we differentiate \eqref{eq:hermitian-quadratic-numerator} with respect to $(a,b,c,d)$ and then substitute
$$
    M_i=\frac{L_i}{f(z_i)f(z_{i+1})},
    \qquad f(x_1,x_2)=r+sx_1+tx_2.
$$
After clearing denominators by $A\prod_{j=0}^5f(z_j)$ and expanding, we obtain
$$
    \begin{aligned}
        \partial_aN_{\mathbf M}=\frac{bN_{\mathbf M}}A,
        \qquad \partial_bN_{\mathbf M}=\frac{aN_{\mathbf M}}A,
        \partial_cN_{\mathbf M}=-\frac{cN_{\mathbf M}}A,
        \qquad \partial_dN_{\mathbf M}=-\frac{dN_{\mathbf M}}A.
    \end{aligned}
$$
Since $\dd A=b\dd a+a\dd b-c\dd c-d\dd d$, we have
\begin{equation}\label{eq:hermitian-quadratic-stationarity}
    \dd (N_{\mathbf M}/A)=0.
\end{equation}

Now we show it is a global minimum.
By \eqref{eq:hermitian-quadratic-numerator}  $N_{\mathbf M}$ is positive definite, whereas the area 
$$
    A=\left(\frac{a+b}{2}\right)^2-\left(\frac{a-b}{2}\right)^2-\frac{c^2+d^2}{2}
$$
has signature $(1,3)$.
By homogeneity, we can normalize $N_{\mathbf M}=1$, so minimizing $N_{\mathbf M}/A$ for $A>0$ is equivalent to maximizing $A$ on this ellipsoid.

We make a linear change of variables that transforms the ellipsoid into the unit sphere.
In these coordinates, $A$ has exactly one positive eigenvalue, so its only critical points with positive value are two antipodal global maximizers, one in each component of $\{A>0\}$.
By \eqref{eq:hermitian-quadratic-stationarity}, the normalized polygon associated with $\mathbf M$ is such a critical point.
Since the admissible polygon parameters lie in the component $a+b>0$, this polygon minimizes $N_{\mathbf M}/A$ among all admissible polygons.
Consequently, it also minimizes $Q_P(\mathbf M)$.

\medskip
\noindent\textbf{Step 2. The remaining polygon families.}

The other four families can be obtained  by collapsing the edges indicated below and setting their fluxes to zero.
All remaining fluxes are positive.
\begin{center}
\begin{tabular}{@{}llll@{}}
\toprule
Polygon & Restrictions & Area & Omitted edges\\
\midrule
Pentagon & $d=0$ & $ab-c^2/2$ & $2$\\
Rectangle & $c=d=0$ & $ab$ & $2,5$\\
Trapezoid & $d=0,\ b=c$ & $ab-b^2/2$ & $2,4$\\
Triangle & $d=0,\ a=b=c$ & $a^2/2$ & $0,2,4$\\
\bottomrule
\end{tabular}
\end{center}

The formula \eqref{eq:hermitian-quadratic-quotient} extends to these restrictions by continuity.
To obtain the stationarity identity, we approximate each reduced polygon by hexagons and use the same positive affine function $f$.
The fluxes of the collapsing edges tend to zero, so \eqref{eq:hermitian-quadratic-stationarity} passes to the limit and holds in the remaining parameter directions.

For the pentagon, rectangle, and trapezoid families, $N_{\mathbf M}$ is positive definite in the independent polygon parameters, while the corresponding area forms have one positive direction and all remaining directions negative.
Thus the same argument in {\bf Step 1} applies.

For the triangle, $N_{\mathbf M}=0$ and
$$
    Q_P(\mathbf M)=\sum_{i<j}M_iM_j,
$$
where the sum is over the remaining edges.
Thus $Q_P(\mathbf M)$ is independent of $P$ so it is easy to see the minimizing property in this case. 
This proves \eqref{eq:hermitian-quadratic-minimum} in all five cases.

\medskip
\noindent\textbf{Step 3. Signature and positivity.}

In each family, we use the independent polygon parameters above as linear coordinates for $\eta\in\R^m$ subject to $\sum_i\eta_i n_i=0$.
The corresponding area formulas therefore imply that $\mathcal A(\eta)$ has signature $(1,m-3)$ on this $(m-2)$-dimensional space.
We write
$$
    \eta=\frac{\mathcal A(P,\eta)}A\mathbf L+\eta_0,
    \qquad \mathcal A(P,\eta_0)=0.
$$
Since $\mathcal A(\mathbf L)=A>0$, the area form is negative definite on the subspace $\mathcal A(P,\eta_0)=0$.
By \eqref{eq:hermitian-quadratic-form}, we have
$$
    Q_P(\mathbf M)=\frac{6}{A}\mathcal A(P,\eta)^2
    +2\mathcal A(\eta_0)-12\int_P(\xi\cdot x)^2\dd x.
$$
The first term is positive in the scaling direction, the second is negative definite in the remaining $m-3$ polygon directions, and the last is negative definite in the two components of $\xi$.
Thus $Q_P$ has signature $(1,m-1)$.

\end{proof}

\begin{proposition}\label{prop:hermitian-volume-concavity}
The function $\sqrt{V(\mathbf M)}$ is strictly concave on the simplex
$$
    \{M_i>0:\ \sum_iM_i=1\}.
$$
\end{proposition}

\begin{proof}
We choose distinct vectors $\mathbf M_0,\mathbf M_1$ in the above simplex and set $\mathbf M_t=(1-t)\mathbf M_0+t\mathbf M_1$ for $0<t<1$.
By Lemma~\ref{lem:hermitian-quadratic-minimum}, we can choose a polygon $P$ attaining the minimum at $\mathbf M_t$.

By the same lemma, $Q_P$ has signature $(1,m-1)$ and is positive on $(0,\infty)^m$, so $(0, \infty)^m$ lies in one component of $\{Q_P>0\}$.
Since $\mathbf M_0$ and $\mathbf M_1$ are distinct and their components sum to one, they are not proportional.
Thus, by the reverse Cauchy--Schwarz inequality and the minimum formula, we obtain
$$
    \begin{aligned}
        &\sqrt{6V(\mathbf M_t)}=\sqrt{Q_P(\mathbf M_t)}\\
        &\qquad>(1-t)\sqrt{Q_P(\mathbf M_0)}+t\sqrt{Q_P(\mathbf M_1)}\\
        &\qquad\geq(1-t)\sqrt{6V(\mathbf M_0)}+t\sqrt{6V(\mathbf M_1)}.
    \end{aligned}
$$
We then obtain the required strict concavity.
\end{proof}

\begin{proof}[Proof of Proposition~\ref{prop:toric-hermitian-uniqueness}]
By the conformally K\"ahler reduction, it suffices to prove uniqueness among the metrics considered above.
We rescale each metric so that $\sum_iM_i=1$.

Since $f$ is affine, the conformal scalar-curvature formula implies that
\begin{equation}\label{eq:hermitian-scalar-equation}
    \Scal_g=-f^5\partial_i\partial_j(f^{-3}H^{ij}).
\end{equation}
Using the facet boundary conditions for $\psi$, we integrate twice by parts and obtain
\begin{equation}\label{eq:hermitian-total-scalar}
    \frac1{(2\pi)^2}\int_{X_{\mR}}\Scal_g\dvol_g=2\sum_iM_i=2.
\end{equation}
For an Einstein metric $g$ with $\Ric_g=\lambda g$, this identity implies that $1=2\lambda V(g)$, so $\lambda>0$.
For every smooth path $g_t$ in the normalized class with $g_0=g$, we apply Einstein stationarity with $\lambda$ fixed to obtain
$$
    0=\left.\frac{\dd}{\dd t}\right|_{t=0}\frac{\EH_\lambda(g_t)}{(2\pi)^2}
    =-2\lambda\left.\frac{\dd}{\dd t}\right|_{t=0}V(g_t).
$$
Thus $g$ is a critical point of $V$, and hence of $\sqrt V$, under these variations.

We can realize every sufficiently small variation of $P$ with the prescribed normals by varying the toric K\"ahler class and choosing invariant K\"ahler representatives smoothly.
Together with variations of $f$, Lemma~\ref{lem:hermitian-flux-data} therefore realizes every nearby flux vector.
Consequently, the flux vector of a normalized Einstein metric is a critical point of $\sqrt{V(\mathbf M)}$ on $\{M_i>0:\sum_iM_i=1\}$.
By Proposition~\ref{prop:hermitian-volume-concavity}, this function is strictly concave, so there is at most one such critical point.

We now consider two normalized Einstein metrics.
Their flux vectors coincide, so by Lemma~\ref{lem:hermitian-flux-data} we can choose representatives with the same pair $(P,f)$.
Their Einstein constants also coincide, since both equal $1/(2V(\mathbf M))$.
We denote their symplectic potentials by $\psi_0,\psi_1$ and their Hessians by $H_0,H_1$.
We divide their scalar-curvature equations \eqref{eq:hermitian-scalar-equation} by $f^5$, subtract, multiply by $\psi_1-\psi_0$, and integrate twice by parts.
Since $\psi_1-\psi_0$ is smooth on $\overline P$ and the two potentials satisfy the same facet boundary conditions, the boundary terms vanish.
Thus we obtain
$$
    \begin{aligned}
        &0=\int_P f^{-3}\operatorname{tr}\bigl((H_1-H_0)(H_1^{-1}-H_0^{-1})\bigr)\dd x\\
        &\qquad=-\int_P f^{-3}\bigl|H_1^{-1/2}(H_1-H_0)H_0^{-1/2}\bigr|^2\dd x.
    \end{aligned}
$$
Therefore $H_1=H_0$, so the potentials differ by an affine function and determine the same K\"ahler metric.
Since $f$ is also the same, the normalized Einstein metrics agree under these identifications.
We conclude that the Einstein metric is unique up to scaling and isometry.
\end{proof}

\begin{proof}[Proof of Corollary~\ref{cor:toric-small-betti-classification}]
We first assume that $b_2^+(X_{\mR})=1$.
By the construction in the proof of Theorem~\ref{thm:toric-weyl-rigidity}, we obtain a $\T^2$-invariant K\"ahler metric conformal to $g$.
Since $H^1(X_{\mR};\R)=0$, the torus action is Hamiltonian, so the associated complex surface is projective and toric.
By Proposition~\ref{lem:rod-self-intersection-bound1} and the adjunction formula, we have
$$
    (-K_{X_{\mR}})\cdot\Sigma_i=2+[\Sigma_i]^2\geq1.
$$
Thus the anticanonical bundle is ample by the toric ampleness criterion.
The classification of smooth toric del Pezzo surfaces therefore places $(X_{\mR},\mathcal J)$ among the surfaces in Proposition~\ref{prop:toric-hermitian-uniqueness}, which implies the desired conclusion.

We now assume that $b_2^+(X_{\mR})=0$.
Since $b_1(X_{\mR})=0$, by Theorem~\ref{thm:improved HT} we have
$$
    2+b_2^-(X_{\mR})\geq3b_2^-(X_{\mR}),
$$
so $b_2^-(X_{\mR})\leq1$.
If $b_2^-(X_{\mR})=1$, we reverse orientation and apply the preceding case.
Since $b_2(X_{\mR})=1$, we obtain the Fubini--Study metric on $\CP^2$.

If $b_2^-(X_{\mR})=0$, we apply Theorem~\ref{thm:toric-weyl-rigidity} in both orientations and obtain $W_g^+=W_g^-=0$.
Thus $g$ has constant positive sectional curvature, and the simple connectivity of $X_{\mR}$ implies that $g$ is the round metric on $\mathbb S^4$, up to scaling and isometry.
\end{proof}

We now prove infinitesimal rigidity.

\begin{proof}[Proof of Theorem~\ref{thm:toric-small-betti-nondegeneracy}]
Let $k\in K_g^{\T^2}$. For the round sphere, we have
$$
    0=\int_{X_{\mR}}\left( |\nabla_g k|_g^2+\frac{2\lambda}{3}|k|_g^2 \right)\dvol_g,
$$
so $k=0$. We consider the remaining metrics in their complex orientations, and write $\mathcal J$ for the corresponding complex structure.

We first represent $k$ by a variation of $(P,f,\psi)$. Set $g_t=g+tk$. Since $k$ solves the linearized Einstein equation with fixed Einstein constant, we have
$$
\Ric(g_t)-\lambda g_t=O(t^2),\qquad E_t^0:=\Ric(g_t)-\frac{\operatorname{Scal}_{g_t}}4g_t=O(t^2).
$$
In the complex orientation, $W_g^+$ has a simple positive eigenvalue and a double negative eigenvalue. By compactness, the simple eigenvalue continues smoothly to a positive eigenvalue $\mu_t$ of $W_{g_t}^+$ for sufficiently small $t$. Following \cite[Section 3.1]{AA}, set $\widehat g_t=\mu_t^{2/3}g_t$ and choose the corresponding normalized self-dual eigenform $\widehat F_t$, continued from the background K\"ahler form. These choices are $\mathbb T^2$-invariant. At $t=0$, $\widehat g_0$ is K\"ahler and $\widehat\nabla\widehat F_0=0$; its constant normalization can be absorbed into the equivalence $(h,f)\sim(s^2h,sf)$.

We now apply the integral identities in \cite[Lemma 3.1]{AA}. Since $X_{\mathfrak R}$ is closed, the boundary integrals in equations (3.12a) and (3.12c) vanish. Comparing the quadratic coefficients in the integrated identity $\widehat\nabla_a\widehat V^a=\widehat A+\widehat B$, equation (3.12b) therefore gives an integral of nonnegative squares equal to zero. In particular,
$$
\int_{X_{\mathfrak R}}
\left|\left.\frac{\dd}{\dd t}\right|_{t=0}
\bigl(\widehat\nabla^{(t)}\widehat F_t\bigr)\right|_{\widehat g_0}^2
\dvol_{\widehat g_0}=0.
$$
Thus the almost-complex structures $\mathcal J_t$ associated with $\widehat F_t$ satisfy
$$
\widehat\nabla^{(t)}\mathcal J_t=O(t^2).
$$
Consequently the conformally K\"ahler condition holds to first order, and $\dot{\mathcal J}$ is an infinitesimal integrable complex-structure deformation.

Since $H^1(X_{\mathfrak R},T_{\mathcal J}^{1,0}X_{\mathfrak R})=0$, there is a real vector field $V$ such that $\dot{\mathcal J}+\mathcal L_V\mathcal J=0$. Averaging $V$ over $\mathbb T^2$ preserves this equation, so its flow is equivariant. After pulling back by this flow, we have $\mathcal J_t=\mathcal J+O(t^2)$. The resulting first variation $\dot\omega$ of the K\"ahler form is therefore a closed real $(1,1)$-form for $\mathcal J$. The path $\omega+t\dot\omega$ realizes the same first-order variation by genuine $\mathbb T^2$-invariant K\"ahler metrics $h_t$ for the fixed complex structure. Choosing positive functions $f_t$ with the prescribed initial value and first derivative, we obtain
$$
\left.\frac{\dd}{\dd t}\right|_{t=0}(f_t^{-2}h_t)
=k+\mathcal L_Vg,
\qquad
\operatorname{Ric}(f_t^{-2}h_t)-\lambda f_t^{-2}h_t=O(t^2).
$$

Under conformal change, we have
$$
\operatorname{Ric}(f_t^{-2}h_t)
=\operatorname{Ric}(h_t)+2f_t^{-1}\nabla_{h_t}^2f_t+q_t h_t
$$
for a scalar function $q_t$. Since $\operatorname{Ric}(h_t)$ and $h_t$ are $\mathcal J$-invariant, its anti-Hermitian part gives $(\nabla_{h_t}^2f_t)^-=O(t^2)$, and hence
$$
\mathcal L_{\mathcal J\nabla_{h_t}f_t}h_t=O(t^2).
$$
In smoothly varying moment coordinates $x_t$, the mixed components of this equation are
$$
H_t^{j\ell}\frac{\partial^2f_t}{\partial x_t^i\partial x_t^\ell}=O(t^2).
$$
Thus $f_t$ agrees to first order with an affine function of $x_t$. Replacing it by this positive affine function preserves the first-order metric variation. The moment polygons $P_t$ and symplectic potentials $\psi_t$ of $h_t$ therefore give a curve $(P_t,f_t,\psi_t)$ of the data used above which represents $k$ modulo equivariant infinitesimal diffeomorphisms.

Let $v=\dot{\mathbf M}$ be its flux variation. By \eqref{eq:hermitian-total-scalar}, the Einstein-Hilbert functional with fixed $\lambda$ is
$$
    \frac1{(2\pi)^2} \int_{X_{\mR}}(\Scal_g-2\lambda)\dvol_g =2M_{\mathrm{tot}}-2\lambda V.
$$
Every flux variation is realized by Lemma~\ref{lem:hermitian-flux-data}, so $2M_{\mathrm{tot}}-2\lambda V$ is critical at $\mathbf M$. Its Hessian vanishes on the variation induced by $k$, because $k$ solves the linearized Einstein equation. Since $M_{\mathrm{tot}}=\sum_iM_i$ is linear,
$$
    0=D^2(2M_{\mathrm{tot}}-2\lambda V)_{\mathbf M}(v,v) =-2\lambda D^2V_{\mathbf M}(v,v).
$$
Also, $\tr_gk=0$ gives $DV_{\mathbf M}[v]=0$.

Fix the polygon $P$ associated with $\mathbf M$. By \eqref{eq:hermitian-quadratic-minimum} and \eqref{eqn:volume equals quadratic form},
$$
    Q_P(\mathbf N)\geq6V(\mathbf N), \qquad Q_P(\mathbf M)=6V(\mathbf M)>0.
$$
Thus $Q_P-6V$ has a minimum at $\mathbf M$, so $DQ_P(\mathbf M)[v]=6DV_{\mathbf M}[v]=0$: the vector $v$ is orthogonal to $\mathbf M$ for the bilinear form associated with $Q_P$. Since $Q_P$ has signature $(1,m-1)$, it is negative definite on this orthogonal complement. On the other hand, differentiating the preceding inequality twice gives
$$
    0\leq D^2(Q_P-6V)_{\mathbf M}(v,v)=2Q_P(v).
$$
These two facts force $v=0$.

By Lemma~\ref{lem:hermitian-flux-data}, we may now choose representatives with $\dot P=\dot f=0$. The linearized scalar-curvature equation becomes
$$
    \partial_i\partial_j \left(f^{-3}H^{ia}\dot\psi_{ab}H^{bj}\right)=0.
$$
Multiplying by $\dot\psi$ and integrating twice by parts gives
$$
    0=\int_P f^{-3}H^{ia}H^{bj} \dot\psi_{ij}\dot\psi_{ab}\dd x.
$$
Here $\dot\psi$ is smooth up to the boundary, and the boundary terms vanish by the same fixed-facet calculation used in the uniqueness proof. Positivity implies $D^2\dot\psi=0$, so this variation does not change the metric.

Consequently, the original deformation is $k=\mathcal L_Zg$. Since $\delta_gk=0$,
$$
    \int_{X_{\mR}}|k|_g^2\dvol_g =2\int_{X_{\mR}}\langle\delta_gk,Z\rangle_g \dvol_g=0.
$$
Thus $k=0$.
\end{proof}

\section{Proof of Theorem \ref{thm:dichotomy} and Theorem \ref{thm:main} }\label{sec:degree}

We use the notations of the gluing construction in Section~\ref{sec:gluing}. Let $g_0$ be a torically nondegenerate smooth compact Einstein metric with $\Ric(g_0)=3g_0$, and let $(Y,h)$ be the equal-angle gravitational instanton associated with $W$, with $A_*=A(h)$. For each sufficiently small $\varepsilon>0$, the approximate metric is $b_\varepsilon=b_{\varepsilon,\beta_*}$, and we write the nearby metrics as 
$$
    g=b_{\varepsilon,\beta}+J_\beta k,\qquad \sigma=\varepsilon(\beta-\beta_*).
$$
Thus $(k,\sigma)=(0,0)$ corresponds to $b_\varepsilon$. The construction in Theorem~\ref{thm:glue} gives a solution $(k_\varepsilon,\sigma_\varepsilon)$ of
$$
    \mathcal F_\varepsilon(k_\varepsilon,\sigma_\varepsilon)=0.
$$
We set $\beta_\varepsilon=\beta_*+\varepsilon^{-1}\sigma_\varepsilon$.
For $\eta=\varepsilon^2A_*$, the resulting Einstein metric satisfies
$$
    \begin{aligned}
g_\eta=b_{\varepsilon,\beta_\varepsilon}+J_{\beta_\varepsilon}k_\varepsilon, 
   \qquad     A(g_\eta)=\eta.
    \end{aligned}
$$
Recall that $\mathcal L_\varepsilon=D\mathcal F_\varepsilon(0,0)$ is the augmented linearization at the approximate metric $b_\varepsilon$, whereas $D\mathcal F_\varepsilon(k_\varepsilon,\sigma_\varepsilon)$ is the augmented linearization at the Einstein metric $g_\eta$. Both are derivatives of the same map in the fixed gluing coordinates and act from $X_\varepsilon\times\mathbb R$ to $Y_\varepsilon\times\mathbb R$.

For sufficiently small $\eta>0$, define the sign of the gluing family to be the boundary sign of $g_\eta\in\partial\mathcal M_\eta$. Since $\mathcal M_\eta$ is defined by the inequality $A\geq\eta$, the outward direction along the gluing family is $-\partial_a$.

For a toric  conical Einstein metric $g$ on $X$ with $\Ric(g)=3g$, we set
$$
    \nu(g)=\dim\ker(\Delta_g-4)^{\mathbb T^2},\qquad \Delta_g=\nabla_g^*\nabla_g.
$$
As before,  the operator acts on invariant functions with regularity measured in local developments. When $g$ is smooth, the Lichnerowicz-Obata equality case gives $\nu(g)=0$ unless $g$ is the round metric on $S^4$, in which case $\nu(g)=1$. We extend this to our toric conical setting.

\begin{lemma}\label{lem:no-invariant-eigenvalue-four}
Let $g$ be a  toric conical Einstein metric on $X$ with $\operatorname{Ric}(g)=3g$. If $X$ has at least 3 torus fixed points, then $\nu(g)=0$.
\end{lemma}

\begin{proof}
Let $f\neq 0$ be a invariant smooth function with $\Delta_gf=4f$. Using the Lichnerowicz estimate we get
$$
    \nabla_g^2f=-fg,\qquad |\dd f|_g^2+f^2=C>0.
$$
Notice that a $\mathbb T^2$ invariant function smooth in local developments is a smooth function of the squared normal radius and the tangential variables along a rod, and of the two squared radial variables at a turning point. It is then easy to see that $f$ is also smooth on the underlying smooth manifold $X$.

At a critical point of $f$, from the above equation we know  $f\ne0$ and the hessian of $f$ is definite. One can see that these also hold in the underlying smooth coordinates on $X$. So $f$ is a Morse function whose critical points have index zero or four. It follows that $f$ has exactly one minimum and one maximum. But every torus fixed point is critical by invariance, contradiction. 
\end{proof}

We will prove the following sign formula in Proposition \ref{prop:gluing-sign} in several steps. By Lemma~\ref{lem:degree-orientation}, the line $D\mathcal F_\varepsilon(k_\varepsilon,\sigma_\varepsilon)^{-1}(\{0\}\times\mathbb R)$ carries an orientation through its identification with the tangent line to the moduli space at $g_a$. Increasing $a$ along the family of Einstein metrics $g_a$ determines the increasing-area direction in this line. We need to determine the sign of this direction with respect to the given orientation. We will first transport this orientation and this direction to the corresponding line at the approximate metric $b_\varepsilon$, which is non-Einstein. Then we construct an auxiliary determinant-line orientation that applies to non-Einstein reference metrics and compute its comparison with the orientation of Section~\ref{sec:moduli} at Einstein metrics. From the uniform estimates for stabilized operators we obtain an orientation-preserving comparison with the compact and bubble pieces, then a finite-dimensional reduction extends this comparison to the angle and area equations and separates the sign from $M$ and the sign from the fixed bubble $Y$. Finally, we compare the auxiliary orientation with the orientation inherited from the moduli space and reverse the increasing-area direction to obtain the boundary sign. The lemmas stated below will follow this order, and at the end we prove  Proposition~\ref{prop:gluing-sign}. 

\begin{proposition}\label{prop:gluing-sign}
There is $\epsilon_W\in\{\pm1\}$, depending only on $W$, such that for any torically nondegenerate Einstein metric $g_0$ on $M$, the metric $g_\eta$ (for $0<\eta\ll1$) given in Theorem~\ref{thm:glue} satisfies
$$
    \operatorname{sign}_{\partial\mathcal M_\eta}(g_\eta)
       =\epsilon_W(-1)^{\mu(g_0)+\nu(g_0)}.
$$
\end{proposition}

\begin{lemma}\label{lem:gluing-area-transport}
For $\varepsilon>0$ small, there is a continuous interpolation between $D\mathcal F_\varepsilon(k_\varepsilon,\sigma_\varepsilon)$ and $\mathcal L_\varepsilon$ through invertible operators. It transports the moduli space orientation to the line $\mathcal L_\varepsilon^{-1}(\{0\}\times\mathbb R)$, and the sign of $\mathcal L_\varepsilon^{-1}(0,1)$ in this inherited orientation equals the sign of the increasing-area direction at $g_\eta$, and its negative gives the induced sign on $\partial \mathcal M_\eta$.
\end{lemma}

\begin{proof}
The estimates in Theorem~\ref{thm:glue} and Proposition~\ref{prop:glueinverse} imply that the straight-line interpolation between $D\mathcal F_\varepsilon(k_\varepsilon,\sigma_\varepsilon)$ and $\mathcal L_\varepsilon$ consists of invertible operators for sufficiently small $\varepsilon$.
The local gauge description identifies $D\mathcal F_\varepsilon(k_\varepsilon,\sigma_\varepsilon)^{-1}(\{0\}\times\mathbb R)$ with $T_{g_\eta}\mathcal M$, equipped with the orientation of Lemma~\ref{lem:degree-orientation}.
Since the area component of $\mathcal F_\varepsilon$ is $\varepsilon^{-1}(A-\eta)$, the vector $D\mathcal F_\varepsilon(k_\varepsilon,\sigma_\varepsilon)^{-1}(0,1)$ has area derivative $\varepsilon>0$.
Along the interpolation, the inverse image of $(0,1)$ varies continuously and never vanishes, so its sign in the transported orientation is unchanged.
The boundary sign for $A\ge \eta$ is its negative, since the outward direction decreases area.
\end{proof}

We first construct an auxiliary determinant-line orientation for the linearized operators.
We consider an invariant metric $g$ on $M$, $X$, or $Y$, with the prescribed conical regularity and, on $Y$, the prescribed asymptotic behavior.
We do not assume that $g$ is Einstein, and we take $\lambda=3$ on $M$ and $X$ and $\lambda=0$ on $Y$.

Notice that even though $\mathscr P_g$ is Fredholm, the associated quadratic form is not coercive after constant shift -it is not bounded below in conformal directions.
To deal with this we define a new operator  $\mathscr A_g=\mathcal J_g\mathscr P_g$, where
$$
    \mathcal J_gk=k-\frac12(\tr_gk)g.
$$
Notice that $\mathcal J_g^2=\Id$. So $\ker\mathscr A_g=\ker\mathscr P_g$, while the map $[f]\mapsto[\mathcal J_gf]$ identifies $\operatorname{coker}\mathscr P_g$ with $\operatorname{coker}\mathscr A_g$. These maps identify their determinant lines:
$$
    \det\mathscr P_g\cong\det\mathscr A_g.
$$
The second-order terms of $\mathscr A_g$ are
$$
    \frac12\nabla_g^*\nabla_gk-\frac12\nabla_g^2(\tr_gk)+\frac12(\delta_g\delta_gk)g.
$$
After integration by parts, the contributions of the last two terms to $\operatorname{Re}\langle\mathscr A_gk,k\rangle_{L^2(g)}$ cancel. The remaining terms in $\mathscr A_g$ have order zero and are bounded by the curvature and the constant $\lambda$ appearing in $\mathscr E_\lambda$. Therefore
\begin{equation}\label{eq:gluing-trace-energy}
    \operatorname{Re}\langle\mathscr A_gk,k\rangle_{L^2(g)}\geq\frac12\|\nabla_gk\|_{L^2(g)}^2-C\int\bigl(|\Rm_g|+|\lambda|\bigr)|k|_g^2\dvol_g.
\end{equation}
Here the boundary terms at the conical strata vanish by regularity in local developments, and on the noncompact $Y$, we consider compactly supported tensors. 

For $g$ on $M$ or $X$, a sufficiently large positive shift makes $\mathscr A_g+sI$ coercive and invertible.
We choose such an $S$ and transport the canonical orientation of $\det(\mathscr A_g+SI)$ to $\det\mathscr A_g$ along this family.
On $Y$, we will use the compactly supported shift defined below.
We call the resulting orientation of $\det\mathscr A_g$ the \emph{auxiliary} orientation. 

The determinant lines form a continuous line bundle even when the kernel dimension changes, so the transport is well defined.
Coercivity for all sufficiently large shifts, locally uniformly in $g$, implies that the orientation is independent of $S$ and varies continuously with $g$.
Also the construction also commutes with equivariant pullbacks.

To compare this construction with Section~\ref{sec:moduli}, suppose that $g$ is Einstein and recall the operator
$$
    \mathcal D_g(t,k)=\mathscr P_gk+tB_g.
$$
We can transfer the auxiliary orientation from $\det\mathscr A_g$ to $\det\mathscr P_g$ using the identification above. Then we tensor it with the orientation of $\Lambda^1\mathbb R$ given by increasing the angle parameter, and use the homotopy
$$
    \mathcal D_{g,s}(t,k)=\mathscr P_gk+stB_g,\qquad 0\leq s\leq1,
$$
as in Lemma~\ref{lem:degree-orientation}. This defines an induced auxiliary orientation on $\det\mathcal D_g$. Thus both this orientation and the orientation of Lemma~\ref{lem:degree-orientation} are orientations of the same determinant line. At a regular point, this line is identified with the tangent line to the moduli space.

\begin{lemma}\label{lem:gluing-auxiliary-orientation}
We consider a compact toric conical Einstein metric $g$ with $\Ric(g)=3g$. The induced auxiliary orientation on $\det\mathcal D_g$ differs from the orientation of Lemma~\ref{lem:degree-orientation} by $(-1)^{1+\nu(g)}$. If $g$ is torically nondegenerate, then $\mathscr A_g$ is invertible and its auxiliary orientation has sign $(-1)^{\mu(g)+1+\nu(g)}$ relative to the canonical orientation of $\det\mathscr A_g$.
\end{lemma}

\begin{proof}
Since $g$ is Einstein, we have $\mathscr A_g=\mathscr P_g=L_g$ on transverse-traceless tensors. Suppose that $(\mathscr A_g+sI)k=0$ for $s\geq0$. Multiplication by $\mathcal J_g$ gives
$$
    \mathscr P_gk=-sk+\frac{s}{2}(\tr_gk)g.
$$
By the trace and linearized Bianchi identities  we get
$$
        (\Delta_g+s-3)\tr_gk=-2\delta_g\delta_gk, \qquad 
        (\delta_g\delta_g^*+s)\delta_gk=-\frac{s}{2}\dd(\tr_gk).
$$
Notice that on 1-forms we have  $\delta_g(\delta_g\delta_g^*)=(\Delta_g-3)\delta_g$, so 
$$
    \bigl((\Delta_g+s-3)^2-s\Delta_g\bigr)\tr_gk=0.
$$
The Lichnerowicz estimate then implies that $\operatorname{Spec}(\Delta_g)\subset\{0\}\cup[4,\infty)$ on invariant functions. For $\lambda>4$, the polynomial $(\lambda+s-3)^2-s\lambda$ has positive leading coefficient and discriminant $-3\lambda(\lambda-4)<0$. Hence $\tr_gk$ has only constant and possibly eigenvalue $4$ components.

Let $E_{\leq0}(L_g)$ denote the sum of the nonpositive eigenspaces of $L_g$. Both $\mathscr A_g$ and its adjoint preserve
$$
    E_{\leq0}(L_g)\oplus\mathbb Rg\oplus\{fg:\Delta_gf=4f\}.
$$
Indeed, $\mathscr P_g$ and $\mathcal J_g$ are self-adjoint, so $\mathscr A_g^*=\mathscr P_g\mathcal J_g$. On the three summands, $\mathscr A_g$ and $\mathscr A_g^*$ act respectively as $L_g$, $-3I$, and $-I$. Their common orthogonal complement is therefore preserved by both operators. If a tensor $k$ in this complement satisfies $(\mathscr A_g+sI)k=0$, the preceding trace calculation gives $\tr_gk=0$. Using the Bianchi identity we obtain $\delta_gk=0$. Thus $k$ is a transverse-traceless eigentensor with eigenvalue $-s\leq0$ hence must vanish. The restriction of $\mathscr A_g+sI$ to the complement has index zero, so is invertible for every $s\geq0$.

The comparison between the auxiliary orientation of $\det\mathcal D_g$ and the orientation of Lemma~\ref{lem:degree-orientation} therefore reduces to the above finite-dimensional space.  On $E_{\leq0}(L_g)$, transport from a positive shift gives the pairing orientation used in Lemma~\ref{lem:degree-orientation}. The other two summands are negative definite of $\mathscr A_g$ and $\mathscr A_g^*$ and have dimensions $1$ and $\nu(g)$, so they contribute $(-1)^{1+\nu(g)}$. The identification induced by $\mathcal J_g$ and the addition of the positively oriented angle parameter carry this same relative sign to $\det\mathcal D_g$. This proves the first assertion.

If $g$ is torically nondegenerate, then $\mathscr P_g$ and $\mathscr A_g$ are invertible. The  contribution from $E_{\leq0}(L_g)$ has sign $(-1)^{\mu(g)}$ relative to the canonical orientation, while the other two summands contribute $(-1)^{1+\nu(g)}$. Thus the auxiliary orientation of $\det\mathscr A_g$ has sign
\begin{equation}\label{eq:gluing-compact-sign}
    (-1)^{\mu(g)+1+\nu(g)}
\end{equation}
relative to its canonical orientation.
\end{proof}

We now compare the auxiliary orientations on $(X,b_\varepsilon)$, $(M,g_0)$, and $(Y,h)$. Choose $R$ sufficiently large such that $r\geq R/2$ lies in the smooth $\Ae$ end of $Y$. Let $\chi$ be a smooth invariant function with $0\leq\chi\leq1$, equal to 1 on the compact part and on $r\leq R$, and equal to 0 on $r\geq2R$. Graft $\chi$ to the scaled copy of $Y$ and extend it by zero, we  obtain a function $\chi_\varepsilon$ on $X$. For $\varepsilon>0$ small, the support of $\chi_\varepsilon$ lies in the region where $b_\varepsilon=\varepsilon^2h$. Consider the three families of operators
\begin{equation}\label{eq:gluing-shifted-families}
        \mathscr A_{g_0}+sI,\qquad \mathscr A_h+s\chi I,\qquad \mathscr A_{b_\varepsilon}+s(1+\varepsilon^{-2}\chi_\varepsilon)I,\qquad 0\leq s\leq S.
\end{equation}
All three families in \eqref{eq:gluing-shifted-families} have Fredholm index zero.

We will first choose $R$ large and then $S$ large so that the three operators at $s=S$ are invertible for all sufficiently small $\varepsilon$. Recall Hardy's inequality in dimension four 
$$
    \int_{\mathbb R^4}\frac{|v|^2}{r^2}\dd x\leq\int_{\mathbb R^4}|\nabla v|^2\dd x,\qquad v\in C_c^\infty(\mathbb R^4).
$$
Using a cutoff function, the $\Ae$ expansion and $|\Rm_h|=O(r^{-4})$  we see that  for compactly supported tensors $u$,
$$
    \begin{aligned}
        \int_{r\geq R}|\Rm_h||u|^2\dvol_h&\leq CR^{-2}\int_{r\geq R}\frac{|u|^2}{r^2}\dvol_h\\
        &\leq CR^{-2}\int_{r\geq R/2}|\nabla_hu|^2\dvol_h+CR^{-4}\int_{R/2<r<R}|u|^2\dvol_h.
    \end{aligned}
$$
For $R$ large, the gradient term is absorbed in \eqref{eq:gluing-trace-energy}. The remaining  terms are controlled by $S\int_Y\chi|u|^2\dvol_h$ when $S$ is  large.  We can apply the same argument on the annular region $\varepsilon R\leq r\leq r_0$, where $r=|x|=\varepsilon|y|$ and $r_0>0$ is fixed independently of $\varepsilon$. Using the fact that
$
    |\Rm_{b_\varepsilon}|\leq C+C\varepsilon^2r^{-4}
$
outside the grafted compact part of $Y$, we get 
$$
    \int_X|\Rm_{b_\varepsilon}||u|^2\dvol_{b_\varepsilon}\leq CR^{-2}\|\nabla u\|_{L^2}^2+C_R\|u\|_{L^2}^2+C_R\varepsilon^{-2}\int_X\chi_\varepsilon|u|^2\dvol_{b_\varepsilon}.
$$
So after choosing $R$ sufficiently large, since we add the terms $SI$ and $S\varepsilon^{-2}\chi_\varepsilon I$, we get coercivity for $S$ large.

The operator $\mathscr A_{g_0}+SI$ on $M$ can be  handled directly by \eqref{eq:gluing-trace-energy}. For the operator $\mathscr A_h+S\chi I$ on $Y$, suppose $u$ is a kernel in the weighted H\"older space, then $\mathscr P_hu=0$ outside $r=2R$. The decay argument in Proposition~\ref{prop:ae-augmented} implies that $|\nabla_h^ju|_h=O(r^{-2-j})$. Then integration by parts we see $u=0$. So we can find $R$ large and $S$ large such that all the three operators at $s=S$ are invertible.

For each family in \eqref{eq:gluing-shifted-families}, we orient the determinant line at $s=0$ by transporting the canonical orientation at the invertible endpoint $s=S$ along the family. We use these orientations in the gluing comparison below. On $X$, we will then verify that this agrees with the auxiliary orientation defined earlier using a constant shift.

We next stabilize the operators on $M$ and $Y$. We first choose an integer $N$ and continuous families of linear maps $j_0(s)$ and $j_h(s)$, with values in smooth invariant tensors supported in fixed compact subsets of $M\setminus\{p\}$ and $Y$, respectively, such that the operators
$$
    \begin{aligned}
        T_0(s):C^{2,\alpha}(M)\oplus\mathbb R^N
        &\longrightarrow C^{0,\alpha}(M),\\
        (u,e_0)&\longmapsto(\mathscr A_{g_0}+sI)u+j_0(s)e_0
    \end{aligned}
$$
and
$$
    \begin{aligned}
        T_h(s):C^{2,\alpha}_{-1}(Y)\oplus\mathbb R^N
        &\longrightarrow C^{0,\alpha}_{-3}(Y),\\
        (v,e_h)&\longmapsto(\mathscr A_h+s\chi I)v+j_h(s)e_h
    \end{aligned}
$$
are surjective for every $s\in[0,S]$. Since the original operators have index zero, the kernels $K_0(s)=\ker T_0(s)$ and $K_h(s)=\ker T_h(s)$ form rank $N$ vector bundles over the interval $[0,S]$. These bundles are trivial, so we choose continuous linear maps $\ell_0(s)$ and $\ell_h(s)$ to $\mathbb R^N$ whose restrictions to the  kernel spaces are isomorphisms. We may choose the components of $\ell_0(s)$ and $\ell_h(s)$ in the following form, for $1\le j\le N$:
$$
    \bigl(\ell_0(s)(u,e_0)\bigr)_j
    =\sum_a c^0_{ja}(s)\int_M\langle u,\gamma_{0,a}\rangle_{g_0}\dvol_{g_0}
    +\sum_{\nu=1}^N b^0_{j\nu}(s)e_{0,\nu},
$$
$$
    \bigl(\ell_h(s)(v,e_h)\bigr)_j
    =\sum_a c^h_{ja}(s)\int_Y\langle v,\gamma_{h,a}\rangle_h \dvol_h
    +\sum_{\nu=1}^N b^h_{j\nu}(s)e_{h,\nu}.
$$
Here $\{\gamma_{0,a}\}$ and $\{\gamma_{h,a}\}$ are finite families of fixed smooth invariant symmetric two-tensors compactly supported in $M\setminus\{p\}$ and $Y$, respectively, and all coefficients depend continuously on $s$. 
Consequently, the operators
$$
    (T_0(s),\ell_0(s)):C^{2,\alpha}(M)\oplus\mathbb R^N
    \longrightarrow C^{0,\alpha}(M)\oplus\mathbb R^N
$$
and
$$
    (T_h(s),\ell_h(s)):C^{2,\alpha}_{-1}(Y)\oplus\mathbb R^N
    \longrightarrow C^{0,\alpha}_{-3}(Y)\oplus\mathbb R^N
$$
are both isomorphisms, with inverses bounded uniformly for $s\in[0,S]$.
Since the original operators are invertible near $s=S$, we may arrange that for $s$ close to $S$, $j_0(s)$ and $j_h(s)$ vanish, and $\ell_0(s)(u,e_0)=e_0$, $\ell_h(s)(v,e_h)=e_h$.

We now transfer these maps to $X$, using the normalizations of Proposition~\ref{prop:glueinverse}.
For $e=(e_0,e_h)\in\mathbb R^{2N}$, we define
$$
    j_\varepsilon(s)e
    =j_0(s)e_0+\varepsilon^{-1}(\iota_\varepsilon^{-1})^*\bigl(j_h(s)e_h\bigr),
$$
and
$$
    \ell_\varepsilon(s)(k,e)
    =
    \bigl(
        \ell_0(s)(k,e_0),
        \ell_h(s)(\varepsilon^{-1}\iota_\varepsilon^*k,e_h)
    \bigr).
$$
Here we use the natural identification under the gluing construction and extend to zero outside. 
We then set
$$
    T_\varepsilon(s)(k,e)
    =
    \bigl(\mathscr A_{b_\varepsilon}
    +s(1+\varepsilon^{-2}\chi_\varepsilon)I\bigr)k
    +j_\varepsilon(s)e
$$
and obtain the enlarged operator
$$
    (T_\varepsilon(s),\ell_\varepsilon(s)):
    X_\varepsilon\oplus\mathbb R^{2N}
    \longrightarrow
    Y_\varepsilon\oplus\mathbb R^{2N}.
$$
The following can be proved by the same argument of Proposition~\ref{prop:glueinverse}, uniformly for $s\in[0,S]$.
\begin{lemma}\label{lem:gluing-stabilized-inverse}
There are constants $C>0$ and $\varepsilon_0>0$ such that
$$
    \|k\|_{X_\varepsilon}+|e|\leq C\left(\|T_\varepsilon(s)(k,e)\|_{Y_\varepsilon}+|\ell_\varepsilon(s)(k,e)|\right)
$$
for every $0<\varepsilon<\varepsilon_0$, $s\in[0,S]$, and $(k,e)\in X_\varepsilon\oplus\mathbb R^{2N}$. In particular, the operators $(T_\varepsilon(s),\ell_\varepsilon(s))$ are invertible, with inverses bounded uniformly in $s$ and $\varepsilon$.
\end{lemma}

Now we denote 
$$
    K_\varepsilon(s)=\ker T_\varepsilon(s).
$$
By Lemma~\ref{lem:gluing-stabilized-inverse}, we can define
$$
    Q_\varepsilon(s)f=(T_\varepsilon(s),\ell_\varepsilon(s))^{-1}(f,0).
$$
Thus $Q_\varepsilon(s):Y_\varepsilon\to X_\varepsilon\oplus\mathbb R^{2N}$ is uniformly bounded and satisfies
$$
    T_\varepsilon(s)Q_\varepsilon(s)=I,\qquad \ell_\varepsilon(s)Q_\varepsilon(s)=0.
$$

By choosing suitable smooth invariant cutoff functions on $M$ and on $Y$, we define a linear map
$$
    C_\varepsilon(s):K_0(s)\oplus K_h(s)\longrightarrow X_\varepsilon\oplus\mathbb R^{2N}
$$
by transferring the cut-off tensors $u$ and $\varepsilon v$ to $X$, adding them, and keeping the auxiliary coordinates $(e_0,e_h)$ unchanged. We choose the cutoffs so that the scalar conditions agree with those prescribed by $\ell_0(s)$ and $\ell_h(s)$. We then define
$$
    \Phi_\varepsilon(s)
    =\bigl(I-Q_\varepsilon(s)T_\varepsilon(s)\bigr)C_\varepsilon(s).
$$
Since $T_\varepsilon(s)Q_\varepsilon(s)=I$ and $\ell_\varepsilon(s)Q_\varepsilon(s)=0$, the map $\Phi_\varepsilon(s)$ takes values in $K_\varepsilon(s)$ and satisfies
$$
    \ell_\varepsilon(s)\Phi_\varepsilon(s)
    =\ell_0(s)\oplus\ell_h(s).
$$
The restrictions of these linear conditions to their respective kernels are isomorphisms, so $\Phi_\varepsilon(s)$ is an isomorphism depending continuously on $s$. By choosing the cutoffs suitably and then taking $\varepsilon$ small, we can  make $\|T_\varepsilon(s)C_\varepsilon(s)\|_{\mathrm{op}}$, and hence $\|\Phi_\varepsilon(s)-C_\varepsilon(s)\|_{\mathrm{op}}$, arbitrarily small uniformly for $s\in[0,S]$.

We next define the resulting comparison of determinant lines. At $s=0$, the stabilization above gives the exact sequence
$$
    0\longrightarrow\ker\mathscr A_{b_\varepsilon}
    \longrightarrow K_\varepsilon(0)
    \longrightarrow\mathbb R^{2N}
    \longrightarrow\operatorname{coker}\mathscr A_{b_\varepsilon}
    \longrightarrow0.
$$
This gives the identification
$$
    \det\mathscr A_{b_\varepsilon}\cong
    \Lambda^{\max}K_\varepsilon(0)\otimes
    \bigl(\Lambda^{\max}\mathbb R^{2N}\bigr)^*.
$$
The same construction applies to the shifted operators on $M$ and on $Y$.  Here  we use the continuous stabilization identifications of Section~\ref{sec:moduli}; see also \cite{Zinger}.

The isomorphism $\Phi_\varepsilon(0)$ identifies $K_\varepsilon(0)$ with $K_0(0)\oplus K_h(0)$. So we have naturally
$$
    \begin{aligned}
        \det\mathscr A_{b_\varepsilon}
        &\cong\Lambda^{\max}K_\varepsilon(0)\otimes\bigl(\Lambda^{\max}\mathbb R^{2N}\bigr)^*\\
        &\cong\Lambda^{\max}\bigl(K_0(0)\oplus K_h(0)\bigr)\otimes\bigl(\Lambda^{\max}(\mathbb R^{N}\oplus\mathbb R^{N})\bigr)^*\\
        &\cong\left(\Lambda^{\max}K_0(0)\otimes\bigl(\Lambda^{\max}\mathbb R^{N}\bigr)^*\right)\otimes\left(\Lambda^{\max}K_h(0)\otimes\bigl(\Lambda^{\max}\mathbb R^{N}\bigr)^*\right)\\
        &\cong\det\mathscr A_{g_0}\otimes\det\mathscr A_h.
    \end{aligned}
$$

\begin{lemma}\label{lem:gluing-determinant-comparison}
The above isomorphism induced by $\Phi_\varepsilon(0)$,
$$
    \det\mathscr A_{b_\varepsilon}\cong
    \det\mathscr A_{g_0}\otimes\det\mathscr A_h,
$$
preserves the auxiliary orientations. 
\end{lemma}

\begin{proof}
At $s=S$, all stabilizing maps vanish and the three unaugmented operators are invertible. Thus $\Phi_\varepsilon(S)$ identifies the two auxiliary vectors with $(e_0,e_h)$, and the determinant-line comparison sends $1$ to $1\otimes1$. It follows that the comparison at $s=0$ preserves the orientations transported from $s=S$.

To identify the transported orientation on $X$ with its auxiliary orientation, we fix $\varepsilon$ and choose $S'$ sufficiently large that $\mathscr A_{b_\varepsilon}+S'I$ is coercive. The convex interpolation with $\mathscr A_{b_\varepsilon}+S(1+\varepsilon^{-2}\chi_\varepsilon)I$ remains coercive and invertible. Hence the two paths defining the orientations are homotopic through Fredholm operators, with invertible initial operators and common terminal operator $\mathscr A_{b_\varepsilon}$, so the orientations agree.
\end{proof}

We now extend Lemma~\ref{lem:gluing-determinant-comparison} to include the angle and area equations. We evaluate all stabilizations and kernel spaces at $s=0$, and we will omit this parameter, and write $t$ for the angle variable $\sigma$.

On $X$, we use $\mathscr A_{b_\varepsilon}=\mathcal J_{b_\varepsilon}\mathscr P_{b_\varepsilon}$ and first solve
$$
    T_\varepsilon(k,e)+t\mathcal J_{b_\varepsilon}B_\varepsilon=0.
$$
Since $T_\varepsilon Q_\varepsilon=I$, every solution is uniquely parametrized by $(z,t)\in K_\varepsilon\oplus\mathbb R$ through
$$
    (k,e)=z-tQ_\varepsilon(\mathcal J_{b_\varepsilon}B_\varepsilon).
$$
We then impose $e=0$ and the area equation using the finite-dimensional map
$$
    \begin{aligned}
        \mathcal R_\varepsilon:K_\varepsilon\oplus\mathbb R
        &\longrightarrow\mathbb R^{2N}\oplus\mathbb R,\\
        \mathcal R_\varepsilon(z,t)
        &=\bigl(e,a_\varepsilon(k)+c_\varepsilon t\bigr).
    \end{aligned}
$$
The invertibility of $\mathcal L_\varepsilon$ implies that $\mathcal R_\varepsilon$ is invertible.
Moreover, $\mathcal R_\varepsilon(z,t)=(0,1)$ is equivalent to $e=0$ and $(k,t)=\mathcal L_\varepsilon^{-1}(0,1)$.

We next carry out the same construction on $M$ and $Y$.
We denote by $B_h$ and $c_0$ the derivatives at $\beta=\beta_*$ of $\mathscr E_0(J_\beta^Yh)$ and $A(J_\beta^Yh)$, respectively.
As in Proposition~\ref{prop:glueinverse}, the rescaling $v=\varepsilon^{-1}\iota_\varepsilon^*k$ identifies the limiting operator with
\begin{equation}\label{eq:gluing-model-operator}
    (u,v,t)\longmapsto
    \bigl(\mathscr A_{g_0}u,\ \mathscr A_hv+t\mathcal J_hB_h,\ \dd A_h(v)+c_0t\bigr),
\end{equation}
where $u\in C^{2,\alpha}(M)$ and $v\in C^{2,\alpha}_{-1}(Y)$.
This operator is invertible by toric nondegeneracy on $M$ and Proposition~\ref{prop:ae-augmented} on $Y$.

We define
$$
    Q_hf=(T_h,\ell_h)^{-1}(f,0).
$$
For $((u,e_0),(v,e_h))\in K_0\oplus K_h$, we set
$$
    (v_t,e_{h,t})=(v,e_h)-tQ_h(\mathcal J_hB_h).
$$
The corresponding finite-dimensional map is
$$
    \begin{aligned}
        \mathcal R_{\mathrm{mod}}:(K_0\oplus K_h)\oplus\mathbb R
        &\longrightarrow\mathbb R^{2N}\oplus\mathbb R,\\
        \mathcal R_{\mathrm{mod}}\bigl((u,e_0),(v,e_h),t\bigr)
        &=\bigl((e_0,e_{h,t}),\dd A_h(v_t)+c_0t\bigr).
    \end{aligned}
$$
As on $X$, this map is invertible, and its inverse at $(0,1)$ determines a triple $(u,v_t,t)$ whose image under \eqref{eq:gluing-model-operator} is $(0,0,1)$.

Both parametrizations respect the angle-parameter convention of Lemma~\ref{lem:degree-orientation}, since scaling the angle terms continuously from zero identifies them with the addition of the positive angle variable.
We use $\Phi_\varepsilon(0)$ to identify the kernel spaces and set
$$
    \widehat{\mathcal R}_\varepsilon
    =\mathcal R_\varepsilon\circ\bigl(\Phi_\varepsilon(0)\oplus I\bigr).
$$
The two maps now have the same domain and range, so we consider
$$
    \mathcal R_\tau
    =(1-\tau)\mathcal R_{\mathrm{mod}}+\tau\widehat{\mathcal R}_\varepsilon,
    \qquad 0\leq\tau\leq1.
$$
Through the stabilization exact sequences, the auxiliary orientations and the standard orientations of the auxiliary coordinates give orientations on $K_0$ and $K_h$.
Together with the positive angle direction, these orient $(K_0\oplus K_h)\oplus\mathbb R$ according to Lemma~\ref{lem:degree-orientation}.
With $\pi(e,a)=e$, whenever $\mathcal R_\tau$ is invertible, the exact sequence
$$
    0\longrightarrow\mathcal R_\tau^{-1}(\{0\}\times\mathbb R)
    \longrightarrow(K_0\oplus K_h)\oplus\mathbb R
    \xrightarrow{\pi\mathcal R_\tau}\mathbb R^{2N}
    \longrightarrow0
$$
therefore orients the line containing $\mathcal R_\tau^{-1}(0,1)$.
We will compare the signs of these vectors at $\tau=0$ and $\tau=1$.

\begin{lemma}\label{lem:gluing-augmented-deformation}
For suitable cutoff functions and all sufficiently small $\varepsilon$, every $\mathcal R_\tau$ is invertible, and the sign of $\mathcal R_\tau^{-1}(0,1)$ is independent of $\tau\in[0,1]$. In particular,  the increasing-area direction for \eqref{eq:gluing-model-operator} and $\mathcal L_\varepsilon^{-1}(0,1)$ have the same sign in the orientations defined above.

\end{lemma}

\begin{proof}
We choose the cutoffs equal to 1 near the finite rod spheres and on the fixed supports of the angle variation, stabilizing tensors, and test tensors.
The cutoff construction then preserves the auxiliary coordinates and the normalized area variation.
Since the area functionals are uniformly bounded, the smallness of $\Phi_\varepsilon(0)-C_\varepsilon(0)$ implies that $\widehat{\mathcal R}_\varepsilon$ and $\mathcal R_{\mathrm{mod}}$ are arbitrarily close on $t=0$.

To compare the terms proportional to $t$, we set
$$
    (v,e_h)=Q_h(\mathcal J_hB_h).
$$
We apply the same cutoff construction to this pair and the zero pair on $M$, obtaining $(k,e)$ on $X$ with $\ell_\varepsilon(k,e)=0$.
The angle variations agree under rescaling, so the same cutoff estimates make $\|T_\varepsilon(k,e)-\mathcal J_{b_\varepsilon}B_\varepsilon\|_{Y_\varepsilon}$ arbitrarily small.
Here we first choose the transitions sufficiently slowly in $\log r$ and then take $\varepsilon$ small.
Lemma~\ref{lem:gluing-stabilized-inverse} therefore implies that
$$
    \left\|(k,e)-Q_\varepsilon(\mathcal J_{b_\varepsilon}B_\varepsilon)\right\|_{X_\varepsilon\oplus\mathbb R^{2N}}
    \leq C\left\|T_\varepsilon(k,e)-\mathcal J_{b_\varepsilon}B_\varepsilon\right\|_{Y_\varepsilon}.
$$
Since the cutoff construction preserves the normalized area variation and $c_\varepsilon=c_0$, this proves the required closeness for the terms proportional to $t$ as well.

Thus we can choose the cutoffs and then $\varepsilon$ so that
$$
    \left\|\mathcal R_{\mathrm{mod}}^{-1}
    \bigl(\widehat{\mathcal R}_\varepsilon-\mathcal R_{\mathrm{mod}}\bigr)\right\|_{\mathrm{op}}<1.
$$
It follows that every $\mathcal R_\tau$ is invertible. The vectors $\mathcal R_\tau^{-1}(0,1)$ vary continuously and never vanish, so their signs in the continuously oriented lines are constant.

At $\tau=1$, applying $\Phi_\varepsilon(0)\oplus I$ and then the reconstruction above takes $\mathcal R_1^{-1}(0,1)$ to $\mathcal L_\varepsilon^{-1}(0,1)$. This identification preserves orientations by Lemma~\ref{lem:gluing-determinant-comparison} and the angle-parameter convention.
\end{proof}

\begin{lemma}\label{lem:gluing-auxiliary-sign}
There is a sign $\epsilon_W\in\{\pm1\}$, depending only on $W$, such that, for all sufficiently small $\varepsilon$, the vector $\mathcal L_\varepsilon^{-1}(0,1)$ has sign
$$
    (-\epsilon_W)(-1)^{\mu(g_0)+1+\nu(g_0)}
$$
in the orientation induced by the auxiliary orientation of $\det\mathscr A_{b_\varepsilon}$ and the positive angle parameter.
\end{lemma}

\begin{proof}
By Proposition~\ref{prop:ae-augmented}, we take the unique pair $(v,t)$ satisfying
$$
    \begin{cases}
        \mathscr A_hv+t\mathcal J_hB_h=0,\\
        \dd A_h(v)+c_0t=1.
    \end{cases}
$$
We define $-\epsilon_W$ to be its sign in the kernel of the first equation, oriented by the auxiliary orientation of $\det\mathscr A_h$ and the positive angle parameter. This sign is unchanged by equivariant pullback and positive rescaling. Thus uniqueness of the metric $(Y,h)$ associated with $W$ implies that $\epsilon_W$ depends only on $W$.

In \eqref{eq:gluing-model-operator}, the metric variation $u$ on $M$ satisfies $\mathscr A_{g_0}u=0$, so $u=0$ by invertibility of $\mathscr A_{g_0}$. By \eqref{eq:gluing-compact-sign}, its auxiliary determinant orientation differs from its canonical orientation by
$$
    (-1)^{\mu(g_0)+1+\nu(g_0)}.
$$
The sign of the increasing-area vector $(0,v,t)$ is therefore the product of this factor and $-\epsilon_W$. Lemma~\ref{lem:gluing-augmented-deformation} transfers this sign to $\mathcal L_\varepsilon^{-1}(0,1)$.
\end{proof}

\begin{lemma}\label{lem:gluing-reference-orientation}
For all sufficiently small $\varepsilon$, the orientation on $\mathcal L_\varepsilon^{-1}(\{0\}\times\mathbb R)$ inherited from the moduli space in Lemma~\ref{lem:gluing-area-transport} differs from the orientation induced by the auxiliary orientation of $\det\mathscr A_{b_\varepsilon}$ and the positive angle parameter by
$$
    (-1)^{1+\nu(g_\eta)}.
$$
\end{lemma}

\begin{proof}
By Lemma~\ref{lem:gluing-auxiliary-orientation}, the moduli-space orientation and the orientation induced by $\det\mathscr A_{g_\eta}$ and the positive angle parameter differ by $(-1)^{1+\nu(g_\eta)}$ on
$$
    D\mathcal F_\varepsilon(k_\varepsilon,\sigma_\varepsilon)^{-1}(\{0\}\times\mathbb R).
$$
Thus we only need to check that the latter orientation transports to the orientation induced by $\det\mathscr A_{b_\varepsilon}$ and the positive angle parameter.

Fix $\varepsilon$.
We first change the metric in the gauge term while keeping $g_\eta$ fixed.
For $0\leq t\leq1$, we consider
$$
    \mathscr E_3(h)-\delta_h^*\delta_h
    \bigl((1-t)J_\beta J_{\beta_\varepsilon}^{-1}g_\eta+t b_{\varepsilon,\beta}\bigr)=0.
$$
At $\beta=\beta_\varepsilon$, the metric in the gauge term is $(1-t)g_\eta+t b_{\varepsilon,\beta_\varepsilon}$.
Since $g_\eta$ is Einstein and both $\delta_{g_\eta}g_\eta$ and $\delta_{g_\eta}b_{\varepsilon,\beta_\varepsilon}$ vanish, $h=g_\eta$ solves this equation for every $t$.

At $\beta=\beta_\varepsilon$ we differentiate the left-hand side with respect to $h$ at $h=g_\eta$, obtaining
$$
    k\longmapsto
    \mathscr P_{g_\eta}k-t\delta_{g_\eta}^*
    \bigl((D\delta)_{g_\eta}[k](b_{\varepsilon,\beta_\varepsilon}-g_\eta)\bigr).
$$
At $t=1$, this operator is obtained by differentiating the first component of $\mathcal F_\varepsilon$ with respect to $k$ at $(k_\varepsilon,\sigma_\varepsilon)$ and conjugating by $J_{\beta_\varepsilon}$.
We next consider the linear interpolation
$$
    k\longmapsto
    (1-t)\left[
        \mathscr P_{g_\eta}k-\delta_{g_\eta}^*
        \bigl((D\delta)_{g_\eta}[k](b_{\varepsilon,\beta_\varepsilon}-g_\eta)\bigr)
    \right]
    +tJ_{\beta_\varepsilon}\mathscr P_{b_\varepsilon}J_{\beta_\varepsilon}^{-1}k,
    \qquad 0\leq t\leq1.
$$
This is the first component of the conjugated interpolation in Lemma~\ref{lem:gluing-area-transport}, with the angle variation set to zero.

We compose the operators obtained by varying the gauge term on the left with $\mathcal J_{g_\eta}$.
For the interpolation to $J_{\beta_\varepsilon}\mathscr P_{b_\varepsilon}J_{\beta_\varepsilon}^{-1}$, we compose on the left with
$$
    (1-t)\mathcal J_{g_\eta}
    +tJ_{\beta_\varepsilon}\mathcal J_{b_\varepsilon}J_{\beta_\varepsilon}^{-1}.
$$
Together, these compositions define a continuous family from $\mathscr A_{g_\eta}$ to $J_{\beta_\varepsilon}\mathscr A_{b_\varepsilon}J_{\beta_\varepsilon}^{-1}$. The gluing estimates bound the changes in the second-order coefficients and the difference $\mathcal J_{g_\eta}-J_{\beta_\varepsilon}\mathcal J_{b_\varepsilon}J_{\beta_\varepsilon}^{-1}$ by $O(\varepsilon^{1/3})$ in rescaled charts. Thus the maps $(1-t)\mathcal J_{g_\eta}+tJ_{\beta_\varepsilon}\mathcal J_{b_\varepsilon}J_{\beta_\varepsilon}^{-1}$ are invertible for all $t\in[0,1]$ and sufficiently small $\varepsilon$. The estimate in \eqref{eq:gluing-trace-energy}, with the small changes in the second-order coefficients absorbed into the gradient term, applies to the composed operators uniformly along the interpolation for each fixed $\varepsilon$.
We therefore choose a common $s>0$ such that adding $sI$ makes every composed operator coercive and invertible.

Now we compare the orientations in the diagram
$$
\begin{tikzcd}[column sep=large,row sep=large]
    \det(\mathscr A_{g_\eta}+sI)
    \arrow[r,"\text{shifted path}"]
    \arrow[d,"\text{shift to zero}"']
    &
    \det\bigl(J_{\beta_\varepsilon}(\mathscr A_{b_\varepsilon}+sI)J_{\beta_\varepsilon}^{-1}\bigr)
    \arrow[d,"\text{shift to zero}"]
    \\
    \det\mathscr A_{g_\eta}
    \arrow[r,"\text{unshifted path}"']
    &
    \det\bigl(J_{\beta_\varepsilon}\mathscr A_{b_\varepsilon}J_{\beta_\varepsilon}^{-1}\bigr).
\end{tikzcd}
$$
Along the upper edge, all operators are invertible, so transport preserves their canonical orientations. The vertical arrows define the auxiliary orientations, with conjugation by $J_{\beta_\varepsilon}$ identifying the right-hand orientation with that of $\det\mathscr A_{b_\varepsilon}$. Since the determinant lines extend over the rectangle of path and shift parameters, the square commutes at the level of orientations. Thus the lower edge transports the auxiliary orientation of $\det\mathscr A_{g_\eta}$ to that of $\det\mathscr A_{b_\varepsilon}$.

The gauge identifications of Section~\ref{sec:moduli} preserve the moduli-space orientation and leave the angle parameter unchanged. Thus, by the convention of Lemma~\ref{lem:degree-orientation}, transporting the orientation induced by $\det\mathscr A_{g_\eta}$ and the positive angle parameter yields the orientation induced by $\det\mathscr A_{b_\varepsilon}$ and the same parameter. So its relative sign with the orientation inherited from $\mathcal M$ on $\mathcal L_\varepsilon^{-1}(\{0\}\times\mathbb R)$ is still $(-1)^{1+\nu(g_\eta)}$.
\end{proof}

\begin{remark}
The possible multiplicity of boundary points is a symmetry issue. Even when every compact limit is isometric to $g_0$, different contraction points with their ordered incident rods need not be related by an isometry compatible with the torus action, even allowing a change of lattice basis. Thus $\partial\mathcal{M}_{0,\eta}$ may a priori contain several points. The oriented count handles this multiplicity: Proposition~\ref{prop:gluing-sign} assigns every such boundary point the same sign $\epsilon_W(-1)^{\mu(g_0)+\nu(g_0)}$.
\end{remark}

\begin{proof}[Proof of Proposition~\ref{prop:gluing-sign}]
Since $X$ has $m_0+q\geq3$ torus fixed points, Lemma~\ref{lem:no-invariant-eigenvalue-four} implies that $\nu(g_\eta)=0$. Thus Lemma~\ref{lem:gluing-reference-orientation} implies that the inherited orientation and the auxiliary orientation on $\mathcal L_\varepsilon^{-1}(\{0\}\times\mathbb R)$ are opposite. The outward direction for $A\geq\eta$ also has the opposite sign to the increasing-area direction. These two signs cancel, so Lemma~\ref{lem:gluing-area-transport} identifies the boundary sign with the sign of $\mathcal L_\varepsilon^{-1}(0,1)$ in the auxiliary orientation. By Lemma~\ref{lem:gluing-auxiliary-sign}, we therefore obtain
$$
    \begin{aligned}
        \operatorname{sign}_{\partial\mathcal M_\eta}(g_\eta)
        =(-\epsilon_W)(-1)^{\mu(g_0)+1+\nu(g_0)}=\epsilon_W(-1)^{\mu(g_0)+\nu(g_0)},
    \end{aligned}
$$
where $\epsilon_W$ depends only on $W$.
\end{proof}

\begin{proof}[Proof of Theorem~\ref{thm:dichotomy}]
We first normalize $\Ric(g_0)=3g_0$, and suppose that the second alternative does not hold. Fix an admissible sequence $W$ in \eqref{eq:alternating-compact-word}, and use the notation of Section~\ref{sec:moduli}. Let $\mathcal M_0$ be the connected component of ${\mathcal M}$ containing the family $g_a$ of Theorem~\ref{thm:glue} for sufficiently small $a>0$.

Suppose that $\mathcal M_0$ contains no point with $\beta=1$. Since $\beta(g_a)\to\beta_*>1$, connectedness gives $\beta>1$ throughout $\mathcal M_0$. We let $\mathcal{M}_{0,\eta}=\mathcal{M}_0\cap \mathcal{M}_\eta$. For every $\eta>0$, Theorems~\ref{thm:upperangle} and~\ref{thm:compactness} imply that $\mathcal M_{0,\eta}$ is compact. Indeed, the cone angles are bounded above and below, and the only degeneration in Theorem~\ref{thm:compactness} has $\mathcal A\to0$. Every other limit belongs to $\mathcal M_{0,\eta}$, since a connected component is closed.

By our assumption, every compact Einstein metric arising in such a degeneration is isometric to $g_0$. The given two-tori are maximal in the respective isometry groups. By conjugacy of maximal tori, each limiting metric admits an isometry to $g_0$ that is equivariant up to an automorphism of $\T^2$. Consequently, all these limits are torically nondegenerate and have the same value $\mu(g_0)$ and $\nu(g_0)$.

Theorem~\ref{thm:glue} and Proposition~\ref{prop:converse} therefore imply that, for sufficiently small $\eta>0$, every point of $\partial\mathcal M_{0,\eta}$ is regular and $\mathcal A$ is a local coordinate there. This boundary is nonempty, since it contains the class of $g_\eta$. Proposition~\ref{prop:parity} and Proposition~\ref{prop:gluing-sign} now give
$$
    0 =\sum_{x\in\partial\mathcal M_{0,\eta}} \operatorname{sign}_{\partial\mathcal M_{0,\eta}}(x) =\epsilon_W(-1)^{\mu(g_0)+\nu(g_0)} \#\partial\mathcal M_{0,\eta} \neq0,
$$
a contradiction. Thus $\mathcal M_0$ contains a metric with $\beta=1$. All its cone angles are $2\pi$, so it is a smooth toric Einstein metric on $X$.
\end{proof}

\begin{proof}[Proof of Theorem~\ref{thm:main}]
We equip
$
    \mathbb S^4$, $\mathbb S^2\times \mathbb S^2$ and $\CP^2\#k\overline{\CP^2}\quad(0\leq k\leq3)
$
with their classical Einstein metrics and standard torus actions, normalized by $\Ric(g)=3g$. Their rod directions belong to $\{A,B,C\}$ after a global change of lattice basis. Corollary~\ref{cor:toric-small-betti-classification} gives uniqueness up to isometry, and Theorem~\ref{thm:toric-small-betti-nondegeneracy} gives toric nondegeneracy. Thus the second alternative of Theorem~\ref{thm:dichotomy} is excluded, and every admissible replacement produces a toric Einstein metric.

By the topological calculations in Section~\ref{sec:moduli}, the alternating sequences applied to $\mathbb S^4$ and $\CP^2\#k\overline{\CP^2}$ give 
$$
    \#_r(\mathbb S^2\times \mathbb S^2),\qquad (\CP^2\#k\overline{\CP^2})\#_r(\mathbb S^2\times \mathbb S^2), \qquad r\geq1,\quad 0\leq k\leq3.
$$
The cases $r=0$ are the classical metrics themselves.

The four sequences in \eqref{eq:hermitian-words}, applied to the six manifolds above, give the finite list of connected sums $X\#Y$. The admissible choices of corner and lattice identification give precisely the listed choices of $X$ and $Y$, with the indicated orientations.
\end{proof}

\section{Proof of Theorem \ref{thm:general}}
\label{sec:toric-existence-criterion}

We prove that the improved Hitchin-Thorpe inequality is sufficient for the existence of a toric Einstein metric on the underlying manifold, thereby establishing Theorem~\ref{thm:general}. We compare the virtual counts of toric Einstein metrics mod 2 before and after replacing a corner of the rod structure. Compared with the proof of Theorem~\ref{thm:main}, we must now carry out the gluing construction simultaneously over all toric Einstein metrics associated with a given rod structure. In particular, we cannot assume toric nondegeneracy of the compact Einstein metric, since we have established this only for algebraically special metrics. We address this difficulty using finite-dimensional stabilization for Fredholm equations, as in Ozuch's gluing modulo obstructions framework \cite{OzuchII}. To make the stabilization compatible with gauge fixing and gluing, we use Delay's localized divergence correction \cite{Delay2012} to choose divergence-free stabilizing tensors that vanish near the chosen fixed point. We recall the statement of Theorem \ref{thm:general}.

\begin{theorem}\label{thm:toric-existence-criterion}
Let $X$ be a  simply connected smooth compact four-manifold admitting an effective $\mathbb T^2$-action. Then $X$ admits a toric Einstein metric if and only if
\begin{equation}\label{eq:toric-existence-criterion}
  \chi(X)\geq3|\tau(X)|.
\end{equation}
\end{theorem}

Notice that the torus action in this theorem is allowed to vary. For prescribed torus actions we prove the following stronger (virtual) counting statement. For a given compact rod structure $\mR$, we denote by $\mathcal M_{\mR}$ the moduli space of toric Einstein metrics $g$ on $X_{\mR}$ with $\Ric(g)=3g$, modulo equivalence relation given by pulling back by an orientation-preserving diffeomorphism commuting with the fixed $\mathbb T^2$ action and preserving each sphere $\Sigma_\alpha$. The number $\mathfrak n(\mR)$  is the virtual modulo two count we define below.

\begin{theorem}\label{thm:prescribed-three-direction-existence}
If the rod directions of $\mR$ belong to $\{A,B,C\}$ after a global change of lattice basis, then
$$
  \mathfrak n(\mR)\equiv1\pmod2.
$$
In particular, $X_{\mR}$ admits a toric Einstein metric for the $\mathbb T^2$ prescribed action.
\end{theorem}
\begin{remark}\label{rmk:more than one}
If a smooth compact four-manifold admits more than one non-equivalent torus action with rod directions in $\{A, B, C\}$, then it follows that it admits more than one toric Einstein metrics. For example, on $2\CP^2\#2\overline {\CP}^2$, there are at least two non-equivalent torus actions, whose corresponding rod vector sequences are $\{ABABAC\}$ and $\{ABACBC\}$. In general it is a combinatorial question to determine the number of non-equivalent $\mathbb T^2$ actions on a given smooth four-manifold.
\end{remark}

\subsection{Definition of the virtual count}

\begin{lemma}\label{lem:three-direction-smooth-compactness}
The moduli space $\mathcal M_{\mR}$ is compact whenever $\mR$ has rod directions in $\{A,B,C\}$.
\end{lemma}

\begin{proof}
By Theorem~\ref{geom:noncollapse} and Myers' theorem, we have  volume noncollapsing and diameter upper bound. The curvature energy is also uniformly bounded by the Gauss-Bonnet theorem. We claim the curvature is uniformly bounded, then Proposition~\ref{prop:bounded-curvature-compactness} gives smooth convergence after orientation-preserving diffeomorphisms commuting with the fixed torus action and preserving each rod sphere.

If the maximum curvature were unbounded, rescaling at a curvature maximum would give a complete smooth nonflat Ricci-flat limit with finite curvature energy and Euclidean volume growth, by similar argument as in Theorem~\ref{thm:compactness}. One can also see that there is an effective limiting torus action with the original primitive rod directions. The limit is $\ALE$ by \cite[Theorem~1.5]{BKN}, contradicting Lemma~\ref{deg:smooth-ae-rigidity}.
\end{proof}

We next define the virtual count $\mathfrak n(\mR)\in\mathbb Z/2\mathbb Z$ for a compact rod structure with directions in $\{A,B,C\}$. We first assume that $\mathcal M_{\mR}$ is nonempty.

At each $[g]\in\mathcal M_{\mR}$, we identify the cokernel of $\mathscr P_g$ with the space of smooth invariant infinitesimal Einstein deformations, by Lemma~\ref{lem:conical-linear}. By compactness and the extension construction in the proof of Proposition~\ref{prop:parity}, we choose finitely many $C^1$ section $\tau_1,\ldots,\tau_N$ whose values span these cokernels along $\mathcal M_{\mR}$. We take these sections to be $C^1$ on a common open neighborhood of invariant $L^2_8$ metrics, with values in $L^2_6$. We require
$$
  \delta_h\tau_\nu(h)=0,\qquad \tau_\nu(F^*h)=F^*\tau_\nu(h)
$$
for every admissible  diffeomorphism $F$. We construct these sections by transporting smooth tensors through the divergence gauges as in the proof of Proposition~\ref{prop:localmoduli},  and applying the projection
$$
  \mathsf Q_h\tau =\tau-\delta_h^*(\delta_h\delta_h^*)_\perp^{-1}\delta_h\tau,
$$
and multiplying by cutoffs on the space of metrics. The inverse is taken on invariant one-forms $L^2(h)$-orthogonal to the metric duals of the two vector fields generating the $\mathbb T^2$-action.

Now we fix a reference Einstein metric $g$ and consider the stabilized equation
\begin{equation}\label{eq:three-direction-stabilized-equation}
  \mathscr E_3(h)-\delta_h^*\delta_h g =\sum_{\nu=1}^N u_\nu\tau_\nu(h), \qquad u\in\mathbb R^N.
\end{equation}
Since the perturbing tensors are divergence-free, we can apply the contracted Bianchi identity and integrate by parts, as in Proposition~\ref{prop:localmoduli}, to see that a solution $(h, u)$ to \eqref{eq:three-direction-stabilized-equation} satisfies that
$$
  \delta_hg=0; \qquad \mathscr E_3(h)=\sum_{\nu=1}^N u_\nu\tau_\nu(h).
$$
At $([g],0)\in\mathcal M_{\mR}\times\{0\}$, we obtain the linearization
$$
  (k,v)\longmapsto \mathscr P_gk-\sum_{\nu=1}^N v_\nu\tau_\nu(g).
$$
Since the values $\tau_\nu(g)$ span the cokernel of $\mathscr P_g$, this map is surjective. It has index $N$, so its kernel has dimension $N$. Then we can apply the implicit-function theorem to obtain a local $C^1$ solution manifold of dimension $N$.

On overlaps, equivariance identifies solutions in different divergence gauges without changing $u$. The argument in the proof of Proposition~\ref{prop:parity} gives $C^1$ coordinate transitions between the local solution manifolds. Since the only admissible diffeomorphisms fixing a metric are torus translations, which act trivially on invariant tensors, we obtain a $C^1$ manifold $\mathcal W$ near $\mathcal M_{\mR}\times\{0\}$ and a $C^1$ map
$$
  p:\mathcal W\longrightarrow\mathbb R^N,\qquad p([h],u)=u,\qquad p^{-1}(0)=\mathcal M_{\mR}.
$$

By compactness, we can choose an open neighborhood $\mathcal V\subset\subset\mathcal W$ of $\mathcal M_{\mR}$ and $c>0$ such that $|p|\geq c$ on $\partial\mathcal V$. For a regular value $u$ with $|u|<c$, we obtain a compact zero-dimensional manifold $p^{-1}(u)\cap\mathcal V$, and hence a finite set.

\begin{definition}
We define the mod 2 virtual count
\begin{equation}\label{eq:three-direction-count-definition}
  \mathfrak n(\mR) =\deg_2(p,\mathcal V,0) =\#\bigl(p^{-1}(u)\cap\mathcal V\bigr)\pmod2,
\end{equation}
where $u$ is any sufficiently small regular value. If $\mathcal M_{\mR}$ is empty, we define $\mathfrak n(\mR)=0$.
\end{definition}

By standard arguments using a common stabilization and homotopy invariance of the degree mod 2, we can see that $\mathfrak n(\mR)$ is independent of the choice of divergence gauges, stabilizing sections, and regular perturbations. We can also see directly that $\mathfrak n(\mR)$ is unchanged if we cyclically relabel the rods, make a global change of lattice basis, or reverse the orientation.

\begin{remark}
If every metric in $\mathcal M_{\mR}$ is torically nondegenerate, then $\mathcal M_{\mR}$ is discrete and compact, hence finite. In this case we have
$$
  \mathfrak n(\mR)=\#\mathcal M_{\mR}\pmod2.
$$
Without nondegeneracy, we can only regard $\mathfrak n(\mR)$ as a virtual count. In general it is possible that $\mathcal M_{\mR}$ may not be discrete. In all cases, we have
$$
  \mathfrak n(\mR)\neq0 \quad\Longrightarrow\quad \mathcal M_{\mR}\neq\varnothing.
$$
\end{remark}

\subsection{The main gluing formula}

Now we fix a replacement $\mR=\mR_0\#_{\mathfrak p}W$ from \eqref{eq:alternating-compact-word}, and write $M=X_{\mR_0}$ and $X=X_{\mR}$. Contracting the intermediate rods of $W$ recovers the corner $\mathfrak p$ of $\mR_0$. The first and last rods of $W$ correspond to the two rods incident to $\mathfrak p$, with directions $B,A$ in boundary order. We identify every other rod with its corresponding rod in $\mR_0$, preserving the cyclic order and lattice directions. Throughout, we use the same torus lattice and the orientations fixed above. Denote by $p_0\in M$ the torus fixed point lying over $\mathfrak p$.

The main result of this section is the following gluing formula

\begin{theorem}\label{thm:gluing formula}
If the rod directions of $\mR_0$ belong to $\{A,B,C\}$, then
\begin{equation}\label{eq:parametrized-parity-law}
  \mathfrak n(\mR_0\#_{\mathfrak p}W) \equiv\mathfrak n(\mR_0)\pmod2.
\end{equation}
\end{theorem}

Assuming Theorem~\ref{thm:gluing formula}, we first prove Theorem~\ref{thm:prescribed-three-direction-existence} and then deduce Theorem~\ref{thm:toric-existence-criterion}.

\begin{proof}[Proof of Theorem~\ref{thm:prescribed-three-direction-existence}]
We prove the count by induction on the number of rods, reducing each rod structure by the inverse of an allowed gluing operation.

For two or three rods, Corollary~\ref{cor:toric-small-betti-classification} identifies every normalized Einstein metric with the corresponding standard metric. Every symmetry of its rod structure is realized by an isometry normalizing the standard torus. So we can  adjust the identifying isometry so that it commutes with the fixed torus action and preserves each rod sphere. Thus $\mathcal M_{\mR}$ consists of a single point, which is torically nondegenerate by Theorem~\ref{thm:toric-small-betti-nondegeneracy}. Consequently, $\mathfrak n(\mR)\equiv1\pmod2$.

We now assume that there are at least four rods. Every permutation of the three unsigned directions is induced by a global lattice change, so we may relabel them in the following reductions. If two labels at distance two in the cyclic list agree, four consecutive labels can be written as either $BABA$ or $BCBA$. We apply the corresponding reduction
$$
  BABA\longmapsto BA
$$
or
$$
  BCBA\longmapsto BA.
$$
Otherwise, every three consecutive labels are distinct, so the cyclic list is $(BCA)^d$ after relabeling. Since there are at least four rods, we have $d\ge2$, and we apply
$$
  BCABCA\longmapsto BA.
$$

Each reduction preserves the endpoint labels, so adjacent labels remain distinct, including across the ends of the replaced block. We therefore obtain a compact rod structure with fewer rods and directions in $\{A,B,C\}$. The inverse operations are precisely replacements listed in \eqref{eq:alternating-compact-word}. Theorem~\ref{thm:gluing formula} therefore preserves $\mathfrak n$ modulo two under each reduction. Induction gives $\mathfrak n(\mR)\equiv1\pmod2$. In particular, $\mathcal M_{\mR}$ is nonempty, which proves existence for the prescribed action.
\end{proof}

\begin{proof}[Proof of Theorem~\ref{thm:toric-existence-criterion}]
The necessity part is Theorem~\ref{thm:improved HT}. For the sufficiency part, the connected-sum classification recalled in Section~\ref{sec:rod structure} reduces the manifold to
$$
  \#_r(\mathbb S^2\times\mathbb S^2),\quad r\geq0, \qquad\text{or}\qquad p\CP^2\#q\overline{\CP^2},\quad p,q\geq0, \quad p+q\geq1.
$$
Here $r=0$ denotes $\mathbb S^4$. We construct a rod structure with directions in $\{A,B,C\}$ on each manifold satisfying \eqref{eq:toric-existence-criterion}.

The alternating list $(AB)^{r+1}$ gives the first family by the topological calculation in Section~\ref{sec:moduli}. For the second family, reverse orientation if necessary so that $q\geq p$, and set $d=q-p$. If $d\geq1$, the inequality is equivalent to $p\geq d-1$. The cyclic list $(ABC)^d$ has every rod sphere of square $-1$, so
$$
  \chi=3d,\qquad \tau=-d,\qquad (b_2^+,b_2^-)=(d-1,2d-1).
$$
Its intersection form is odd. The connected-sum classification in Section~\ref{sec:rod structure} therefore identifies its oriented smooth manifold as
\begin{equation}\label{eq:three-direction-cycle-topology}
  (d-1)\CP^2\#(2d-1)\overline{\CP^2}.
\end{equation}
Apply the alternating replacement $p-d+1$ times. Each adds $\mathbb S^2\times\mathbb S^2$, increasing both $b_2^+$ and $b_2^-$ by one. The intersection form remains odd, so the resulting smooth manifold is $p\CP^2\#q\overline{\CP^2}$. All rod directions remain in $\{A,B,C\}$.

If $d=0$, then $p=q\geq1$. We start with $ABAC$, whose underlying manifold is $\CP^2\#\overline{\CP^2}$, and apply the alternating replacement $p-1$ times. Theorem~\ref{thm:prescribed-three-direction-existence} gives the desired toric Einstein metric in all cases.
\end{proof}

\begin{remark}
Theorem~\ref{thm:toric-existence-criterion} concerns existence for some torus action on $M$. The inequality alone does not give existence for every prescribed action. For example, the action on the second Hirzebruch surface has rod spheres of square $2$ and $-2$, so it is excluded by Proposition~\ref{lem:rod-self-intersection-bound1}. Its underlying manifold is $\mathbb S^2\times\mathbb S^2$, which admits the round product metric with a different action.
\end{remark}

The rest of this section is devoted to the proof of Theorem \ref{thm:gluing formula}.

\

Let us explain the main ideas first. We will compare the two boundaries of the conical moduli space with $\beta\geq1$ and $A\geq\eta$, where $\eta>0$ is small. At $\beta=1$, the metrics are smooth on $X$. At $A=\eta$, the rod spheres with small area lie in a small region obtained by gluing a rescaled gravitational instanton into $M$ near $p_0$. We therefore expect these boundaries to contribute $\mathfrak n(\mR)$ and $\mathfrak n(\mR_0)$, respectively.

The main difficulty is that the Einstein metrics on $M$ may be degenerate. We stabilize the equation by divergence-free tensors supported away from $p_0$, so the perturbed metrics remain Einstein near the gluing point. This allows us to carry out gluing over the entire stabilized solution space and establish the small-area count. The compactness and angle bounds then give a compact region between the two boundaries. A regular perturbation gives a compact one-dimensional manifold, whose boundary counts agree mod 2.

\subsection{Localization of stabilized sections}

We begin with the required localization of the stabilizing sections.

\begin{proposition}\label{prop:localized-compact-stabilization}
We can choose finitely many equivariant sections $\tau_1,\ldots,\tau_N$, defined on a common $L^2_8$ neighborhood of the invariant metrics representing $\mathcal M_{\mR_0}$, and a constant $\rho>0$. These sections are continuous with values in $L^2_7$ and $C^1$ with values in $L^2_6$, and they satisfy
$$
  \delta_{h}\tau_\nu(h)=0, \qquad \tau_\nu(h)\big|_{B_{h}(p_0,\rho)}=0
$$
for every invariant metric $h$ in this neighborhood and every $\nu=1,\ldots,N$. Moreover, at each Einstein metric $g$ representing a point of $\mathcal M_{\mR_0}$, the $L^2(g)$-orthogonal projections of $\tau_1(g),\ldots,\tau_N(g)$ onto $\ker\mathscr P_g$ span this kernel.
\end{proposition}

\begin{remark}
Our arguments below do not give $C^1$ dependence with values in $L^2_7$. This is a common issue regarding derivative loss for the action of diffeomorphism group. Differentiating transport by a varying diffeomorphism introduces a Lie derivative of $\tau_\nu(h)$. Using the $L^2_7$ regularity of these tensors, we control this term in $L^2_6$. Thus our estimates give continuity as maps $L^2_8\to L^2_7$ and continuous differentiability as maps $L^2_8\to L^2_6$.
\end{remark}

\begin{proof}
We need to modify the divergence-free sections constructed in the definition of $\mathfrak n(\mR)$. We fix a reference Einstein metric $g$ and work in its divergence gauge. For each tensor $\tau_0$ in a basis of $\ker\mathscr P_g$, we define
$$
  \tau(h)=\mathsf Q_h\tau_0,
$$
where $\mathsf Q_h$ is the divergence-free projection defined above. Since $\delta_g\tau_0=0$, we have $\tau(g)=\tau_0$, while $\delta_h\tau(h)=0$ for every nearby invariant metric $h$. We localize this section in the divergence gauge before transporting it to other representatives.

Choose $r$ sufficiently small and an invariant cutoff $\psi$ equal to zero on $B_g(p_0,r)$ and one outside $B_g(p_0,2r)$, with its transition compactly contained in $\{5r/4<d_g(p_0,\cdot)<7r/4\}$. Set
$$
  \Omega=B_g(p_0,2r)\setminus\overline{B_g(p_0,r)}.
$$
Let $Z_1,Z_2$ be the two vector fields generating the $\mathbb T^2$-action. Since $\delta_{h}\tau(h)=0$ and $(\psi-1)\tau(h)$ vanishes near the boundary of $B_g(p_0,2r)$, integration by parts gives, for $i=1,2$,
$$
  \begin{aligned}
    \int_\Omega \delta_{h}(\psi\tau(h))(Z_i) \dvol_{h} &=\frac12\int_{B_g(p_0,2r)} \bigl\langle(\psi-1)\tau(h), \mathcal L_{Z_i}h\bigr\rangle_{h} \dvol_{h}\\
                                                        &=0.
  \end{aligned}
$$
Here $\mathcal L_{Z_i}h=0$ because $h$ is invariant. Thus we can apply Lemma~\ref{lem:annular-divergence-correction} below to get an invariant symmetric two-tensor $s(h)$, supported in $\overline\Omega$ and extended by zero to $M$, such that
$$
  \delta_{h}s(h)=\delta_{h}(\psi\tau(h)), \qquad \|s(h)\|_{L^2(h)} \le Cr\|\delta_{h}(\psi\tau(h))\|_{L^2(h)}.
$$
Choose $\psi$ by rescaling a fixed cutoff. Since $|\dd\psi|_{h}\le C/r$, $\delta_{h}\tau(h)=0$, and $\operatorname{Vol}_{h}(B_g(p_0,2r))\le Cr^4$, we obtain
\begin{equation}\label{eq:localized-source-error}
  \|\psi\tau(h)-s(h)-\tau(h)\|_{L^2(h)} \le Cr^2\|\tau(h)\|_{C^0(B_g(p_0,2r),h)}.
\end{equation}
The constants are uniform for $h$ in a sufficiently small neighborhood. The modified section $\psi\tau(h)-s(h)$ is divergence-free and vanishes on $B_g(p_0,r)$.

For $k=7,8$, the divergence-free projection gives $C^1$ dependence of $\tau(h)$ on $h\in L^2_k$ with values in $L^2_{k-1}$. Since $\delta_h\tau(h)=0$, the expression $\delta_h(\psi\tau(h))$ involves only the contraction of $\tau(h)$ with $\dd\psi$ and has the same regularity. Lemma~\ref{lem:annular-divergence-correction} therefore gives $C^1$ dependence of the zero extension of $s(h)$ in $L^2_{k-1}$. The modified local section $\psi\tau(h)-s(h)$ has the same regularity.

We apply \eqref{eq:localized-source-error} to each basis element $\tau_0\in\ker\mathscr P_g$. Since $\tau(g)=\tau_0$, the $L^2(g)$-orthogonal projections of the corresponding tensors $\psi\tau(g)-s(g)$ onto $\ker\mathscr P_g$ still span this space when $r$ is sufficiently small. This also holds at every nearby Einstein metric, after shrinking the neighborhood.

By compactness, we choose finitely many Einstein metrics $g_\alpha$ and open neighborhoods $\mathcal V_\alpha\subset\mathcal U_\alpha$ of $[g_\alpha]$ in the space of invariant $L^2_8$ metrics modulo equivalence, such that
$$
  \overline{\mathcal V_\alpha}\subset\mathcal U_\alpha
$$
and
$$
  \mathcal M_{\mR_0}\subset\bigcup_\alpha\mathcal V_\alpha.
$$
On $\mathcal U_\alpha$, we use representatives $h$ satisfying $\delta_hg_\alpha=0$. We choose cutoffs $\chi_\alpha$ equal to one on $\mathcal V_\alpha$ and supported in $\mathcal U_\alpha$. We multiply the local sections constructed near $g_\alpha$ by $\chi_\alpha$, transport them to other representatives through the divergence gauges, and extend them by zero outside $\mathcal U_\alpha$.

The map assigning to a metric its representative satisfying $\delta_hg_\alpha=0$ is continuous in $L^2_8$ and $C^1$ with values in $L^2_7$. The estimates at $k=8$ and $k=7$ therefore give continuity of the local sections in $L^2_7$ and $C^1$ dependence in $L^2_6$ after this substitution. Differentiating the pullback uses one spatial derivative of the section, so the transported sections have the same regularity. Multiplication by $\chi_\alpha$ preserves these properties. Since $\chi_\alpha=1$ on $\mathcal V_\alpha$, the resulting sections have $L^2(g)$-orthogonal projections spanning $\ker\mathscr P_g$ at every $[g]\in\mathcal M_{\mR_0}$.

We denote by $r_\alpha$ the radius used in the localization near $g_\alpha$. We choose the neighborhoods $\mathcal U_\alpha$ small enough that their representatives in the chosen gauges are uniformly equivalent to $g_\alpha$. Since there are finitely many $\alpha$, we can choose $\rho>0$ such that
$$
  B_h(p_0,\rho)\subset B_{g_\alpha}(p_0,r_\alpha)
$$
for every such representative $h$ in $\mathcal U_\alpha$. The local sections vanish on these balls. The admissible diffeomorphisms fix $p_0$, so the transported sections also vanish on $B_h(p_0,\rho)$.

\end{proof}

The following special case of Delay's construction is used in the above proof; see \cite[Section~9.3]{Delay2012}. We give a self-contained proof here for the convenience of readers.

\begin{lemma}\label{lem:annular-divergence-correction}
We use the notation of the preceding proposition. For $k=7,8$, let $h$ be an invariant metric in a sufficiently small $L^2_k$-neighborhood of $g$. For all sufficiently small $r>0$ and every invariant one-form $f\in L^2_{k-2}$ supported in $\operatorname{supp}(\dd\psi)$ and satisfying
\begin{equation}\label{eq:compatibility}
  \int_\Omega f(Z_i) \dvol_h=0,\qquad i=1,2,
\end{equation}
we can choose an invariant symmetric two-tensor $s\in L^2_{k-1}(M)$ such that
\begin{equation}\label{eq:conclusion}
  \begin{aligned}
    \operatorname{supp}s&\subset\overline\Omega,\\
    \delta_hs           &=f,\\
    \|s\|_{L^2(h)}      &\le Cr\|f\|_{L^2(h)}.
  \end{aligned}
\end{equation}
For fixed $r$, $s$ depends $C^1$ on $(h,f)\in L^2_k\times L^2_{k-2}$, and  $s$ is smooth on $M$ whenever $h$ and $f$ are smooth.
\end{lemma}

\begin{proof}
We use equivariant normal coordinates for $g$ and rescale by $r$. The annulus becomes $A=\{1<|x|<2\}$, the action is given by independent rotations in two orthogonal planes in $\R^4$, and the rescaled metric, still denoted by $h$, is $C^1$-close to $g_{\mathrm{Euc}}$. We set $K=\{5/4\le|x|\le7/4\}$, which contains the rescaled support of $f$.

For $A_t=\{1+t<|x|<2-t\}$, $0\le t\le1/8$, Korn's inequality gives
$$
  \|V\|_{L^2_1(A_t,g_{\mathrm{Euc}})} \le C_t\bigl(\|\mathcal L_Vg_{\mathrm{Euc}}\|_{L^2(A_t,g_{\mathrm{Euc}})} +\|V\|_{L^2(A_t,g_{\mathrm{Euc}})}\bigr).
$$
Here $\mathcal L_Vg_{\mathrm{Euc}}$ has components $\partial_iV_j+\partial_jV_i$, so the inequality controls the full first derivative by its symmetric part and the $L^2$ norm of $V$. For invariant vector fields we impose the condition
\begin{equation}\label{eq:annular-normalization}
  \int_K\langle V,Z_i\rangle_{g_{\mathrm{Euc}}} \dd x=0, \qquad i=1,2.
\end{equation}
Then we can remove $\|V\|_{L^2}$ in the above inequality. Indeed, otherwise a sequence with unit $L^2_1$ norm and symmetric derivatives tending to zero would converge weakly in $L^2_1$ and strongly in $L^2$ to an invariant Euclidean Killing vector field satisfying \eqref{eq:annular-normalization}. Every Euclidean Killing vector field is affine, and invariance under the two independent rotations makes it a constant linear combination of $Z_1,Z_2$. The normalization makes this limit zero, contradicting Korn's inequality. Thus
\begin{equation}\label{eq:annular-korn}
  \|V\|_{L^2_1(A_t,g_{\mathrm{Euc}})} \le C\|\mathcal L_Vg_{\mathrm{Euc}}\|_{L^2(A_t,g_{\mathrm{Euc}})},
\end{equation}
where $C$ is uniform in $t$.

Since
$$
  |\mathcal L_Vh-\mathcal L_Vg_{\mathrm{Euc}}|_{g_{\mathrm{Euc}}} \le C\|h-g_{\mathrm{Euc}}\|_{C^1} (|V|_{g_{\mathrm{Euc}}}+|\partial V|_{g_{\mathrm{Euc}}}),
$$
we can absorb the error when $\|h-g_{\mathrm{Euc}}\|_{C^1}$ is sufficiently small, and obtain
$$
  \|V\|_{L^2_1(A_t,g_{\mathrm{Euc}})} \le C\|\mathcal L_Vh\|_{L^2(A_t,h)}
$$
for every invariant vector field satisfying \eqref{eq:annular-normalization}, with $C$ independent of $t\in[0,1/8]$ and of $h$ in the chosen neighborhood.

We choose a smooth nondecreasing function $\chi:[0,\infty)\to[0,\infty)$ with $\chi(0)=0$, $\chi(t)=\exp(-2/t)$ for $t>0$ small, and $\chi(t)$ constant for $t\ge1/8$. On $A$, we denote
$$
  w(x)=\chi(\min\{|x|-1,2-|x|\}\bigr).
$$
It is smooth and positive on $A$, constant near $K$, and vanishes to infinite order at $\partial A$. For  $x\in A$, we have
$$
  w(x)=\int_0^{1/8}\chi'(t)\mathbf 1_{A_t}(x) \dd t.
$$
We square the previous estimate on $A_t$, integrate it against $\chi'(t) \dd t$, and apply Fubini's theorem to obtain
\begin{equation}\label{eq:weighted-korn}
  \begin{aligned}
    \|V\|_H^2 &:=\int_Aw (|V|_{g_{\mathrm{Euc}}}^2+ |\partial V|_{g_{\mathrm{Euc}}}^2) \dd x\\
              &\le C\int_Aw|\mathcal L_Vh|_h^2 \dvol_h.
  \end{aligned}
\end{equation}
We define $H$ as the completion in this norm of the invariant vector fields smooth on $\overline A$ and satisfying \eqref{eq:annular-normalization}. Notice that $H$ is independent of $h$.

The 1-form $f$ is supported in $K$, on which we know $w$ is bounded below. The Lax-Milgram theorem then yields a unique $V\in H$ such that
\begin{equation}\label{eq:variational}
  \frac12\int_Aw\langle\mathcal L_Vh,\mathcal L_Wh\rangle_{h}  \dvol_{h} =\int_Af(W) \dvol_{h}\quad \text{for}\  W\in H,
\end{equation}
and
$$
  \|V\|_H\le C\|f\|_2.
$$
Set $s=w\mathcal L_Vh$ on $A$ and $s=0$ outside $A$. Then we have
$$
  \frac12\int\langle s,\mathcal L_Wh\rangle_{h} \dvol_{h} =\int f(W) \dvol_{h}
$$
for every compactly supported smooth $W$. This implies $\delta_{h}s=f$ in the sense of distributions, including across $\partial A$. Taking $W=V$ in \eqref{eq:variational}, and using that $w$ is bounded, we see that $\|s\|_2\le C\|f\|_2$.

Interior elliptic regularity gives $V\in C^\infty(A)$, and $\delta_{h}(w\mathcal L_Vh)=0$ near $\partial A$. For $p\in\partial A_t$ with $t>0$ small, apply interior estimates on balls of radius $ct^2$, where $c>0$ is fixed and small. On these balls, $w/w(p)$ is bounded above and below, with uniformly bounded derivatives after rescaling. Since $w(p)=\exp(-2/t)$, these estimates and $\|V\|_H\le C\|f\|_2$ give
\begin{equation}\label{eq:flatness}
  |\partial^js(p)| \le C_jt^{-2j-6}\exp(-1/t)\|f\|_2,\qquad j\ge0.
\end{equation}
So every derivative of $s$ tends to zero at $\partial A$, and $s$, defined to be zero outside $A$, is smooth across $\partial A$.

We fix $k\in\{7,8\}$. Since $H$ is fixed, uniform coercivity of \eqref{eq:variational} and $C^1$ dependence of its bilinear form and right-hand side imply that $V$ depends $C^1$ in $H$ on $(h,f)\in L^2_k\times L^2_{k-2}$. The interior and weighted boundary estimates applied to $V$ and its first parameter derivative then imply that the zero extension of $s$ depends $C^1$ in $L^2_{k-1}(M)$.
\end{proof}

\subsection{The enlarged moduli space}

Now we choose the sections $\tau_1,\ldots,\tau_N$ from Proposition~\ref{prop:localized-compact-stabilization} and consider the stabilized equation for $(h, u)$
\begin{equation}\label{eq:compact-universal-equation}
  \mathscr E_3(h)=\sum_{\nu=1}^N u^\nu\tau_\nu(h), \qquad u\in\mathbb R^N.
\end{equation}
We denote by $\mathcal W$ the set of solutions modulo equivalence near $\mathcal M_{\mR_0}\times\{0\}$, with the quotient topology induced by $L^2_8\times\mathbb R^N$.

We now define local coordinates on $\mathcal W$. Fix an Einstein representative $g$ and use the divergence gauge $\delta_{h}g=0$. We consider the gauged map $\widetilde{\mathscr E}$ given by
$$
  (h,u)\longmapsto \mathscr E_3(h) -\delta_{h}^*\delta_{h}g -\sum_{\nu=1}^N u^\nu\tau_\nu(h).
$$
Since the sections $\tau_\alpha$ are divergence-free, as before, an integration by parts argument identifies ther zeroes of $\widetilde{\mathscr E}$ with the solutions of \eqref{eq:compact-universal-equation} in this divergence gauge.

At $(g,0)$, our choice of $\tau_\alpha$'s imply that the linearization of $\widetilde{\mathscr E}$ is surjective and has  an $N$-dimensional kernel. Now we choose smooth invariant test tensors $\gamma_1,\ldots,\gamma_N$, supported away from $p_0$, and constants $b_{j\nu}$ such that
$$
  \ell(k,v)= \left( \int_M\langle k,\gamma_j\rangle_g \dvol_g +\sum_{\nu=1}^N b_{j\nu}v^\nu \right)_{j=1}^N
$$
restricts to an isomorphism on $\ker D\widetilde{\mathscr E}_{(g,0)}$. The support condition can be achieved by using a cut-off near $p_0$ and the fact that the kernel is finite-dimensional.

For a nearby solution, we choose its representative $h$ satisfying $\delta_hg=0$ and then define local coordinates by
$$
  \xi([h],u)=\ell(h-g,u).
$$
After shrinking the  neighborhood, each nearby equivalence class in $\mathcal W$ has a unique representative $h$ in this neighborhood satisfying $\delta_{h}g=0$. Hence the coordinate map $\xi$ is well defined.

\begin{lemma}
After restricting to a sufficiently small neighborhood of $\mathcal M_{\mR_0}\times\{0\}$, the maps defined above give $\mathcal W$ the structure of an $N$-dimensional $C^1$ manifold.
\end{lemma}

\begin{proof}
Since $(h,u)\mapsto(\widetilde{\mathscr E}(h,u),\ell(h-g,u))$ has invertible differential at $(g,0)$, the implicit-function theorem implies that $\xi$ defines local $C^1$ coordinates. The transition maps are $C^1$ because the gauge changes depend $C^1$ with values in $L^2_7$ and the functionals defining $\ell$ are continuous at this order.
\end{proof}

We define the parameter projection $p:\mathcal W\to\mathbb R^N$ by
$$
  p([h],u)=u.
$$
It is a $C^1$ map, and $p^{-1}(0)$ is naturally identified with $\mathcal M_{\mR_0}$. We choose a small open neighborhood $\mathcal V$ of $p^{-1}(0)$ with compact closure in $\mathcal W$. We may assume that any $([h], u)$ in $\mathcal V$ has the property that $\tau_\nu(h)=0$ vanishes on $B_h(p_0, \rho)$ (for $\rho$ given in Proposition \ref{prop:localized-compact-stabilization}). In each local chart, we write $\xi=([g_{\xi}],u_{\xi})$, where $g_{\xi}$ is the unique nearby representative satisfying $\delta_{g_{\xi}}g=0$ for the fixed reference metric $g$ of that chart.

Denote by $\mathcal L_{\xi}:L^2_8\times\mathbb R^N\to L^2_6$ the linearization of $\widetilde{\mathscr E}$.  After possibly shrinking $\mathcal V$, we may assume $\mathcal L_{\xi}$ is surjective with $N$-dimensional kernel for $\xi\in \mathcal V$. Using a finite cover and a $C^1$ partition of unity, we may combine the local test functionals $\ell$ chosen above into an equivariant $C^1$ family
$$
  F_{\xi}:L^2_8\times\mathbb R^N\longrightarrow\mathbb R^d
$$
which is injective on every $\ker\mathcal L_{\xi}$.

Now we have orthogonal projections $\pi_{\xi}$ onto $F_{\xi}(\ker\mathcal L_{\xi})$ which vary continuously on $\xi$. Then we choose $C^1$ rank-$N$ orthogonal projections $\Pi_{\xi}$ which are uniformly close to $\pi_{\xi}$. Set
$$
  \mathcal E_{\xi}=\operatorname{Im}\Pi_{\xi}, \qquad \ell_{\xi}=\Pi_{\xi}F_{\xi}.
$$
Then $\ell_{\xi}$ restricts to an isomorphism on $\ker\mathcal L_{\xi}$. Under the identification of $T_{\xi}\mathcal W$ with $\ker\mathcal L_{\xi}$ by divergence-free representatives, we obtain
\begin{equation}\label{eq:smooth-test-rows}
  \ell_{\xi}|_{T_{\xi}\mathcal W}:T_{\xi}\mathcal W\longrightarrow\mathcal E_{\xi} \quad\text{is an isomorphism}.
\end{equation}
Therefore,
$$
  (\mathcal L_{\xi},\ell_{\xi}): L^2_8\times\mathbb R^N \longrightarrow L^2_6\times\mathcal E_{\xi}
$$
is an isomorphism with uniformly bounded inverses over $\overline{\mathcal V}$. The maps $\ell_{\xi}$ are equivariant and depend $C^1$ on $\xi$, and they are defined using test tensor supported uniformly away from $p_0$.

\subsection{The uniform gluing construction}

Let $(Y,h)$ be the $\Ae$ toric conical  gravitational instanton associated with $W$ in Section~\ref{sec:gluing}. Its two smooth infinite rods correspond, in order, to the $B$- and $A$-rods at $\mathfrak p$, and its finite rods have common cone angle $2\pi\beta_*$, where $\beta_*>1$. Fix a scale for $h$ and denote by  $A_*=A(h)$ the total area of its finite rod spheres.

For $\xi=([g_{\xi}],u_{\xi})\in\overline{\mathcal V}$, we use the  representative $g_{\xi}$ on $M$ chosen in the preceding subsection. We construct $b_{\varepsilon,\xi}$ as in Section~\ref{sec:gluing}, interpolating between $\varepsilon^2h$ and $g_{\xi}$ on the annulus $\varepsilon^{2/3}<r<2\varepsilon^{2/3}$ in the $g_{\xi}$-normal coordinates at $p_0$. We then set
$$
  b_{\varepsilon,\xi,\beta}=J_\beta b_{\varepsilon,\xi},
$$
where we recall that $J_\beta$ changes the cone angle along $\mathcal S$ from $2\pi\beta_*$ to $2\pi\beta$ and is the identity outside the chosen neighborhood of $\mathcal S$. In particular,
$$
  A(b_{\varepsilon,\xi,\beta_*})=\varepsilon^2A_*.
$$

We fix a large integer $m$ and $0<\alpha<1$, and denote by $X_\varepsilon$ and $Y_\varepsilon$ the completions of the space of invariant smooth tensors (with respect to cone angle $2\pi \beta_*$) in the norms of $C^{m+2,\alpha}_{-1}$ and $C^{m,\alpha}_{-3}$, respectively, as in Section \ref{sec:gluing}. Notice that the weights hence the norm depend on $\varepsilon$ whereas the underlying vector spaces are independent of $\varepsilon>0$..

We choose $R>0$ uniformly small so that the test tensors defining $\ell_\xi$ and the tensors $\tau_\nu(h)$ vanish on $\overline{B_{g_\xi}(p_0,2R)}$ for every $\nu$ and every metric $h$ sufficiently close to $g_\xi$ on $M$. We choose a $\mathbb T^2$-invariant cutoff $\chi_{\xi}$, depending $C^1$ on $\xi$, which is zero on $B_{g_{\xi}}(p_0, R)$ and one outside $B_{g_{\xi}}(p_0, 2R)$. For $\varepsilon>0$ sufficiently small, our gluing construction will only be carried out inside $B_{g_{\xi}}(p_0, R)$. We can also naturally identify $M\setminus B_{g_{\xi}}(p_0, R)$ with its image in the glued manifold $X$. This allows us to freely graft tensors on $M$ and on $X$ which are supported in this region. We will use this identification throughout.

We write a nearby metric on $X$ as
$$
  g=b_{\varepsilon,\xi,\beta}+J_\beta k,
$$
where $k\in X_\varepsilon$ is the unknown invariant symmetric two-tensor at cone angle $2\pi\beta_*$. Now we regard $\chi_{\xi}k$ as a tensor on $M$ by extending it by zero across $B_{g_{\xi}}(p_0, R)$, and set
\begin{equation}\label{eq:filled-compact-metric}
  \widehat g_{\xi}=g_{\xi}+\chi_{\xi}k.
\end{equation}
This metric equals $g_{\xi}$ near $p_0$ and agrees with $g$ outside $B_{g_{\xi}}(p_0, 2R)$. Then we have the tensor $\tau_\nu(\widehat g_{\xi})$ on $M$ and we can naturally view it as a tensor $X$, which we denote by $\widehat\tau_{\nu,\varepsilon,\xi}(g)$. Since $\widehat g_{\xi}$ agrees with $g$ on a neighborhood of the support of $\tau_\nu(\widehat g_{\xi})$, the divergence-free property of $\tau_\nu$ implies that
\begin{equation}\label{eq:exact-reference-divergence}
  \begin{aligned}
    \delta_g\widehat\tau_{\nu,\varepsilon,\xi}(g)                &=0,\\
    \widehat\tau_{\nu,\varepsilon,\xi}(b_{\varepsilon,\xi,\beta})&=\tau_\nu(g_{\xi}).
  \end{aligned}
\end{equation}
For fixed $\varepsilon$ and $\xi$, the tensors $\widehat\tau_{\nu,\varepsilon,\xi}(g)$ depend $C^1$ on $k$ and are independent of $\beta$. The first derivatives with respect to $k$ are uniformly bounded and uniformly continuous as maps from $X_\varepsilon$ to $Y_\varepsilon$, for $\xi\in\overline{\mathcal V}$ and sufficiently small $\varepsilon$.

We denote
$$
  \sigma=\varepsilon(\beta-\beta_*).
$$
For fixed $\varepsilon$ and $\xi$, we choose a sufficiently small open neighborhood $\mathcal U_{\varepsilon,\xi}$ of $(0,0,0)$ in $X_\varepsilon\times\mathbb R\times\mathbb R^N$ and define a map
$$
  \mathcal F_{\varepsilon,\xi}: \mathcal U_{\varepsilon,\xi} \longrightarrow Y_\varepsilon\times\mathbb R\times\mathcal E_{\xi}
$$
by
\begin{equation}\label{eq:augmented-family-map}
  \mathcal F_{\varepsilon,\xi}(k,\sigma,v)=
  \begin{pmatrix}
    J_\beta^{-1}\left\{ \mathscr E_3(g)-\delta_g^*\delta_gb_{\varepsilon,\xi,\beta} -\displaystyle\sum_{\nu=1}^N (u_{\xi}^\nu+v^\nu)\widehat\tau_{\nu,\varepsilon,\xi}(g) \right\}\\[2mm]
    \varepsilon^{-1}\bigl(A(g)-\varepsilon^2A_*\bigr)\\[2mm]
    \ell_{\xi}(k,v)
  \end{pmatrix}.
\end{equation}
The first two components of $\mathcal F_{\varepsilon,\xi}=0$ express the gauged stabilized Einstein equation in the divergence gauge with $A(g)=\varepsilon^2A_*$. The last component selects the auxiliary parameter $\xi$. We will show below that, near the constructed solutions with $u=0$, this condition determines $\xi$ uniquely from $([g],u)$, where $u=u_{\xi}+v$.

We denote
$$
  N_\varepsilon(k,\sigma,v) =\|k\|_{X_\varepsilon}+|\sigma|+|v|.
$$

\begin{proposition}\label{prop:stabilized-gluing}
After shrinking $\mathcal V$, there exist $\varepsilon_0,C>0$ such that, for every $\xi\in\overline{\mathcal V}$ and $0<\varepsilon<\varepsilon_0$, the equation $\mathcal F_{\varepsilon,\xi}(k,\sigma,v)=0$ has a unique solution in the ball $N_\varepsilon(k,\sigma,v)\le C\varepsilon^{4/3}$, which depends $C^1$ on $(\varepsilon, \xi)$, and satisfies
\begin{equation}\label{eq:family-correction}
  N_\varepsilon(k_\varepsilon,\sigma_\varepsilon, v_\varepsilon)\le C\varepsilon^{4/3}.
\end{equation}
Denote by $g_\varepsilon(\xi)$ the metric associated to the data $(k_\varepsilon, \sigma_\varepsilon, v_\varepsilon)$ and denote $p_\varepsilon(\xi)=p(\xi)+v_\varepsilon(\xi)$. Then $\xi\mapsto[g_\varepsilon(\xi)]$ gives a bijection from $p_\varepsilon^{-1}(0)$ onto the set of toric conical Einstein metrics on $X$ with common cone angle $2\pi\beta$ along $\mathcal S$ for some $\beta\ge1$, smooth along the other rod spheres, and satisfying $A(g)=\varepsilon^2A_*$, modulo equivalence.
\end{proposition}

\begin{proof}
The proof is a bit long. We divide into three steps.

\medskip
\noindent\textbf{Step 1. Uniform invertibility of the linearized operator.}
\medskip

Fix finitely many charts of $\mathcal W$ centered at points $([g],0)\in p^{-1}(0)$, using representatives satisfying $\delta_{g_\xi}g=0$, and choose local frames of $\mathcal E$. The estimates for $\mathsf Q_h$ and the correction in Lemma~\ref{lem:annular-divergence-correction} imply that $h\mapsto\tau_\nu(h)$ is $C^1$ from $C^{m+2,\alpha}$ metrics to $C^{m,\alpha}$ tensors. It is also easy to see that  $\xi\mapsto g_\xi$ is $C^1$ with values in $C^{m+2,\alpha}(M)$.

We set
$$
  \mathcal L_{\varepsilon,\xi} =D\mathcal F_{\varepsilon,\xi}(0)
$$
and denote
$$
  \mathcal L^0_{\varepsilon,\xi}(k,\sigma,v) = \mathcal L_{\varepsilon,\xi}(k,\sigma,v) +
  \begin{pmatrix}
    \displaystyle\sum_\nu u_{\xi}^\nu D\widehat\tau_{\nu,\varepsilon,\xi} (b_{\varepsilon,\xi,\beta_*})[k]\\
    0\\
    0
  \end{pmatrix}.
$$
As $\varepsilon\to0$, convergence on compact subsets of $M\setminus\{p_0\}$ identifies the limiting operator as
$$
  (k,v)\longmapsto \left( \mathscr P_{g_{\xi}}k-\sum_\nu v^\nu\tau_\nu(g_{\xi}), \ell_{\xi}(k,v) \right).
$$
At $\xi\in p^{-1}(0)$, this operator equals $(\mathcal L_\xi,\ell_\xi)$ and is therefore invertible by the construction of $\ell_\xi$. Then we obtain a relatively compact neighborhood $\mathcal V_0$ of $p^{-1}(0)$ on which these inverses are uniformly bounded.

The proof of Proposition~\ref{prop:glueinverse} now applies to $\mathcal L^0_{\varepsilon,\xi}$, uniformly for $\xi\in\overline{\mathcal V_0}$. Indeed, in a contradicting sequence, the limit on $Y$ has no stabilizing terms. Proposition~\ref{prop:ae-augmented} therefore eliminates the limiting tensor and angle variable. On $M\setminus\{p_0\}$, the limiting tensor has growth $O(r^{-1})$ and extends across $p_0$ by removal of singularities. The inverse above on $M$ then eliminates this tensor and the limiting variable $v$. The remaining possible limit lies on a flat annulus and is excluded by the weight $-1$ argument in Proposition~\ref{prop:glueinverse}.

For fixed $\varepsilon>0$, the operator has Fredholm index zero, hence
$$
  \|(\mathcal L^0_{\varepsilon,\xi})^{-1}\|\le C_0
$$
uniformly on $\overline{\mathcal V_0}$ for small $\varepsilon$.

Now since $\chi_\xi$ vanishes on $B_{g_\xi}(p_0,R)$, the uniform $C^1$ bounds for $\tau_\nu$ imply that $\mathcal L^0_{\varepsilon,\xi}-\mathcal L_{\varepsilon,\xi}$ has operator norm at most $C_1|u_\xi|$, uniformly in $\varepsilon$ and $\xi$. We shrink $\mathcal V\Subset\mathcal V_0$ so that
$$
  C_0C_1\sup_{\overline{\mathcal V}}|u_{\xi}|<\frac12.
$$
Then we get a uniform inverse for $\mathcal L_{\varepsilon,\xi}$.

\medskip
\noindent\textbf{Step 2. Construction of the solutions.}
\medskip

At $(k,\sigma,v)=0$, the area and $\ell_\xi$ components of $\mathcal F_{\varepsilon,\xi}$ vanish, while the first component vanishes outside the gluing region by \eqref{eq:compact-universal-equation} and \eqref{eq:exact-reference-divergence}. Since $g_\xi$ is Einstein near $p_0$, the estimate in the proof of Theorem~\ref{thm:glue} implies that
$$
  \|\mathcal F_{\varepsilon,\xi}(0)\|\le C\varepsilon^{4/3},
$$
uniformly in $\xi$. For $z=(k,\sigma,v)$, the uniform inverse from {\bf Step 1} and the nonlinear estimates in that proof, together with the uniform $C^1$ dependence of $\widehat\tau_{\nu,\varepsilon,\xi}$ on $k$, imply that
$$
  z\longmapsto z-\mathcal L_{\varepsilon,\xi}^{-1}\mathcal F_{\varepsilon,\xi}(z)
$$
maps the ball $N_\varepsilon(z)\le C\varepsilon^{4/3}$ into itself and is a contraction for sufficiently large $C$ and sufficiently small $\varepsilon$. Thus there is a unique solution in this ball satisfying \eqref{eq:family-correction}, and the implicit-function theorem implies the  solution depends $C^1$ on $(\varepsilon, \xi)$.

By \eqref{eq:exact-reference-divergence}, the Bianchi identity, and the integration-by-parts argument in Theorem~\ref{thm:glue}, the gauge term vanishes, so
$$
  \mathscr E_3(g_\varepsilon(\xi))=\sum_\nu p_\varepsilon^\nu(\xi)\widehat\tau_{\nu,\varepsilon,\xi}(g_\varepsilon(\xi)).
$$
The area equation implies that $A(g_\varepsilon(\xi))=\varepsilon^2A_*$, while \eqref{eq:family-correction} implies that
$$
  |\beta_\varepsilon(\xi)-\beta_*|\le C\varepsilon^{1/3}.
$$
Since $\beta_*>1$, a zero of $p_\varepsilon$ therefore determines a conical  Einstein metric with $\beta>1$ for sufficiently small $\varepsilon$. Equivariance and uniqueness identify the local constructions on the overlaps of different charts.

\medskip
\medskip
\noindent\textbf{Step 3. The converse to gluing.}
\medskip

We follow closely the proof of Proposition~\ref{prop:converse}, with the additional choice of $\xi$. For a sequence $g_j$ as in the statement with $\beta_j\ge1$ and $A(g_j)\to0$, the same compactness and uniqueness arguments imply, after passing to a subsequence, convergence to an Einstein metric $g_0$ on $M$ and to $(Y,h)$ after area normalization, with $\beta_j\to\beta_*$. We set $\varepsilon_j^2=A(g_j)/A_*$ and work near $\xi_0=([g_0],0)$ in the chart satisfying $\delta_{g_\xi}g_0=0$.

The proof of Lemma~\ref{scd:gauge} does not require nondegeneracy of $g_0$ and applies uniformly near $\xi_0$. For each $\xi$ near $\xi_0$, we can choose a representative $g_j(\xi)$ of the fixed class $[g_j]$, by an equivariant pullback preserving each $\Sigma_\alpha$, such that
$$
  \delta_{g_j(\xi)}b_{\varepsilon_j,\xi,\beta_j}=0.
$$
We then define
$$
  k_j(\xi)=J_{\beta_j}^{-1}\bigl(g_j(\xi)-b_{\varepsilon_j,\xi,\beta_j}\bigr).
$$
On the supports of the test tensors defining $\ell_\xi$, $k_j(\xi)$ converges to $\bar g_\xi-g_\xi$, where $\bar g_\xi$ is the representative of $g_0$ satisfying $ \delta_{\bar g_\xi}g_\xi=0. $

In a local frame of $\mathcal E$, we define
$$
  \Phi_j(\xi)=\ell_\xi(k_j(\xi),-u_\xi).
$$
Applying the estimates in Lemma~\ref{scd:gauge} also to the $\xi$-derivatives of the gauge equation, and using removability and uniqueness on $M$, we obtain $C^1$ convergence of $\Phi_j$, on a fixed neighborhood of $\xi_0$, to
$$
  \Phi_0(\xi)=\ell_\xi(\bar g_\xi-g_\xi,-u_\xi).
$$
At $\xi_0$, we have $\bar g_\xi=g_\xi=g_0$, so $\Phi_0(\xi_0)=0$.

To solve $\Phi_j(\xi)=0$ nearby, we now check that $D\Phi_0(\xi_0)$ is invertible. Differentiating $\delta_{g_\xi}g_0=0$ and $\delta_{\bar g_\xi}g_\xi=0$ at $\xi_0$ shows that the derivatives of $\bar g_\xi$ with respect to $\xi$ are divergence-free. Since they are also Lie derivatives of $g_0$, they vanish by integration by parts. Thus
$$
  D\Phi_0(\xi_0)(\dot g,\dot u) =-\ell_{\xi_0}(\dot g,\dot u).
$$
This is an isomorphism by \eqref{eq:smooth-test-rows}, so the implicit-function theorem implies that there is a unique nearby zero $\xi_j$ of $\Phi_j$ for all sufficiently large $j$, with $\xi_j\to\xi_0$.

We now set
$$
  z_j=\bigl(k_j(\xi_j),\varepsilon_j(\beta_j-\beta_*),-u_{\xi_j}\bigr).
$$
By construction, $\mathcal F_{\varepsilon_j,\xi_j}(z_j)=0$. As in the proof of Proposition~\ref{prop:converse}, elliptic bootstrapping, \eqref{scd:tame}, and the uniform inverse estimate imply that
$$
  N_{\varepsilon_j}(z_j) \le C\varepsilon_j^{4/3}+o(1)N_{\varepsilon_j}(z_j),
$$
where the additional stabilizing terms are controlled by uniform continuity of their first derivatives. After absorbing the last term, we obtain $N_{\varepsilon_j}(z_j)\le C\varepsilon_j^{4/3}$. Uniqueness in {\bf Step 2} then identifies $g_j$ with $g_{\varepsilon_j}(\xi_j)$ modulo equivalence, and $v_{\varepsilon_j}(\xi_j)=-u_{\xi_j}$ implies that $p_{\varepsilon_j}(\xi_j)=0$.

We finally prove injectivity on $p_\varepsilon^{-1}(0)$ for sufficiently small $\varepsilon$. Otherwise, there are $\varepsilon_j\to0$ and distinct $\xi_j,\xi'_j\in p_{\varepsilon_j}^{-1}(0)$ such that
$$
  [g_{\varepsilon_j}(\xi_j)]=[g_{\varepsilon_j}(\xi'_j)].
$$
By compactness and \eqref{eq:family-correction}, both parameter sequences have subsequential limits in $p^{-1}(0)$. The corresponding isometries converge away from $p_0$, and their limits extend across $p_0$ and identify the limiting metrics. Thus both parameter sequences converge to the same $\xi_0\in p^{-1}(0)$. Applying the preceding construction to the common metric classes, both $\xi_j$ and $\xi'_j$ are zeros of $\Phi_j$ in a fixed neighborhood of $\xi_0$, contradicting the uniqueness established above.
\end{proof}

\begin{remark}\label{remark:C1 family}
Near the pairs $([g_\varepsilon(\xi)],0)$ with $p_\varepsilon(\xi)=0$, the parameter selection in the proof of {\bf Step 3} above extends to a $C^1$ map $(\varepsilon,[g],u)\mapsto\xi_\varepsilon([g],u)$. By {\bf Step 2} and the divergence identity, the equations
$$
  \mathscr E_3(g)=\sum_\nu u^\nu\widehat\tau_{\nu,\varepsilon,\xi_\varepsilon([g],u)}(g),\qquad A(g)=\varepsilon^2A_*
$$
have surjective differential on equivalence classes, and their solutions are locally parametrized by
$$
  \xi\longmapsto([g_\varepsilon(\xi)],p_\varepsilon(\xi)),
$$
with inverse $([g],u)\mapsto\xi_\varepsilon([g],u)$.
\end{remark}

\subsection{Proof of Theorem~\ref{thm:gluing formula}}

The idea is to compare the perturbed count at small fixed area with the perturbed count at cone angle $2\pi$. We will construct a compact one-dimensional manifold whose boundary points are precisely the solutions counted at these two ends. The gluing construction identifies the first count with $\mathfrak n(\mR_0)$, while the second count is $\mathfrak n(\mR)$.

We will use the equivalence relation of Section~\ref{sec:moduli} and we abbreviate the class of a pair $(\beta,g)$ by $[g]$, where $g$ is a toric conical metric on $X$ with common cone angle $2\pi\beta$ along $\mathcal S$ and is smooth along the other rod spheres. We write $\beta(g)=\beta$. The space $\mathcal M$ consists of the classes satisfying $\Ric(g)=3g$, and $\mathcal M_{\mR}$ is a subset of $\mathcal M$ with $\beta(g)=1$.

For $\eta>0$, we set
$$
  \mathcal K_\eta = \{[g]\in\mathcal M: \beta(g)\ge1,\ A(g)\ge\eta\}.
$$
Theorems~\ref{thm:compactness} and \ref{thm:upperangle} imply that $\mathcal K_\eta$ is compact: the cone angle is bounded, and the only possible degeneration must have  $A\to0$. Moreover, Lemma~\ref{lem:three-direction-smooth-compactness} gives compactness of $\mathcal M_{\mR}$. Thus, for sufficiently small $\eta$, no Einstein metric satisfies both $A(g)=\eta$ and $\beta(g)=1$.

We first assume that $\mathcal M_{\mR_0}$ is nonempty. We recall the $N$-dimensional manifold $\mathcal W$ used to define $\mathfrak n(\mR_0)$ and the map
$$
  p:\mathcal W\longrightarrow\mathbb R^N, \qquad ([g_{\xi}],u_{\xi})\longmapsto u_{\xi}.
$$
Its zero set is identified with $\mathcal M_{\mR_0}$. We choose $\mathcal V\Subset\mathcal W$ containing $p^{-1}(0)$ and $c>0$ such that $|p|\ge c$ on $\partial\mathcal V$. Then we have
$$
  \deg_2(p,\mathcal V,0)=\mathfrak n(\mR_0).
$$
By \eqref{eq:family-correction}, we have $|p_\varepsilon-p|<c/2$ on $\overline{\mathcal V}$ for sufficiently small $\varepsilon$. The straight-line homotopy between $p$ and $p_\varepsilon$ is therefore nonzero on $\partial\mathcal V$, and hence
\begin{equation}\label{eq:boundary-degree-identification}
  \deg_2(p_\varepsilon,\mathcal V,0) =\mathfrak n(\mR_0).
\end{equation}
We choose such an $\varepsilon$ and set $\eta=\varepsilon^2A_*$. We denote
$$
  \mathcal K_A=\{[g]\in\mathcal K_\eta:A(g)=\eta\}, \qquad \mathcal K_\beta=\{[g]\in\mathcal K_\eta:\beta=1\}.
$$
For sufficiently small $\varepsilon$, these compact sets are disjoint and $\mathcal K_\beta=\mathcal M_{\mR}$.

We now construct a perturbed Einstein equation connecting the two counts. We work in an open neighborhood $\mathcal U$ of $\mathcal K_\eta$ in the space of invariant conical metrics on $X$ modulo equivalence, allowing the common cone angle along $\mathcal S$ to vary. We choose disjoint open neighborhoods $\mathcal U_A,\mathcal U_\beta\subset\mathcal U$ of $\mathcal K_A,\mathcal K_\beta$, respectively.

For $[g]\in\mathcal U_A$, we denote
$$
  \varepsilon(g)=\sqrt{A(g)/A_*}.
$$
By the remark following Proposition~\ref{prop:stabilized-gluing}, after shrinking $\mathcal U_A$, the map $([g],u)\mapsto\xi_{\varepsilon(g)}([g],u)$ is well defined and $C^1$ for $[g]\in\mathcal U_A$ and sufficiently small $u\in\mathbb R^N$. We choose smaller open neighborhoods $\mathcal V_A,\mathcal V_\beta$ of $\mathcal K_A,\mathcal K_\beta$, respectively, with
$$
  \overline{\mathcal V_A}\subset\mathcal U_A, \qquad \overline{\mathcal V_\beta}\subset\mathcal U_\beta.
$$
We choose  cutoff functions $\chi_A,\chi_\beta$ equal to 1 on $\mathcal V_A,\mathcal V_\beta$ and supported in $\mathcal U_A,\mathcal U_\beta$, respectively. For $[g]\in\mathcal U_A$, we define
$$
  T_\nu(g,u)=\chi_A([g])\widehat\tau_{\nu,\varepsilon(g),\xi_{\varepsilon(g)}([g],u)}(g), \qquad 1\le\nu\le N,
$$
and extend $T_\nu$ by zero to $\mathcal U$.

On $\mathcal U_\beta$, we extend the tensors defining $\mathfrak n(\mR)$ to varying cone angle as in the proof of Proposition~\ref{prop:parity}. We denote by $S_1(g),\ldots,S_L(g)$ the resulting tensors multiplied by $\chi_\beta([g])$ and extended by zero to $\mathcal U$. These extensions are equivariant and divergence-free. Since the cutoff functions depend only on $[g]$ and are constant on $X$, multiplication by them preserves the divergence identity.

We also choose finitely many equivariant $C^1$ assignments $g\mapsto H_j(g)$ for $[g]\in\mathcal U$, where each $H_j(g)$ is an invariant symmetric two-tensor on $X$ satisfying
$$
  \delta_gH_j(g)=0.
$$
We require $H_j(g)=0$ whenever $[g]\in\mathcal V_A\cup\mathcal V_\beta$. For $u\in\mathbb R^N$, $v\in\mathbb R^L$, and $w\in\mathbb R^K$, we consider the system
\begin{equation}\label{eq:parity-comparison-equation}
  \mathscr E_3(g) = \sum_{\nu=1}^N u^\nu T_\nu(g,u) +\sum_{\mu=1}^L v^\mu S_\mu(g) +\sum_{j=1}^K w^j H_j(g).
\end{equation}
By the preceding remark and the construction of $\mathfrak n(\mR)$, the tensors $T_\nu(g,0)$ and $S_\mu(g)$ already ensure surjectivity at Einstein metrics in $\overline{\mathcal V_A}$ and $\overline{\mathcal V_\beta}$, respectively, after shrinking these neighborhoods. As in the proof of Proposition~\ref{prop:parity}, we choose $H_j$ so that the differential of \eqref{eq:parity-comparison-equation} in $([g],u,v,w)$ is onto $\ker\delta_g$ at every $([g],0,0,0)$ with $[g]\in\mathcal K_\eta$. Compactness of $\mathcal K_\eta$ ensures that finitely many $H_j$ suffice. At $u=0$, differentiating $T_\nu(g,u)$ contributes no additional term, since it is multiplied by $u^\nu$.

We write $\mathcal Z$ for the solutions of \eqref{eq:parity-comparison-equation} near $\mathcal K_\eta\times\{(0,0,0)\}$ satisfying
$$
  \beta(g)\ge1,\qquad A(g)\ge\eta.
$$
The implicit-function theorem gives a $C^1$ manifold of dimension $N+L+K+1$ before imposing these inequalities. Near $A=\eta$, its solutions are parametrized by
$$
  (r,\xi,v,w) \longmapsto ([g_r(\xi)],p_r(\xi),v,w),
$$
where $r$ varies near $\varepsilon$ and $A(g_r(\xi))=r^2A_*$. Thus $A$ is a submersion there. Near $\beta=1$, the surjectivity of the stabilized equation with cone angle fixed makes $\beta$ a submersion. After shrinking the neighborhood, $\mathcal Z$ is therefore a $C^1$ manifold with two disjoint boundary parts, given by $A=\eta$ and $\beta=1$.

We define the map
$$
  q:\mathcal Z\longrightarrow \mathbb R^N\times\mathbb R^L\times\mathbb R^K, \qquad ([g],u,v,w)\longmapsto(u,v,w).
$$
On the boundary $A=\eta$, the gluing parametrization gives the explicit formula
$$
  q\bigl([g_\varepsilon(\xi)],p_\varepsilon(\xi),v,w\bigr) = (p_\varepsilon(\xi),v,w).
$$
On the boundary $\beta=1$, the metric equation is
$$
  \mathscr E_3(g) = \sum_{\mu=1}^L v^\mu S_\mu(g),
$$
while $u$ and $w$ are free coordinates. Consequently, prescribing a sufficiently small regular value $(a_A,a_\beta,a_I)$ of $q$ on these two boundary parts gives respectively $\mathfrak n(\mR_0)$ and $\mathfrak n(\mR)$ solutions modulo two.

Now choose a relatively compact neighborhood $\mathcal U\subset\mathcal Z$ containing $\mathcal K_\eta\times\{(0,0,0)\}$. Since these are all the zeros of $q$, there is a $\delta>0$ such that
$$
  |(u,v,w)|\ge\delta \quad\text{on }\partial_{\mathcal Z}\mathcal U,
$$
where $\partial_{\mathcal Z}\mathcal U$ denotes the boundary of $\mathcal U$ relative to $\mathcal Z$. We choose $a=(a_A,a_\beta,a_I)$ with $|a|<\delta/2$, sufficiently small for the preceding count identifications and regular for the restrictions of $q$ to both $A=\eta$ and $\beta=1$.

By relative transversality, we perturb $q$ inside $\mathcal U$ to a $C^1$ map $\widetilde q$ for which $a$ is a regular value, keeping $\widetilde q=q$ near $A=\eta$, $\beta=1$, and $\partial_{\mathcal Z}\mathcal U$. Then
$$
  \mathcal C = \{z\in\overline{\mathcal U}: \widetilde q(z)=a\}
$$
is a compact one-dimensional manifold. It does not meet $\partial_{\mathcal Z}\mathcal U$. Its boundary points at $A=\eta$ correspond exactly to $p_\varepsilon^{-1}(a_A)\cap\mathcal V$; its boundary points at $\beta=1$ are the solutions counted by $\mathfrak n(\mR)$. Since a compact one-dimensional manifold has an even number of boundary points, \eqref{eq:parametrized-parity-law} follows.

Finally, suppose that $\mathcal M_{\mR_0}$ is empty. There is no sequence in $\mathcal M$ with $\beta\ge1$ and $A\to0$: smooth convergence would preserve positive area, while degeneration would produce a metric in $\mathcal M_{\mR_0}$. We may therefore choose $\eta>0$ such that $\mathcal M$ has no points with $\beta\ge1$ and $A\le\eta$. We repeat the construction above, omitting $u$ and the tensors $T_\nu$, and using
$$
  q([g],v,w)=(v,w).
$$
The resulting compact one-dimensional manifold has boundary only at $\beta=1$. Its boundary count is even, so $\mathfrak n(\mR)=0$. Since $\mathfrak n(\mR_0)=0$, this proves \eqref{eq:parametrized-parity-law} also in this case.

\section{Discussion}\label{sec:further-questions}
Our results suggest several future directions for exploration.

\begin{question}
For a fixed compact rod structure,
\begin{itemize}
\item(Uniqueness) Is a toric Einstein metric unique up to scaling and equivariant isometry?
\item(Non-degeneracy) Is every compact toric Einstein metric torically nondegenerate?
\item(Stability) Which compact toric Einstein metrics are stable? 
\end{itemize}
\end{question}
It would be more interesting if the answers to the first two questions are NO.  Notice that by Theorem~\ref{thm:fixed-rod-smooth-compactness},  the  toric Einstein moduli space $\mathcal M_{\mR}$ is compact for an arbitrary compact rod structure $\mR$.
Thus the construction in Section~\ref{sec:toric-existence-criterion} defines a virtual count $\mathfrak n(\mR)\in\mathbb Z/2\mathbb Z$ for every such $\mR$, even when $\mathcal M_{\mR}$ is not discrete. It is also intersting to try to define an oriented virtual count. We leave these directions for future investigation.

The toric Einstein theory developed in this paper are likely useful in studying toric shrinking Ricci solitons in  dimensions. In particular, a natural question is

\begin{question}
Does there exist a compact simply connected toric gradient shrinking Ricci soliton in dimension four that is neither Einstein nor K\"ahler?
\end{question}

There is also a natural Ricci-flow problem. Liu proved that integrability of the horizontal distribution is preserved under smooth Ricci flow on closed simply connected four-manifolds; see \cite[Theorem~8.1]{Liu}. Thus metrics with this property can be studied through their evolving quotient metrics and orbit matrices.

\begin{question}
Can we have a theory of Ricci flow with surgery for $\mathbb T^2$-invariant Ricci flow with integrable horizontal distribution?
\end{question}

This could serve as an experiment towards understanding of general singularity formation of 4d Ricci flow.

\appendix

\section{Asymptotics of toric gravitational instantons}
\label{app:asymptotics of toric gravitational instantons}

In this appendix we establish the classification of asymptotic geometries for toric conical gravitational instantons and discuss some applications.

\subsection{Definitions and statements}
As in \cite{LiSun}, a gravitational instanton is a complete connected oriented Ricci-flat four-manifold with quadratic curvature decay
\begin{equation}\label{eq:quadratic decay}
    |\Rm_g|_g(p)\leq C(1+d_g(o,p))^{-2},
\end{equation}
where $o$ is a fixed base point.

We first recall relevant notations in \cite{LiSun}. We say a Ricci-flat four-manifold $(M,g)$ is \emph{toric}, if every point on $M$ has a connected stabilizer group. Denote the quotient by $(N,\bg)$, which is a smooth manifold with boundary, and the interior of $N$ by $N^\circ$. On $M$ there is the horizontal distribution that is orthogonal to the $\T^2$-orbits. The toric Ricci-flat four-manifold is called \emph{integrable} if the horizontal distribution is integrable, which automatically holds when the $\T^2$-action has a fixed point. After fixing a set of generators $\{\partial_{\phi_1},\partial_{\phi_2}\}$, there is the \emph{Gram matrix}
$$
    G_{ij}:=g(\partial_{\phi_i},\partial_{\phi_j}).
$$
Write $\rho=\sqrt{\det G}$ which is a function on the quotient $N$ and define
$$
    \Phi:=\rho^{-1}G,
$$
which is a map from $N^\circ$ into the space of symmetric positive $(2\times 2)$ matrices with determinant one. The function $\rho$ is smoothly defined on $N^\circ$ and continuously extends to $N$ and vanishes on $\partial N$. We say $(M,g)$ is \emph{regular} if $\dd\rho\neq0$ on $N^\circ$.
\begin{definition}\label{def:app gravitational instanton}
A toric gravitational instanton is a complete Ricci-flat oriented four-manifold $(M,g)$ that is integrable and regular, with quadratic curvature decay \eqref{eq:quadratic decay}.
\end{definition}

For a toric gravitational instanton, the function $\rho$ on $(N,\bg)$ is harmonic on $N^\circ$, and locally there is the harmonic conjugate $z$ that satisfies $\dd z=\underline{\ast} \dd\rho$.
\begin{definition}\label{def: app gravitational instanton proper}
A toric gravitational instanton is proper if there is a globally defined harmonic conjugate $z$ of $\rho$ on $N$ such that the map $(\rho,z)$ induces a homeomorphism from $N$ to $\mathbb H$, the two-dimensional right half-plane
\end{definition}

\begin{remark}
The properness assumption is a genuine restriction. The periodic quotients of the four-dimensional examples in \cite[Theorems~2 and~4]{KRWY24} are complete toric gravitational instantons with Kasner asymptotics; see also \cite[Section~2.3]{OzuchStability}. Their orbit spaces are half-cylinders $[0,\infty)\times(\mathbb R/L\mathbb Z)$, where $L>0$, and the local harmonic conjugate $z$ is periodic. Consequently, $\underline{\ast}\dd\rho$ has a nonzero period around the cylinder and admits no globally defined real-valued primitive. These gravitational instantons are therefore not proper.
The corresponding infinite periodic covers have half-plane orbit spaces, but do not satisfy the quadratic curvature decay condition \eqref{eq:quadratic decay}.
\end{remark}
For a proper toric gravitational instanton, we have the global coordinates $(\rho,z)$ for $N$,  and the Ricci-flat metric can be written as
\begin{equation}\label{eq:integrable toric metric}
    g=e^{2\nu}(\dd\rho^2+\dd z^2)+\rho\Phi_{ij}\dd\phi_i \dd\phi_j.
\end{equation}
It has been shown in \cite{LiSun} that the Ricci-flat equation for a metric in the form of \eqref{eq:integrable toric metric} is equivalent to the following system:
\begin{align}
    &(\rho \Phi^{-1}\Phi_\rho)_\rho+(\rho\Phi^{-1}\Phi_z)_z=0,\label{eq:axisymmetric harmonic map equation}\\
    &\partial_z\nu=\frac{1}{4\rho}\tr (\mathcal{U}\mathcal{V}),\\
    &\partial_\rho\nu=-\frac{1}{2\rho}+\frac{1}{8\rho}\tr (\mathcal{U}^2-\mathcal{V}^2),
\end{align}
where
\begin{equation}
    \mathcal{U}=\mathrm{Id}+\rho\Phi^{-1}\Phi_\rho, \qquad \mathcal{V}=\rho\Phi^{-1}\Phi_z.
\end{equation}
The target of $\Phi$, which is the space of symmetric positive $(2\times 2)$ matrices with determinant one, can be identified with the two-dimensional hyperbolic upper-half-plane with curvature $-2$
$$
    \mathcal{H}=\{(X,Y)\mid X>0\}, \quad g_{\mathcal{H}}=\frac{\dd X^2+\dd Y^2}{2X^2}
$$
via
\begin{equation}\label{eq:identification of hyperbolic space and matrix}
    (X,Y)\mapsto \left(\begin{matrix}
        X+X^{-1}Y^2 & X^{-1}Y \\
        X^{-1}Y & X^{-1}
    \end{matrix}\right).
\end{equation}

In the Euclidean $\mathbb R^3$, denote the $z$-axis as $\Gamma$. Let $\rho$ be the distance function to $\Gamma$ in $\mathbb R^3$. Then \eqref{eq:axisymmetric harmonic map equation} is equivalent to that
$$
    \Phi:\mathbb{R}^3\setminus \Gamma \to \mathcal{H}
$$
is an axisymmetric harmonic map. For each toric gravitational instanton $(M,g)$, there is an associated \emph{rod structure}
$$
    \mR =(\{z_j\}_{j=1}^n,\{\mv_j\}_{j=0}^n),
$$
which consists of
\begin{itemize}
\item A collection of \emph{turning points} $0=z_1<\ldots <z_n$ on $\Gamma=\R$. We also write $z_0=-\infty$ and $z_{n+1}=\infty$. For $0\leq j\leq n$, the open interval $\mI_j:=(z_{j},z_{j+1})\subset\Gamma$ is called a rod.
\item A collection of \emph{rod vectors} $\mv_0,\ldots, \mv_n\in\mathbb{RP}^1$ with $\mv_j\neq \mv_{j+1}$.
\end{itemize}
Note that the definition of rod structures we use in this appendix is actually the \emph{untyped rod structure} under the terminology of \cite{LiSun}. Under the coordinate $(\rho,z)$ on $N$, the boundary $\partial N$ is identified with $\Gamma$. Each fixed point of $\T^2$ in $M$ descends to a point in $\partial N$, and corresponds to a turning point $z_j$ on $\Gamma$. Each rod $\mI_j$ corresponds to a totally geodesic $\mathbb S^2$ inside $M$, whose isotropy group is $S^1$. Under our choice of generator $\{\partial_{\phi_1},\partial_{\phi_2}\}$ where we identify them as $(1,0)^\top$ and $(0,1)^\top$, the isotropy of the rod $\mI_j$ is described by the rod vector $\mv_j$. Under our terminology of \cite{LiSun}, the axisymmetric harmonic map $\Phi$ is locally tamed by the rod structure $\mR$. Along each rod, there is the notion of \emph{cone angles}. See \cite[Definition 4.19]{LiSun}. In the smooth case the cone angles are all $2\pi$. For toric \emph{conical} gravitational instantons, cone angles are allowed to be different from $2\pi$. Throughout this appendix, we take the lattice to be $\Lambda=(2\pi\bZ)^2$, which means the periods of $\phi_1,\phi_2$ are $2\pi$.

In this appendix, we work with proper toric conical gravitational instantons, which also corresponds to axisymmetric harmonic maps as above. The only difference is the cone angles along each rod is allowed to be different from $2\pi$. For a precise formulation, see \cite{LiSun}. We will prove
\begin{theorem}\label{thm:appendix classification of asymptotics}
A proper toric conical gravitational instanton $(M,g)$ with finite rods is of Asymptotic Type $\ALE$, $\ALF^\pm$, or $\AF_\beta$. In the $\ALF^\pm$ and $\AF_\beta$ cases, there are a corresponding model end $(E_0,g_0)$ of Asymptotic Type $\ALE,\ALF^\pm$, or $\AF_\beta$, and a diffeomorphism $\varpi:E_0\to M\setminus K$ such that
\begin{equation}
    |\nabla_{g_0}^k(\varpi^* g-g_0)|_{g_0}=O(\mathfrak r^{-1-k}), \quad k\geq0,
\end{equation}
where $\mathfrak r$ is the distance to a fixed base point.
\end{theorem}

\subsection{Maximal volume growth}

In the maximal volume growth case, the gravitational instanton is $\ALE$.

\begin{proposition}\label{app:maximal-volume-growth}
Let $(M,g)$ be a proper toric conical gravitational instanton with maximal volume growth. Then the two infinite rods have distinct directions. After changing the angular lattice on the complement of a compact set, the resulting smooth Ricci-flat end is asymptotically Euclidean. With the original angular lattice, the end is asymptotic to the corresponding flat toric cone, allowing conical singularities along its two infinite rods and a finite angular quotient.
\end{proposition}
\begin{proof}
We first enlarge the lattice so that all cone angles are strictly less than $2\pi$, and continue to denote the resulting space by $(M,g)$.
As in \cite{LiSun}, Bishop--Gromov comparison applies.
Together with maximal volume growth, this implies that
$$
    vt^4\leq\operatorname{Vol}_g B_g(p,t)\leq\omega_4t^4
$$
for some $v>0$, every $p\in M$, and every $t>0$, where $\omega_4$ is the volume of the Euclidean unit ball.

We choose a sequence $R_j\to\infty$ and set $g_{R_j}=R_j^{-2}g$.
After passing to a subsequence, we obtain a pointed limit $(Z,g_\infty,o_\infty)$.
The quadratic curvature decay and volume lower bound imply smooth convergence on annuli in local developments.
Thus $g_\infty$ is Ricci-flat on its regular set.
Writing $v_\infty$ for the asymptotic volume ratio, volume convergence implies that
$$
    \operatorname{Vol}_{g_\infty}B_{g_\infty}(o_\infty,t)
    =v_\infty\omega_4t^4
$$
for every $t>0$.

We denote the distance from $o_\infty$ by $r$.
Weak Laplacian comparison passes to the limit, so
$$
    \Delta_{g_\infty}\frac{r^2}{2}\leq4.
$$
Testing this inequality against nonnegative radial cutoff functions and using the volume identity, we see that equality holds on every annulus, across the conical strata. In particular,  $r^2/2$ is smooth on the regular set.
Since $|\nabla(r^2/2)|^2=r^2$, the Bochner formula implies that
$$
    \left|\operatorname{Hess}_{g_\infty}\frac{r^2}{2}\right|^2
    =\frac12\Delta_{g_\infty}r^2
    =4.
$$
This implies that
$$
    \operatorname{Hess}_{g_\infty}\frac{r^2}{2}=g_\infty.
$$
So the limit is  a metric cone
$$
    (C(Y),g_C),\qquad g_C=\dd r^2+r^2h.
$$
It follows that $g_C$ is flat on its regular set.
The limiting effective torus action fixes the cone vertex and, in flat developments, acts by rotations in two orthogonal two-planes.
Thus, up to a finite angular quotient, the cone is
$$
    \mathbb R^2_{\alpha_-}\times\mathbb R^2_{\alpha_+},
$$
where $\mathbb R^2_\alpha$ denotes the flat cone of angle $\alpha$.
Its two rods have distinct stabilizer circles.

If the two original infinite rods had the same direction, their common primitive generator $X$ would vanish on both.
Its normal rotational derivative along either rod is fixed and nonzero.
Smooth convergence in local developments would therefore produce a nonzero limiting Killing field vanishing on both limiting rods, contradicting the two-plane description above.
Hence the original infinite rod directions are distinct.

We now return to the original angular lattice.
We choose primitive period vectors $\ell_-,\ell_+\in\Lambda$ representing the two infinite rod directions and denote the corresponding cone angles by $2\pi\beta_-$ and $2\pi\beta_+$.
The associated $2\pi$-periodic infinitesimal generators are $\ell_\pm/(2\pi)$.
On an end containing no finite rods, we replace $\Lambda$ by
$$
    \Lambda_{\mathrm{sm}}
    =\mathbb Z\frac{\ell_-}{\beta_-}
    +\mathbb Z\frac{\ell_+}{\beta_+}.
$$
The two generators are independent and form a basis of this lattice.
Both infinite rods now have cone angle $2\pi$, so we obtain a smooth Ricci-flat end $(E_{\mathrm{sm}},g_{\mathrm{sm}})$ with cross-section $S^3$.
This change of lattice is made only on the end.

The modified end still has quadratic curvature decay and maximal volume growth.
After this change of lattice, the annular limits above are Euclidean.
The Chern--Gauss--Bonnet formula on truncations of the end therefore implies finite curvature energy, since the outer boundary terms remain bounded.
We can now apply the  argument of Bando--Kasue--Nakajima \cite{BKN}; its application to a Ricci-flat end is also explained in \cite[proof of Proposition~4.38]{LiSun}.
Thus the modified end is ALE.
Since its cross-section is $S^3$, the asymptotic group is trivial, so the end is AE.
As in that argument, we can choose the asymptotic identification equivariantly.
Restoring the original angular periods then yields the required flat conical asymptotic model, with the differentiated decay supplied by the exterior estimates.
\end{proof}

\subsection{Estimates on the quotient}

In the following we consider the case where the gravitational instanton has non-maximal volume growth.  For a proper toric gravitational instanton $(M,g)$, the quotient $(N,\bg)$ is a smooth manifold with boundary and corners. The rods are the boundaries and turning points are the corners. The rods $\mI_j$ are totally geodesic in $N$.
\begin{lemma}
The quotient $(N,\bg)$ has non-positive sectional curvature $K_{\bg}\leq 0$ and quadratic curvature decay. Moreover, $(N,\bg)$ is a $\mathrm{CAT}(0)$ space.
\end{lemma}
\begin{proof}
Because $(M,g)$ is integrable, the sectional curvature of $(N,\bg)$ is equal to the sectional curvature of $(M,g)$ along $\partial_\rho\wedge\partial_z$, which therefore decays quadratically.  In particular $K_{\bg}\to 0$ at infinity, so any positive global maximum is attained. The proof of Lemma \ref{geom:scalar} and Proposition \ref{geom:Kupper} go through as well for the Ricci-flat case, hence we know there is no interior maximum. It follows $K_{\bg}\leq0$. Since it has geodesic boundaries and right angle corners, it is $\mathrm{CAT}(0)$ \cite[Theorem II.4.1]{BH}.
\end{proof}

In $(N,\bg)$, we set $o\in \partial N$ to be the turning point $z_1$. If there is no turning point in the rod structure, we just set $o$ to be an arbitrary base point in $\partial N$ and normalize $z$ to be zero at $o$. Write $r=\sqrt{\rho^2+z^2}$ and $d=d_{\bg}(o,\cdot)$. Denote
$$
    A_d(a,b)=\{p\in N\mid a<d(p)<b\},\qquad A_r(a,b)=\{p\in N\mid a<r(p)<b\}.
$$
By the next several lemmas, we prove as $R\to\infty$, the scaled down sequence  $(A_d(aR,bR),R^{-2}\bg)$ is noncollapsed. We denote the metric annulus in $(M,g)$ centered at $o$ as $A_g(a,b)$, which covers $A_d(a,b)$.

Consider any point $q$ that is on the two infinite rods that are sufficiently far away from $o$. Write $O=\T^2\cdot q$ as the the orbits at $q$, which is an $S^1$, and $F$ as the fixed surface corresponding to the infinite rod that contains $q$. Set $R=d(q)$ and the length of orbit $O$ at $q$ as $\sigma(q)=|O|$. Let $X$ be a primitive Killing vector field with period $2\pi$ that vanishes along $F$.
\begin{lemma}\label{lem:injectivity radius lower bound for quotient}
There are constants $\epsilon,\delta>0$, such that the normal exponential map of the orbit $O$ is injective on its radius $\epsilon R$ normal disk bundle. Moreover,
\begin{equation}\label{eq:volume cubic lower bound}
    \Vol \left(B_g(o,(1+\epsilon)R)\right)\geq \delta \sigma(q)R^3.
\end{equation}
\end{lemma}
\begin{proof}
The fixed surface $F$ is totally geodesic and
$$
    g_F=\dd l^2+\sigma(l)^2\dd\psi^2,\qquad \sigma''+K_F\sigma=0,
$$
where $\partial_\psi$ is a primitive integral Killing field in $\T^2$ with period $2\pi$ and $l$ is the arclength from the torus fixed point contained in $F$. The Ricci-flat version of \eqref{eqn:relation ambient and sigma curvature} and $K_{\bg}\leq 0$ give $\sigma''=2K_{\bg}\sigma\leq0$. When there is no turning point for $\mR$, the function $\sigma$ is nowhere vanishing and concave, hence constant. When there are turning points, since $\sigma(0)=0$ and $\sigma$ stays positive on the whole ray, concavity gives $0\leq \sigma'\leq \sigma/l$. Since the ray is a geodesic in $(N,\bg)$, so $l(q)=R+O(1)$. On the ray when $|l-l(q)|<R/4$ we consequently have
\begin{equation}\label{eq:sigma'/sigma}
    |\sigma'|/\sigma\leq C/R.
\end{equation}

The orbit circle has second fundamental form of the norm $|\sigma'|/\sigma$. Normal Jacobi comparison gives a focal radius at least $c_0R$. If normal injectivity radius were smaller than a sufficiently small $\epsilon R< c_0R$, then there would be a geodesic $\gamma$ from $O$ to $O$, normal to both ends, of length $L<2\epsilon R$. The Jacobi field $J_{\beta_\varepsilon}=X|_\gamma$ vanishes at both ends. From the Jacobi field equation we obtain
$$
    \int_0^L|D_t J_{\beta_\varepsilon}|^2\dd t\leq CR^{-2}\int_0^L|J_{\beta_\varepsilon}|^2\dd t.
$$
If $\epsilon$ is small enough, the Poincar\'e inequality forces $J_{\beta_\varepsilon}=0$. Thus $\gamma$ lies in $F$. The velocity $\dot\gamma$ is orthogonal to $\partial_\psi$, since $\frac{d}{dt}g(\dot\gamma, \partial_\psi)=0$ and $\gamma$ starts orthogonally to $O$. So it is a radial $l$-geodesic, which cannot go back to $O$, a contradiction.

To see the lower bound on the volume, the $\epsilon R$ lower bound on injectivity radius gives
$$
    \Vol \left(B_g(o,(1+\epsilon)R)\right)\geq \sigma(q)\delta R^3
$$
for a small constant $\delta$. Since $\sigma$ is concave, this in particular implies $(M,g)$ has at least cubic volume growth.
\end{proof}

\begin{lemma}\label{lem:quotient annuli converge}
For every $0<a<b<\infty$, as $R\to\infty$, the sequence of scaled down annuli $(A_d(aR,bR),R^{-2}\bg)$ has uniformly bounded curvature and is noncollapsed. Each annulus can be covered by a uniform number of balls with radius $\epsilon$.
\end{lemma}
\begin{proof}
This is clear from Lemma \ref{lem:injectivity radius lower bound for quotient}, where we have a uniform injectivity radius lower bound $\epsilon R$ at points on the infinite rods. The scaled down annulus has uniformly bounded curvature since $(N,\bg)$ has quadratic curvature decay. The number of balls to cover $A_d(aR,bR)$ can be controlled by the number of balls to cover $A_g(aR,bR)$ in the Ricci-flat $(M,g)$, which is controlled by volume comparison.
\end{proof}

Next we control $\nu$ and $\rho$. Recall $\bg=e^{2\nu}(\dd\rho^2+\dd z^2)$. Recall that $q$ is any point on the two infinite rods that are far away from $o$.
\begin{lemma}\label{lem:control on weyl annulus}
There is a constant $C>0$ such that, outside a compact set,
\begin{equation}\label{eq:bound nu}
    re^{\nu}\leq Cd.
\end{equation}
Moreover, on $A_d(R/2,2R)$,
\begin{equation}
    C^{-1}\sigma(q)\leq|\dd\rho|_{\bg}\leq C\sigma(q),\quad C^{-1}R\sigma(q)\leq|z(q)|\leq CR\sigma(q).
\end{equation}
\end{lemma}
\begin{proof}
Consider annuli $A_g(aR,bR)$ of the gravitational instanton $(M,g)$ centered at $o$. Since $(M,g)$ has quadratic curvature decay, in $g_R=R^{-2}g$ the curvature of $(A_g(aR,bR),g_R)$ is uniformly bounded. By Cheeger-Fukaya-Gromov for any $p\in(A_g(aR,bR),g_R)$ there is a small ball $(B_{\delta}(p),g_R)$ such that its universal cover $(\tilde B_{\delta}(p),\tilde g_R)$ is noncollapsed. We pick the point $p\in A_g(aR,bR)$ to be above $q$. Then the universal cover $(\tilde B_{\delta}(p),\tilde g_R)$ is noncollapsed, and by the regularity of Einstein equation, it converges smoothly to a Ricci-flat limit after passing to subsequences. By Lemma \ref{lem:quotient annuli converge}, we know $(A_d(aR,bR),\bg_R)$ converges smoothly to a limit $(A_{d,\infty},\bg_\infty)$ with geodesic boundaries.

As before, we let $X$ be a primitive Killing field with period $2\pi$ that vanishes along $F$. Let $Y$ be another Killing field in $\T^2$, normalized such that $|Y(p)|_{\bg}=\sigma(p)/2\pi$. Such a $Y$ restricts to $\partial_\psi$ on $F$. We then normalize $X,Y$ in $A_g(aR,bR)$ to
$$
    X,\quad V_R=\frac{R}{\sigma(q)}(Y-b_RX),
$$
where we take the constant $b_R$ to cancel the normal rotational part of $Y$ at $p$. In the noncollapsed limit of the universal cover $(\tilde{B}_{g_\infty},\tilde{g}_\infty)$, the normalized Killing fields $X,V_R$ converge to nontrivial limit Killing fields $X_\infty,V_\infty$. From the Killing field equation, this implies the normalized Killing fields $X,V_R$ converge on the local noncollapsed cover for any point in $A_d(aR,bR)$. This implies
$$
    \rho_R:=\frac{\rho}{R\sigma(q)},
$$
which is the square root of the determinant of the Gram matrix formed by $X,V_R$ in $g_R$, converges smoothly to a limit harmonic function $\rho_\infty$ on $A_{d,\infty}$. The conjugate of $\rho_R$, which is related to $z$ via
$$
    z_R=\frac{z-z(q)}{R\sigma(q)},
$$
converges smoothly to a limit harmonic function $z_\infty$ as well on $A_{d,\infty}$. The holomorphic functions $w_R=z_R+\sqrt{-1}\rho_R$ have nowhere vanishing derivative on $A_{d}(aR,bR)$, so Hurwitz's theorem implies that $w_\infty=z_\infty+\sqrt{-1}\rho_\infty$ has no critical points on $A_{d,\infty}$. Therefore, we have uniform upper and lower bound on $|\dd \rho_R|$. This proves $C^{-1}\sigma(q)\leq|\dd\rho|_{\bg}\leq C\sigma(q)$.

From the convergence of $z_R$ we know for any two points $q_1,q_2\in A_{d}(aR,bR)$, we have
$$
    |z(q_1)-z(q_2)|\leq CR\sigma(q).
$$
Since $A_d(aR,bR)$ contains part of the two infinite rods, on which the values of $z$ are positive and negative respectively, there exists points in $A_d(aR,bR)$ with $z=0$. Hence on $A_d(aR,bR)$ we obtain the bound
$$
    \rho+|z|\leq CR\sigma(q).
$$
On a rod segment of intrinsic length $cR$, \eqref{eq:sigma'/sigma} gives that $\sigma$ is comparable to $\sigma(q)$. Harmonic conjugacy gives the identity $|\dd z/\dd l|=a$ on the rod. Since the base point $o$ is certainly not contained in $A_R(aR,bR)$, integration gives $|z(q)|\geq cR\sigma(q)$.

Finally, since $|\dd \rho|_{\bg}=e^{-\nu}$ and $d(q)=R$, we have obtained $re^{\nu}\leq Cd$.
\end{proof}

\subsection{Limits of the rescaled harmonic maps}
Now we justify the existence of scale-downs of $\Phi$. Recall $\Phi$ is an axisymmetric harmonic map from $\mathbb{H}^\circ=\{(X,Y)\mid X>0\}$ into $\mathcal{H}$. For each $R\gg1$, there is the scaled down harmonic map $\Phi(R\rho,Rz)$. Set $p_*=(1,0)$ in $\mathbb H$. For $R\to\infty$, we choose an isometry $T_R$ of the hyperbolic space with $T_R\Phi(Rp_*)=(1,0)$. As $R\to\infty$, we define
$$
    \Phi_R(\rho,z)=T_R\Phi(R\rho,Rz).
$$
\begin{lemma}
As $R\to\infty$, after passing to subsequences, $\Phi_R$ converges smoothly on any compact subset of $\mathbb{H}^\circ$ to a limit harmonic map $\Phi_\infty$.
\end{lemma}
\begin{proof}
The distance function $d$ is $1$-Lipschitz for $\bg$, and \eqref{eq:bound nu} therefore gives that in the Euclidean metric $\dd \rho^2+\dd z^2$ for $\mathbb{H}$,
$$
    |\nabla_{\mathrm{Eucl}}\log d|\leq C/r
$$
almost everywhere. Integrating along paths of bounded length in an annulus $A_r(aR,bR)$ shows that in such an annulus the values of the distance function $d$ are uniformly comparable. Hence, if we set $d_R=d((R,0))$ and enlarge the annulus so that $a<1<b$, $A_r(aR,bR)$ is contained in $A_d(a'd_R,b'd_R)$ for some $0<a'<b'<\infty$, independent of $R$.

Here $R$ is the Weyl coordinate scale, whereas the geometric rescaling is $g_{d_R}=d_R^{-2}g$. Choose $q=(0,R)$ on the positive infinite rod and use the normalization
$$
    V_R=\frac{d_R}{\sigma(q)}(Y-b_RX),
$$
where $X,Y$ and $b_R$ are chosen as in the proof of Lemma \ref{lem:control on weyl annulus}. Since $d(q)$ is uniformly comparable to $d_R$, that proof gives smooth convergence of the normalized Killing fields $X,V_R$ on the local noncollapsed limits around any point in $(A_g(a'd_R,b'd_R),g_{d_R})$. Hence the associated Gram matrix converges on each local noncollapsed cover. After division by the square root of its determinant, this Gram matrix differs from $\Phi$ by an isometry of the hyperbolic space via \eqref{eq:identification of hyperbolic space and matrix}; the scalar metric rescaling cancels in this determinant normalization. Hence after a correct isometry, the rescaled harmonic map converges on any compact subset of $\mathbb H^\circ$. Since the isometry $T_R$ normalizes $\Phi(Rp_*)$ to $(1,0)$ and the stabilizer of $(1,0)$ is compact, we know $\Phi_R$ converges to a limit harmonic map $\Phi_\infty$ after passing to a subsequence.
\end{proof}

\begin{lemma}\label{lem:limit harmonic map and limit rod structure with same rod vectors}
Any limit harmonic map $\Phi_\infty$ has a rod structure $\mR_\infty$ that only consists of two infinite rods and one turning point. The two infinite rods have the same rod vector. The harmonic map $\Phi_\infty$ is locally tamed by $\mR_\infty$ in the interior of the infinite rods.
\end{lemma}
Here, after applying an isometry of the hyperbolic space, we may assume the rod vectors is $[1,0]^\top$, then $\Phi_\infty$ being locally tamed in the interior of infinite rods means that for any point in the interior of infinite rods, there is a neighborhood such that
$$
    d_{\mathcal{H}}(\Phi_\infty,\Phi_{[1,0]^\top})\text{ is bounded},\qquad \Phi_{[1,0]^\top}=\left(\begin{matrix}
        \rho & 0 \\
        0 & \rho^{-1}
    \end{matrix}\right).
$$
\begin{proof}
We keep $R$ as the Weyl coordinate scale and set $d_R=d((R,0))$. We take $q=(0,R)$ on the positive infinite rod, and let $F$ be the fixed surface in $M$ corresponding to this infinite rod as before. Let $X$ be a primitive Killing field with period $2\pi$ that vanishes along $F$, and let $Y$ be another Killing field normalized such that $|Y(q)|_g=\sigma(q)/2\pi$, which restricts to $\partial_\psi$ on $F$. From \eqref{eq:volume cubic lower bound}, applied at the geometric scale $d(q)$, since $(M,g)$ has non-maximal volume growth, we know $\sigma(q)/d(q)\to0$ along the infinite rod. This implies that $|Y(q)|_g/d(q)\to0$.

We also pick a point $q'$ on the negative infinite rods, such that $q,q'\in A_r(aR,bR)$ for fixed $0<a<1<b$. Let $F'$ be the fixed surface, and let $X'$ be a primitive Killing field with period $2\pi$ that vanishes along the rod and $Y'$ be another Killing field. Then $X=mX'+nY'$ for some constants $m,n$. Similarly, $|Y'(q')|_g/d(q')\to0$. By the preceding lemma, both $d(q)/d_R$ and $d(q')/d_R$ are bounded above and below independently of $R$.

Now we use $g_{d_R}=d_R^{-2}g$ and $V_R=\frac{d_R}{\sigma(q)}(Y-b_RX)$, with $b_R$ chosen as above, and consider the Gram matrix for $g_{d_R}$ under the basis $X,V_R$. We express this matrix in the rescaled Weyl coordinates and denote it by
$$
    G_R=\left(\begin{matrix}
        G_{R,11} & G_{R,12}\\
        G_{R,21} & G_{R,22}
    \end{matrix}\right).
$$
Since we have the smooth convergence of $X,V_R$ in any local noncollapsed cover and $A_r(aR,bR)\subset A_d(a'd_R,b'd_R)$, $G_R$ converges smoothly to a limit $G_\infty$ on $\mathbb H^\circ$. The component $G_{R,11}$ vanishes along the positive infinite rod. On the negative infinite rod, since $X'=0$ there, we have
$$
    \begin{aligned}
        G_{R,11}(0,z(q')/R)&=\frac{|nY'(q')|_g^2}{d_R^2}\\
        &=\frac{|nY'(q')|_g^2}{d(q')^2}\frac{d(q')^2}{d_R^2}\to0.
    \end{aligned}
$$
This holds for any $q'\in A_r(aR,bR)$ on the other infinite rod,  hence $G_{\infty,11}=0$ along the other infinite rod as well. This indicates the limit Gram matrix $G_\infty$ has the same one-dimensional kernel along the two infinite rods, hence the two infinite rods share the same rod vector.

Since $G_R$ is defined using the basis $X,V_R$ and $X,V_R$ converge smoothly, it is clear that $\Phi_\infty$ is locally tamed in the interior of the positive infinite rod. As for the negative infinite rod, a similar normalization of Killing fields $X',V_R'$ converge smoothly in any local noncollapsed cover so we obtain a limit harmonic map $\Phi_\infty'$ that is locally tamed in the interior of the negative infinite rod. Comparing the two normalization shows $\Phi_\infty$ is also locally tamed in the interior of the negative infinite rod.
\end{proof}

For a limit harmonic map $\Phi_\infty$, we always apply an isometry of the hyperbolic space so that the rod vector is $[1,0]^\top$.
\begin{proposition}
The convergence $\Phi_R\to\Phi_\infty$ is $C^{1,\alpha}$ on any compact subset of $\mathbb H\setminus\{o\}$. There exists a constant matrix $P\in \SL(2,\mathbb R)$ such that
$$
    \Phi_\infty=(P^{-1})^\top \Phi_{[1,0]^\top}P^{-1}.
$$
\end{proposition}
For the precise definition of $C^{1,\alpha}$ convergence of harmonic maps, see \cite[Definition 4.20]{LiSun}.
\begin{proof}
We write $\rho=r\sin\theta$ and $z=r\cos\theta$. Let $X=\rho e^u$ and $Y=v$ where we parametrize $\mathcal{H}$ via $(X,Y)$. From the smooth convergence to $\Phi_\infty$ and the local tameness of $\Phi_\infty$ by $\mR_\infty$, from the regularity of \cite[Proposition 4.12]{LiSun}, we have the $C^{1,\alpha}$ convergence of $\Phi_R$ to $\Phi_\infty$.

Now we consider a limit harmonic map $\Phi_\infty=(\rho e^{u_\infty},v_\infty)$, and choose a sequence $R_j\to\infty$ such that $\Phi_{R_j}\to\Phi_\infty$. The annular estimates and rod regularity used above also give uniform estimates for the normalized secondary rescalings at every scale $R>0$. Indeed, fix $R>0$, let $\xi_j$ be the ideal endpoint associated with the positive rod of $\Phi_{R_j}$, and choose the orientation-preserving hyperbolic isometry $S_{j,R}$ that sends $\xi_j$ to $0$ and $\Phi_{R_j}(R,0)$ to $(1,0)$. Since $\xi_j\to0$, these isometries converge to the isometry $S_R$ satisfying
$$
    S_R(0)=0,\qquad S_R\Phi_\infty(R,0)=(1,0).
$$
The map $S_{j,R}\Phi_{R_j}(R\rho,Rz)$ is a normalized rescaling of the original map at scale $R_jR\to\infty$. Consequently, its annular interior estimates and rod-adapted estimates pass to $S_R\Phi_\infty(R\rho,Rz)$. Their constants are independent of $R$, since they are independent of the original scale $R_jR$. The limiting rod values on both sides are $0$ by Lemma \ref{lem:limit harmonic map and limit rod structure with same rod vectors}. In particular, for every $R>0$,
$$
    d_{\mathcal{H}}\Big(S_R\Phi_\infty(R\sin\theta,R\cos\theta),(\sin\theta,0)\Big)\leq C,\qquad 0<\theta<\pi,
$$
and the family $S_R\Phi_\infty(R\rho,Rz)$ has uniform interior derivative bounds on compact subsets of $\mathbb H^\circ$. Writing
$$
    S_R\Phi_\infty(R\rho,Rz)=(\rho e^{u_{\infty,R}},v_{\infty,R}),
$$
the same rod regularity gives uniform $C^{1,\alpha}$ bounds for $u_{\infty,R}$ and the corresponding weighted derivative bounds for $v_{\infty,R}$ on compact subsets of $\mathbb H\setminus\{o\}$; see \cite[Proposition 4.12]{LiSun}. Thus these estimates hold at both arbitrarily large and arbitrarily small secondary scales.

The curves $S_R^{-1}(\sin\theta,0)$ and $(R\sin\theta,0)$, parametrized by $-\log\sin\theta$, are both geodesics in $\mathcal{H}$ with the same constant speed $1/\sqrt{2}$, sharing the same ideal endpoint. Their distance is non-increasing in the parameter $-\log\sin\theta$ due to hyperbolic geometry, where $-\log\sin\theta\geq0$ and $-\log\sin\theta=0$ at $\theta=\frac{\pi}{2}$. So
\begin{align*}
    d_{\mathcal{H}}\Big(S_R^{-1}(\sin\theta,0),(R\sin\theta,0)\Big)&\leq d_{\mathcal{H}}\Big(S_R^{-1}(1,0),(R,0)\Big)\\
    &= d_{\mathcal{H}}\Big(\Phi_\infty(R,0),(R,0)\Big).
\end{align*}
From this we obtain that at the point $(R\sin\theta,R\cos\theta)$,
\begin{align*}
    d_{\mathcal{H}}(\Phi_\infty,\Phi_{[1,0]^\top})&= d_{\mathcal{H}}\Big(\Phi_\infty,(R\sin\theta,0)\Big)\\
    &\leq d_{\mathcal{H}}\Big(\Phi_\infty,S_R^{-1}(\sin\theta,0)\Big)+d_{\mathcal{H}}\Big(S_R^{-1}(\sin\theta,0),(R\sin\theta,0)\Big)\\
    &\leq C+d_{\mathcal{H}}\Big(\Phi_\infty(R,0),(R,0)\Big).
\end{align*}
The uniform interior derivative bound for $S_R\Phi_\infty(R\rho,Rz)$ at $(\rho,z)=(1,0)$, valid for every $R>0$, gives
$$
    |\partial_\rho\Phi_\infty(\rho,0)|_{\mathcal{H}}\leq \frac{C}{\rho},\qquad \rho>0.
$$
Integrating from $1$ to $\rho$, in either direction, gives
$$
    d_{\mathcal{H}}\Big(\Phi_\infty(\rho,0),(1,0)\Big)\leq C(|\log \rho|+1).
$$
Since the hyperbolic distance from $(\rho,0)$ to $(1,0)$ is $|\log\rho|/\sqrt{2}$, we conclude
$$
    d_{\mathcal{H}}\Big(\Phi_\infty(R,0),(R,0)\Big)\leq C(|\log R|+1),\qquad R>0,
$$
and therefore
\begin{equation}\label{eq:limit harmonic map to model log bound}
    d_{\mathcal{H}}(\Phi_\infty,\Phi_{[1,0]^\top})\leq C(1+|\log r|),\qquad r>0.
\end{equation}

Finally we compare $\Phi_\infty$ to itself after a $z$-translation by $t$
$$
    D_t(\rho,z):=d_{\mathcal{H}}\Big(\Phi_\infty(\rho,z+t),\Phi_\infty(\rho,z)\Big).
$$
Both are harmonic maps from $\mathbb{R}^3\setminus \Gamma$ to $\mathcal{H}$, so $D_t$ is subharmonic. The uniform interior and rod estimates for the secondary rescalings at every scale $R>0$ imply the following uniform bound on $\mathbb{H}^\circ$
\begin{equation}\label{eq:limit harmonic map z derivative bound}
    |\partial_z\Phi_\infty|_{\mathcal{H}}\leq \frac{C}{r}.
\end{equation}
Indeed, writing $S_R\Phi_\infty(R\rho,Rz)=(\rho e^{u_{\infty,R}},v_{\infty,R})$, the rod estimates give uniform bounds for $u_{\infty,R}$, $\partial_z u_{\infty,R}$ and $\rho^{-1}\partial_zv_{\infty,R}$ on the unit semicircle. Thus the hyperbolic norm of the $z$-derivative of the normalized map is uniformly bounded. Rescaling and using that $S_R$ is an isometry give \eqref{eq:limit harmonic map z derivative bound} for every $r>0$. Integration along the segment from $(\rho,z)$ to $(\rho,z+t)$ gives
$$
    D_t(\rho,z)\leq C\frac{|t|}{r},\qquad \text{if }r\gg1+|t|.
$$
The function $D_t$ is locally bounded near the rods except possibly at $(0,0)$ and $(0,-t)$, since $\Phi_\infty$ is locally tamed in the interior of the two rods. Hence it extends subharmonically across the rods except the two points. As $\Phi_{[1,0]^\top}$ is independent of $z$, the bound \eqref{eq:limit harmonic map to model log bound} gives at most logarithmic growth of $D_t$ at these two exceptional points. For $\epsilon>0$, apply maximum principle to
$$
    D_t(\rho,z)-\epsilon\left( r^{-1}+\frac{1}{\sqrt{\rho^2+(z+t)^2}}\right)
$$
on large balls with small balls about the two exceptional points deleted. The logarithmic growth makes the inner boundary values negative; the outer supremum tends to zero. For fixed $t$ and $\epsilon>0$, let the radii of the deleted balls tend to zero and the outer radius tend to infinity. The maximum principle gives
$$
    0\leq D_t(\rho,z)\leq\epsilon\left(r^{-1}+\frac{1}{\sqrt{\rho^2+(z+t)^2}}\right).
$$
Letting $\epsilon\to0$ at each fixed point with $\rho>0$ gives $D_t=0$, which holds for any $t$.

Thus $\Phi_\infty$ depends only on $\rho$. The axisymmetric harmonic map equation then says $t\mapsto \Phi_\infty(e^t,z)$ is a geodesic in $\mathcal{H}$. After an isometry $P\in \SL(2,\mathbb R)$ it is identified to the geodesic $t\mapsto (e^t,0)$ that corresponds to the harmonic map $\Phi_{[1,0]^\top}$. This gives $\Phi_\infty=(P^{-1})^\top\Phi_{[1,0]^\top} P^{-1}$, as desired.
\end{proof}
\subsection{Convergence at infinity}
So far we have shown that for the normalized harmonic map $\Phi_R(\rho,z)=T_R\Phi(R\rho,Rz)$, any subsequential limit of it is equivalent to $\Phi_{[1,0]^\top}$. In the following, let $s=\log r$. We first study the behavior of $\Phi$ along the ray $(R,0)$ in $\mathbb H$ as $R\to\infty$.

\begin{lemma}
There is a constant matrix $Q_1^{-1}\in\SL(2,\mathbb R)$ and a constant $b$, such that for $\Phi_1=(Q_1^{-1})^\top\Phi Q_1^{-1}$, as $R\to\infty$,
$$
    Y(\Phi_1(R,0))\to b,\quad X(\Phi_1(R,0))\to0,\quad \frac{|Y(\Phi_1(R,0))-b|}{X(\Phi_1(R,0))}\to0.
$$
\end{lemma}
\begin{proof}
We first perform a M\"obius transformation
$$
    Y+\sqrt{-1}X\mapsto-\frac{1}{Y+\sqrt{-1}X-c}.
$$
This is a fixed hyperbolic isometry, which corresponds to conjugation by a matrix $Q_1^{-1}\in\SL(2,\mathbb R)$. The value of $Y$ is a constant along the positive infinite rod, and we take $c$ to be that value.  Then the boundary value of $\Phi_1:=(Q_1^{-1})^\top\Phi Q_1^{-1}$ on the positive rod, which is a point in the ideal boundary, is infinity.

We follow the ray $(R,0)$ in $\mathbb H$ to infinity and consider the values $\Phi_1(R,0)$. For $R>0$, define the real functions
$$
    f(R)=Y(\Phi_1(R,0)),\quad h(R)=X(\Phi_1(R,0)).
$$
Let $A_R$ be the hyperbolic isometry that is defined by
$$
    A_R:(X,Y)\mapsto\left(\frac{X}{h(R)},\frac{Y-f(R)}{h(R)}\right).
$$
Thus the normalized harmonic map
$$
    \Phi_{1,R}(\rho,z)=A_R\Big(\Phi_1(R\rho,Rz)\Big)
$$
has value $(X,Y)=(1,0)$ at $(\rho,z)=(1,0)$, while its positive rod value remains infinity.

Write
$$
    F_R=(Y+\sqrt{-1}X)\circ\Phi_{1,R},
$$
so $F_R$ is the complex coordinate expression of $\Phi_{1,R}$. The preceding blow-down result shows $F_R$ has a limit, which corresponds to a geodesic in $\mathcal{H}$. The boundary value of $Y(\Phi_{1,R})$ is infinity on the positive rod, and $F_R(1,0)=(1,0)$. Since the limit geodesic has to pass through $(1,0)$ and diverge to infinity, it is the vertical line $(X,Y)=(\frac{1}{\rho},0)$, and
$$
    F_R(\rho,z)\to\frac{\sqrt{-1}}{\rho}.
$$
The convergence is smooth near $(1,0)$. Differentiating in the domain variable $\rho$, with $R$ fixed, therefore gives
$$
    \partial_\rho F_R(1,0)=\frac{R f'(R)+\sqrt{-1}R h'(R)}{h(R)}\to-\sqrt{-1},\qquad \text{as }R\to\infty.
$$
Consequently,
$$
    \frac{R f'(R)}{h(R)}=o(1),\quad\frac{R h'(R)}{h(R)}=-1+o(1),\qquad \text{as }R\to\infty.
$$
Integration of the second relation yields
$$
    h(R)=R^{-1+o(1)}.
$$
Moreover, for all sufficiently large $R$, we have $h'(R)/h(R)<-\frac{1}{2R}$. So integration gives
$$
    \frac{h(t)}{h(R)}\leq \left(\frac{t}{R}\right)^{-1/2},\quad \text{for }t\geq R.
$$
It follows that $|f'|\leq \frac{h(t)}{t}$, so $f'$ is absolutely integrable at infinity. Hence $f(R)\to b$ for some $b\in\mathbb R$, and
$$
    \frac{|f(R)-b|}{h(R)}\leq\left(\sup_{t\geq R}\frac{t|f'(t)|}{h(t)}\right)\int_R^\infty\left(\frac{t}{R}\right)^{-1/2}\frac{dt}{t}=2\sup_{t\geq R}\frac{t|f'(t)|}{h(t)}\to0.
$$
\end{proof}

\begin{lemma}\label{lem:rough estimate on u and v}
There is a constant matrix $Q\in\SL(2,\mathbb R)$ and constant $a$, such that for
$$
    (Q^{-1})^\top\Phi Q^{-1}=(\rho e^u,v),
$$
as $s\to\infty$,
$$
    u=o(s),\quad u_s=o(1),\quad u_\theta=o(1),
$$
uniformly in $\theta$. Moreover, for every $\epsilon>0$,
$$
    \frac{|v-a\cos\theta|}{\sin^2\theta} \leq C_\epsilon e^{(1+\varepsilon)s}
$$
for all sufficiently large $s$.
\end{lemma}
\begin{proof}
Introduce the M\"obius transformation
$$
    Y+\sqrt{-1}X\mapsto-\frac{1}{Y+\sqrt{-1}X-b},
$$
which corresponds to conjugation by another matrix $Q_2^{-1}\in\SL(2,\mathbb R)$. Then we set
$$
    \Phi_2:=(Q_2^{-1})^\top\Phi_1 Q_2^{-1}.
$$
The M\"obius transformation sends $(b,0)$ to infinity and the positive rod value of $\Phi_1$ to $(0,0)$. We claim boundary values of $\Phi_2$ on both infinite rods are finite. For the positive rod, the value of $\Phi_1$ is infinite, so the value of $\Phi_2$ is $(0,0)$, which is finite. For the negative rod, if the value of $\Phi_2$ was not finite, then $\Phi_1$ had negative rod value $b$. Then for the normalized maps
$$
    \Phi_{1,R}(\rho,z)=A_R(\Phi_1(R\rho,Rz)),
$$
the positive rod value is infinity, whereas the negative rod value is
$$
    \frac{b-f(R)}{h(R)}\to0.
$$
This is a contradiction since $F_R\to\frac{\sqrt{-1}}{\rho}$, where the negative rod value should diverge to infinity.

Let $c_+,c_-$ be the positive and negative finite rod values of $Y(\Phi_2)$, so $c_+=0$. Define a third hyperbolic isometry via
$$
    Y+\sqrt{-1}X\mapsto Y+\sqrt{-1}X-\frac{c_++c_-}{2},\quad a=\frac{c_+-c_-}{2}.
$$
The hyperbolic isometry corresponds to a conjugation by $Q_3^{-1}\in \SL(2,\mathbb R)$. Set
$$
    \Phi_3=(Q_3^{-1})^\top\Phi_2 Q_3^{-1}.
$$
Then the boundary values of $\Phi_3$ are $a$ and $-a$. Set
$$
    H_R=X(\Phi_3(R,0)),\quad K_R=Y(\Phi_3(R,0)),\quad\alpha_R=\frac{f(R)-b}{h(R)}.
$$
Direct calculation gives
$$
    H_R=\frac{1}{h(R)(1+\alpha_R^2)}=R^{1+o(1)}.
$$
Define the hyperbolic isometry
$$
    B_R:(X,Y)\mapsto\left(\frac{X}{H_R},\frac{Y-K_R}{H_R}\right),
$$
where $B_R(\Phi_3(R\rho,Rz))$ has value $(1,0)$ at $(\rho,z)=(1,0)$. This indicates that $B_R(\Phi_3(R\rho,Rz))$ converges to a limit. The complex coordinate expression
$$
    \widehat F_R(\rho,z)=\frac{(Y+\sqrt{-1}X)(\Phi_3(R\rho,Rz))-K_R}{H_R}
$$
satisfies
$$
    \widehat F_R=\frac{\alpha_R F_R-1}{F_R+\alpha_R}.
$$
Since $\alpha_R\to0$ and $F_R\to \sqrt{-1}/\rho$, we obtain
$$
    \widehat F_R\longrightarrow \sqrt{-1}\rho.
$$

Define $Q=Q_3Q_2Q_1$ and $(Q^{-1})^\top\Phi Q^{-1}=(\rho e^u,v)$, so $\Phi_3=(Q^{-1})^\top\Phi Q^{-1}$. Let
$$
    U_R(\rho,z)=\log\frac{\operatorname{Im}\widehat F_R}{\rho}=u(R\rho,Rz)-\log\frac{H_R}{R},\quad V_R(\rho,z)=\operatorname{Re}\widehat F_R=\frac{1}{H_R}(v(R\rho,Rz)-K_R).
$$
The established rod compactness gives $C^{1,\alpha}$ convergence of these functions to zero on closed annuli away from $o$.  For each fixed $R$, the definition of $U_R$ gives,
$$
    u(s,\theta)=U_R\left(\frac{e^s}{R}\sin\theta,\frac{e^s}{R}\cos\theta\right)+\log\frac{H_R}{R}
$$
for $s$ near $\log R$. In particular,
$$
    u(\log R,\theta)=U_R(\sin\theta,\cos\theta)+\log\frac{H_R}{R}.
$$
Since $U_R\to0$ uniformly on the unit semicircle and $H_R=R^{1+o(1)}$, we obtain
$$
    u=o(s)
$$
uniformly in $\theta$. To estimate the derivatives, keep $R$ fixed in the preceding identity, differentiate with respect to $s$ and $\theta$, and then evaluate at $s=\log R$. This gives
$$
    \begin{aligned}
        u_s(\log R,\theta)&=\sin\theta \partial_\rho U_R(\sin\theta,\cos\theta)+\cos\theta\partial_z U_R(\sin\theta,\cos\theta),\\
        u_\theta(\log R,\theta)&=\cos\theta\partial_\rho U_R(\sin\theta,\cos\theta)-\sin\theta\partial_z U_R(\sin\theta,\cos\theta).
    \end{aligned}
$$
The $C^1$ convergence of $U_R$ therefore implies
$$
    u_s=o(1),\qquad u_\theta=o(1),
$$
uniformly in $\theta$.

It remains to estimate $v$. On the portions of the outer rods contained in the fixed annulus, the boundary values of $V_R$ are
$$
    V_R(0,z)= \begin{cases}
        (a-K_R)/H_R,& z>0,\\
        (-a-K_R)/H_R,& z<0.
    \end{cases}
$$
The uniform $C^{1,\alpha}$ convergence and the regularity of \cite[Proposition~4.12]{LiSun} consequently gives
$$
    \left|V_R(\rho,z)-\frac{a-K_R}{H_R}\right|\leq C\rho^2\qquad\text{near }(0,1),
$$
and
$$
    \left|V_R(\rho,z)-\frac{-a-K_R}{H_R}\right|\leq C\rho^2\qquad\text{near }(0,-1),
$$
with neighborhoods and constants independent of $R$. For $\theta$ sufficiently close to $0$, we therefore have
$$
    \begin{aligned}
        \left|V_R(\sin\theta,\cos\theta)-\frac{a\cos\theta-K_R}{H_R}\right|&\leq\left|V_R(\sin\theta,\cos\theta)-\frac{a-K_R}{H_R}\right|+\frac{|a|}{H_R}(1-\cos\theta)\\
        &\leq C\sin^2\theta.
    \end{aligned}
$$
Here we used the boundedness of $|a|/H_R$. The analogous estimate near $\theta=\pi$ follows from the negative rod value.

On the remaining compact angular interval, $\sin\theta$ is bounded below. Moreover, $V_R$ and its two rod values are uniformly bounded. Their interpolation
$$
    \frac{a\cos\theta-K_R}{H_R}=\frac{1+\cos\theta}{2}\frac{a-K_R}{H_R}+\frac{1-\cos\theta}{2}\frac{-a-K_R}{H_R}
$$
is therefore uniformly bounded as well. Increasing $C$ if necessary, we obtain on the entire unit semicircle
$$
    \left|V_R(\sin\theta,\cos\theta)-\frac{a\cos\theta-K_R}{H_R}\right|\leq C\sin^2\theta.
$$
By the definition of $V_R$, multiplication by $H_R$ gives
$$
    |v(R\sin\theta,R\cos\theta)-a\cos\theta|\leq C H_R\sin^2\theta.
$$
Finally, $H_R=R^{1+o(1)}$ implies that, for every $\epsilon>0$,
$$
    \frac{|v(s,\theta)-a\cos\theta|}{\sin^2\theta}\leq C_\epsilon e^{(1+\epsilon)s}
$$
for all sufficiently large $s$ and $0<\theta<\pi$. This proves the stated estimates.
\end{proof}

Next we improve the estimate that we obtain in Lemma \ref{lem:rough estimate on u and v}. In the follwing, scalar functions are regarded as axisymmetric functions on Euclidean $\mathbb R^3$, and $\Delta$ and $\operatorname{div}$ denote its Euclidean operators. For $f=f(\rho,z)$, the gradient norm is $|\nabla f|^2=f_\rho^2+f_z^2$, which also equals its flat half-plane gradient norm.
\begin{lemma}
We have the following bounds on $v$ when $r$ is sufficently large
$$
    |v-a\cos\theta|\leq C\sin^{3/2}\theta,\quad |\nabla v|\leq \frac{C}{r}\sin^{1/2}\theta.
$$
\end{lemma}

\begin{proof}
Recall that our notation is
$$
    X\circ\Phi_3=\rho e^u,\qquad Y\circ\Phi_3=v,
$$
and that the two outer rod values of $v$ are $a$ and $-a$. The harmonic map equations are
\begin{equation}\label{eq:strong-tameness-hm}
    \Delta_{\mathbb R^3}u=-\rho^{-2}e^{-2u}|\nabla v|^2, \qquad \operatorname{div}_{\mathbb R^3}\bigl(\rho^{-2}e^{-2u}\nabla v\bigr)=0.
\end{equation}
By Lemma \ref{lem:rough estimate on u and v},
\begin{equation}\label{eq:strong-tameness-preliminary}
    u=o(s),\qquad u_s=o(1),\qquad u_\theta=o(1),\qquad s=\log r,
\end{equation}
as $s\to\infty$, uniformly in $\theta$. Moreover, for every $\epsilon>0$,
\begin{equation}\label{eq:strong-tameness-growth}
    |v-a\cos\theta|\leq C_\varepsilon e^{(1+\epsilon)s}\sin^2\theta
\end{equation}
for all sufficiently large $s$.

We first improve \eqref{eq:strong-tameness-growth} to a bound independent of $s$. In logarithmic spherical coordinates, the second equation in \eqref{eq:strong-tameness-hm} is $Lv=0$, where
$$
    L=\partial_s^2-(1+2u_s)\partial_s+\partial_\theta^2-\cot\theta\partial_\theta-2u_\theta\partial_\theta.
$$
Consequently,
$$
    L(v-a\cos\theta)=-2a u_\theta\sin\theta.
$$
For $\lambda=0$ and $\lambda=\frac32$, direct computation gives
$$
    \frac{L(e^{\lambda s}\sin^{3/2}\theta)}{e^{\lambda s}\sin^{3/2}\theta}=\lambda^2-\lambda-\frac32-\frac34\cot^2\theta-2\lambda u_s-3u_\theta\cot\theta.
$$
Choose $s_0$ sufficiently large that
$$
    |u_s|,|u_\theta|\leq\frac1{12},\qquad\text{for }s\ge s_0.
$$
Using $|\cot\theta|\leq(1+\cot^2\theta)/2$, we obtain
\begin{equation}\label{eq:strong-tameness-barrier}
    L(e^{\lambda s}\sin^{3/2}\theta)\leq-\frac38 e^{\lambda s}\sin^{-1/2}\theta,\qquad\lambda\in\left\{0,\frac32\right\}.
\end{equation}

Choose $C$ sufficiently large to control $2a u_\theta\sin\theta$ and the boundary values at $s=s_0$. For $\delta>0$, consider
$$
    \pm(v-a\cos\theta)-C\sin^{3/2}\theta-\delta e^{3s/2}\sin^{3/2}\theta
$$
on $[s_0,T]\times[0,\pi]$. By \eqref{eq:strong-tameness-barrier}, these functions have nonnegative $L$-image and are non-positive at $s=s_0$. Fixing $0<\epsilon<\frac12$ in \eqref{eq:strong-tameness-growth}, we see that they are also non-positive at $s=T$ for all sufficiently large $T$. They vanish continuously at $\theta=0,\pi$. Any positive maximum would therefore occur strictly inside the rectangle, where $L$ is regular, contradicting the maximum principle. Letting $T\to\infty$ and then $\delta\to0$, we obtain
\begin{equation}\label{eq:strong-tameness-v-bound}
    |v-a\cos\theta|\leq C\sin^{3/2}\theta.
\end{equation}

We next estimate the source in the equation for $u$. Away from the axis, the equation for $v$ can be written as
$$
    v_{\rho\rho}+v_{zz}-\left(\rho^{-1}+2u_\rho\right)v_\rho-2u_zv_z=0.
$$
By \eqref{eq:strong-tameness-preliminary},
$$
    |\nabla u|=\frac{\sqrt{u_s^2+u_\theta^2}}{r}\leq\frac Cr.
$$
Fix $p=(\rho_0,z_0)$, set $r_0=(\rho_0^2+z_0^2)^{1/2}$, and consider the Euclidean disk
$$
    B=B_{\rho_0/8}(p)\subset\{\rho>0\}.
$$
On $B$, we have $\rho$ is comparable to $\rho_0$ and $r$ is comparable to $r_0$. Hence
$$
    \rho_0\sup_B\bigl(\rho^{-1}+2|\nabla u|\bigr)\leq C.
$$
Thus, after rescaling $B$ to a unit disk, the equation for $v$ has uniformly bounded first-order coefficients. The standard interior gradient estimate gives
$$
    |\nabla v(p)|\leq \frac{C}{\rho_0}\operatorname{osc}_B v,
$$
where $C$ is independent of $p$. This follows from the interior $W^{2,q}$ estimate and Sobolev embedding with $q>2$; only a uniform bound on the first-order coefficients is required.

To estimate the oscillation, \eqref{eq:strong-tameness-v-bound} gives
$$
    \sup_B|v-a\cos\theta|\leq C\left(\frac{\rho_0}{r_0}\right)^{3/2}.
$$
Moreover, since
$$
    |\nabla\cos\theta|=\frac{\rho}{r^2},
$$
we have
$$
    \operatorname{osc}_B(a\cos\theta)\leq C\frac{\rho_0^2}{r_0^2}.
$$
Consequently,
$$
    \begin{aligned}
        \operatorname{osc}_B v&\leq2\sup_B|v-a\cos\theta|+\operatorname{osc}_B(a\cos\theta)\\
        &\leq C\left(\frac{\rho_0}{r_0}\right)^{3/2},
    \end{aligned}
$$
where we used $0<\rho_0/r_0\le1$. Therefore
$$
    |\nabla v(p)|\leq \frac{C}{\rho_0}\left(\frac{\rho_0}{r_0}\right)^{3/2}.
$$
Since $p$ was arbitrary, we obtain
\begin{equation}\label{eq:strong-tameness-v-gradient}
    |\nabla v|\leq \frac C\rho\left(\frac{\rho}{r}\right)^{3/2} =\frac Cr\sin^{1/2}\theta
\end{equation}
for all sufficiently large $r$.
\end{proof}

Finally we have the following strong tameness.

\begin{proposition}
There exists $c>0$ such that the harmonic model $\Phi_0$, defined by
$$
    \Phi_0(\rho,z)=\left(\rho\Big(c+\frac{|a|}{r}\Big),a\cos\theta\right)
$$
satisfies
$$
    d_{\mathcal H}(\Phi_3,\Phi_0)=O(r^{-1}).
$$
\end{proposition}

\begin{proof}
Since $u=o(\log r)$, equations \eqref{eq:strong-tameness-hm} and \eqref{eq:strong-tameness-v-gradient} imply, for every $\epsilon>0$,
$$
    |\Delta_{\mathbb R^3}u| \leq C_\epsilon r^{-4+\epsilon}\sin^{-1}\theta.
$$
The angular factor is integrable on $\mathbb S^2$. Moreover, $|\nabla u|\le C/r$ implies that the flux through a shrinking cylindrical tube around the axis tends to zero on every compact exterior annulus. Thus the first equation in \eqref{eq:strong-tameness-hm} holds weakly across the exterior axis. Consider the spherical mean
$$
    m(r)=\frac1{4\pi}\int_{\mathbb S^2}u(r,\omega)\dd\omega.
$$
Averaging the equation for $u$ gives
$$
    (r^2m'(r))'=\frac{r^2}{4\pi}\int_{\mathbb S^2}\Delta_{\mathbb R^3}u(r,\omega)\dd\omega=O(r^{-2+\epsilon}).
$$
Taking $0<\epsilon<1$, the right-hand side is integrable at infinity. Hence $r^2m'(r)$ has a finite limit, $m'(r)=O(r^{-2})$, and $m(r)$ converges to a finite constant $u_\infty$. Furthermore,
$$
    \sup_\theta|u(r,\theta)-m(r)|\leq\pi\sup_\theta|u_\theta(r,\theta)|=o(1).
$$
Therefore
\begin{equation}\label{eq:strong-tameness-u-limit}
    u\to u_\infty,\quad\text{uniformly in }\theta.
\end{equation}

Set
$$
    c=e^{u_\infty}>0,\qquad H=c+\frac{|a|}{r},
$$
and define $\Phi_0$ by
$$
    X\circ\Phi_0=\rho H,\qquad Y\circ\Phi_0=a\cos\theta.
$$
This map is harmonic. Write
$$
    D=d_{\mathcal H}(\Phi_3,\Phi_0).
$$
The hyperbolic distance formula gives
$$
    \cosh(\sqrt2D)-1=\frac{(e^u-H)^2}{2e^uH}+\frac{(v-a\cos\theta)^2}{2\rho^2e^uH}.
$$
The first term tends uniformly to zero by \eqref{eq:strong-tameness-u-limit}. Since $e^u$ and $H$ are bounded above and below for large $r$, estimate \eqref{eq:strong-tameness-v-bound} bounds the second term by
$$
    Cr^{-2}\sin\theta.
$$
Thus $D\to0$ uniformly for $0<\theta<\pi$.

Finally, the distance between harmonic maps into $\mathcal H$ is subharmonic away from the axis. Fix a sufficiently large $R_0$ and set
$$
    C=R_0\sup_{r=R_0}D.
$$
Given $\eta>0$, take $T>R_0$ sufficiently large that $D\le\eta$ on $r=T$. For $\delta>0$, the function
$$
    D-\frac Cr-\eta+\delta\log\frac{\rho}{2T}
$$
is subharmonic on
$$
    \{R_0<r<T,\ \rho>0\},
$$
because both $r^{-1}$ and $\log\rho$ are harmonic there. It is non-positive on the spherical boundaries and tends to $-\infty$ at the axis, since $D$ is bounded there. Applying the maximum principle, first on the domain with a sufficiently thin cylindrical tube removed, gives
$$
    D-\frac Cr-\eta+\delta\log\frac{\rho}{2T}\le0.
$$
Letting $\delta\to0$, then $T\to\infty$ and $\eta\to0$, we conclude that
$$
    D\leq\frac Cr.
$$
This proves strong tameness.
\end{proof}

Starting from the strong tameness of the harmonic map, it is a direct application of \cite[Proposition A.1]{LiSun} to conclude Theorem \ref{thm:appendix classification of asymptotics}. Although \cite[Proposition A.1]{LiSun} is stated for smooth ends, the same argument applies to conical ends by working changing the lattice for the end. Accordingly, the asymptotic model in Theorem \ref{thm:appendix classification of asymptotics} is allowed to have cone angles different from $2\pi$ along its infinite rods.

\subsection{Properness and ALE rigidity}

Finally, we record a criterion for properness and a rigidity statement for toric ALE orbifolds. We continue to write $(N,\bg)=M/\T^2$, with $G$, $\rho=\sqrt{\det G}$, and $\Phi=\rho^{-1}G$ as above.

\begin{lemma}\label{geom:proper-ricciflat}
Let $(M,g)$ be a complete connected oriented toric conical Ricci-flat four-manifold that is integrable. Assume that $N^\circ$ is simply connected, $\partial N\neq\varnothing$, and $K_{\bg}\leq0$. Then $\dd\rho\neq0$ on $N^\circ$, and there is a globally defined harmonic conjugate $z$, satisfying $\dd z=\underline{\ast}\dd\rho$, such that
$$
    F=(\rho,z):N\longrightarrow\mathbb H=\{(\rho,z)\in\mathbb R^2:\rho\geq0\}
$$
is a homeomorphism. Its restriction to the interiors is a conformal diffeomorphism. If $g$ also has quadratic curvature decay \eqref{eq:quadratic decay}, then $(M,g)$ is a proper toric conical gravitational instanton in the sense of Definitions \ref{def:app gravitational instanton} and \ref{def: app gravitational instanton proper}.
\end{lemma}

\begin{proof}
Since $(M,g)$ is complete and the torus is compact, the quotient $(N,\bg)$ is complete as a metric space, including its boundary. The reduced identities \eqref{eq:rho} and \eqref{geom:hessu}, applied with $\lambda=0$, give $\Delta_{\bg}\rho=0$ and
\begin{equation}\label{eq:prelim-rho-concavity}
    \Hess_{\bg}\sqrt{\rho}=\frac{\sqrt{\rho}}{2}(K_{\bg}\bg-T)\leq0,\qquad T=\frac14\tr(\Phi^{-1}\dd\Phi\otimes\Phi^{-1}\dd\Phi).
\end{equation}
Here $T$ is the nonnegative tensor introduced before \eqref{geom:hessu}. These are local identities for an integrable Einstein metric, so they apply in the present noncompact Ricci-flat setting. The function $\rho$ is positive on $N^\circ$ and vanishes on $\partial N$, and is therefore nonconstant. If $\dd\rho(p)=0$ at some $p\in N^\circ$, then along every sufficiently short geodesic $\gamma$ issuing from $p$,
$$
    \left.\frac{\dd}{\dd t}\sqrt{\rho(\gamma(t))}\right|_{t=0}=0,\qquad \frac{\dd^2}{\dd t^2}\sqrt{\rho(\gamma(t))}\leq0.
$$
Thus $p$ would be a local maximum of $\rho$, contradicting the strong maximum principle for the nonconstant harmonic function $\rho$. Hence $\dd\rho\neq0$ on $N^\circ$.

Since $N^\circ$ is simply connected and $\dd(\underline{\ast}\dd\rho)=0$, there is a global harmonic conjugate $z$ on $N^\circ$. The map $F=(\rho,z)$ is a local conformal diffeomorphism. Set $\nu=-\log|\dd\rho|_{\bg}$, so that in each such coordinate neighborhood
$$
    \bg=e^{2\nu}(\dd\rho^2+\dd z^2).
$$
Equation \eqref{eq:prelim-rho-concavity} implies
$$
    \begin{pmatrix}
        \nu_\rho+\dfrac{1}{2\rho}&\nu_z\\[3pt]
        \nu_z&-\nu_\rho
    \end{pmatrix}\geq0.
$$
The diagonal entries and determinant therefore give
\begin{equation}\label{eq:prelim-weyl-gradient}
    -\frac{1}{2\rho}\leq\nu_\rho\leq0,\qquad \nu_\rho^2+\nu_z^2\leq-\frac{\nu_\rho}{2\rho}\leq\frac{1}{4\rho^2}.
\end{equation}

We claim that $F:N^\circ\to\mathbb H^\circ$ is a diffeomorphism. Let $\gamma:[0,1]\to\mathbb H^\circ$ be a piecewise smooth path starting at $F(p)$, and lift it locally starting at $p\in N^\circ$. Since $\gamma([0,1])$ is compact, its $\rho$ coordinate is bounded below by some $\varepsilon>0$. Along the lift, \eqref{eq:prelim-weyl-gradient} gives
$$
    \left|\frac{\dd}{\dd t}\nu(\widetilde\gamma(t))\right|\leq\frac{|\dot\gamma(t)|}{2\varepsilon},\qquad \nu(\widetilde\gamma(t))\leq\nu(p)+\frac{L_{\mathrm{Euc}}(\gamma)}{2\varepsilon}.
$$
Consequently the $\bg$-length of every lifted subarc is bounded by a fixed constant times the Euclidean length of the corresponding subarc of $\gamma$. If the lift stopped before reaching the endpoint, its tail would be Cauchy. Completeness would give a limit in $N$, and the lower bound $\rho\geq\varepsilon$ would place this limit in $N^\circ$, where the local inverse extends the lift. Thus every such path lifts throughout, and in particular $F$ is onto. It follows that $F$ is a covering map. Since $N^\circ$ is connected and $\mathbb H^\circ$ is simply connected, $F$ is a diffeomorphism.

We next extend $F$ across $\partial N$, using the local expansions in Lemma \ref{lem:rodasymptotics}. Near an interior point of a rod, let $t\geq0$ be the inward distance to that rod in the quotient metric $\bg$, and let $s$ be arclength along the rod, extended constantly along the normal geodesics. Thus $(t,s)$ are local Fermi coordinates. By \eqref{eq:rhoends},
$$
    \begin{aligned}
        \bg&=\dd t^2+f(t^2,s)^2\dd s^2,\qquad f(0,s)=1,\\
        \rho&=t\kappa(s)+O(t^3),\qquad \kappa(s):=\partial_t\rho(0,s)>0.
    \end{aligned}
$$
In particular, $\kappa(s)$ is the inward normal derivative of $\rho$ along the rod. The one-form $\underline{\ast}\dd\rho$ extends smoothly to $t=0$, so the same is true of $z$. Choose the direction of $s$ so that $\underline{\ast}\dd t=\dd s$ on the rod. Then
$$
    \left.\dd\rho\right|_{t=0}=\kappa(s)\dd t,\qquad \left.\dd z\right|_{t=0}=\kappa(s)\dd s.
$$
Hence $F$ has nonzero Jacobian there and is a local diffeomorphism of surfaces with boundary.

At a turning point, use the quadrant coordinates $(x,y)$ of \eqref{eq:cornerexpansion}, rescaled so that $\bg(0)=\dd x^2+\dd y^2$ and ordered consistently with the orientation defining $\underline{\ast}$. The corner expansion and $\dd z=\underline{\ast}\dd\rho$ give
$$
    \begin{aligned}
        \rho&=cxy+O((x^2+y^2)^2),\\
        z&=z_0+\frac{c}{2}(y^2-x^2)+O((x^2+y^2)^2),\qquad c>0.
    \end{aligned}
$$
Here $z$ again extends smoothly to the closed quadrant, since $\underline{\ast}\dd\rho$ is smooth and closed there. For a sufficiently small quadrant neighborhood, $z$ is strictly decreasing on $y=0$ and strictly increasing on $x=0$, with values on opposite sides of $z_0$. Together with the injectivity of $F$ on $N^\circ$ and the positivity of $\rho$ away from the rods, this shows that $F$ is injective on a small closed quarter-disc. It is therefore a homeomorphism from the quarter-disc onto its image. The boundary of the image consists of an interval on $\{\rho=0\}$ and a simple arc whose interior lies in $\{\rho>0\}$. By the Jordan curve theorem, the image contains a neighborhood of $(0,z_0)$ relative to $\mathbb H$. Thus $F$ is a local homeomorphism at the turning point as well.

It remains to show that the extended boundary map is bijective. For any $z_0\in\mathbb R$, integrating the lower bound for $\nu_\rho$ in \eqref{eq:prelim-weyl-gradient} gives
$$
    e^{\nu(\rho,z_0)}\leq e^{\nu(1,z_0)}\rho^{-1/2},\qquad 0<\rho\leq1.
$$
The curve $\rho\mapsto F^{-1}(\rho,z_0)$ therefore has finite $\bg$-length as $\rho\to0$. By completeness it converges to a point of $\partial N$, whose image is $(0,z_0)$. Thus the boundary map is onto. If two distinct boundary points had the same image, disjoint local homeomorphism neighborhoods would give two preimages of nearby interior points, contradicting interior injectivity. Hence $F:N\to\mathbb H$ is a bijective local homeomorphism, and therefore a homeomorphism. We have proved regularity and the existence of global $(\rho,z)$ coordinates. With quadratic curvature decay, the final assertion follows from Definitions \ref{def:app gravitational instanton} and \ref{def: app gravitational instanton proper}.
\end{proof}

\begin{lemma}\label{geom:ale-rigidity}
Let $(M,g)$ be a complete connected Ricci-flat ALE four-orbifold with at most isolated orbifold singularities. Suppose it carries an effective isometric $\T^2$-action with connected stabilizers. In a fixed basis of $2\pi$-periodic angular generators, set
$$
    A=(1,0)^\top,\qquad B=(0,1)^\top,\qquad C=(1,1)^\top.
$$
If every circle stabilizer has one of the unoriented directions represented by $A$, $B$, and $C$, then $(M,g)$ is isometric to Euclidean four-space.
\end{lemma}

\begin{proof}
We first show that $M$ is smooth. Every isolated orbifold point is fixed by the connected torus, since the torus preserves the discrete singular set. Let $v,w$ be the primitive rod vectors incident to such a point. In the local toric model, the orbifold group is the quotient of $\mathbb Z^2$ by the sublattice generated by $v,w$, and hence has order $|\det(v,w)|$; see \cite[Section 3]{LiSun}. Any two distinct directions among $A,B,C$ have determinant $\pm1$. Thus the local group is trivial at every torus fixed point, and there are no orbifold singularities. Since the only singularities allowed in the statement were isolated orbifold points, the metric is now smooth everywhere, including along all rods.

Let $\Gamma$ be the finite group at infinity. A sufficiently large invariant cross-section $L$ of the ALE end is diffeomorphic to $S^3/\Gamma$. To obtain such a cross-section, extend the squared ALE radius smoothly over $M$ and average it over $\T^2$. The averaged function has Hessian $2g+o(1)$ and is strictly increasing along radial rays sufficiently far out, so a large level has the required properties. The torus action on $L$ is effective: an element fixing $L$ pointwise also fixes its normal vector, and is therefore the identity isometry. Its stabilizers are connected, and it has no fixed points, since $\T^2$ has no faithful orthogonal representation in dimension three. Consequently $L/\T^2$ is a circle or an interval. The circle case would make $L$ a principal $\T^2$-bundle over a circle and give an infinite fundamental group, contrary to $\pi_1(L)\cong\Gamma$. Thus $L/\T^2$ is an interval.

Let $v_-,v_+\in\{A,B,C\}$ represent the two endpoint stabilizers. The corresponding circles bound meridian discs in the two solid tori making up $L$, so
\begin{equation}\label{eq:prelim-ale-end-group}
    \Gamma\cong\pi_1(L)\cong\mathbb Z^2/(\mathbb Z v_-+\mathbb Z v_+).
\end{equation}
Here $\mathbb Z^2$ records integer coefficients in the angular basis, whose periods are $2\pi$. Finiteness forces the endpoint directions to be distinct, so the same determinant calculation gives $\Gamma=\{1\}$. The end is therefore AE.

We can now apply ordinary Bishop-Gromov comparison, as in the volume comparison argument of \cite[proof of Proposition 4.42]{LiSun}. For any $p\in M$, the function
$$
    R\longmapsto\frac{\operatorname{Vol}_g(B_g(p,R))}{\omega_4R^4}
$$
is nonincreasing, where $\omega_4$ is the volume of the Euclidean unit four-ball. Its limit is one both as $R\downarrow0$, by smoothness, and as $R\to\infty$, by the AE asymptotics. Thus it is identically one, and the equality case of Bishop-Gromov comparison gives $(M,g)\cong\mathbb R^4$.
\end{proof}

\begin{remark}
Lemma \ref{geom:proper-ricciflat} allows arbitrary positive cone angles. In Lemma \ref{geom:ale-rigidity}, however, all rod interiors are smooth and have cone angle $2\pi$. The orbifold formulation is useful for the intermediate rescalings in the proof of Theorem \ref{thm:compactness}, where the conical surfaces have collapsed to an isolated point; the determinant calculation above removes that possible singularity. The main conical bubble in Proposition \ref{prop:bubble} retains its finite conical rods and is instead treated by the angle comparison argument of \cite[proof of Proposition 4.42]{LiSun}.
\end{remark}

\subsection{Topological bounds for toric gravitational instantons}
\label{subsec:instanton-rod-obstructions}
We prove some topological bounds for toric gravitational instantons similar to the case of compact toric Einstein metrics in Section \ref{ss:improved HT}, which may be of independent interest. For AF ends, the sharper bound below strengthens the restrictions obtained by Nilsson \cite[Section~4.1 and Remark~5]{Nilsson}.

Let $(M,g)$ be a smooth, simply connected, proper toric gravitational instanton with $n\geq1$ turning points. The finite rods $\mI_j$, $1\leq j<n$, give spheres $\Sigma_j=\pi^{-1}(\overline{\mI_j})$; set $e_j=[\Sigma_j]^2$. Choose representatives $\mv_j=2\pi u_j$ with $\det(u_{j-1},u_j)=1$. Then
\begin{equation}\label{eq:instanton-rod-recurrence}
    e_j=-\det(u_{j-1},u_{j+1}),\qquad u_{j+1}=-e_j u_j-u_{j-1}.
\end{equation}
The sign differs from the compact recurrence because increasing $z$ is opposite to the positive boundary orientation of the half-plane. The plumbing description gives
\begin{equation}\label{eq:instanton-chain-topology}
    \chi(M)=n,\qquad \tau(M)=\operatorname{sign}Q,\qquad Q= \begin{pmatrix}
        e_1&1&&0\\
        1&e_2&\ddots&\\
        &\ddots&\ddots&1\\
        0&&1&e_{n-1}
    \end{pmatrix};
\end{equation}
see \cite[Theorem~1]{Nilsson}. The matrix represents the intersection form in the basis of finite rod spheres, oriented so adjacent intersections are one. Its zero eigenvalues do not contribute to the signature. For $n=1$, the matrix is empty.

\begin{lemma}\label{lem:instanton-rod-self-intersections}
Every finite rod sphere has $|e_j|\leq2$. Unless $g$ is hyperk\"ahler in one of the two orientations, we must have
\begin{equation}\label{eq:instanton-strict-rod-bound}
    e_j\in\{-1,0,1\}.
\end{equation}
\end{lemma}

\begin{proof}
We apply the Weyl identities of Section~\ref{subsec:toric-weyl-rigidity} in both orientations, with $\lambda=0$. By quadratic curvature decay, the functions $w_\pm=-K\mp J$ tend to zero at infinity, and we have
\begin{equation}\label{eq:instanton-w-equation}
    \Delta_{\bg}w_\pm+3\langle\dd\log\rho,\dd w_\pm\rangle_{\bg}=-3w_\pm^2-4|q_\pm|^2.
\end{equation}
The even extensions across rods and turning points used in the proof of Proposition~\ref{prop:weyl-positive} allow us to apply the maximum principle throughout the quotient. A negative minimum is impossible, so $w_\pm\geq0$. The strong maximum principle then implies that each $w_\pm$ is either positive everywhere or identically zero. In the latter case, \eqref{eq:instanton-w-equation} implies that $q_\pm=0$, and hence $W^\pm=0$.

For each finite rod sphere $\Sigma_j$, the rod curvature identities with $\lambda=0$ imply that
\begin{equation}\label{eq:instanton-rod-curvatures}
    K_j=w_++w_-,\qquad |K_j^\perp|=|w_+-w_-|\leq K_j,
\end{equation}
where $K_j$ and $K_j^\perp$ denote its tangent and normal curvatures. Integrating these identities, we obtain
$$
    2\pi|e_j|\leq\int_{\Sigma_j}|K_j^\perp| \dd A_g
    \leq\int_{\Sigma_j}K_j \dd A_g=4\pi.
$$
If both $W^\pm$ are nonzero, then $w_+,w_->0$, so $|K_j^\perp|<K_j$ and consequently $|e_j|<2$. Since $e_j$ is an integer, we obtain \eqref{eq:instanton-strict-rod-bound}. If one of $W^\pm$ vanishes, Ricci-flatness implies that the induced connection on the corresponding bundle $\Lambda^\pm$ is flat. By simple connectivity, we obtain a global parallel triple, so the metric is hyperk\"ahler in that orientation.
\end{proof}

\begin{proposition}\label{prop:instanton-improved-hitchin-thorpe}
If $g$ is not hyperk\"ahler in either orientation, then $M$ is non-spin unless it is of type $\AF$, and
\begin{equation}\label{eq:instanton-improved-signature}
    |\tau(M)|\leq \left\lfloor\frac{\chi(M)+1}{3}\right\rfloor,
\end{equation}
while in the $\AF$ case one always has
\begin{equation}\label{eq:af-improved-hitchin-thorpe}
    3|\tau(M)|\leq\chi(M).
\end{equation}
\end{proposition}
\begin{remark}
Notice that, it is essential that we exclude hyperk\"ahler metrics in \eqref{eq:instanton-improved-signature}. The $\ALE$ and $\ALF$ Gibbons-Hawking metrics with three collinear centres have $\chi(M)=3$ and $|\tau(M)|=2$, and hence violate this bound.
\end{remark}
\begin{proof}
{\bf Step 1. Bound the signature of the intersection matrix}. Set $m=n-1$; the case $m=0$ is immediate. Let $A_m$ have entries one immediately above and below the diagonal and zero elsewhere. Lemma~\ref{lem:instanton-rod-self-intersections} gives
$$
    A_m-\operatorname{Id}\leq Q\leq A_m+\operatorname{Id}.
$$
Eigenvalue monotonicity bounds the signature by those of the two outer matrices. These have opposite signatures, since conjugating by $\operatorname{diag}(1,-1,1,-1,\ldots)$ changes $A_m$ to $-A_m$. The eigenvalues of $A_m+\operatorname{Id}$ are
$$
    1+2\cos\frac{k\pi}{m+1},\qquad 1\leq k\leq m.
$$
Counting their signs gives $|\tau(M)|\leq\lfloor(m+2)/3\rfloor$, proving \eqref{eq:instanton-improved-signature}.

\medskip\noindent{\bf Step 2. The strengthening in the AF case.} The two infinite rods have the same direction, so $u_n=\pm u_0$ and $n\geq2$. The recurrence gives $\det Q=(-1)^{n-1}\det(u_0,u_n)=0$. Since all adjacent entries are one, a vector in $\ker Q$ is determined by its first component. Thus $\dim\ker Q=1$.

If one of the $W^\pm$ vanished, the signed version of \eqref{eq:instanton-rod-curvatures} would give $e_j=2$ for every finite rod, or $e_j=-2$ for every finite rod. The resulting matrix $\pm2\operatorname{Id}+A_m$ is definite, a contradiction. The flat case is also excluded by Gauss-Bonnet on a finite rod sphere. Hence the first bound applies. Finally,
$$
    \tau(M)\equiv\operatorname{rank}Q=n-2 \equiv\chi(M)\pmod2.
$$
Thus $\chi(M)-3|\tau(M)|$ is even and at least $-1$, so it is nonnegative.

\medskip\noindent{\bf Step 3. Non-Spin.}
Suppose $(M, g)$ is spin and  of $\ALE$ or $\ALF^\pm$ type, we claim that it is hyperk\"ahler for one choice of orientation.  Otherwise, by  Lemma~\ref{lem:instanton-rod-self-intersections} and the spin condition we must have $e_j=0$ for all $j$.
By \eqref{eq:instanton-rod-recurrence}, all primitive rod vectors therefore belong to $\{\pm u_0,\pm u_1\}$, and the two infinite rod vectors form a lattice basis.
In the ALE case, the end must then by  asymptotically Euclidean, so Bishop-Gromov  implies that $(M,g)$ is Euclidean.
In the $\ALF^\pm$ case, the asymptotic fibre circle is generated by $\xi=u_0\pm u_1$ for a suitable choice of sign, and its action is semi-free since $|\det(\xi,u_j)|=1$ along every rod.
By \cite[Proposition~B.1]{LiSun}, the non-hyperkähler possibility must be Taub-Bolt, which is not spin. 

\end{proof}

\begin{remark}\label{rem:instanton-classical-hitchin-thorpe}
The self-intersection bound is proved without using the Hitchin-Thorpe inequality. It implies the signature estimate above, but also gives restrictions on individual rod spheres that neither signature inequality detects. For example, the $\AF$ chains with finite self-intersections $(a,0,-a)$ in \cite[Figure~5]{Nilsson} have $\chi(M)=4$ and $\tau(M)=0$ for every integer $a$. Both the classical and improved signature inequalities hold for every $a$, whereas Lemma~\ref{lem:instanton-rod-self-intersections} requires $|a|\leq1$.

For AF ends, \eqref{eq:af-improved-hitchin-thorpe} strengthens $2\chi(M)\geq3|\tau(M)|$ to $\chi(M)\geq3|\tau(M)|$. Nilsson observes that the former inequality excludes none of the AF rod structures with four turning points \cite[Remark~5]{Nilsson}. The latter excludes all four exceptional structures in \cite[Figure~6]{Nilsson}: their finite self-intersections are
$$
    \pm(1,2,1)\qquad\text{or}\qquad\pm(2,1,2),
$$
so $\chi(M)=4$ and $|\tau(M)|=2$.

For $\ALE$ and $\ALF$ ends, \eqref{eq:instanton-improved-signature} provides an additional restriction outside the hyperk\"ahler case; the classical inequalities also involve boundary terms at infinity \cite[Section~4.1]{Nilsson}.
\end{remark}

For example, the $\AF$ rod sequence with primitive representatives
$$
    (u_0,u_1,u_2,u_3,u_4) =\bigl((0,1),(-1,0),(1,-1),(-1,2),(0,-1)\bigr)
$$
has finite self-intersections $(1,2,1)$ and intersection eigenvalues $0,1,3$. It has $\chi=4$ and $\tau=2$, so it satisfies the classical inequality but violates \eqref{eq:af-improved-hitchin-thorpe}. Likewise, the $\AF$ chains $(a,0,-a)$ have $\chi=4$ and $\tau=0$ for every integer $a$, whereas the individual rod bound requires $|a|\leq1$.

\

An easy consequence of the above discussion is the following uniqueness theorem for Eguchi-Hanson metrics. 

\begin{proposition}\label{prop:toric-eguchi-hanson}
A smooth toric $\ALE$ gravitational instanton with asymptotic cone
$\R^4/\mathbb Z_2$ is isometric, up to scaling, to the Eguchi-Hanson
metric.
\end{proposition}

\begin{proof}
By Anderson's fundamental-group estimate  and  Bishop-Gromov volume comparison, we know $M$ is simply connected. This implies $M/\mathbb T^2$ is also simply connected, then  by Lemma~\ref{geom:proper-ricciflat} we know $(M,g)$ is a proper $\ALE$ gravitational instanton. 

Denote $\Sigma=\pi^{-1}(\overline{\mI_0})$. It is  a totally geodesic surface diffeomorphic to $\R^2$. The residual circle collapses smoothly at the endpoint, so $\Sigma\cong\R^2$. We denote its tangent and normal curvatures by $K_\Sigma$ and $K_\Sigma^\perp$. The local identity \eqref{eq:instanton-rod-curvatures} also applies to $\Sigma$.

We use an equivariant $\ALE$ identification as in the proof of Proposition~\ref{app:maximal-volume-growth}. On the covering end, $\Sigma$ corresponds to a plane $P\subset\R^4$, and the nontrivial deck transformation acts as $-I$ on both $P$ and $P^\perp$. In particular, the induced metric on $\Sigma$ has asymptotic cone angle $\pi$. If $D_R\subset\Sigma$ is the disc bounded by a large radial cross-section $\gamma_R$, then
$$
    \int_{\gamma_R}k_g \dd s_g\longrightarrow\pi.
$$
By Gauss-Bonnet, we get
\begin{equation}\label{eq:eh-tangent-integral}
    \int_\Sigma K_\Sigma \dd A_g=\pi.
\end{equation}

We next compute the limiting normal holonomy. A lift of $\gamma_R$ joins $x$ to $-x$ in $P$. In the normal frame obtained by projecting and orthonormalizing a fixed Cartesian frame of $P^\perp$, parallel transport tends to the identity: the connection is $O(R^{-1-\delta})$ along a curve of length $O(R)$, for some $\delta>0$. Hence we know
$$
    \operatorname{Hol}_{N\Sigma}(\gamma_R)
    \longrightarrow-1\in\SO(2).
$$
Since $\Sigma\cong\R^2$, its oriented normal bundle is trivial. Stokes' theorem  then implies that 
$$
    \exp\left(i\int_{D_R}K_\Sigma^\perp \dd A_g\right)
    =\operatorname{Hol}_{N\Sigma}(\gamma_R)^{\pm1}.
$$
Taking limit, we obtain
\begin{equation}\label{eq:eh-normal-integral}
    \int_\Sigma K_\Sigma^\perp \dd A_g\in\pi+2\pi\Z.
\end{equation}
Combining \eqref{eq:instanton-rod-curvatures}, \eqref{eq:eh-tangent-integral}, and \eqref{eq:eh-normal-integral}, we find
$$
    \pi\leq\left|\int_\Sigma K_\Sigma^\perp \dd A_g\right|\leq\int_\Sigma|K_\Sigma^\perp| \dd A_g\leq\int_\Sigma K_\Sigma \dd A_g=\pi.
$$
If both $W^\pm$ were nonzero, the strict positivity established in the proof of Lemma~\ref{lem:instanton-rod-self-intersections} would give
$$
    K_\Sigma-|K_\Sigma^\perp|=2\min(w_+,w_-)>0,
$$
contradicting equality of the integrals. Thus one of $W^\pm$ vanishes identically, and $g$ is hyperk\"ahler by simple connectivity. It then follows easily that $g$ must be the Eguchi-Hanson metric, up to scaling and isometry.
\end{proof}

\section{Toric \texorpdfstring{$\Ae$}{AE} gravitational instantons  with equal cone angles}
\label{app:ae-models}

In this appendix we study toric conical equal cone angle $\Ae$ gravitational instantons that we need for the gluing construction. See Appendix \ref{app:asymptotics of toric gravitational instantons} for the reduction of proper toric conical gravitational instantons to axisymmetric harmonic maps. Theorem \ref{thm:appendix classification of asymptotics} has already proved that for a proper toric conical gravitational instanton, it is of Asymptotic Type $\ALE,\ALF^\pm$, or $\AF_\beta$. As it proved in \cite{LiSun}, given a rod structure $\mR$ with an Asymptotic Type, there is a unique axisymmetric harmonic map that is strongly tamed by $\mR$ and the Asymptotic Type.  Note that in \cite{LiSun}, the definition of rod structures includes the Asymptotic Types, which is different from the notion of rod structures we use here. Here $\Ae$ means that the end is $\ALE$ and has trivial group at infinity.

The sequences $W$ used in our construction admit two explicit descriptions. For $(B,(A,B)^n,A)$, diagonality forces the harmonic map to take values in a geodesic of the hyperbolic plane. For $BCA$, $BACA$, $BCACA$, and $BCABCA$, the Hermitian structure allows us to use the Tod ansatz, with its description by a convex piecewise linear function on the axis. The first alternating case $BABA$ is also Hermitian. In both families, the resulting formulas give control of the cone-angle map. We use this to prove existence and uniqueness of the equal-angle metric up to homothety, and to show that its decaying $\mathbb T^2$-invariant infinitesimal deformations at fixed cone angles consist only of scaling modulo coordinate changes. These properties provide the linear analysis required for gluing.

\subsection{Toric conical \texorpdfstring{$\Ae$}{Ae} gravitational instantons with two rod directions}
\label{subsec:app AE with two rods}

In this section, we consider rod structures that only consist of two rod directions, $e_1=[1,0]^\top$ and $e_2=[0,1]^\top$.  We prove

\begin{theorem}\label{thm:alternating-ae}
For every $n\geq1$ there is a proper toric conical $\Ae$ gravitational instanton, whose rod vectors are
$$
    (e_2,e_1,e_2,\ldots,e_2,e_1),
$$
with $2n$ finite rods, smooth cone angles along infinite rods, and the same cone angle $2\pi\beta_n$ on all finite rods. It is unique up to homothety and equivariant isometry among proper toric conical AE gravitational instantons with these rod vectors, whose finite rods all have the same cone angle and infinite rods have cone angle $2\pi$. Moreover, $\beta_n>1$.
\end{theorem}

These $\Ae$ gravitational instantons are actually explicit. The rod lengths are invariant under reflection of $z$ at the middle turning point, and the metric admits an orientation-reversing isometry.

We claim that the condition that there are only two rod directions forces the metric to be diagonal. For any proper toric conical $\Ae$ gravitational instanton with rod directions $e_1,e_2$ and smooth infinite rods, let $\Phi$ be its associated harmonic map. The target isometry $\Phi\mapsto S^\top\Phi S$, where $S=\operatorname{diag}(1,-1)$, fixes both rod lines and preserves the normalized Euclidean asymptotics. Thus $S^\top\Phi S$ is strongly tamed by the same rod structure and has the same asymptotic data. By the uniqueness theorem \cite[Theorem~4.24]{LiSun}, we have
$$
    S^\top\Phi S=\Phi,
$$
which forces $\Phi_{12}=0$. Consequently, the metric has the diagonal form \eqref{eq:alternating-metric}, and the harmonic-map equation reduces to the scalar axisymmetric Laplace equation. This argument applies for arbitrary positive rod lengths, independently of the equal-angle condition.

Set $m=2n+1$ and choose
$$
    a_1<\cdots<a_m,\qquad \epsilon_i=(-1)^{i+1}.
$$
For $\rho>0$, set
\begin{gather*}
    R_i=\sqrt{\rho^2+(z-a_i)^2},\qquad\mu_i=R_i-(z-a_i),\qquad E_{ij}=R_iR_j+(z-a_i)(z-a_j)+\rho^2,\\
    U=\prod_{i=1}^m\mu_i^{\epsilon_i},\qquad V=\frac{\prod_{i<j}E_{ij}^{-\epsilon_i\epsilon_j}}{2^{n+1}\prod_iR_i}.
\end{gather*}
Then
\begin{equation}\label{eq:alternating-metric}
    g_a=V(\dd\rho^2+\dd z^2)+U\dd\phi_1^2+\rho^2U^{-1}\dd\phi_2^2
\end{equation}
is positive definite and Ricci-flat. Indeed, with $u=\frac12\log U$ and $\nu=\frac12\log V$, the Ricci-flat equations are
\begin{align*}
    &u_{\rho\rho}+\rho^{-1}u_\rho+u_{zz}=0,\\
    &\nu_\rho=\rho(u_\rho^2-u_z^2)-u_\rho,\\
    &\nu_z=2\rho u_\rho u_z-u_z.
\end{align*}
The first follows from the same equation for $\log\mu_i$. Differentiating the displayed product for $V$ gives the other two. All factors are positive when $\rho>0$.

The outer rods have cone angle $2\pi$, and the cone angle $2\pi\alpha_k$ of $(a_k,a_{k+1})$ is
\begin{equation}\label{eq:alternating-angles}
    \alpha_k(a)=\prod_{i\leq k<j}(a_j-a_i)^{\epsilon_i\epsilon_j},\qquad 1\leq k\leq 2n.
\end{equation}
For example, on a finite rod the expansions
$$
    \mu_i=\frac{\rho^2}{2(z-a_i)}(1+O(\rho^2))\quad(i\leq k),\qquad \mu_j=2(a_j-z)(1+O(\rho^2))\quad(j>k)
$$
give $\alpha_k^2=\lim_{\rho\to0} U/(\rho^2V)$ for odd $k$, and $\alpha_k^2=\lim_{\rho\to0} 1/(UV)$ for even $k$. Substitution gives \eqref{eq:alternating-angles}; all factors involving $z$ cancel. The same calculation on either outer rod gives one.

At infinity let $R=\sqrt{\rho^2+z^2}$. Since $\sum_i\epsilon_i=1$ and $\sum_{i<j}\epsilon_i\epsilon_j=-n$,
$$
    U=(R-z)(1+O_\infty(R^{-1})),\qquad V=(2R)^{-1}(1+O_\infty(R^{-1})).
$$
The quotients in these formulas have the same differentiated expansions at the outer axes: one can use $\mu_i=\rho^2/(R_i+z-a_i)$ near the positive axis. Setting $x^2=R+z$ and $y^2=R-z$ therefore gives
$$
    g_a=\dd x^2+x^2\dd\phi_2^2+\dd y^2+y^2\dd\phi_1^2+O_\infty(r^{-2}),\qquad r^2=x^2+y^2.
$$
The error is understood in the Euclidean tensor norm, with its corresponding differentiated bounds. This is an AE end with smooth outer axes.

For every choice of turning points $a_1<\cdots<a_{2n+1}$, the formulas above define a harmonic map
$$
    \Phi_a=\rho^{-1}\operatorname{diag}(U,\rho^2U^{-1})
$$
strongly tamed by the corresponding alternating rod structure of $\Ae$ type. By the existence and uniqueness theorem in \cite{LiSun}, this is precisely the harmonic map associated with these data. The regularity theorem in \cite{LiSun} then gives the conical extension of $g_a$ across the fixed surfaces and torus fixed points. Since the turning points are arbitrary, the metrics $g_a$ exhaust all proper toric conical $\Ae$ gravitational instantons with these alternating rod vectors, smooth infinite rods, and the prescribed angular lattice. Here the augmentation is fixed by the Euclidean normalization at infinity.

By a translation and homothety we may normalize so that $a_1=0$, $a_m=1$. Write
$$
    \ell_k=a_{k+1}-a_k>0,\qquad \sum_{k=1}^{2n}\ell_k=1,\qquad d_{ij}=a_j-a_i.
$$
On the simplex of $\bl$, define
\begin{equation}\label{eq:alternating-functional}
    F(\bl)=\sum_{i<j}\epsilon_i\epsilon_j d_{ij}\log d_{ij}.
\end{equation}
The convention $0\log0=0$ extends $F$ continuously to the closed simplex. The sum $s_k=\sum_{i\leq k}\epsilon_i$ is zero or one. Hence
$$
    \sum_{i\leq k<j}\epsilon_i\epsilon_j=s_k(1-s_k)=0,
$$
and differentiation gives
\begin{equation}\label{eq:alternating-gradient}
    \dd F_{\bl}(\dot{\bl})=\sum_{k=1}^{2n}\log\alpha_k(\bl) \dot\ell_k, \qquad \sum_{k=1}^{2n}\dot\ell_k=0.
\end{equation}
Thus equal angles are precisely the critical points of $F$ under the normalization of total length.

\begin{lemma}\label{lem:alternating-hessian}
On the open simplex
$$
    \left\{\bl=(\ell_1,\ldots,\ell_{2n})\in(0,\infty)^{2n}:\sum_{k=1}^{2n}\ell_k=1\right\},
$$
the Hessian of $F$ is negative definite. More precisely, for every nonzero variation $\dot{\bl}=(\dot\ell_1,\ldots,\dot\ell_{2n})$ satisfying $\sum_k\dot\ell_k=0$, one has
$$
    D^2F(\bl)[\dot\bl,\dot\bl]<0.
$$
Consequently, the map
$$
    \bl\longmapsto\bigl(\log\alpha_1(\bl),\ldots,\log\alpha_{2n}(\bl)\bigr)
$$
has differential of rank $2n-1$ everywhere on this simplex.
\end{lemma}

\begin{proof}
Fix $\bl$ and a variation $\dot\bl$ satisfying $\sum_k\dot\ell_k=0$. Set $m=2n+1$ and define the corresponding variations of the turning points by
$$
    \dot a_1=0,\quad \dot a_i=\sum_{k=1}^{i-1}\dot\ell_k,\qquad 2\leq i\leq m.
$$
Thus $\dot a_m=0$, consistently with the normalization $a_1=0$, $a_m=1$. Differentiating twice along $\bl+t\dot\bl$ gives
$$
    D^2F(\bl)[\dot\bl,\dot\bl]=\sum_{i<j}\epsilon_i\epsilon_j\frac{(\dot a_j-\dot a_i)^2}{a_j-a_i}.
$$

For $\tau\in\mathbb R$, define
$$
    H(\tau)=\sum_{i<j}\epsilon_i\epsilon_j\left|\tau(a_j-a_i)+\dot a_j-\dot a_i\right|.
$$
Order the numbers $\tau a_i+\dot a_i$ as $y_1\leq\cdots\leq y_m$, and let $q_j\in\{1,-1\}$ be the sign $\epsilon_i$ attached to $y_j$. Set
$$
    S_k=\sum_{j=1}^k q_j.
$$
Since $\sum_jq_j=\sum_i\epsilon_i=1$, writing each difference $y_j-y_i$ as a sum of consecutive gaps yields
$$
    H(\tau)=\sum_{k=1}^{m-1}S_k(1-S_k)(y_{k+1}-y_k)\leq0.
$$
Indeed, $S_k$ is an integer, so $S_k(1-S_k)\leq0$. For sufficiently large positive or negative $\tau$, the signs occur in the original or reversed alternating order. In either case every $S_k$ is zero or one, so $H(\tau)=0$.

For $d>0$ and $c\in\mathbb R$, direct integration gives
$$
    \int_{\mathbb R}\left(|d\tau+c|-d|\tau|-c \operatorname{sgn}\tau\right) \dd\tau=\frac{c^2}{d}.
$$
The identity
$$
    \sum_{i\leq k<j}\epsilon_i\epsilon_j=0
$$
implies both
$$
    \sum_{i<j}\epsilon_i\epsilon_j(a_j-a_i)=0,\qquad \sum_{i<j}\epsilon_i\epsilon_j(\dot a_j-\dot a_i)=0.
$$
Applying the integral formula to each pair therefore gives
$$
    D^2F(\bl)[\dot\bl,\dot\bl]=\int_{\mathbb R}H(\tau)\dd\tau\leq0.
$$

We now determine when equality holds. If the integral vanishes, then $H\equiv0$, since $H$ is continuous and non-positive. Whenever the numbers $\tau a_i+\dot a_i$ are distinct, all their consecutive gaps are positive. Hence every $S_k$ must be zero or one, which forces the ordered signs to alternate, beginning and ending with $+1$.

We claim that the points $(a_i,\dot a_i)$ are collinear. We use Kelly's elementary minimum-distance proof of the Sylvester-Gallai theorem. Suppose otherwise. Among all positive distances from one of these points to a line through two others, choose a smallest one, realized by a point $P$ and a line $L$. The line $L$ contains exactly two of the points. Indeed, if it contained at least three, two of them, say $A$ and $B$, would lie on the same closed side of the perpendicular foot from $P$, with $B$ farther from that foot. Then
$$
    \operatorname{dist}(A,PB)=\operatorname{dist}(P,L)\frac{|AB|}{|PB|}<\operatorname{dist}(P,L),
$$
contradicting the choice of $P$ and $L$.

Since the $a_i$ are distinct, $L$ is not vertical. Thus there is a value $\tau_0$ at which the two affine functions $\tau a_i+\dot a_i$ corresponding to the points on $L$ have the same value, and no third function has that value. These two functions occupy adjacent positions in the ordered list near $\tau_0$, and exchange positions as $\tau$ passes through $\tau_0$. Their signs must be opposite for the ordering to alternate on one side, but exchanging them prevents the ordering from alternating on both sides. Any simultaneous crossings at other values do not affect these two adjacent positions. This contradicts $H\equiv0$.

The points $(a_i,\dot a_i)$ are therefore collinear. Since they include $(a_1,\dot a_1)=(0,0)$ and $(a_m,\dot a_m)=(1,0)$, we obtain $\dot a_i=0$ for every $i$, and consequently $\dot\ell_k=\dot a_{k+1}-\dot a_k=0$ for every $k$. Thus the Hessian is negative definite on the tangent space $\{\dot\bl:\sum_k\dot\ell_k=0\}$.

Finally, the identity
$$
    \dd F=\sum_{k=1}^{2n}\log\alpha_k \dd\ell_k
$$
gives
$$
    D^2F(\bl)[\dot\bl,\dot\bl]=\sum_{k=1}^{2n}\dot\ell_k D_{\dot\bl}\log\alpha_k.
$$
If $\dot\bl$ lies in the kernel of the differential of the logarithmic cone-angle map, the right-hand side vanishes. Negative definiteness then forces $\dot\bl=0$. The differential is therefore injective on the $(2n-1)$-dimensional tangent space, proving the rank statement.
\end{proof}

The continuous function $F$ attains a maximum on the closed simplex. We show that no boundary point is a maximum. On the boundary of the simplex, some $\ell_k=a_{k+1}-a_k=0$. This means some consecutive turning points coincide. We show that slightly separating the coincident turning points increases $F$, so such a configuration cannot maximize it.

At a boundary point, partition the turning points into maximal consecutive blocks having the same position. Since the coefficients $\epsilon_i$ alternate between $+1$ and $-1$, their sum within each group is zero if the number of turning points is even, and is $+1$ or $-1$ if that number is odd. At least one group contains more than one turning point.

Suppose first that a group contains an even number of turning points. Divide it into two consecutive groups, each containing an odd number of turning points. The sums of their respective coefficients $\epsilon_i$ are opposite, so their product is $-1$. Separate the positions of these two groups by a small distance $\delta>0$, keeping all turning points in $[0,1]$ and preserving $a_1=0$ and $a_m=1$. The contribution to $F$ from pairs with one turning point in each of these two groups changes from zero to $-\delta\log\delta$. Contributions involving turning points outside the original group change by $O(\delta)$, since their distances from that group remain bounded away from zero. Consequently,
$$
    F_{\mathrm{new}}-F_{\mathrm{old}}=-\delta\log\delta+O(\delta)>0
$$
for sufficiently small $\delta>0$.

It remains to consider a group containing an odd number, at least three, of coincident turning points. Let $c$ be their common position and define
$$
    \sigma=\sum_{\{i \mid a_i=c\}}\epsilon_i\in\{-1,1\}.
$$
Divide this group into three consecutive groups, each containing an odd number of turning points. The sums of their respective coefficients $\epsilon_i$ are $\sigma,-\sigma,\sigma$.

Suppose first that $0<c<1$. Move these three groups to $c-\delta,c,c+\delta$, respectively, choosing $\delta>0$ sufficiently small to preserve the ordering of all turning points. The contribution to $F$ from pairs within the original group increases by
$$
    -\delta\log\delta-\delta\log\delta+2\delta\log(2\delta)=2\delta\log2.
$$
For each turning point outside this group, its distances from the three new positions remain bounded away from zero. The terms linear in $\delta$ cancel between the groups at $c-\delta$ and $c+\delta$, since both have coefficient sum $\sigma$. Taylor expansion therefore shows that all contributions involving turning points outside the original group change by $O(\delta^2)$. Consequently,
$$
    F_{\mathrm{new}}-F_{\mathrm{old}} =2\delta\log2+O(\delta^2)>0
$$
for sufficiently small $\delta>0$.

Suppose next that $c=0$. Place the three groups at $0,p\delta,\delta$, respectively, where $0<p<1$. This preserves $a_1=0$. Their contribution to $F$ increases by
$$
    \delta h(p),\qquad h(p)=-p\log p-(1-p)\log(1-p).
$$
For each turning point $a_j$ outside the original group, we have $a_j>0$, and the change in its contribution with the original group is exactly
$$
    \sigma\epsilon_j\left[(a_j-\delta)\log(a_j-\delta)-(a_j-p\delta)\log(a_j-p\delta)\right].
$$
Choose $\delta$ smaller than half the distance from zero to every turning point outside the original group. The mean value theorem then bounds the absolute value of the sum of these changes by
$$
    C(1-p)\delta,
$$
where $C$ is independent of $p$ and $\delta$. Thus
$$
    F_{\mathrm{new}}-F_{\mathrm{old}}\geq\delta\bigl(h(p)-C(1-p)\bigr).
$$
Since
$$
    \frac{h(p)}{1-p}\longrightarrow+\infty\qquad\text{as }p\longrightarrow1,
$$
we can first choose $p$ sufficiently close to one and then choose $\delta>0$ sufficiently small. The resulting change in $F$ is positive. The case $c=1$ follows by reflecting the construction about $1/2$, preserving $a_m=1$.

All these perturbations remain in the closed simplex. Therefore no boundary point can maximize $F$. Its maximum lies in the open simplex and is unique by Lemma \ref{lem:alternating-hessian}. Since
$$
    dF=\sum_{k=1}^{2n}\log\alpha_k \dd\ell_k,
$$
the critical-point condition under $\sum_k\ell_k=1$ is
$$
    \log\alpha_1=\cdots=\log\alpha_{2n}.
$$
Hence there is a unique choice of normalized rod lengths for which all finite cone angles are equal. The identification of all metrics with the family $g_a$ then gives uniqueness of the corresponding metric up to homothety and equivariant isometry.

\ 

Finally, after normalizing the total rod length to one, the Jacobian of the equations $\log\alpha_i-\log\alpha_{2n}=0$, $1\leq i\leq 2n-1$, with respect to the rod lengths is non-singular at the equal-angle solution. We finish by verifying this assertion. Consider the map
$$
    \begin{aligned}
        \left\{ \bl\in(0,\infty)^{2n}: \sum_{k=1}^{2n}\ell_k=1 \right\}&\longrightarrow \mathbb R^{2n-1},\\
        \bl &\longmapsto \left( \log\frac{\alpha_1(\dot\bl)}{\alpha_{2n}(\dot\bl)},\ldots, \log\frac{\alpha_{2n-1}(\dot\bl)}{\alpha_{2n}(\dot\bl)} \right).
    \end{aligned}
$$
At the equal-angle solution, its differential is an isomorphism from the tangent space
$$
    \left\{ \dot\bl\in\mathbb R^{2n}: \sum_{k=1}^{2n}\dot\ell_k=0 \right\}
$$
onto $\mathbb R^{2n-1}$. Indeed, if $\dot\bl$ lies in its kernel, then
$$
    D_{\dot\bl}\log\alpha_k=c
$$
for some common value $c$ and every $k$. Therefore
$$
    D^2F(\dot\bl)[\dot\bl,\dot\bl]=\sum_k\dot\ell_kD_{\dot\bl}\log\alpha_k=c\sum_k\dot\ell_k=0.
$$
Lemma \ref{lem:alternating-hessian} implies $\dot\bl=0$. The differential is thus injective, and the domain and target both have dimension $2n-1$.

Denote the common finite cone angle by $\beta_n$. The formula for the first cone angle gives the common cone angle is
$$
    \beta_n=\alpha_1=\ldots=\alpha_{2n}=\prod_{j=1}^{n}\frac{a_{2j+1}-a_1}{a_{2j}-a_1}>1.
$$
This justifies $\beta_n>1$. For the toric conical AE gravitational instanton with equal cone angle on finite rods, the configuration of turning points on $[0,1]$ must be invariant under the reflection $z\mapsto 1-z$, which follows from the uniqueness of the maximum of $F$. The map $(\rho,z,\phi_1,\phi_2)\mapsto(\rho,z,-\phi_1,\phi_2)$ induces an isometry of the gravitational instanton, swapping its two orientations.

We also record an identity used below. Let $\dot\bl$ be any variation tangent to the normalized simplex, and let $\dot\alpha_k$ and $\dot d_{ij}$ denote the corresponding variations of $\alpha_k$ and $d_{ij}$. Since
$$
    d_{ij}=\sum_{k=i}^{j-1}\ell_k, \qquad \log\alpha_k =\sum_{i\leq k<j}\epsilon_i\epsilon_j\log d_{ij},
$$
we have
$$
    \begin{aligned}
        \sum_k\ell_k\frac{\dot\alpha_k}{\alpha_k}&= \sum_{i<j}\epsilon_i\epsilon_j\frac{\dot d_{ij}}{d_{ij}}\sum_{k=i}^{j-1}\ell_k\\
        &=\sum_{i<j}\epsilon_i\epsilon_j\dot d_{ij}\\
        &=\sum_k\left(\sum_{i\leq k<j}\epsilon_i\epsilon_j\right)\dot\ell_k=0.
    \end{aligned}
$$
Thus, on the normalized simplex,
\begin{equation}\label{eq:alternating-balance}
    \sum_k\ell_k \dd\log\alpha_k=0, \qquad\sum_k A_k\dd\alpha_k=0,\qquad A_k=\frac{2\pi\ell_k}{\alpha_k}.
\end{equation}
Here $A_k$ is the area of the fixed surface corresponding to the $k$-th finite rod. Indeed, its induced area element is $\alpha_k^{-1} \dd z \dd\phi$, where $\phi$ is the noncollapsing angular coordinate of period $2\pi$. In particular, if a variation in this family satisfies $\dot\alpha_1=\cdots=\dot\alpha_{2n}=c$, then
$$
    0=\sum_k A_k\dot\alpha_k=c\sum_k A_k.
$$
Since every $A_k>0$, necessarily $c=0$.

\subsection{Four Hermitian \texorpdfstring{$\Ae$}{Ae} metrics}\label{app:hermitian-ae}

In this section we study four Hermitian toric conical $\Ae$ gravitational instantons using the piecewise linear description of the Tod ansatz in \cite{BG,AL}.

Set $A=e_1$, $B=e_2$, and $C=e_1+e_2$. The rod vectors below include the two infinite rods. We write $2\pi\alpha_j$ for the cone angle along the $j$-th finite rod, and $\beta_*$ for their common cone angle when all finite cone angles are equal. The infinite rods are smooth. Let $\phi_1,\phi_2$ be the $2\pi$-periodic angular coordinates associated with $e_1,e_2$, and set $o_0=\dd\rho\wedge\dd z\wedge\dd\phi_1\wedge\dd\phi_2$.

\begin{theorem}\label{thm:hermitian-ae}
For each sequence of rod vectors in the following table there is a unique proper toric conical $\Ae$ gravitational instanton with smooth infinite rods and common finite cone angle $2\pi\beta_*$, up to homothety and equivariant isometry:
$$
    \begin{array}{c|c|c}
        \text{rod vectors}&\beta_*&\text{rod lengths, up to scale}\\ \hline
        BCA&2&(1)\\
        BACA&1+2/\sqrt3&(1,2)\\
        BCACA&4/\sqrt3&(3,2\sqrt3-2,\sqrt3+2)\\
        BCABCA&1+\sqrt2&(1,\sqrt2-1,\sqrt2-1,1).
    \end{array}
$$
\end{theorem}

The first two metrics are Hermitian in orientation $o_0$, and the last two are Hermitian in $-o_0$. The first is hyperk\"ahler in $-o_0$, and the second is Hermitian in both orientations. This follows from couting the $\ALE$ topological degree of \cite{LiSun}. Combining this with the representation of \cite[Section~3]{BG} and \cite[Propositions~2.4 and~4.5]{AL}, every such instanton, in the Hermitian non-K\"ahler orientation, takes the form of the ansatz below.

Other presentations in the theorem give equivalent metrics. For example, reversing the order of the rods and interchanging $A$ and $B$ takes $BACA$ to $BCBA$. This identifies the metrics after the corresponding change of angular basis, and gives the two more sequences that appear in \eqref{eq:alternating-compact-word}.

\subsubsection{The ansatz and its rod data}

Choose $N\geq2$, turning points $z_1<\cdots<z_N$, and constants $a_i>0$ with $\sum_{i=1}^N a_i=1$. The ansatz is parametrized by the positive convex piecewise linear function
$$
    f(z)=\sum_{i=1}^N a_i|z-z_i|.
$$
Its slopes on the two infinite rods are $-1$ and $1$. The absence of an additive constant in this formula is the $\Ae$ condition. We recall the reconstruction of the metric to fix the normalization.

The axisymmetric harmonic function determined by $f$ is
\begin{equation}\label{eq:hermitian-potential}
    \begin{aligned}
        \mathcal U(\rho,z)&=\sum_{i=1}^N a_i\mathcal U_0(\rho,z-z_i),\\
        \mathcal U_0(\rho,t)&=2\sqrt{\rho^2+t^2}-t\log\frac{\sqrt{\rho^2+t^2}+t}{\sqrt{\rho^2+t^2}-t}.
    \end{aligned}
\end{equation}
In particular, $\mathcal U(\rho,z)=f(z)\log\rho^2+O(1)$ as $\rho\to0$ on the interior of each rod. Set
$$
    \mu=\sum_{i=1}^N a_i z_i,\qquad c=-2\sum_{i=1}^N a_i(z_i-\mu)^2<0,
$$
and let $H$ satisfy
$$
    H_\rho=-\rho\mathcal U_z,\qquad H_z=\rho\mathcal U_\rho.
$$
Tod's reconstruction is
\begin{equation}\label{eq:hermitian-ae-explicit}
    g_f=W^{-1}(\dd\tau+F\dd y)^2+W\rho^2\dd y^2+e^{2\nu}(\dd\rho^2+\dd z^2),
\end{equation}
where
$$
    \begin{aligned}
        W&=\frac1{2c}\left(\rho\mathcal U_\rho+\frac{\mathcal U_\rho^2\mathcal U_{zz}}{\mathcal U_{zz}^2+\mathcal U_{\rho z}^2}\right),\\
        e^{2\nu}&=\frac{W\rho^2}{4}\left(\mathcal U_{zz}^2+\mathcal U_{\rho z}^2\right),\\
        F&=\frac1{2c}\left(\frac{\rho\mathcal U_\rho^2\mathcal U_{\rho z}}{\mathcal U_{zz}^2+\mathcal U_{\rho z}^2}-\rho^2\mathcal U_z-2H\right).
    \end{aligned}
$$
Choose the additive constant in $H$ so that $F=0$ on the right infinite rod. The choice of $c$ then gives $F=-1$ on the left infinite rod. These are the standard Tod formulas of \cite[Proposition~2.4]{AL}; all subsequent calculations use only the function $f$ and its rod data.

For $BCA$ and $BACA$ take $(\tau,y)=(-\phi_2,\phi_1-\phi_2)$, while for $BCACA$ and $BCABCA$ take $(\tau,y)=(\phi_2,\phi_1+\phi_2)$. These choices identify the two infinite rod lines with $B$ and $A$ and fix the angular lattice. The construction gives the normalized Euclidean end. It depends smoothly on the parameters while the turning points remain distinct and the $a_i$ remain positive. In the corresponding normalized end charts, its asymptotic expansion gives
$$
    |\nabla_{g_{\mathrm{Euc}}}^k(g_f-g_{\mathrm{Euc}})|_{g_{\mathrm{Euc}}}=O(r^{-2-k}),\qquad k\geq0,
$$
where $r$ is the Euclidean radius. The same estimates hold after differentiation in the parameters, locally uniformly in the parameter space. When the collapsing fields agree with the prescribed primitive rod lines, as verified below for the metrics and nearby families under consideration, the cited regularity results give the conical extension along the fixed surfaces and at the torus fixed points.

There are $q=N-1$ finite rods. Write
$$
    \begin{aligned}
        I_0&=(-\infty,z_1),\qquad I_j=(z_j,z_{j+1})\quad(1\leq j\leq q),\qquad I_N=(z_N,\infty),\\
        s_j&=f'|_{I_j},\qquad f_i=f(z_i),\qquad \ell_j=z_{j+1}-z_j.
    \end{aligned}
$$
Thus $s_0=-1$, $s_N=1$, $a_i=(s_i-s_{i-1})/2$, and
$$
    f_{j+1}-f_j=s_j\ell_j.
$$
The condition that $f$ have no additive constant is equivalent to
\begin{equation}\label{eq:tod-ae-height-identity}
    f_1+f_N=z_N-z_1=\sum_{j=1}^q\ell_j.
\end{equation}

We next recall the rod formulas from \cite[Lemmas~5.1-5.2]{AL}. A normalized rod generator means a collapsing Killing field for which period $2\pi$ would give cone angle $2\pi$. Up to sign, it is
$$
    u_j=\begin{cases}
        s_j(\partial_y-F_j\partial_\tau),&s_j\ne0,\\[2pt]
        -c^{-1}f_j^2\partial_\tau,&s_j=0,
    \end{cases}
$$
where $F_j$ is the constant value of $F$ on $I_j$ when $s_j\ne0$. If $v_j$ is the primitive lattice vector of the rod, then $v_j=\pm\alpha_j u_j$, where we set $\alpha_0=\alpha_N=1$ for the two infinite rods. Thus its cone angle is $2\pi\alpha_j$.

If $s_{i-1}$ and $s_i$ are both nonzero, then
\begin{equation}\label{eq:tod-normalized-rods}
    F_i-F_{i-1}=\frac{f_i^2}{c}\left(\frac1{s_i}-\frac1{s_{i-1}}\right).
\end{equation}
If $s_j=0$, then $f_j=f_{j+1}$ and
\begin{equation}\label{eq:tod-zero-slope}
    F_{j+1}-F_{j-1}=\frac1c\left[f_j^2\left(\frac1{s_{j+1}}-\frac1{s_{j-1}}\right)-2\ell_j f_j\right].
\end{equation}
Consecutive primitive rod vectors in our four configurations have determinant of absolute value one. Taking determinants of their normalized generators gives
$$
    f_i^2(s_i-s_{i-1})=\frac{|c|}{\alpha_{i-1}\alpha_i},\qquad 1\leq i\leq N.
$$
This identity also holds at an endpoint of a zero-slope rod.

\subsubsection{The equal-angle metrics}

For the first two configurations, choose the signed rod lists $(B,-C,A)$ and $(B,-A,-C,A)$. The coefficient of $\partial_y$ is $-1$ on $B$, $1$ on $A$, and $0$ on $C$. Consequently, the slopes for $BCA$ are $(-1,0,1)$. Up to translation, its function is
$$
    f(z)=\frac12|z|+\frac12|z-L|,\qquad L>0.
$$
Its finite rod has height $L/2$, and $c=-L^2/2$. Since its primitive rod vector is $C=-\partial_\tau$, the zero-slope generator gives
$$
    \alpha_1=\frac{|c|}{(L/2)^2}=2.
$$
There is therefore exactly one such metric up to scale.

For $BACA$, set the turning points at $0,L,L+M$, where $L,M>0$, and write $b=\alpha_1$. The slopes are
$$
    (-1,-1/b,0,1).
$$
In particular, $b>1$. Set $h=f_2=f_3$ and $x=L/(bh)$, so that $f_1=h+L/b=h(1+x)$. The zero-slope generator gives $\alpha_2=|c|/h^2$. The jump at the first turning point, the jump across the zero-slope rod, and \eqref{eq:tod-ae-height-identity} give, respectively,
$$
    \alpha_2=(b-1)(1+x)^2,\qquad M=\frac{b+1}{2}h,\qquad L+M=2h+\frac Lb.
$$
Hence
\begin{equation}\label{eq:hermitian-four-angle-curve}
    \alpha_2=\frac{(b+1)^2}{4(b-1)},\qquad \frac ML=\frac{b^2-1}{b(3-b)},\qquad 1<b<3.
\end{equation}
Conversely, these formulas determine a positive convex function $f$ with the prescribed rods for each $1<b<3$ and $L>0$. The derivative of the length ratio is
$$
    \frac{\dd}{\dd b}\left(\frac{b^2-1}{b(3-b)}\right)=\frac{3b^2-2b+3}{b^2(3-b)^2}>0,
$$
and the ratio ranges from zero to infinity. Thus this family realizes every pair of positive finite rod lengths. The equal-angle equation $\alpha_2=b$ has the unique solution
$$
    b=1+\frac2{\sqrt3},\qquad \frac ML=2.
$$

For the last two configurations, the two infinite rod lines again have slope coefficients of absolute value one. Opposite coefficients on $A$ and $B$ would give zero coefficient on both occurrences of $C$, contrary to the strict increase of the slopes. The signed rod lists are therefore
$$
    (-B,-C,-A,C,A),\qquad (-B,-C,-A,B,C,A),
$$
respectively. At common finite cone angle $b$, their slopes are
$$
    \begin{aligned}
        BCACA:&\qquad (-1,-2/b,-1/b,2/b,1),\\
        BCABCA:&\qquad (-1,-2/b,-1/b,1/b,2/b,1).
    \end{aligned}
$$
In particular, $b>2$. A homothety allows us to set $|c|=1$. The determinant identity then gives the heights
$$
    \begin{aligned}
        BCACA:&\qquad (f_1,f_2,f_3,f_4)=\left(\frac1{\sqrt{b-2}},\frac1{\sqrt b},\frac1{\sqrt{3b}},\frac1{\sqrt{b-2}}\right),\\
        BCABCA:&\qquad (f_1,f_2,f_3,f_4,f_5)=\left(\frac1{\sqrt{b-2}},\frac1{\sqrt b},\frac1{\sqrt{2b}},\frac1{\sqrt b},\frac1{\sqrt{b-2}}\right).
    \end{aligned}
$$
Recover the lengths from $\ell_j=(f_{j+1}-f_j)/s_j$. The identity \eqref{eq:tod-ae-height-identity} becomes, respectively,
$$
    \begin{aligned}
        2\sqrt{b-2}+(1-\sqrt3)\sqrt b&=0,\\
        2\sqrt{b-2}+2(1-\sqrt2)\sqrt b&=0.
    \end{aligned}
$$
Their unique solutions for $b>2$ are $4/\sqrt3$ and $1+\sqrt2$. Substitution gives the positive length ratios in the theorem.

Conversely, these heights and lengths determine $z_i$ up to translation and $a_i=(s_i-s_{i-1})/2>0$. The identity \eqref{eq:tod-ae-height-identity} ensures that the resulting function is $f(z)=\sum_i a_i|z-z_i|$. The rod formulas then give the required primitive rod vectors and common cone angles. The reconstruction, including the normalization of $c$, is also verified in the nearby-angle calculation below. This proves existence. Every competitor has the ansatz already described, so the preceding calculations force its finite lengths to have the same ratios. After fixing scale and translation, the uniqueness theorem for harmonic maps with prescribed rod data and normalized end, \cite[Theorem~4.24]{LiSun}, identifies the metrics. This proves Theorem~\ref{thm:hermitian-ae}.

\subsubsection{Nearby angles}

We describe the nearby families for use in the next subsection. For $BCA$ the finite cone angle remains $4\pi$, and the only parameter is scale. For $BACA$, equation \eqref{eq:hermitian-four-angle-curve} describes the nearby cone angles, with scale providing the second parameter.

For $BCACA$, normalize $|c|=1$ and set
$$
    \begin{aligned}
        (s_0,s_1,s_2,s_3,s_4)&=\left(-1,-\frac2{\alpha_1},-\frac1{\alpha_2},\frac2{\alpha_3},1\right),\\
        (d_1,d_2,d_3,d_4)&=(\alpha_1-2,2\alpha_2-\alpha_1,2\alpha_2+\alpha_3,\alpha_3-2),\qquad f_i=d_i^{-1/2}.
    \end{aligned}
$$
The lengths are
\begin{equation}\label{eq:hermitian-five-lengths}
    \ell_1=\frac{\alpha_1}{2}(f_1-f_2),\qquad \ell_2=\alpha_2(f_2-f_3),\qquad \ell_3=\frac{\alpha_3}{2}(f_4-f_3),
\end{equation}
and the condition \eqref{eq:tod-ae-height-identity} is
\begin{equation}\label{eq:hermitian-five-angle}
    H_5(\alpha):=\sqrt{d_1}+\sqrt{d_2}-\sqrt{d_3}+\sqrt{d_4}=0.
\end{equation}

For $BCABCA$, set
$$
    \begin{aligned}
        (s_0,s_1,s_2,s_3,s_4,s_5)&=\left(-1,-\frac2{\alpha_1},-\frac1{\alpha_2},\frac1{\alpha_3},\frac2{\alpha_4},1\right),\\
        (d_1,d_2,d_3,d_4,d_5)&=(\alpha_1-2,2\alpha_2-\alpha_1,\alpha_2+\alpha_3,2\alpha_3-\alpha_4,\alpha_4-2),\qquad f_i=d_i^{-1/2}.
    \end{aligned}
$$
Then
\begin{equation}\label{eq:hermitian-six-lengths}
    \begin{aligned}
        \ell_1&=\frac{\alpha_1}{2}(f_1-f_2),&\qquad \ell_2&=\alpha_2(f_2-f_3),\\
        \ell_3&=\alpha_3(f_4-f_3),&\qquad \ell_4&=\frac{\alpha_4}{2}(f_5-f_4),
    \end{aligned}
\end{equation}
and the corresponding condition is
\begin{equation}\label{eq:hermitian-six-angle}
    H_6(\alpha):=\sqrt{d_1}+\sqrt{d_2}-2\sqrt{d_3}+\sqrt{d_4}+\sqrt{d_5}=0.
\end{equation}

All $d_i$ and $\ell_j$ are positive near the equal-angle solutions. Conversely, a nearby solution of $H_5=0$ or $H_6=0$ reconstructs $f$: set $z_1=0$, $z_{j+1}=z_j+\ell_j$, and $a_i=(s_i-s_{i-1})/2$. The slopes are strictly increasing, and \eqref{eq:tod-ae-height-identity} ensures that $f(z)=\sum_i a_i|z-z_i|$ has the prescribed heights.

For completeness, the normalization of the end imposes no additional equation. Write $t_j=\alpha_j s_j$, where $\alpha_0=\alpha_N=1$. These are the $\partial_y$ coefficients of the signed primitive rod vectors, namely $(-1,-2,-1,2,1)$ and $(-1,-2,-1,1,2,1)$. The determinant identity gives
$$
    \begin{aligned}
        \sum_{i=1}^N f_i^2\left(\frac1{s_i}-\frac1{s_{i-1}}\right)&=-\sum_{i=1}^N\frac1{t_{i-1}t_i}=-1,\\
        \sum_{i=1}^N f_i^2\left(\frac1{s_i}-\frac1{s_{i-1}}\right)&=f_1^2+f_N^2-\sum_{j=1}^q\ell_j(f_j+f_{j+1})=-2\sum_{i=1}^N a_i(z_i-\mu)^2.
    \end{aligned}
$$
Thus the variance formula for $c$ indeed gives $c=-1$. Moreover, each jump of $F$ is $1/(t_{i-1}t_i)$, which is exactly the jump prescribed by the consecutive signed primitive rod vectors. Starting from $F_0=-1$ therefore recovers every rod and gives $F_N=0$. A common rescaling of the lengths supplies the scale parameter.

Differentiating the defining functions gives
\begin{equation}\label{eq:hermitian-angle-gradient}
    \frac{\partial H_5}{\partial\alpha_j}=\frac{\ell_j}{\alpha_j}>0\quad(1\leq j\leq3),\qquad \frac{\partial H_6}{\partial\alpha_j}=\frac{\ell_j}{\alpha_j}>0\quad(1\leq j\leq4).
\end{equation}
For $BACA$, use
$$
    H_4(\alpha_1,\alpha_2)=\alpha_2-\frac{(\alpha_1+1)^2}{4(\alpha_1-1)}.
$$
Its partial derivatives satisfy
$$
    \frac{\partial H_4}{\partial\alpha_1}=\frac{(\alpha_1+1)(3-\alpha_1)}{4(\alpha_1-1)^2}>0,\qquad \frac{\partial H_4}{\partial\alpha_2}=1,
$$
for $1<\alpha_1<3$. Consequently, in each case the admissible cone angles form a smooth hypersurface near the equal-angle point, and the common-angle direction $(1,\ldots,1)$ is not tangent to it.

More explicitly, on this hypersurface the map
$$
    \alpha\longmapsto\left(\log\alpha_1-\log\alpha_q,\ldots,\log\alpha_{q-1}-\log\alpha_q\right)
$$
has invertible differential at the equal-angle point. Indeed, a tangent vector in its kernel has all components equal, and the positive partial derivatives of the defining function force this common component to vanish. For $q=1$ there are no angle-difference equations. Thus, including scale, each metric belongs to a smooth $q$-parameter family. Proposition~\ref{prop:ae-deformations} below proves that the finite Weyl lengths are local coordinates on this family; after fixing their total length, the same calculation gives a nonsingular Jacobian for the $q-1$ equal-angle equations.

\subsection{Infinitesimal deformations}

Let $(X,g)$ be one of the equal-angle metrics constructed above. There are $q$ finite rods, with lengths $\ell_j=z_{j+1}-z_j$ and cone angles $2\pi\alpha_j$. The two infinite rods have cone angle $2\pi$. We keep the torus action, the primitive rod vectors, and the normalized Euclidean end fixed.

In this section, we will prove that every decaying torus-invariant solution of the linearized Ricci-flat equation in the class specified below comes from varying the rod lengths in the families above, up to a change of coordinates. If all finite cone angles have the same first variation, that variation must be zero and only scaling remains.

We first specify the symmetry and regularity of the fixed-angle tensors, using the coordinates in Lemma \ref{lem:rodasymptotics}. Variations of the cone angles will be specified separately. Write $(N,\bg)=X/\T^2$ and let $\pi:X\to N$ be the quotient map. All coordinates below are chosen for the background metric $g$ and are kept fixed when describing tensor regularity. We consider $\T^2$-invariant symmetric two-tensors $k$; their coefficients in the quotient and angular coordinates are therefore independent of the angular variables.

Near an interior point of the $j$-th rod $\mI_j$, take the coordinates of Lemma \ref{lem:rodasymptotics}, denoted here by $(\xi,\eta,\phi,\psi)$, with $(\xi,\eta)$ the quotient coordinates and $\phi,\psi$ the angular coordinates of period $2\pi$. Here the rod is given by $\xi=0$, $\eta$ parametrizes it, and $\partial_\phi$ is the primitive angular generator that vanishes on the rod surface $\pi^{-1}(\mI_j)$. Thus $(\eta,\psi)$ are local coordinates along that surface. If the cone angle is $2\pi\alpha_j$, set
$$
    x=\xi\cos(\alpha_j\phi),\qquad y=\xi\sin(\alpha_j\phi).
$$
We require that $k$, when expressed in $(x,y,\eta,\psi)$, extend as a smooth symmetric two-tensor across $x=y=0$.

Near the torus fixed point over a turning point between rods $\mI_i$ and $\mI_j$, take the coordinates supplied by the corner description in Lemma \ref{lem:rodasymptotics}, denoted here by $(\xi_1,\xi_2,\phi_1,\phi_2)$. The quotient is locally the quadrant $\xi_1,\xi_2\geq0$, with
$$
    \mI_i=\{\xi_1=0\},\qquad \mI_j=\{\xi_2=0\},
$$
and the primitive angular generators $\partial_{\phi_1},\partial_{\phi_2}$ collapse on the corresponding rod surfaces. No distance interpretation of $\xi_1,\xi_2$ is required. If the two cone angles are $2\pi\alpha_i,2\pi\alpha_j$, set
$$
    \begin{aligned}
        x_1&=\xi_1\cos(\alpha_i\phi_1),& y_1&=\xi_1\sin(\alpha_i\phi_1),\\
        x_2&=\xi_2\cos(\alpha_j\phi_2),& y_2&=\xi_2\sin(\alpha_j\phi_2).
    \end{aligned}
$$
We require that $k$, when expressed in $(x_1,y_1,x_2,y_2)$, extend as a smooth symmetric two-tensor across both $\{x_1=y_1=0\}$ and $\{x_2=y_2=0\}$, including their intersection.

To allow changing cone angles, for each finite rod $\mI_j$ choose a smooth one-parameter family of torus-invariant conical metrics $g^{(j)}(t)$ with $g^{(j)}(0)=g$, keeping the torus action and rod vectors fixed. Such families can be constructed by modifying the Gram matrix. The family can be chosen to agree with $g$ outside a compact neighborhood of the finite rod surfaces. Its cone angle along $\mI_i$ is $2\pi\alpha_i^{(j)}(t)$, where
$$
    \alpha_i^{(j)}(t)=\alpha_i+t\delta_{ij},\qquad 1\leq i\leq q,
$$
while both infinite cone angles remain $2\pi$. Define the reference variations by
$$
    b_j=\left.\frac{\dd}{\dd t}\right|_{t=0}g^{(j)}(t).
$$
Thus $b_j$ is a compactly supported invariant symmetric two-tensor which changes the cone angle along $\mI_j$ at speed $2\pi$ and leaves every other cone angle unchanged.

A general variation considered here is
\begin{equation}\label{eq:ae-deformation-class}
    h=k+\sum_{j=1}^q c_jb_j,
\end{equation}
where $k$ is a torus-invariant fixed-angle tensor with the regularity just described and $c_j\in\mathbb R$. The variation of its cone angle along $\mI_j$ is $2\pi c_j$, equivalently $\dot\alpha_j(h)=c_j$, and its infinite cone angles have zero variation.

For the deformation statement below, $k$ is smooth in these coordinates and satisfies
\begin{equation}\label{eq:decay of the variation tensor}
    |\nabla_{g_{\mathrm{Euc}}}^m k|_{g_{\mathrm{Euc}}}=O(r_{\mathrm{Euc}}^{-1-m}),\qquad m\geq0,
\end{equation}
on the Euclidean end. Here $g_{\mathrm{Euc}}$ is the Euclidean metric and $r_{\mathrm{Euc}}$ is its distance to the origin; all end quantities are understood using a fixed equivariant Euclidean end chart. Since the $b_j$ have compact support, the same decay holds for $h$. The equation $(D\Ric)_g h=0$ is imposed on the complement of the rod surfaces. The coordinate changes we quotient out are generated by torus-invariant vector fields $Y$ which are smooth in the same coordinates, preserve the rod surfaces and torus fixed points, and satisfy
\begin{equation}\label{eq:decay of the gauge vector fields}
    |\nabla_{g_{\mathrm{Euc}}}^mY|_{g_{\mathrm{Euc}}}=O(r_{\mathrm{Euc}}^{-m}),\qquad m\geq0.
\end{equation}
In particular, $Y$ is bounded at infinity.

Scaling gives the family $e^{2t}g$, whose derivative at $t=0$ is $2g$. This tensor is smooth in the fixed-angle coordinates described above and preserves all cone angles, but it does not decay at infinity. To obtain a decaying variation representing scaling, choose a smooth torus-invariant vector field $Z_0$ which vanishes near all finite rod surfaces and is equal to $r_{\mathrm{Euc}}\partial_{r_{\mathrm{Euc}}}$ on the end, transported to $X$ by the fixed equivariant Euclidean end chart. Let $\Psi_t$ be the flow of $Z_0$, and define
\begin{equation}\label{eq:ae-scale-deformation}
    h_{\mathrm{sc}}=\left.\frac{\dd}{\dd t}\right|_{t=0}\left(e^{2t}\Psi_{-t}^*g\right)=2g-\mathcal L_{Z_0}g.
\end{equation}
The family $e^{2t}\Psi_{-t}^*g$ is torus-invariant and Ricci-flat, with the same cone angles as $g$. Consequently, $h_{\mathrm{sc}}$ is a smooth fixed-angle tensor satisfying $(D\Ric)_g h_{\mathrm{sc}}=0$. This tensor represents scaling and decays at infinity.

Let $E_j$ be the compact surface corresponding to the $j$-th finite rod. Its area and the total rod area are
\begin{equation}\label{eq:ae-total-area}
    A_j=\Area_g(E_j)=\frac{2\pi\ell_j}{\alpha_j},\qquad A_*=\sum_{j=1}^q A_j>0.
\end{equation}
The area formula follows from $\det G=\rho^2$. The induced area element on the rod is $\alpha_j^{-1}\dd z\dd\psi$, with $\psi$ of period $2\pi$. Scaling gives
$$
    \dot\ell_j(h_{\mathrm{sc}})=2\ell_j,\qquad \dot\alpha_j(h_{\mathrm{sc}})=0,\qquad D A_*(h_{\mathrm{sc}})=2A_*.
$$

We study tensors $h$ satisfying $(D\Ric)_gh=0$ and taking the form \eqref{eq:ae-deformation-class}, with the decay \eqref{eq:decay of the variation tensor}. Recall that the reduction of the toric Ricci-flat metric gives
$$
    g=\bg+\sum_{i,j=1}^2G_{ij}\dd\phi_i\dd\phi_j.
$$
We first prove that, up to a gauge term, we can remove the mixed components of $h$ between the quotient and angular directions.

\begin{lemma}
There is a torus-invariant vector field $Y'$, which is tangent to the torus orbits and satisfying \eqref{eq:decay of the gauge vector fields}, such that
$$
    h+\mathcal L_{Y'}g=\dot{\bg}+\sum_{i,j=1}^2\dot G_{ij}\dd\phi_i\dd\phi_j.
$$
\end{lemma}
\begin{proof}
Let $K_i=\partial_{\phi_i}$. On $N^\circ$, use the background coordinates $(\rho,z)$ and define one-forms by
$$
    a^i=\sum_{j=1}^2G^{ij}\left(h(\partial_\rho,K_j)\dd\rho+h(\partial_z,K_j)\dd z\right),\qquad i=1,2,
$$
where $(G^{ij})=G^{-1}$. These one-forms encode precisely the mixed components:
$$
    h=\dot{\bg}+2\sum_{i,j=1}^2G_{ij}a^j\dd\phi_i+\sum_{i,j=1}^2\dot G_{ij}\dd\phi_i\dd\phi_j.
$$
Here $\dot{\bg}$ and $\dot G$ denote the quotient and angular components of $h$. The reference variations $b_j$ preserve the block form, so the mixed components come entirely from the fixed-angle tensor $k$. Although $G^{-1}$ is singular at a rod, the regularity of $k$ supplies the vanishing needed in this formula. Indeed, in the adapted coordinates of Lemma \ref{lem:rodasymptotics}, smoothness and $\mathbb T^2$ invariance give
$$
    a^i=\xi p_i(\xi^2,\eta)\dd \xi+q_i(\xi^2,\eta)\dd \eta
$$
near an interior rod point, and
$$
    a^i=\xi_1p_{i1}(\xi_1^2,\xi_2^2)\dd \xi_1+\xi_2p_{i2}(\xi_1^2,\xi_2^2)\dd \xi_2
$$
near a turning point, with smooth coefficient functions. Thus these one-forms are smooth in the variables $(\xi^2,\eta)$ and $(\xi_1^2,\xi_2^2)$, respectively.

Linearizing the identities used in the proof of Lemma \ref{lem:polar}, the mixed components of $(D\Ric)_g h=0$ give
$$
    \dd\left(\rho\sum_{j=1}^2G_{ij}*_{\bg}\dd a^j\right)=0,\qquad i=1,2,
$$
where $*_{\bg}$ is the Hodge star on the quotient. Each expression inside the parentheses is therefore a constant function. The regularity just described and the vanishing of $\rho$ on the rods imply that these constants are zero. Since $\rho>0$ and $G$ is invertible in the interior, it follows that $\dd a^i=0$.

The quotient interior is simply connected, so there are functions $f^i$ satisfying $a^i=\dd f^i$. Set
$$
    Y'=-\sum_{i=1}^2f^iK_i.
$$
This vector field changes only the angular coordinates. Its Lie derivative cancels the mixed components:
$$
    h+\mathcal L_{Y'}g=\dot{\bg}+\sum_{i,j=1}^2\dot G_{ij}\dd\phi_i\dd\phi_j.
$$

We check that this is an allowed coordinate change. The local forms of $a^i$ above imply that their primitives are smooth in $(\xi^2,\eta)$ near a rod and in $(\xi_1^2,\xi_2^2)$ near a turning point. Consequently, $Y'$ extends smoothly in the fixed-angle Cartesian coordinates and preserves the rod surfaces and torus fixed points.

For the end estimates, use the fixed equivariant Euclidean end chart. Applying the same division argument on rescaled annuli, the differentiated decay of $h$ gives
$$
    a^i=O_\infty(r_{\mathrm{Euc}}^{-2}),
$$
uniformly up to the smooth infinite rods. The primitives have a common limit along the end; subtracting their constant limits gives $f^i=O_\infty(r_{\mathrm{Euc}}^{-1})$. Since the torus generators have size $O(r_{\mathrm{Euc}})$, we obtain
$$
    |\nabla_{g_{\mathrm{Euc}}}^mY'|_{g_{\mathrm{Euc}}}=O(r_{\mathrm{Euc}}^{-m}),\qquad m\geq0.
$$
In particular, $Y'$ is bounded.
\end{proof}

\begin{proposition}\label{prop:ae-deformations}
The tensor $h$ can be written as
$$
h=\sum_{j=1}^q\dot\ell_j\frac{\partial g_{\dot\bl}}{\partial\ell_j}+\mathcal L_Yg,
$$
where the coefficients $\dot\ell_j$ are uniquely determined and $Y$ is a bounded vector field of the kind specified above.

At the equal-angle metric with rod-length vector $\dot\bl=(\ell_1,\ldots,\ell_q)$,
\begin{equation}\label{eq:common-length-rank}
    \ker D_{\dot\bl}(\log\alpha_1,\ldots,\log\alpha_q)=\mathbb R\dot\bl.
\end{equation}
If $c_1=\cdots=c_q$ in \eqref{eq:ae-deformation-class}, then every $c_j$ is zero and
$$
    h=a h_{\mathrm{sc}}+\mathcal L_Yg
$$
for some $a\in\mathbb R$. In particular, fixing all cone angles leaves only scaling modulo coordinate changes. If also $D A_*(h)=0$, then $h$ is pure gauge.
\end{proposition}

\begin{proof}
By the preceding lemma, we may replace $h$ by $h+\mathcal L_{Y'}g$ and assume that it has no mixed components between the quotient and angular directions. Write
$$
    h=\dot{\bg}+\sum_{i,j=1}^2\dot G_{ij}\dd\phi_i\dd\phi_j.
$$
We first prove that if all rod-length variations vanish, then $h$ is pure gauge. We then use this fact to identify the general deformation with a tangent to the families constructed above.

Write $\underline{\ast}$ for the Hodge star of the quotient metric $\bg$, so that $\dd z=\underline{\ast}\dd\rho$, and let $\dot{\underline{\ast}}$ be its first variation induced by $\dot{\bg}$. The first variations of $\rho=\sqrt{\det G}$ and its harmonic conjugate satisfy
$$
    \begin{aligned}
        \dot\rho&=\frac{\rho}{2}\operatorname{tr}(G^{-1}\dot G),\\
        \dd\dot z&=\dot{\underline{\ast}}\dd\rho+\underline{\ast}\dd\dot\rho.
    \end{aligned}
$$
Indeed, linearizing \eqref{eq:rho} with $\lambda=0$ shows that the one-form in the second line is closed. Since $N^\circ$ is simply connected, $\dot z$ exists globally and is unique up to an additive constant. Its values $\dot z_i$ at the turning points determine the rod-length variations
$$
    \dot\ell_j=\dot z_{j+1}-\dot z_j.
$$
These differences are unchanged by the allowed coordinate changes.

Suppose now that every $\dot\ell_j$ is zero. All the $\dot z_i$ are then equal, so we choose the additive constant in $\dot z$ to make them zero. On $N^\circ$, define
$$
    V=-\dot\rho \partial_\rho-\dot z \partial_z,
$$
and lift $V$ horizontally to $X$. Replacing $h$ by $h+\mathcal L_Vg$ makes the first variations of both $\rho$ and $z$ zero. We check that $V$ is an allowed vector field, including along the rod surfaces and at the torus fixed points.

Near a turning point, use the quotient coordinates $(\xi_1,\xi_2)$ of Lemma \ref{lem:rodasymptotics}. By reflecting four times, we can glue four copies of a small disc along its two rods around the turning point, and we obtain a reflected disc, where $\underline{g}$ extends as a smooth metric. The original disk in the quotient is a quadrant of the reflected disc. The reference variations $b_j$ change only the Gram matrix, so $\dot{\bg}$ is the quotient component of the fixed-angle tensor $k$. Its smoothness and $\mathbb T^2$ invariance give
$$
    \dot{\bg}=\dot a \dd\xi_1^2+2\xi_1\xi_2\dot b \dd\xi_1\dd\xi_2+\dot c \dd\xi_2^2,
$$
where $\dot a,\dot b,\dot c$ are smooth functions of $(\xi_1^2,\xi_2^2)$. This formula extends $\dot{\bg}$ smoothly to either sign of $\xi_1,\xi_2$ and makes it invariant under both coordinate reflections. Together with the corresponding extension of $\bg$, it gives a smooth metric and a smooth symmetric two-tensor on the full reflected disc.

Choose a fixed isothermal coordinate $\zeta=x+\sqrt{-1}y$ on the reflected disc, centered at the turning point and compatible with both reflections, such that $\bg=e^{2\sigma}|\dd\zeta|^2$. Set $w=z+\sqrt{-1}\rho$ and $\dot w=\dot z+\sqrt{-1}\dot\rho$, so that the equations defining $V$ become $\dd w(V)=-\dot w$. Choose the complex orientation of $\zeta$ so that $w$ is holomorphic. The map $w$ takes the original quadrant to the right half-plane, and the corner expansion gives
$$
    w=z_i+a_i\zeta^2+O(\zeta^4),\qquad \partial_\zeta w=2a_i\zeta+O(\zeta^3),\qquad a_i\ne0.
$$
Thus $\zeta$ is a regular coordinate at the corner, whereas $w$ has vanishing differential there. The corner expansions of $G$ and $\dot G$, together with $\rho=\sqrt{\det G}$, give
$$
    \dot\rho=\xi_1\xi_2 f(\xi_1^2,\xi_2^2)=O(|\zeta|^2)
$$
for a smooth function $f$. In particular, $\dd\rho$ and $\dd\dot\rho$ are $O(|\zeta|)$. Differentiating $\dd z=\underline{\ast}\dd\rho$ gives
$$
    \dd\dot z=\dot{\underline{\ast}}\dd\rho+\underline{\ast}\dd\dot\rho.
$$
The right-hand side extends smoothly to the reflected disc and is $O(|\zeta|)$. Since $\dot z_i=0$, it follows that $\dot z=O(|\zeta|^2)$, and hence $\dot w=O(|\zeta|^2)$. In the fixed coordinate $\zeta$, the differentiated conjugacy equation reads
$$
    \partial_{\bar\zeta}\dot w=\mu \partial_\zeta w,\qquad \mu=\frac{\dot{\bg}_{xx}-\dot{\bg}_{yy}+2\sqrt{-1}\dot{\bg}_{xy}}{4e^{2\sigma}},
$$
where $\mu$ is smooth on the reflected disc. Since $w$ is holomorphic and $\dd w(V)=-\dot w$, the complex component of $V$ satisfies
$$
    v:=\dd\zeta(V)=-\frac{\dot w}{\partial_\zeta w}=O(|\zeta|),\qquad \partial_{\bar\zeta}v=-\mu
$$
away from the corner. Extending $v(0)=0$, the latter equation holds distributionally on the full disc: integration by parts around a circle of radius $\varepsilon$ produces a boundary term of order $O(\varepsilon^2)$. We now apply interior elliptic regularity to $\partial_{\bar\zeta}v=-\mu$. Thus $v$ is smooth across the corner.

The reflection symmetries then imply
$$
    V=\xi_1v_1(\xi_1^2,\xi_2^2)\partial_{\xi_1} +\xi_2v_2(\xi_1^2,\xi_2^2)\partial_{\xi_2},
$$
with smooth $v_1,v_2$. Its lift is therefore smooth in the fixed-angle Cartesian coordinates and preserves both rod surfaces and the torus fixed point. At an interior rod point, $\partial_\zeta w\ne0$, and the same conclusion follows directly from one reflection.

On the end, the decay of $h$ gives $\dot\rho,\dot z=O_\infty(r_{\mathrm{Euc}})$. The Weyl-coordinate differential has size comparable to $r_{\mathrm{Euc}}$, so
$$
    |\nabla_{g_{\mathrm{Euc}}}^mV|_{g_{\mathrm{Euc}}}=O(r_{\mathrm{Euc}}^{-m}),\qquad m\geq0.
$$
Thus the torus-invariant lift of $V$ is an allowed vector field. Replacing $h$ by $h+\mathcal L_Vg$, we may henceforth assume $\dot\rho=\dot z=0$.

\ 

Linearizing $\dd z=\underline{\ast}\dd\rho$ now gives $\dot{\underline{\ast}}\dd\rho=0$. Since $\underline{\ast}^2=-1$ on one-forms, we also have $\dot{\underline{\ast}}\dd z=0$. Thus the Hodge star is fixed to first order, which in two dimensions means that the quotient-metric variation is conformal. Write
$$
    \bg=e^{2\nu}(\dd\rho^2+\dd z^2),\qquad \dot{\bg}=2\dot\nu \bg.
$$
With the coordinates $(\rho,z)$ fixed, linearizing \eqref{eq:axisymmetric harmonic map equation} shows that $Z=\dot\Phi$, where $\Phi=\rho^{-1}G$, is a Jacobi field along the axisymmetric harmonic map $\Phi$. The target $\mathcal H=SL(2,\mathbb R)/SO(2)$ has negative curvature. Consequently,
\begin{equation}\label{eq:ae-jacobi-bochner}
    \Delta_{\mathbb R^3}|Z|^2=2|DZ|^2-2\sum_{a=1}^3\left\langle\mathcal R(Z,\dd\Phi(e_a))\dd\Phi(e_a),Z\right\rangle\geq0.
\end{equation}
Here the norms and inner product use the target metric, $D$ is the pullback connection, $\mathcal R$ is the target curvature, and $e_a$ is an orthonormal frame for the Euclidean $\mathbb R^3$. Its axisymmetric Laplacian is $\partial_\rho^2+\rho^{-1}\partial_\rho+\partial_z^2$.

To apply the maximum principle, we check that $|Z|$ is bounded near the auxiliary axis and tends to zero at infinity. At a turning point, rescale the rows and columns of $G$ and $\dot G$ by $\xi_1^{-1},\xi_2^{-1}$ in the adjacent primitive rod basis. The prescribed regularity gives smooth bounded matrices after this rescaling, and the rescaled background matrix is positive definite. Hence the relative Gram variation, and therefore $|Z|$, is bounded. At an interior rod point, the same argument rescales only the row and column corresponding to the collapsing circle. On the end, $|Z|=O(r_{\mathrm{Euc}}^{-1})$ by the decay of $h$.

The axis has locally zero Newtonian capacity in $\mathbb R^3$, so the bounded subharmonic function $|Z|^2$ extends across it. Since it tends to zero at infinity, the maximum principle gives $Z=0$. We already have $\dot\rho=0$, so $\dot G=0$ as well. The remaining reduced equations for $\nu$ then give $\dd\dot\nu=0$, and the end normalization forces this constant to be zero. Thus the tensor vanishes after the allowed coordinate changes. This proves that every deformation with zero rod-length variations is pure gauge.

\ 

We next show that the rod lengths are smooth local coordinates on the families constructed in Sections \ref{subsec:app AE with two rods}-\ref{app:hermitian-ae}. For the alternating families, they are already independent parameters. For a Hermitian family, the parameters $\alpha_j$ determining the finite cone angles lie on the hypersurface described in the preceding subsection. Together with scale, they give $q$ parameters. If a tangent to this family has zero rod-length variations, the result just proved makes it pure gauge. It therefore has zero variation of every cone angle, so its component tangent to the angle hypersurface is zero. The remaining scale component changes the length vector by a nonzero multiple of $\dot\bl$, unless that component is zero. Thus the length differential is injective between two $q$-dimensional spaces, and the inverse function theorem makes the rod lengths local coordinates.

For a general solution $h$, choose the unique family tangent with the same rod-length variations and subtract it. The difference has zero rod-length variations and is therefore pure gauge. This gives
$$
    h=\sum_{j=1}^q\dot\ell_j\frac{\partial g_{\dot\bl}}{\partial\ell_j}+\mathcal L_Yg.
$$
The coefficients are uniquely determined because the allowed coordinate changes do not change the rod-length variations.

\ 

It remains to consider the cone angles. In an alternating family, split a length variation into a multiple of $\dot\bl$ and a tangent to the simplex of fixed total rod length. The differential of the angle differences on that simplex is nonsingular by Lemma \ref{lem:alternating-hessian} and the calculation following it. If all finite cone angles have the same first variation, the angle differences have zero variation, so the tangent component vanishes. Only scale remains, and scale preserves all cone angles.

For a Hermitian family with $q\geq2$, differentiating its angle relation gives
$$
    \sum_{j=1}^q\frac{\partial H}{\partial\alpha_j}c_j=0,
$$
where $H=H_4,H_5$ or $H_6$ and $c_j=\dot\alpha_j(h)$. All the partial derivatives are positive by \eqref{eq:hermitian-angle-gradient} and the formula for $H_4$. If the $c_j$ have a common value, that value must be zero. The tangent component along the angle hypersurface then vanishes, leaving only scale. For $BCA$, the single finite cone angle is identically $4\pi$, and the family has only the scale parameter. This proves \eqref{eq:common-length-rank} and gives, in the common-angle case,
$$
    h=a h_{\mathrm{sc}}+\mathcal L_Yg.
$$
Finally, the allowed coordinate changes preserve the areas of the compact rod surfaces, while scaling gives $D A_*(h_{\mathrm{sc}})=2A_*$. Therefore
$$
    D A_*(a h_{\mathrm{sc}}+\mathcal L_Yg)=2aA_*.
$$
Since $A_*>0$, the condition $D A_*(h)=0$ forces $a=0$, and $h$ is pure gauge.
\end{proof}
\subsection{Angle variation and invertibility at fixed area}
\label{subsec:app angle variation and invertibility}

We consider the equal-angle metric $(X,g)$ and keep the fixed-angle tensor regularity from the preceding subsection. Using the notation of Section~\ref{sec:moduli}, with Einstein constant zero, set
$$
    \mathscr E_0(g)=\Ric_g-\frac12\Scal_g g,\qquad \mathscr P_g k=(D\mathscr E_0)_g k+\delta_g^*\delta_g k,
$$
where $\delta_g=-\operatorname{div}_g$ and $\delta_g^*$ is its formal adjoint. The domain of $\mathscr P_g$ consists of torus-invariant fixed-angle symmetric two-tensors. All components are allowed, including the mixed components between quotient and angular directions.

Let $b_j$ be the reference variations chosen in the preceding subsection, and set
$$
    b=\sum_{j=1}^q b_j,\qquad B=(D\mathscr E_0)_g b.
$$
Thus $\dot\alpha_j(b)=1$ for every finite rod, so $b$ changes each finite cone angle at speed $2\pi$. It has compact support and leaves both infinite cone angles unchanged.

Although $b$ itself need not be smooth in the fixed-angle Cartesian coordinates, $B$ is smooth. We explain this in the following. Let $g_t$ be the reference family inducing $b$, with $g_0=g$ and $\dot g_0=b$. The quotient metric is fixed, while the Gram matrix varies.

Near an interior point of a finite rod, set $\alpha_j(t)=\alpha_j+t$ and consider, on an angular sector, the map
$$
    F_t(\xi,\eta,\phi,\psi) =\bigl(\xi\cos(\alpha_j(t)\phi),\xi\sin(\alpha_j(t)\phi),\eta,\psi\bigr).
$$
The regularity of the reference family means that
$$
    g_t=F_t^*\widetilde g_t,
$$
where $\widetilde g_t$ is a smooth family of smooth metrics on a fixed Cartesian domain with coordinates $(x,y,\eta,\psi)$, including $x=y=0$. These Cartesian coordinates are held fixed when differentiating $\widetilde g_t$. At $t=0$, the map $F_0$ is exactly the map used to define fixed-angle regularity for the background metric. Set $S_t=\mathscr E_0(\widetilde g_t)$. This is a smooth family of smooth tensors on the Cartesian domain. Naturality of the Einstein tensor gives
$$
    \mathscr E_0(g_t)=F_t^*S_t.
$$
Since $g$ is Ricci-flat, $S_0=0$. Differentiating on the regular locus therefore gives
$$
    \begin{aligned}
        B =\left.\frac{\dd}{\dd t}\right|_{t=0}\bigl(F_t^*S_t\bigr)=F_0^*\dot S_0 +\left.\frac{\dd}{\dd t}\right|_{t=0}\bigl(F_t^*S_0\bigr) =F_0^*\dot S_0.
    \end{aligned}
$$
The term coming from the changing coordinates vanishes because $F_t^*S_0$ is identically zero. Thus the coefficients of $B$ in the background Cartesian coordinates are precisely those of the smooth tensor $\dot S_0$. This proves smooth extension across the rod.

At a torus fixed point, apply the same argument using both angle-adapted Cartesian coordinate pairs from the preceding subsection. The reference metrics are jointly smooth in $t$ and these four Cartesian variables, so the same calculation proves smoothness of $B$ across both rod surfaces and their intersection. Hence $B$ is a fixed-angle smooth tensor. It has compact support because $b$ does, and differentiating the Bianchi identity at the Ricci-flat background gives $\delta_g B=0$.

\begin{lemma}\label{lem:pair}
Choose the vector field $Z_0$ in \eqref{eq:ae-scale-deformation} to vanish on a neighborhood of $\operatorname{supp}b$. Then
$$
    h_{\mathrm{sc}}=2g-\mathcal L_{Z_0}g=O_\infty(r_{\mathrm{Euc}}^{-2})
$$
satisfies $(D\mathscr E_0)_g h_{\mathrm{sc}}=0$, $D A_*(h_{\mathrm{sc}})=2A_*$, and
\begin{equation}\label{eq:ae-angle-scale-pairing}
    \int_X\langle B,h_{\mathrm{sc}}\rangle_g \dvol_g=-8\pi A_*\ne0.
\end{equation}
\end{lemma}

\begin{proof}
The equation and the area variation follow from the scaling construction in the preceding subsection. Since $Z_0$ agrees with the Euclidean dilation field on the end and $g-g_{\mathrm{Euc}}=O_\infty(r_{\mathrm{Euc}}^{-2})$, its leading Euclidean contribution cancels that of $2g$, giving the stated decay.

We compute the pairing by integrating the scalar-curvature variation. Near an interior point of $E_j$, use the coordinates $(\xi,\eta,\phi,\psi)$ from the preceding subsection. Here $\xi$ is normal distance to the rod, and the leading normal metric is
$$
    \dd\xi^2+\alpha_j^2\xi^2\dd\phi^2.
$$
The reference variation keeps the quotient metric fixed and satisfies $\dot\alpha_j(b)=1$. Differentiating the angular term therefore gives the leading normal part of $b$:
$$
    2\alpha_j\xi^2 \dd\phi\otimes\dd\phi.
$$
At a Ricci-flat metric, the scalar-curvature variation is
$$
    (D\Scal)_g b=\operatorname{div}_g(J^\sharp),\qquad J=\operatorname{div}_g b-\dd\operatorname{tr}_g b.
$$
Here $J$ is a one-form, and $J^\sharp$ is its metric dual. The rod expansions give
$$
    \begin{aligned}
        (\operatorname{div}_g b)(\partial_\xi) =-\frac{2}{\alpha_j\xi}+O(1),\qquad \dd(\operatorname{tr}_g b)(\partial_\xi)=O(1).
    \end{aligned}
$$
Indeed, for the leading normal metric, $\Gamma^\phi_{\phi\xi}=\xi^{-1}$ and the angular component of $b$ with one index raised is $2/\alpha_j$. These give the term $-2/(\alpha_j\xi)$ in the divergence. The trace of the leading normal variation is the constant $2/\alpha_j$, and the remaining terms in the rod expansions contribute only bounded errors.

Delete tubes of radius $\varepsilon$ about the finite rod surfaces. Along the boundary near an interior rod point, the outward unit normal of the remaining region is $-\partial_\xi$. Hence
$$
    J(-\partial_\xi)=\frac{2}{\alpha_j\varepsilon}+O(1).
$$
After integration in $\phi$, whose period is $2\pi$, the boundary measure is
$$
    \bigl(2\pi\alpha_j\varepsilon+O(\varepsilon^2)\bigr)\dd A_g.
$$
The limiting flux along $E_j$ is therefore
$$
    \int_{E_j}\frac{2}{\alpha_j\varepsilon} \bigl(2\pi\alpha_j\varepsilon\bigr)\dd A_g =4\pi A_j.
$$
The boundary pieces obtained by rounding the tubes near the torus fixed points contribute $O(\varepsilon^2)$: their measure is $O(\varepsilon^3)$ and the current is $O(\varepsilon^{-1})$. There is no contribution at infinity because $b$ has compact support. Applying the divergence theorem and letting $\varepsilon$ tend to zero gives
$$
    \int_X(D\Scal)_g b\dvol_g=4\pi A_*.
$$

In dimension four, at the Ricci-flat background,
$$
    \operatorname{tr}_g\bigl((D\mathscr E_0)_g b\bigr)=-(D\Scal)_g b.
$$
Moreover, $\operatorname{supp}B\subset\operatorname{supp}b$, and $Z_0$ vanishes on a neighborhood of this support. Thus $h_{\mathrm{sc}}=2g$ on $\operatorname{supp}B$, and consequently
$$
    \begin{aligned}
        \int_X\langle B,h_{\mathrm{sc}}\rangle_g\dvol_g &=2\int_X\operatorname{tr}_g B\dvol_g\\
        &=-2\int_X(D\Scal)_g b\dvol_g =-8\pi A_*.
    \end{aligned}
$$
\end{proof}

We now specify the weighted spaces. Choose a smooth positive torus-invariant function $\omega$ on $X$ which equals $1+r_{\mathrm{Euc}}$ on the end. For an integer $m\geq0$ and a weight $\tau\in\mathbb R$, set
\begin{equation}\label{eq:ae-sobolev-norm}
    \|k\|_{H^m_\tau}^2=\sum_{j=0}^m\int_X\omega^{-2\tau-4+2j}|\nabla_g^j k|_g^2 \dvol_g.
\end{equation}
The space $H^m_\tau$ is the completion in this norm of compactly supported torus-invariant tensors with the fixed-angle Cartesian regularity specified above. Write $L^2_\tau=H^0_\tau$. On the compact part, the corresponding H\"older spaces are defined using the same Cartesian coordinates and an invariant partition of unity. On the end, the $C^{m,\gamma}_\tau$ norm uses $\omega^{j-\tau}|\nabla_g^j k|_g$ and the scaled H\"older seminorm on dyadic annuli, where $0<\gamma<1$. We use the same conventions for functions and one-forms. All spaces below are $\mathbb T^2$ invariant.

\begin{proposition}\label{prop:ae-augmented}
The operator
$$
    \mathscr P_g:H^2_{-1}\longrightarrow L^2_{-3}
$$
is Fredholm of index zero. Its kernel and its cokernel, the latter represented using the unweighted integral pairing, are generated by a same divergence-free, trace-free tensor $s_T$ representing scaling variation. It satisfies $D A_*(s_T)=2A_*$. For any $c_0\in\mathbb R$, the map
\begin{equation}\label{eq:ae-augmented-isomorphism}
    (k,t)\longmapsto\bigl(\mathscr P_g k+tB,\ D A_*(k)+c_0t\bigr)
\end{equation}
is an isomorphism from $H^2_{-1}\oplus\mathbb R$ to $L^2_{-3}\oplus\mathbb R$. It is also an isomorphism from $C^{2,\gamma}_{-1}\oplus\mathbb R$ to $C^{0,\gamma}_{-3}\oplus\mathbb R$ for $0<\gamma<1$.
\end{proposition}

In the rest part of this subsection, we prove this proposition.

We first establish the local and weighted estimates. In the fixed-angle Cartesian coordinates from the preceding subsection, the metric and the operator have smooth coefficients, and the tensors under consideration are invariant under rotations in the corresponding coordinate pairs. Ordinary interior elliptic estimates therefore apply across a rod surface and across a torus fixed point.

Let $\mathscr P_0$ be the operator for $g_{\mathrm{Euc}}$. The standard fundamental solution construction for constant-coefficient elliptic operators gives a kernel $K$, smooth away from the origin, satisfying
$$
    |\partial^jK(x)|\leq C_j|x|^{-2-j},\qquad x\ne0,\quad j\geq0.
$$
For every smooth compactly supported symmetric two-tensor $f$ on $\mathbb R^4$, the convolution $k=K*f$ satisfies $\mathscr P_0k=f$. We use this kernel to establish the weighted estimates below. The kernel bound, together with rescaling on annuli and interior elliptic estimates, gives
$$
    \|K*f\|_{H^2_\tau(\mathbb R^4)} \leq C_\tau\|f\|_{L^2_{\tau-2}(\mathbb R^4)}, \qquad -2<\tau<0,
$$
for smooth compactly supported $f$. By density, convolution with $K$ therefore defines a bounded right inverse in these weighted spaces. The homogeneous equation has no nonzero solution in $H^2_\tau$: applying scaled interior estimates on expanding Euclidean balls and using $\tau<0$ forces such a solution to vanish. Consequently,
$$
    \mathscr P_0:H^2_\tau(\mathbb R^4)\longrightarrow L^2_{\tau-2}(\mathbb R^4)
$$
is an isomorphism for $-2<\tau<0$.

Since $g-g_{\mathrm{Euc}}=O_\infty(r_{\mathrm{Euc}}^{-2})$, the coefficients of $\mathscr P_g-\mathscr P_0$ in orders two, one, and zero are respectively $O(r_{\mathrm{Euc}}^{-2})$, $O(r_{\mathrm{Euc}}^{-3})$, and $O(r_{\mathrm{Euc}}^{-4})$. We may therefore absorb these terms outside a sufficiently large compact set. Combining the resulting end estimate with the local estimates gives
\begin{equation}\label{eq:ae-fredholm-estimate}
    \|k\|_{H^2_\tau}\leq C\bigl(\|\mathscr P_g k\|_{L^2_{\tau-2}}+\|k\|_{L^2(K_0)}\bigr),\qquad -2<\tau<0,
\end{equation}
where $K_0\subset X$ is compact. Rellich compactness implies that the kernel is finite dimensional and the range is closed.

At a Ricci-flat metric, $(D\mathscr E_0)_g$ is, up to an overall sign, the Hessian of the Einstein-Hilbert functional. Hence $\mathscr P_g$ is formally self-adjoint. The Cartesian regularity makes the inner boundary terms in Green's formula vanish. Under the unweighted integral pairing, the dual of $L^2_{\tau-2}$ is $L^2_{-\tau-2}$. If a tensor in this dual space annihilates the range, it satisfies the adjoint equation across the rod surfaces and torus fixed points as well as in the regular locus. Indeed, in a Cartesian chart, average an arbitrary smooth test tensor over the rotations. Its pairing with an invariant tensor is unchanged, and the averaged test tensor belongs to the specified domain. Distributional elliptic regularity on the whole chart therefore applies.  End estimates then place the adjoint solution in $H^2_{-\tau-2}$. This proves the Fredholm property. At $\tau=-1$, the complementary weight is again $-1$, so formal self-adjointness gives index zero.

We will also need a decay improvement. If $\mathscr P_g k=0$ and $k\in H^2_\tau$ for some $-2<\tau<0$, scaled elliptic estimates first give $k=O_\infty(r_{\mathrm{Euc}}^\tau)$. Multiply $k$ by a cutoff supported on the end and extend the resulting tensor by zero to $\mathbb R^4$. Its flat equation has a right-hand side which is compactly supported plus $O_\infty(r_{\mathrm{Euc}}^{\tau-4})$. This right-hand side is integrable because $\tau<0$. Convolution with $K$ gives $O(r_{\mathrm{Euc}}^{-2})$; the difference from the cutoff tensor is a decaying entire flat solution and is therefore zero. Differentiated estimates give
$$
    k=O_\infty(r_{\mathrm{Euc}}^{-2}).
$$

\ 

To set the scaling variation in divergence gauge, consider
$$
    Q=2\delta_g\delta_g^*
$$
on invariant one-forms with the fixed-angle regularity specified above. Since $g$ is Ricci-flat,
$$
    Q\vartheta=\nabla_g^*\nabla_g\vartheta+\dd(\delta_g\vartheta).
$$
This is a formally self-adjoint elliptic operator. Its Euclidean fundamental kernel satisfies the same bounds used above. The preceding weighted estimates and decay argument therefore apply to $Q$: the map
$$
    Q:H^2_\tau\longrightarrow L^2_{\tau-2}
$$
is Fredholm for $-2<\tau<0$, and every element of its kernel is $O_\infty(r_{\mathrm{Euc}}^{-2})$. If $Q\vartheta=0$, integration by parts gives
$$
    0=\int_X\langle Q\vartheta,\vartheta\rangle_g\dvol_g =2\int_X|\delta_g^*\vartheta|_g^2\dvol_g.
$$
The boundary terms vanish by the fixed-angle regularity at the rods and the decay at infinity. Hence $\vartheta^\sharp$ is a decaying Killing field. Ricci-flatness and the Bochner maximum principle imply that $\vartheta=0$. By formal self-adjointness, the cokernel is represented by the kernel of $Q$ at the complementary weight $-\tau-2$. This weight also lies in $(-2,0)$, so the same argument shows that the cokernel is zero. Consequently,
$$
    Q:H^2_\tau\longrightarrow L^2_{\tau-2}
$$
is an isomorphism for every $-2<\tau<0$.

Fix $0<\epsilon<1$ and solve
$$
    Q\vartheta=-\delta_g h_{\mathrm{sc}}
$$
at weight $\tau=-1+\epsilon$. This small shift avoids the borderline weight: $\delta_g h_{\mathrm{sc}}=O_\infty(r_{\mathrm{Euc}}^{-3})$ belongs to $L^2_{-3+\epsilon}$. Define
$$
    s_T=h_{\mathrm{sc}}+2\delta_g^*\vartheta=h_{\mathrm{sc}}+\mathcal L_{\vartheta^\sharp}g.
$$
The vector field $\vartheta^\sharp$ is invariant, smooth in the fixed-angle coordinates, and tangent to the rod surfaces; it vanishes at the torus fixed points. By construction and diffeomorphism invariance,
$$
    \delta_g s_T=0,\qquad (D\Ric)_g s_T=0.
$$
Initially $s_T=O_\infty(r_{\mathrm{Euc}}^{-2+\epsilon})$, so $s_T\in H^2_{-1}$. Its trace is a decaying harmonic function and hence is zero by the maximum principle, applied across the rods using the Cartesian coordinates. It follows that $\mathscr P_gs_T=0$, and the decay improvement gives $s_T=O_\infty(r_{\mathrm{Euc}}^{-2})$. Since the gauge field preserves every rod surface, it does not change its area. Thus $D A_*(s_T)=2A_*$.

Now suppose $\mathscr P_g k=0$ in $H^2_{-1}$. The linearized Bianchi identity gives
$$
    Q(\delta_g k)=0.
$$
The decay improvement for $k$ gives $\delta_g k=O_\infty(r_{\mathrm{Euc}}^{-3})$, so the kernel calculation for $Q$ implies $\delta_g k=0$. Taking the trace of $\mathscr P_g k=0$ then shows that $\operatorname{tr}_g k$ is harmonic. It decays, so it vanishes. Consequently $(D\Ric)_g k=0$, and Proposition~\ref{prop:ae-deformations} applies. Since $k$ is a fixed-angle variation, it is a multiple of $h_{\mathrm{sc}}$ plus an allowed bounded gauge term.

Subtracting the corresponding multiple of $s_T$, we obtain a tensor $\mathcal L_Yg$ with zero divergence and trace. The vector field $Y$ has the regularity and differentiated bounds specified in \eqref{eq:decay of the gauge vector fields}. The divergence equation and Ricci-flatness give $\nabla_g^*\nabla_gY=0$. Comparing this equation with the Euclidean operator on the end gives
$$
    \Delta_{\mathbb R^4}Y=O_\infty(r_{\mathrm{Euc}}^{-4}).
$$
The Newton kernel estimate gives
$$
    Y=v_0+O_\infty(r_{\mathrm{Euc}}^{-2}\log r_{\mathrm{Euc}})
$$
for a constant vector $v_0$. The standard torus action on the Euclidean end has no nonzero fixed vector, so invariance forces $v_0=0$. The Bochner identity and the maximum principle now give $Y=0$. Hence the kernel and cokernel of $\mathscr P_g$ are both represented by $\mathbb R s_T$. The same conclusion holds at every weight $-2<\tau<0$, because every homogeneous solution in these weights has the decay just proved.

Since $B$ is divergence-free and compactly supported, integration by parts gives
$$
    \int_X\langle B,2\delta_g^*\vartheta\rangle_g \dvol_g=0.
$$
Lemma~\ref{lem:pair} therefore yields
$$
    \int_X\langle B,s_T\rangle_g \dvol_g=-8\pi A_*\ne0.
$$
Given a first component $f$, pairing $\mathscr P_g k+tB=f$ with $s_T$ determines $t$ uniquely. The Fredholm alternative then gives $k$ modulo $\mathbb R s_T$. The second equation in \eqref{eq:ae-augmented-isomorphism} determines the remaining coefficient because $D A_*(s_T)=2A_*>0$. The area functional is bounded on $H^2$: in the Cartesian charts, the codimension-two trace theorem gives a bounded restriction to each compact rod surface. Thus the bounded inverse theorem proves the Sobolev assertion, for every $c_0\in\mathbb R$.

Finally, let $f\in C^{0,\gamma}_{-3}$. Then it belongs to $L^2_{-3+\epsilon}$ for every $\epsilon>0$. Solve the augmented problem first at weight $-1+\epsilon$, with $0<\epsilon<1$, using the same kernel and cokernel calculation. By local regularity and scaled estimates we get $k=O(r_{\mathrm{Euc}}^{-1+\epsilon})$ with the corresponding derivative bounds through order two. On the end, its flat equation has right-hand side
$$
    O(r_{\mathrm{Euc}}^{-3})+O(r_{\mathrm{Euc}}^{-5+\epsilon}).
$$
Convolution with $K$ improves $k$ to $O(r_{\mathrm{Euc}}^{-1})$. By scaled Schauder estimates we then obtain $k\in C^{2,\gamma}_{-1}$ and the bounded inverse estimate in these spaces. Uniqueness follows from the kernel and area calculation above. By the same argument we also obtain the corresponding higher regularity statements.

\section{Smooth compactness for arbitrary rod structure}
\label{app:fixed-rod-smooth-compactness}

We fix a smooth compact rod structure $\mR$ and denote by $\mathcal M_{\mR}$ the space of smooth toric Einstein metrics $g$ on $X_{\mR}$ with $\Ric(g)=3g$, modulo orientation-preserving diffeomorphisms commuting with the fixed torus action and preserving each rod sphere. As in Section~\ref{sec:moduli}, we equip $\mathcal M_{\mR}$ with the quotient topology induced by the $C^\infty$ topology on metrics. The purpose of this appendix is to establish the following compactness of $\mathcal M_{\mR}$, extending Lemma \ref{lem:three-direction-smooth-compactness} to the case of arbitrary rod structure.

\begin{theorem}\label{thm:fixed-rod-smooth-compactness}
The space $\mathcal M_{\mR}$ is compact.
\end{theorem}

\begin{corollary}
Given $N>0$, the moduli space of simply-connected compact toric Einstein four-manifolds $(X,g)$ with $\Ric_g=3g$ satisfying $\chi(X)\leq N$ is compact.
\end{corollary}

\begin{proof}
By Theorem~\ref{thm:fixed-rod-smooth-compactness}, it suffices to show that only finitely many equivalence classes of rod structures occur. The number of rods is $m=\chi(X)\leq N$, and Proposition~\ref{lem:rod-self-intersection-bound1} implies that every rod self-intersection belongs to $\{-1,0,1\}$.

We choose representative vectors $\mv_1,\ldots,\mv_m$ so that $(\mv_i,\mv_{i+1})$ is a positively oriented lattice basis for $1\leq i<m$. After an orientation-preserving lattice change, we may fix $\mv_1=(2\pi,0)$ and $\mv_2=(0,2\pi)$. By the self-intersection formula in Section~\ref{sec:rod structure}, we have
$$
    \mv_{i+1}=e_i\mv_i-\mv_{i-1}, \qquad 2\leq i\leq m-1.
$$
Thus the integers $e_2,\ldots,e_{m-1}\in\{-1,0,1\}$ determine all remaining rod vectors, so there are only finitely many possibilities. When $m=2$, the first two vectors already determine the rod structure. The conclusion then follows from Theorem~\ref{thm:fixed-rod-smooth-compactness}.
\end{proof}

The rest of this appendix is devoted to the proof of Theorem \ref{thm:fixed-rod-smooth-compactness}. By the geometric bounds in Section~\ref{sec:end structure}, it suffices to exclude curvature concentration along a sequence in $\mathcal M_{\mR}$. We will use the convergence and gauge estimates of \cite{OzuchI,OzuchII} and the ALE dilation obstruction of \cite{OzuchIntegrability}, and the new input in our setting is the curvature identities and the self-intersection bounds.

We first heuristically explain this obstruction for a smooth gluing family, before justifying the calculation for a general degenerating sequence.

We consider a smooth toric Ricci-flat ALE metric $(Y,h)$ attached at a fixed point $p$ of an Einstein orbifold $(M,g)$, where $\Ric_g=\lambda g$. We choose coordinates near $p$ and at infinity on $Y$ that identify the two torus actions with the same linear action on $\R^4/\Gamma$, and we may write
$$
    h=g_{\mathrm{Euc}}+H^4+O_\infty(R^{-5}), \qquad g=g_{\mathrm{Euc}}+H_2+O_\infty(r^3),
$$
where $H^4$ and $H_2$ have components homogeneous of degrees $-4$ and $2$, respectively.

We suppose that a family of smooth Einstein metrics $g_t$ is obtained by gluing $(Y,th)$ into $(M,g)$ near $p$, with $\Ric_{g_t}=\lambda_tg_t$ and $\lambda_t\to\lambda$ as $t\rightarrow 0$. Then we obtain metrics $\gamma_t$ on expanding subsets of $Y$ satisfying
$$
    \Ric_{\gamma_t}=t\lambda_t\gamma_t.
$$
We assume the first-order expansion $ \gamma_t=h+tk+o(t) $ in $C^\infty$ on every compact subset of $Y$. The first-order matching condition requires
$$
    k=H_2+o_\infty(R^2) \qquad\text{as }R\to\infty.
$$
Differentiating the Einstein equation at $t=0$, we obtain
\begin{equation}\label{eqn:Ricci deformation}
    (D\Ric)_h k=\lambda h.
\end{equation}
The relevant infinitesimal deformation is generated by  a change of scale of $h$. Up to sign, this deformation is represented by the tensor $\zeta_h$ constructed below in \eqref{eq:fixed-rod-dilation-tensor} which, in suitably normalized coordinates,  satisfies  that
$$
    (D\Ric)_h \zeta_h=0, \qquad \delta_h \zeta_h=0, \qquad \tr_h \zeta_h=0, \qquad \zeta_h=-4H^4+O_\infty(R^{-5}).
$$
We contract the equation  \eqref{eqn:Ricci deformation} with $\zeta_h$ and integrate over the region $Y_R$ bounded by the coordinate sphere of radius $R$. Integrating by parts, we obtain
\begin{equation}\label{eq:fixed-rod-heuristic-green}
    \int_{\partial Y_R} \left[ \frac12\bigl(\langle k,\nabla_\nu \zeta_h\rangle_h -\langle\nabla_\nu k,\zeta_h\rangle_h\bigr) -\zeta_h(\nu,B_hk) \right] \dd A_h=0,
\end{equation}
where $\nu$ is the outward unit normal and $B_hk=\delta_hk+\frac12\dd\tr_hk$ is the Bianchi operator. Letting $R$ tend to infinity, we obtain an identity
$$
    4\int_{\mathbb S^3/\Gamma} \left(3\langle H_2,H^4\rangle_{g_{\mathrm{Euc}}} +H^4(\partial_ R,B_{g_{\mathrm{Euc}}}H_2)\right) \dd A_{g_{\mathrm{Euc}}}=0.
$$
This is the boundary obstruction computed in \cite[proof of Theorem~5.6, equation~(56)]{OzuchIntegrability}, following the earlier calculation in \cite[Propositions~2.1--2.2]{BiquardI}. For the limits arising in our compactness argument, we will show that the toric curvature identities and the  self-intersection bounds force this integral to be strictly positive, so a contradiction.

Below we will make the above argument rigorous. One difficulty, as explained in \cite{OzuchIntegrability}, is that a priori we do not know if the bubble $(Y, h)$ has integrable deformations.  We avoid this assumption in our setting by allowing the finite cone angles to vary and constructing a real-analytic family of conical Ricci-flat metrics with prescribed nearby rod lengths.

We consider a sequence $[g_i]\in\mathcal M_{\mR}$.  Notice that we have uniform volume lower bound, diameter bound, and curvature energy bound from Section~\ref{sec:end structure}. Then passing to a subsequence we obtain Gromov-Hausdorff convergence to an Einstein orbifold, and a bubble tree  structure \cite{Bando}. Each bubble is an orbifold toric $\ALE$ gravitational instanton.

Suppose now that the curvature of $g_i$ is not uniformly bounded, then in the deepest bubble we obtain a non-flat smooth toric $\ALE$ gravitational instanton $(Y,h)$. In the next scale, we obtain an non-flat Einstein orbifold $(M, g)$, compact or complete. We denote the corresponding orbifold point by $p$ and write
$$
    \Ric_g=\lambda g, \qquad \lambda=3\ \text{or}\ 0.
$$
By Proposition~\ref{lem:rod-self-intersection-bound1}, every finite rod of $Y$ has self-intersection in $\{-1,0,1\}$.

\subsection{Positivity of the curvature pairing}
\label{subsec:fixed-rod-positive-pairing}

We use the convention of Section~\ref{subsec:toric-weyl-rigidity} and for a compact toric Einstein metric with $\Ric_g=\lambda g$ by Proposition~\ref{prop:weyl-positive} we have
$$
    w_\pm=\frac{\lambda}{3}-K\mp J \geq0.
$$
Clearly these inequalities also pass to the (rescaled) limits and across the orbifold points. By the boundary curvature calculation leading to \eqref{eq:weyl-boundary}, applied in both orientations, we can choose orthonormal bases of $\Lambda_p^\pm$ so that
\begin{equation}\label{eq:fixed-rod-curvature-at-vertex}
    \Rm_g^\pm(p)=\frac{\lambda}{3}\Id+w_\pm(p)D_0, \qquad D_0=\operatorname{diag}(2,-1,-1).
\end{equation}
We first consider a deepest bubble, which is a smooth nonflat toric Ricci-flat ALE gravitational instanton $(Y,h)$ with end $\R^4/\Gamma$. We choose equivariant coordinates at infinity in the volume gauge of \cite[Definition~4.12]{OzuchIntegrability}. In these coordinates, we have
\begin{equation}\label{eq:fixed-rod-ale-expansion}
    h=g_{\mathrm{Euc}}+H^4+O_\infty(r^{-5}),
\end{equation}
where the Cartesian components of $H^4$ are homogeneous of degree $-4$. Denote the triple of parallel two-forms on the Euclidean $\R^4$ by $\omega_a^\pm$, with $\omega_1^\pm=\dd x_1\wedge\dd x_2\pm\dd x_3\wedge\dd x_4$. We define
$$
    \theta_a^\mp =\omega_a^\pm-2\dd r\wedge\iota_{\partial_r}\omega_a^\pm,
$$
where $\iota_{\partial_r}$ denotes contraction with the Euclidean radial unit vector. By \cite[Proposition~4.13]{OzuchIntegrability}, expressed in the bases chosen above, there are constant real symmetric $3\times3$ matrices $S^+$ and $S^-$ such that
\begin{equation}\label{eq:fixed-rod-ale-matrices}
    H^4=-r^{-4}\left(\sum_{a,b}S^+_{ab}\theta_a^-\circ\omega_b^+ +\sum_{a,b}S^-_{ab}\theta_a^+\circ\omega_b^-\right).
\end{equation}
Here we identify a 2-form as a skew symmetric endomorphism and a symmetric 2-tensor as a symmetric endomorphism, and recall the isomorphism
$$
    \Lambda^-\otimes\Lambda^+ \simeq \operatorname{Sym}^2_0(T^*\mathbb R^4), \qquad \alpha^-\otimes\beta^+\longmapsto\alpha^-\circ\beta^+=\beta^+\circ\alpha^-.
$$
Notice that $(S^+,S^-)$ are uniquely determined up to addition of $(c\Id, -c\Id)$, so the number
$$
    A_h=\tr S^++\tr S^-
$$
is well defined. In particular, we have
\begin{equation}\label{eq:fixed-rod-radial-coefficient}
    H^4(\partial_r,\partial_r)=-A_hr^{-4}.
\end{equation}

We now construct the  infinitesimal deformation of $h$ corresponding to rescaling, which we will  use to test the linearized Einstein equation. By \cite[Proposition~5.1 and Lemma~5.4]{OzuchIntegrability}, there is a unique invariant function $f$ on $Y$ satisfying
$$
    \Delta_hf=8, \qquad f=r^2+A_hr^{-2}+o_\infty(r^{-2}), \qquad \nabla_hf=2r\partial_r+o_\infty(r^{-3}).
$$
Then we define
\begin{equation}\label{eq:fixed-rod-dilation-tensor}
    \zeta_h=\mathcal L_{\frac12\nabla_hf}h-2h=\nabla_h^2f-2h.
\end{equation}
It follows from the Ricci-flat equation and the volume-gauge expansion that
\begin{equation}\label{eq:fixed-rod-exact-dilation}
    (D\Ric)_h \zeta_h=0, \qquad \tr_h \zeta_h=0, \qquad \delta_h \zeta_h=0, \qquad \zeta_h=-4H^4+O_\infty(r^{-5}).
\end{equation}

\begin{lemma}\label{lem:fixed-rod-positive-trace}
We have $ A_h>0.$

\end{lemma}

\begin{proof}
Since $\tr_h \zeta_h=0$ and $\delta_h \zeta_h=0$,  integrating by parts we obtain
$$
    \int_Y|\zeta_h|_h^2 \dvol_h =\lim_{R\to\infty}\int_{\{r=R\}}\zeta_h(\nabla_hf,\nu) \dd A_h =8\Vol(\mathbb S^3/\Gamma)A_h.
$$
Thus $A_h\geq0$. If $A_h=0$, then $\nabla_h^2f=2h$. In particular, $h$ is a metric cone hence is flat. Contradiction.
\end{proof}
Using the above chosen bases, we can regard the symmetric matrices $S^\pm$ as self-adjoint endomorphisms of $\Lambda^\pm$ by setting
$$
    S^\pm\omega_b^\pm=\sum_a S^\pm_{ab}\omega_a^\pm.
$$
It follows from the $\mathbb T^2$ invariance that $S^\pm$ commute with the rotations fixing $\omega_1^\pm$. So we can write
\begin{equation}\label{eq:fixed-rod-ale-diagonal}
    S^\pm=\frac{\tr S^\pm}{3}\Id+a_\pm D_0.
\end{equation}

\begin{lemma}\label{lem:fixed-rod-strict-ale-coefficients}
We have  $a_\pm>0$.
\end{lemma}

\begin{proof}
By \cite[Proposition~4.13]{OzuchIntegrability}, the coefficients $S^\pm_{ab}$ defined in the basis $\{\omega_a^\pm\}$ determine the leading curvature term in the reflected basis $\{\theta_a^\mp\}$ by
$$
    W_h^\mp(\theta_b^\mp) =24r^{-6}\sum_{a=1}^3 \left(S^\pm_{ab}-\frac{\tr S^\pm}{3}\delta_{ab}\right)\theta_a^\mp +O_\infty(r^{-7}).
$$
Thus \eqref{eq:fixed-rod-ale-diagonal} implies that the leading curvature eigenvalue is $48a_\pm r^{-6}$ along $\theta_1^\mp$ and $-24a_\pm r^{-6}$ on its orthogonal complement.

We now identify the coefficient in the expansion of $w_\mp$. By the curvature calculation in Section~\ref{subsec:toric-weyl-rigidity}, the sum or difference of the horizontal and vertical area forms is an eigenform of $W_h^\mp$ with eigenvalue $-w_\mp$. At leading order, $\theta_1^\mp$ is a sum of two-forms with one horizontal and one vertical factor, so it is orthogonal to this eigenform. Consequently, we obtain
\begin{equation}\label{eq:fixed-rod-w-infinity}
    w_\mp=24a_\pm r^{-6}+O_\infty(r^{-7}).
\end{equation}
By \eqref{eq:instanton-w-equation}, we have
$$
    \Delta_{\underline h}w_\pm+3\langle\dd\log\rho,\dd w_\pm\rangle_{\underline h} =-3w_\pm^2-4|q_\pm|^2.
$$
Multiplying by $\rho^3$, we obtain the divergence identity
$$
    \diver_{\underline h}\bigl(\rho^3\nabla_{\underline h}w_\pm\bigr) =-\rho^3\bigl(3w_\pm^2+4|q_\pm|^2\bigr).
$$
We integrate this identity over the quotient $Y/\mathbb T^2$ truncated at $r=R$, with small neighborhoods of the rods and turning points removed. The boundary integrals along the rods vanish in the limit because $\rho$ vanishes linearly and $w_\pm$ are smooth and even in the normal variable. At a turning point, $\rho=O(r^2)$ and the curvature is smooth on the local cover, so the boundary integral over a shrinking circular arc also vanishes. On the end, we use polar coordinates $r,\theta$, with $0\leq\theta\leq\pi/2$, and write
$$
    \rho=\kappa r^2\sin\theta\cos\theta+O_\infty(r^{-2}), \qquad \underline h=\dd r^2+r^2\dd\theta^2+O_\infty(r^{-4}),
$$
where $\kappa>0$ is a constant. By \eqref{eq:fixed-rod-w-infinity}, so we obtain
\begin{equation}\label{eq:fixed-rod-strict-coefficient-identity}
    12\kappa^3a_\mp =\int_{Y/\mathbb T^2}\rho^3\bigl(3w_\pm^2+4|q_\pm|^2\bigr) \dd A_{\underline h}.
\end{equation}
In particular, $a_\mp\geq0$.

To show the strict inequality holds, we notice that if $a_\mp=0$, then  $w_\pm=q_\pm=0$, hence $W_h^\pm=0$. It follows easily from the proof of Lemma~\ref{lem:instanton-rod-self-intersections}  that every finite rod sphere has self-intersection $\mp2$. This contradicts our self-intersection bound $\{-1,0,1\}$.
\end{proof}

We now turn the orbifold fixed point $p$ of the limit $(M,g)$ at the next larger scale. We have $\Ric_g=\lambda g$ for $\lambda=3$ or $0$. Using the identification of the asymptotic cone of $h$ with $T_pM/\Gamma$, in equivariant coordinates on the local cover at $p$ we can write
$$
    g=g_{\mathrm{Euc}}+H_2+O_\infty(r^3)
$$
where the coefficients of $H_2$ are homogeneous of degree two. For a symmetric two-tensor $\xi$ whose coefficients are homogeneous of degree two, we define the linear functional
$$
    \mathcal B(\xi) =4\int_{\mathbb S^3/\Gamma} \left( 3\langle\xi,H^4\rangle_{g_{\mathrm{Euc}}} +H^4(\partial_r,B_{g_{\mathrm{Euc}}}\xi) \right) \dd A_{g_{\mathrm{Euc}}},
$$
where $B_{g_{\mathrm{Euc}}}$ is the Bianchi operator.

\begin{lemma}\label{lem:fixed-rod-dilation-pairing}
We have
\begin{equation}\label{eq:fixed-rod-dilation-pairing}
    \mathcal B(H_2) =4\Vol(\mathbb S^3/\Gamma) \mathfrak d(h,g,p),
\end{equation}
where
\begin{equation}\label{eq:fixed-rod-curvature-pairing}
    \mathfrak d(h,g,p) =\tr\bigl(S^+\Rm_g^+(p)\bigr) +\tr\bigl(S^-\Rm_g^-(p)\bigr).
\end{equation}
\end{lemma}

\begin{proof}
We first prove that $\mathcal B(\mathcal L_Vg_{\mathrm{Euc}})=0$ for every homogeneous cubic equivariant vector field $V$. Since $H^4$ is divergence-free, we have
$$
    \operatorname{div}_{g_{\mathrm{Euc}}} \bigl(H^4(V,\cdot)^\sharp\bigr) =\frac12\langle\mathcal L_Vg_{\mathrm{Euc}},H^4\rangle_{g_{\mathrm{Euc}}}.
$$
Integrating by parts over $\{\varepsilon<r<1\}$ and let $\varepsilon\to0$ we obtain
$$
    \int_{\mathbb S^3/\Gamma} \langle\mathcal L_Vg_{\mathrm{Euc}},H^4\rangle_{g_{\mathrm{Euc}}} \dd A_{g_{\mathrm{Euc}}} =4\int_{\mathbb S^3/\Gamma} H^4(\partial_r,V) \dd A_{g_{\mathrm{Euc}}}.
$$
Since $H^4$ is harmonic and divergence-free, the vector-valued function $x\mapsto H^4_x(x,\cdot)^\sharp$ is harmonic and homogeneous of degree $-3$. Thus on $\mathbb S^3$ we have
$$
    \begin{aligned} &\Delta_{\mathbb S^3}\bigl(H^4(\partial_r,\cdot)^\sharp\bigr) =-3H^4(\partial_r,\cdot)^\sharp,\\
                    &\Delta_{g_{\mathrm{Euc}}}V =\Delta_{\mathbb S^3}V+15V. \end{aligned}
$$
Integration by parts therefore implies that
$$
    \int_{\mathbb S^3/\Gamma} H^4(\partial_r,\Delta_{g_{\mathrm{Euc}}}V)\dd A_{g_{\mathrm{Euc}}} =12\int_{\mathbb S^3/\Gamma} H^4(\partial_r,V)\dd A_{g_{\mathrm{Euc}}}.
$$
Now noticing that
$$
    \bigl(B_{g_{\mathrm{Euc}}}(\mathcal L_Vg_{\mathrm{Euc}})\bigr)^\sharp =-\Delta_{g_{\mathrm{Euc}}}V.
$$
So we have
$$
    \mathcal B(\mathcal L_Vg_{\mathrm{Euc}}) =4\int_{\mathbb S^3/\Gamma} \left( 3\langle\mathcal L_Vg_{\mathrm{Euc}},H^4\rangle_{g_{\mathrm{Euc}}} -H^4(\partial_r,\Delta_{g_{\mathrm{Euc}}}V) \right)\dd A_{g_{\mathrm{Euc}}} =0. 
$$

Now we choose geodesic normal coordinates at $p$. Since the Christoffel symbols vanish at $p$ in both coordinate systems, the coordinate transition has an equivariant homogeneous cubic leading correction $V$. The quadratic term therefore changes by $\mathcal L_Vg_{\mathrm{Euc}}$, whose boundary integral vanishes by the above calculation. Thus we may compute $\mathcal B(H_2)$ in normal coordinates, where
$$
    (H_2)_{ij}=-\frac13R_{iajb}(p)x_ax_b.
$$
On the unit sphere, we then have
$$
    H_2(\partial_r,\cdot)=0,\qquad \tr_{g_{\mathrm{Euc}}}H_2=-\frac{\lambda}{3},\qquad B_{g_{\mathrm{Euc}}}H_2=-\frac{2\lambda}{3}\dd r.
$$
Since $\mathcal B$ is linear in $H^4$, we compute separately the contributions of the two sums in \eqref{eq:fixed-rod-ale-matrices}. We first consider the sum involving $S^+$; the contribution of $S^-$ follows from the same calculation with the opposite orientation. We set $s=\tr S^+$ and $u_a=\omega_a^+\partial_r$. On the unit sphere, the summand involving $S^+$ in \eqref{eq:fixed-rod-ale-matrices} is given by
$$
    -\sum_{a,b}S^+_{ab}\theta_a^-\circ\omega_b^+ = s g_{\mathrm{Euc}} -2s \dd r\otimes\dd r -2\sum_{a,b}S^+_{ab}u_a^\flat\otimes u_b^\flat.
$$
The contribution to $\mathcal B(H_2)$ is then
$$
    8\int_{\mathbb S^3/\Gamma} \sum_{a,b}S^+_{ab}R(u_a,\partial_r,u_b,\partial_r) \dd A_{g_{\mathrm{Euc}}} -\frac{4\lambda s}{3}\Vol(\mathbb S^3/\Gamma).
$$
By the curvature decomposition and averaging over the unite sphere we get
\begin{equation}\label{eqn:spherical average}
    \frac{1}{\Vol(\mathbb S^3)} \int_{\mathbb S^3} R(u_a,\partial_r,u_b,\partial_r) \dd A_{g_{\mathrm{Euc}}} =\frac12\left((\Rm_g^+(p))_{ab}+\frac{\lambda}{3}\delta_{ab}\right).
\end{equation}
Indeed, we have
$$
    u_a^\flat\wedge\dd r=\frac12(\theta_a^--\omega_a^+),
$$
so the projections to $\Lambda^+$ are independent of the radial direction. Up to fixed orthogonal changes of basis, the reflected frame varies by the adjoint action of the unit quaternions. Thus spherical averaging agrees with averaging over rotations, which replaces $\Rm_g^-(p)$ by $\frac13\tr(\Rm_g^-(p))\Id=\frac{\lambda}{3}\Id$. Since  $|\omega_a^+|^2=|\theta_a^-|^2=2$, we obtain \eqref{eqn:spherical average}.

It follows that the contribution to $\mathcal B$ by $S^+$ is given by
$$
    4\Vol(\mathbb S^3/\Gamma)\tr\bigl(S^+\Rm_g^+(p)\bigr).
$$
The same calculation in the opposite orientation yields the contribution $4\Vol(\mathbb S^3/\Gamma)\tr(S^-\Rm_g^-(p))$, so summing the two terms proves \eqref{eq:fixed-rod-dilation-pairing}.
\end{proof}

\begin{proposition}\label{prop:positivity of pairing}
We have $\mathcal B(H_2)>0$.
\end{proposition}
\begin{proof}
By \eqref{eq:fixed-rod-curvature-at-vertex} and \eqref{eq:fixed-rod-ale-diagonal}, we have
\begin{equation}\label{eq:fixed-rod-positive-pairing}
    \mathfrak d(h,g,p) =\frac{\lambda}{3}A_h+6a_+w_+(p)+6a_-w_-(p).
\end{equation}
If $\lambda=3$, Lemma~\ref{lem:fixed-rod-positive-trace} implies that $\mathfrak d(h,g,p)>0$. If $\lambda=0$, the metric $g$ is nonflat, so $W_g^+\not\equiv0$ or $W_g^-\not\equiv0$. By the strong maximum principle in the proof of Lemma~\ref{lem:instanton-rod-self-intersections}, applied on the local orbifold covers, the corresponding function $w_+$ or $w_-$ is positive at $p$. Lemma~\ref{lem:fixed-rod-strict-ale-coefficients} therefore implies that $\mathfrak d(h,g,p)>0$ in this case as well.
\end{proof}
\subsection{Conical deformations of the $\ALE$ metric}
\label{subsec:fixed-rod-angle-family}

We fix a smooth  $\ALE$ toric gravitational instanton $(Y,h)$ with end $\mathbb R^4/\Gamma$ and finite rod spheres $E_1,\ldots,E_q$. By Lemma~\ref{geom:proper-ricciflat} we have natural coordinates $(\rho, z)$ on $Y/\mathbb T^2$. If $z_1<\cdots<z_{q+1}$ are the turning points, the finite rod lengths are $\ell_j=z_{j+1}-z_j$. We write $\bl^*=(\ell_1, \cdots, \ell_q)$ for the length vector of $h$ and $\mathbf1=(1,\ldots,1)$. When we vary the rod lengths $\bl=(\ell_1, \cdots, \ell_q)$ by \cite{LiSun}, one can still uniquely solve a toric conical $\ALE$ gravitational instantons with cone angles $2\pi\beta_j$ along $E_j$. In \cite{LiSun} this is obtained by solving the axisymmetric harmonic map. For our purpose here we need real analytic dependence on $\bl$, so we will proceed with an implicit function theorem argument. Notice that, for a toric conical gravitational instanton, the rod length is related to the area of rod spheres by (see \eqref{eq:ae-total-area})
$$
    \beta_j\Area(E_j)=2\pi\ell_j.
$$

We will use the angle identifications $J_{\boldsymbol{\beta}}$ from Section~\ref{sec:moduli}, with one independent angle for each finite rod and $J_{\mathbf1}=\Id$. We can choose fixed charts and cutoffs so that $J_{\boldsymbol{\beta}}$ depends real-analytically on $\boldsymbol{\beta}$ and equals the identity outside a compact neighborhood of $E_1\cup\cdots\cup E_q$. For a family $h(\bl)$ with angle vector $\boldsymbol{\beta}(\bl)$, we will  measure the dependence through the map
$$
    \bl\longmapsto \bigl(J_{\boldsymbol{\beta}(\bl)}^{-1}h(\bl)-h,\boldsymbol{\beta}(\bl)\bigr) \in C^{2,\alpha}_{-1}\times(0,\infty)^q.
$$

We will use the notations in Section~\ref{sec:moduli}:
$$
    \mathscr E_0(h)=\Ric_h-\frac12\Scal_h h, \qquad \mathscr P_h=(D\mathscr E_0)_h+\delta_h^*\delta_h, \qquad \delta_h=-\operatorname{div}_h.
$$

\begin{lemma}\label{lem:fixed-rod-angle-family}
On a neighborhood $U$ of $\bl^*$ in $(0,\infty)^q$, there is a real-analytic family of toric conical gravitational instantons $h(\bl)$ with angle vector $\boldsymbol{\beta}(\bl)$ such that
$$
    h(\bl^*)=h,\qquad \boldsymbol{\beta}(\bl^*)=\mathbf1,\qquad \beta_j(\bl)\Area_{h(\bl)}(E_j)=2\pi\ell_j \quad (1\leq j\leq q).
$$
\end{lemma}

\begin{proof}
For $k\in C^{2,\alpha}_{-1}$ sufficiently small, we set $g(k,\boldsymbol{\beta})=J_{\boldsymbol{\beta}}(h+k)$ and consider the map
\begin{equation}\label{eq:fixed-rod-analytic-map}
    (k,\boldsymbol{\beta})\longmapsto \left( J_{\boldsymbol{\beta}}^{-1}\left\{ \mathscr E_0(g(k,\boldsymbol{\beta})) -\delta_{g(k,\boldsymbol{\beta})}^*\delta_{g(k,\boldsymbol{\beta})}(J_{\boldsymbol{\beta}} h) \right\}, \left( \frac{\beta_j\Area_{g(k,\boldsymbol{\beta})}(E_j)}{2\pi} \right)_{j=1}^q \right).
\end{equation}
Its coefficients in the fixed angular developments are analytic functions of the metric, its first two derivatives, and $\boldsymbol{\beta}$. Since $J_{\boldsymbol{\beta}}$ is the identity on the end and the area functions are analytic, the map \eqref{eq:fixed-rod-analytic-map} is real-analytic with values in $C^{0,\alpha}_{-3}\times\mathbb R^q$.

We set
$$
    b_j=\left.\partial_{\beta_j}(J_{\boldsymbol{\beta}} h)\right|_{\boldsymbol{\beta}=\mathbf1}, \qquad B_j=(D\mathscr E_0)_h b_j.
$$
By the calculation preceding Lemma~\ref{lem:pair}, the tensors $B_j$ are smooth, compactly supported, and divergence-free. The derivative of \eqref{eq:fixed-rod-analytic-map} at $(0,\mathbf1)$ is
\begin{equation}\label{eq:fixed-rod-all-angle-derivative}
    (k,c)\longmapsto \left( \mathscr P_hk+\sum_{j=1}^q c_jB_j, \left( \frac{D\Area_h(E_j)[k+\sum_{a=1}^q c_ab_a] +c_j\Area_h(E_j)}{2\pi} \right)_{j=1}^q \right).
\end{equation}
The Fredholm argument in the proof of Proposition~\ref{prop:ae-augmented} uses only the ALE end estimates and local regularity, so it applies to $h$ and implies that $\mathscr P_h:C^{2,\alpha}_{-1}\to C^{0,\alpha}_{-3}$ has index zero. Adding $q$ angle variables and $q$ length equations preserves we still have index $0$.

For a pair $(k,c)$ in the kernel, the Bianchi identity and the decaying vector-field estimate in the proof of Proposition~\ref{prop:ae-augmented} imply that $\delta_hk=0$. Consequently,
$$
    \dot g=k+\sum_{j=1}^q c_jb_j
$$
satisfies $(D\Ric)_h\dot g=0$ and has zero variation of every rod length. The argument  in the proof of Proposition~\ref{prop:ae-deformations} also applies to our situation, and we see that $\dot g=\mathcal L_Vh$ for a bounded invariant vector field $V$ preserving the rods. Since such a variation preserves every cone angle, we have $c=0$ and $k=\mathcal L_Vh$. The equation $\delta_hk=0$ then implies that $\delta_h\delta_h^*V^\flat=0$. As in the proof of in Proposition~\ref{prop:ae-augmented}, the $\ALE $ expansion and torus invariance imply that $V$ decays, and integration by parts implies that $k=0$. Therefore \eqref{eq:fixed-rod-all-angle-derivative} is injective and hence invertible.

By the analytic implicit-function theorem, we can solve \eqref{eq:fixed-rod-analytic-map} with prescribed value $(0,\bl)$. The gauge term can be eliminated by the Bianchi identity and integration by parts, as in Proposition~\ref{prop:localmoduli}. Then the conclusion follows.
\end{proof}

Now the coordinate construction of \cite[Section~4.2.1]{OzuchIntegrability} can be applied uniformly to the family $h(\bl)$. Since the conical singularities lie in a fixed compact set, they do not affect the choice of volume gauge on the end.  After shrinking $U$, we choose equivariant end coordinates, followed by the radial adjustment defining the volume gauge, such that
\begin{equation}\label{eq:fixed-rod-family-end}
    h(\bl)=g_{\mathrm{Euc}}+H^4(\bl)+O_\infty(r^{-5}),
\end{equation}
where the coefficients of $H^4(\bl)$ are homogeneous of degree $-4$. We fix the normalization in the exterior gauge equations, so their solutions depend smoothly on $\bl$. Differentiating these equations and applying the exterior elliptic estimates, we obtain remainder estimates uniform on a smaller neighborhood of $\bl^*$ and after any fixed number of parameter derivatives.

The construction of $f$ and $\zeta_h$ in \eqref{eq:fixed-rod-dilation-tensor} also extends to the family $h(\bl)$, using the volume normalization in \cite[proof of Lemma~5.4]{OzuchIntegrability}. After shrinking $U$, the resulting tensors $\zeta_{h(\bl)}$ satisfy \eqref{eq:fixed-rod-exact-dilation} with remainder estimates uniform in $\bl$.

We now write $A_j(\bl)=\Area_{h(\bl)}(E_j)$ and consider
$$
    \mathcal S(\bl) =4\pi\sum_{j=1}^q (1-\beta_j(\bl))A_j(\bl).
$$
This has a natural interpretation as the Einstein-Hilbert functional-- the  metric $h(\bl)$ is Ricci-flat away from the rod spheres, and the curvature contribution from the conical singularities yields the above sum.

\begin{lemma}\label{lem:fixed-rod-angle-cancellation}
After shrinking $U$, there are constants $C>0$ and $0<\vartheta\leq1$ such that
\begin{equation}\label{eq:fixed-rod-angle-cancellation}
    \left|\sum_{j=1}^q A_j(\bl)(1-\beta_j(\bl))\right| \leq C|\boldsymbol{\beta}(\bl)-\mathbf1|^{1+\vartheta}.
\end{equation}
\end{lemma}

\begin{proof}
We consider a variation of $\bl$. Since $h(\bl)$ is Ricci-flat, we have $(D\Scal)_{h(\bl)}(\dot h)=0$. Integrating this and using the scalar-curvature flux calculation in the proof of Lemma~\ref{lem:pair}, we obtain
\begin{equation}
    \sum_{j=1}^q A_j \dd\beta_j=0.
\end{equation}
So we have
\begin{equation}\label{eq:fixed-rod-angle-potential}
    \dd\mathcal S =4\pi\sum_{j=1}^q(1-\beta_j) \dd A_j.
\end{equation}
The functions $A_j$ and $\beta_j$ are real analytic, so $\mathcal S$ is real analytic. Since $\boldsymbol{\beta}(\bl^*)=\mathbf1$, we have
$$
    \mathcal S(\bl^*)=0, \qquad \dd\mathcal S(\bl^*)=0.
$$
After shrinking $U$, the derivatives of $A_j$ are uniformly bounded, and hence
$$
    |\nabla_{\bl}\mathcal S(\bl)| \leq C|\boldsymbol{\beta}(\bl)-\mathbf1|,
$$
where the gradient is taken with respect to the Euclidean metric in the length coordinates. By the \L{}ojasiewicz gradient inequality at $\bl^*$, we obtain
$$
    |\mathcal S(\bl)| \leq C|\nabla_{\bl}\mathcal S(\bl)|^{1+\vartheta} \leq C|\boldsymbol{\beta}(\bl)-\mathbf1|^{1+\vartheta}
$$
for some $0<\vartheta\leq1$, which proves \eqref{eq:fixed-rod-angle-cancellation}.
\end{proof}

\subsection{Comparison with the conical ALE family}

We consider a smooth nonflat toric Ricci-flat ALE space $(Y,h)$ attached at an orbifold point $p$ of an Einstein orbifold $(M,g)$, where $M$ may be compact or ALE. We write $g_i$ for the metrics on $X_\mR$ at the scale of $M$ and assume that
$$
    \Ric(g_i)=\lambda_i g_i,\qquad \lambda_i\longrightarrow\lambda,\qquad g_i\longrightarrow g\quad\text{smoothly away from }p.
$$
We denote the relative squared scale by $t_i\to0$, so that $t_i^{-1}g_i\to h$ smoothly on every compact subset of $Y$. In equivariant coordinates with common background metric $g_{\mathrm{Euc}}$, we write $r$ and $R$ for the radii at the two scales, related by $r=\sqrt{t_i}R$.

We denote the finite rod spheres by $E_1,\ldots,E_q$ and choose the comparison metrics from Lemma~\ref{lem:fixed-rod-angle-family} by setting
$$
    \ell_{i,j}=\frac{\Area_{t_i^{-1}g_i}(E_j)}{2\pi},\qquad h_i=h(\bl_i),\qquad \boldsymbol{\beta}_i=\boldsymbol{\beta}(\bl_i),\qquad c_i=\mathbf1-\boldsymbol{\beta}_i.
$$
Thus $\bl_i\to\bl^*$, and the defining identities for the family imply that
\begin{equation}\label{eq:fixed-rod-local-area-matching}
    \beta_{i,j}\Area_{h_i}(E_j) =\Area_{t_i^{-1}g_i}(E_j), \qquad 1\leq j\leq q.
\end{equation}
To obtain a smooth comparison metric, we then change the cone angle along $E_j$ by $2\pi c_{i,j}=2\pi(1-\beta_{i,j})$.

We choose equivariant identifications $\Phi_i$ from expanding domains in $Y$ onto the corresponding regions in $(X_\mR, g_i)$, preserving each labelled rod sphere, and set
$$
    \gamma_i=t_i^{-1}\Phi_i^*g_i.
$$
We write $J_{i,c}$ for the angle identification from $\boldsymbol{\beta}_i$ to $\boldsymbol{\beta}_i+c$, with $J_{i,0}=\Id$, and choose these identifications to equal the identity outside a fixed compact set. Since $\boldsymbol{\beta}_i+c_i=\mathbf1$, the identification $J_i=J_{i,c_i}$ changes the finite cone angles of $h_i$ to $2\pi$. We therefore compare $\gamma_i$ with $b_i=J_i h_i$ and define $k_i$ by
$$
    \gamma_i=J_i(h_i+k_i)=b_i+J_i k_i.
$$
We measure $k_i$ in the fixed-angle developments of $h_i$ and denote by $|k_i|_{2,\alpha;1+R}$ its local $C^{2,\alpha}$ norm, with derivatives of order $m$ multiplied by $(1+R)^m$ and the H\"older seminorm rescaled accordingly. For $r_0>0$ small, we set
$$
    R_i=\frac{r_0}{\sqrt{t_i}},\qquad Y_{R_i}=\{R\leq R_i\}.
$$
We will choose the identifications $\Phi_i$ on one fixed neighborhood and restrict them as $r_0$ decreases. The estimates below hold with the same constants for all sufficiently small $r_0$, which will allow us to let $r_0\to0$ in the later argument.

We denote by $y$ the ALE end coordinates on the cover and define $x_i$ by $x_i\circ\Phi_i=\sqrt{t_i} y$, using lifts of $\Phi_i$ to the corresponding covers.

\begin{lemma}\label{lem:fixed-rod-local-estimate}
For every sufficiently small $\delta>0$, we can choose the identifications $\Phi_i$ and pass to a subsequence so that $x_i$ converge smoothly away from $p$ to coordinates $x$ extending smoothly across $p$ on the local cover, in which $g(p)=g_{\mathrm{Euc}}$ and $\partial g(p)=0$, and
\begin{align*}
    &\delta_{\gamma_i}b_i=0,\\
    &\limsup_{i\to\infty}\sup_{Y_{R_i}}|k_i|_{2,\alpha;1+R}\leq Cr_0^2,\\
    &\limsup_{i\to\infty}t_i^{-1+\delta/2} \left(\|k_i\|_{C^{2,\alpha}_{2-\delta}(Y_{R_i})}+|c_i|\right) \leq C_\delta r_0^\delta.
\end{align*}
\end{lemma}

\begin{proof}
We first establish an inverse estimate using \eqref{eq:fixed-rod-local-area-matching} and the two coordinate rescalings. We then apply this estimate twice: first to identify the coordinates near $p$, and then to obtain the required quatitative rate.

\medskip
\noindent\textbf{Step 1. Linearization of the equation and the area constraints.}

We define
$$
    B_{i,j} =\left.\frac{\partial}{\partial c_j}\right|_{c=0} \mathscr E_0(J_{i,c}h_i).
$$
After transporting the gauged equation to the fixed-angle developments of $h_i$, its value at $(k,c)=(0,0)$ is $\lambda_i t_i h_i$, and its linearization is
$$
    \mathscr P_{h_i}k+\sum_{j=1}^q c_jB_{i,j} +\lambda_i t_i k.
$$
We first establish an inverse estimate for the operator
$$
    (k,c)\longmapsto \mathscr P_{h_i}k+\sum_{j=1}^q c_jB_{i,j},
$$
We denote by $L_i(k,c)$ the variation of the weighted areas $\beta_j\Area(E_j)/(2\pi)$,  namely
$$
    (L_i(k,c))_j =\frac1{2\pi}\left.\frac{\dd}{\dd u}\right|_{u=0} \left[(\beta_{i,j}+uc_j) \Area_{J_{i,uc}(h_i+uk)}(E_j)\right].
$$

By the proof of Lemma~\ref{lem:fixed-rod-angle-family}, the map
$$
    (k,c)\longmapsto \left(\mathscr P_{h_i}k+\sum_jc_jB_{i,j},\,L_i(k,c)\right)
$$
is invertible at weight $-1$, uniform for large $i$.

\medskip
\noindent\textbf{Step 2. The inverse estimate at weight $2-\delta$.}

Increasing the weight from $-1$ to $2-\delta$ introduces two coordinate rescalings, generated at infinity by
$$
    X_1=x_1\partial_{x_1}+x_2\partial_{x_2},\qquad X_2=x_3\partial_{x_3}+x_4\partial_{x_4}.
$$
We extend $X_1$ and $X_2$ to invariant vector fields on $Y$ using a cutoff on the ALE end. By the vector-field inverse established in the proof of Proposition~\ref{prop:ae-augmented}, we can add decaying corrections to obtain invariant fields $W_{i,1}$ and $W_{i,2}$ satisfying
$$
    \delta_{h_i}\delta_{h_i}^*W_{i,\mu}=0, \qquad W_{i,\mu}=X_\mu+O_\infty(R^{1-\nu}), \qquad \mu=1,2,
$$
where $0<\nu<1$ is smaller than the decay exponent in the annular coordinate estimate. Here we use $h_i$ to identify vector fields with one-forms. Since the flows of $W_{i,\mu}$ preserve the rod spheres and cone angles, we have
$$
    \mathscr P_{h_i}(\mathcal L_{W_{i,\mu}}h_i)=0, \qquad L_i(\mathcal L_{W_{i,\mu}}h_i,0)=0.
$$
We denote the corresponding fields for $h$ by $W_1,W_2$ and choose compactly supported invariant test tensors $\psi_1,\psi_2$ such that
$$
    s(k)=\left( \int_Y\langle k,\psi_1\rangle_h\dvol_h,\, \int_Y\langle k,\psi_2\rangle_h\dvol_h \right)
$$
defines an isomorphism on the span of $\mathcal L_{W_1}h$ and $\mathcal L_{W_2}h$. We use the angle identifications to define the same pairings for tensors with the fixed angles of $h_i$. Thus the conditions $s(k)=0$ fix the two coordinate rescalings.

We now classify the pairs $(k,c)$ at weight $2-\delta$ satisfying the equation
$$
    \mathscr P_{h_i}k+\sum_jc_jB_{i,j}=0, \qquad L_i(k,c)=0.
$$
Since the tensors $B_{i,j}$ have compact support, we have $\mathscr P_{h_i}k=0$ on the ALE end. The limiting operator on the Euclidean orbifold cover is
$$
    \mathscr P_0 =(D\mathscr E_0)_{g_{\mathrm{Euc}}} +\delta_{g_{\mathrm{Euc}}}^*\delta_{g_{\mathrm{Euc}}},
$$
so the possible additional leading terms are determined by homogeneous solutions of $\mathscr P_0k=0$ on $\mathbb R^4\setminus\{0\}$.

We consider a torus-invariant solution whose  coefficients are homogeneous of degree $\sigma\in(-1,2)$. The singularity at the origin is removable, and elliptic estimates together with the growth bound imply that the coefficients are polynomials of degree at most one. Invariance under the simultaneous half-turn excludes the linear terms. Thus the only additional leading terms are linear combinations of
$$
    \dd x_1^2+\dd x_2^2, \qquad \dd x_3^2+\dd x_4^2,
$$
which correspond to rescaling the two rotation planes.  Consequently, the expansion for the linearized equation on the ALE end implies that every solution of
$$
    \mathscr P_{h_i}k+\sum_jc_jB_{i,j}=0,\qquad L_i(k,c)=0
$$
at weight $2-\delta$ differs from a linear combination of $\mathcal L_{W_{i,1}}h_i$ and $\mathcal L_{W_{i,2}}h_i$ by a solution at negative weight. By Step 1, the latter solution vanishes.

The weighted Fredholm argument in the proof of Proposition~\ref{prop:ae-augmented} shows that the index increases by two when the weight passes from $-1$ to $2-\delta$. Thus adding the two conditions $s(k)=0$ gives us back the invertibility, and we obtain
\begin{equation}\label{eq:fixed-rod-local-inverse}
    \|k\|_{C^{2,\alpha}_{2-\delta}}+|c| \leq C_\delta\left( \left\|\mathscr P_{h_i}k+\sum_jc_jB_{i,j}\right\|_{C^{0,\alpha}_{-\delta}} +|L_i(k,c)|+|s(k)|\right).
\end{equation}
Since $h_i\to h$ in the fixed angle identifications, the constant $C_\delta$ can be chosen uniformly for large $i$.

\medskip
\noindent\textbf{Step 3. Choice of coordinates on a fixed neighborhood.}

We choose $r_*>0$ sufficiently small and construct the identifications $\Phi_i$ on
$$
    Y_{4r_*/\sqrt{t_i}}.
$$
We will keep these identifications fixed when we later decrease $r_0$. We use the annular coordinates and weighted norms from the proof of Lemma~\ref{scd:gauge}, with neck weight
$$
    \omega_i(r)=\left(\frac r{r_*}\right)^\nu +\left(\frac{\sqrt{t_i}}r\right)^\nu.
$$
As in the proof of Lemma~\ref{scd:gauge}, we subtract cutoffs of invariant constant tensors and linear vector fields before applying the weighted norms, and measure their coefficients separately. We start with equivariant normal coordinates for $g$ on the local orbifold cover at $p$, so that
$$
    g-g_{\mathrm{Euc}}=O_\infty(r^2).
$$
Since $\omega_i(r)\geq(r/r_*)^\nu$, the contribution of the error term to the weighted norm is bounded by
$$
    \frac{r^2}{\omega_i(r)} \leq r_*^\nu r^{2-\nu} \leq Cr_*^2.
$$
Combining this bound with convergence to $h$ and $g$ at the two scales and the annular coordinate estimate, we obtain a bound of $Cr_*^2+o_i(1)$ for the difference from the reference metric $t_i b_i$.

We now use this metric comparison to choose the identifications $\Phi_i$ so that the condition $\delta_{\gamma_i}b_i=0$ is satisfies. Pulling back the reference metric $b_i$ by the flow of an invariant vector field $W$, we obtain the linearized operator $2\delta_{b_i}\delta_{b_i}^*W$ at $\gamma_i=b_i$. Thus, to construct the coordinate correction, we first establish a uniform inverse for $\delta_{b_i}\delta_{b_i}^*$.

We first require the coordinate correction to fix the outer coordinate sphere, so we impose zero Dirichlet boundary values on $W$. On the flat orbifold ball $B_{4r_*}/\Gamma$, the operator $\delta_{g_{\mathrm{Euc}}}\delta_{g_{\mathrm{Euc}}}^*$ with zero Dirichlet boundary values has trivial kernel by integration by parts. On $(Y,h)$, the equation $\delta_h\delta_h^*W=0$ has no nonzero invariant solution with $W=O(R^{1-\nu})$, by the argument in Lemma~\ref{scd:gauge}. Together with the annular estimate used there, these facts imply a uniform Dirichlet inverse for $\delta_{b_i}\delta_{b_i}^*$ on $Y_{4r_*/\sqrt{t_i}}$.

We now vary the boundary values of $\Phi_i$ by the two linear rescalings generated by $X_1$ and $X_2$. For each sufficiently small choice of the two parameters, we extend the boundary map equivariantly to the interior and use the Dirichlet inverse and the nonlinear estimates from Lemma~\ref{scd:gauge} to construct a further coordinate correction fixing the outer boundary and satisfying
$$
    \delta_{\gamma_i}b_i=0.
$$
Thus the divergence gauge condition holds for every choice of the boundary parameters, which remain available to impose $s(k_i)=0$.

At zero boundary parameters, the metric comparison and the bounds for the coordinate correction imply that $s(k_i)=O(r_*^2)+o_i(1)$. The derivative of $s(k_i)$ with respect to the two parameters differs by $O(r_*^2)+o_i(1)$ from the matrix with columns
$$
    s(\mathcal L_{W_1}h),\qquad s(\mathcal L_{W_2}h).
$$
Since this matrix is invertible, for $r_*$ sufficiently small and $i$ sufficiently large, the implicit function theorem allows us to choose boundary parameters of size $O(r_*^2)+o_i(1)$ such that $s(k_i)=0$. The bounds for these rescalings and the Dirichlet coordinate correction, together with the initial metric comparison, imply
\begin{equation}\label{eq:fixed-rod-local-rough}
    \sup_{R\leq2r_*/\sqrt{t_i}}|k_i|_{2,\alpha;1+R} \leq Cr_*^2+o_i(1).
\end{equation}
All coordinate changes preserve each labelled rod sphere, so the identities \eqref{eq:fixed-rod-local-area-matching} continue to hold.

\medskip
\noindent\textbf{Step 4. An estimate on the truncated ALE region.}

We now combine the Einstein equation with the initial smallness from \eqref{eq:fixed-rod-local-rough} to obtain a quantitative estimate. For $r_0\leq r_*$, we set $R_i=r_0/\sqrt{t_i}$ and choose an invariant cutoff $\chi_i$ equal to one on $Y_{R_i}$ and zero for $R\geq2R_i$, satisfying
$$
    |\nabla\chi_i|\leq CR_i^{-1}, \qquad |\nabla^2\chi_i|\leq CR_i^{-2},
$$
with the corresponding scaled H\"older bounds. We set
$$
    \begin{aligned} &\widetilde k_i=\chi_i k_i,\\
                    &E_i(r_0)=\|\widetilde k_i\|_{C^{2,\alpha}_{2-\delta}}+|c_i|,\\
                    &M_i(r_0)=\sup_{R_i\leq R\leq2R_i}|k_i|_{2,\alpha;1+R}. \end{aligned}
$$
For each fixed $r_0>0$ and sufficiently large $i$, the cutoff equals one on the fixed compact set containing the finite rods, the angle modifications, and the supports of $B_{i,j}$ and $\psi_1,\psi_2$. Thus
$$
    L_i(\widetilde k_i,c_i)=L_i(k_i,c_i), \qquad s(\widetilde k_i)=0.
$$

On $Y_{2R_i}$, the Einstein equation and the gauge condition imply that
$$
    \mathscr P_{h_i}k_i+\sum_jc_{i,j}B_{i,j} =-\lambda_i t_i(h_i+k_i)+\mathcal Q_i(k_i,c_i),
$$
where $\mathcal Q_i$ consists of terms at least quadratic in $(k_i,c_i)$. Moreover, the area identities \eqref{eq:fixed-rod-local-area-matching} read
$$
    (\beta_{i,j}+c_{i,j}) \Area_{J_{i,c_i}(h_i+k_i)}(E_j) =\beta_{i,j}\Area_{h_i}(E_j).
$$
Hence $L_i(k_i,c_i)$ consists of quadratic error terms.

We set
$$
    \eta_i=C\left( r_*^2+ \sup_{R\leq2r_*/\sqrt{t_i}}|k_i|_{2,\alpha;1+R} +|c_i| \right).
$$
By the product estimate used in \eqref{scd:tame}, we have
$$
    \|\mathcal Q_i(\widetilde k_i,c_i)\|_{C^{0,\alpha}_{-\delta}} +|L_i(\widetilde k_i,c_i)| \leq C\eta_iE_i(r_0).
$$
Here we estimate one factor in the scaled unweighted norm, bounded by $\eta_i$, and the other in the weighted norm defining $E_i(r_0)$.

Since $\lambda_i$ is bounded and $1+R\leq CR_i$ on the support of $\chi_i$, the Einstein constant terms satisfy
$$
    \begin{aligned} &\|\chi_i\lambda_i t_i h_i\|_{C^{0,\alpha}_{-\delta}} \leq Ct_iR_i^\delta,\\
                    &\|\lambda_i t_i\widetilde k_i\|_{C^{0,\alpha}_{-\delta}} \leq Ct_iR_i^2E_i(r_0) \leq Cr_*^2E_i(r_0). \end{aligned}
$$
The error from differentiating the cutoff is supported on $R_i\leq R\leq2R_i$. By the product rule and the definition of $M_i(r_0)$, we obtain
$$
    \|\mathscr P_{h_i}(\chi_i k_i) -\chi_i\mathscr P_{h_i}k_i\|_{C^{0,\alpha}_{-\delta}} \leq CR_i^{-2+\delta}M_i(r_0).
$$
Since the angle modifications are supported where $\chi_i=1$, the difference between $\chi_i\mathcal Q_i(k_i,c_i)$ and $\mathcal Q_i(\widetilde k_i,c_i)$ is supported on the same annulus. Estimating one factor in each quadratic term by $\eta_i$, we similarly obtain
$$
    \|\chi_i\mathcal Q_i(k_i,c_i) -\mathcal Q_i(\widetilde k_i,c_i)\|_{C^{0,\alpha}_{-\delta}} \leq C\eta_iR_i^{-2+\delta}M_i(r_0).
$$

We now apply \eqref{eq:fixed-rod-local-inverse} to $(\widetilde k_i,c_i)$. Since $\chi_iB_{i,j}=B_{i,j}$ and $s(\widetilde k_i)=0$, the preceding estimates imply that
$$
    E_i(r_0) \leq C_\delta\left( t_iR_i^\delta +R_i^{-2+\delta}M_i(r_0) +\eta_iE_i(r_0) \right).
$$
By \eqref{eq:fixed-rod-local-rough} and $c_i\to0$, we first choose $r_*$ sufficiently small and then $i$ sufficiently large so that $C_\delta\eta_i\leq\frac12$. Absorbing the last term into the left-hand side, we obtain
\begin{equation}\label{eq:fixed-rod-local-cutoff}
    E_i(r_0) \leq C_\delta\left( t_iR_i^\delta +R_i^{-2+\delta}M_i(r_0) \right).
\end{equation}
\medskip
\noindent\textbf{Step 5. Extension of the coordinates across $p$.}

We first apply \eqref{eq:fixed-rod-local-cutoff} with $r_0=r_*$. By \eqref{eq:fixed-rod-local-rough}, the quantities $M_i(r_*)$ are bounded. Since $R_i=r_*/\sqrt{t_i}\to\infty$, both terms on the right tend to zero, so $E_i(r_*)\to0$ and $k_i\to0$ on every fixed compact subset of $Y$.

We now consider the coordinates induced on the punctured neighborhood of $p$ in $M$. In the annular decomposition used in {\bf Step 3}, we denote the common invariant constant tensor by $T_i$ and pass to a subsequence with $T_i\to T$. On the ALE end, the limiting difference from $h$ is $T+O(R^{-\nu})$. Since this difference vanishes by the above convergence, we obtain $T=0$. At the scale of $M$, the same decomposition implies that
$$
    g-g_{\mathrm{Euc}}=O((r/r_*)^\nu) \qquad\text{as }r\to0.
$$

On the local orbifold cover near $p$, the vector-field estimates express the limiting coordinate change as a linear map plus $O_\infty(r^{1+\nu})$. Thus the coordinate change extends in $C^{1,\nu}$ across the point above $p$, and elliptic regularity for the limiting coordinate equation implies smooth extension. In these coordinates, we have $g(p)=g_{\mathrm{Euc}}$. Equivariance under the simultaneous half-turn implies that $\partial g(p)=0$, so
$$
    g=g_{\mathrm{Euc}}+H_2+O_\infty(r^3).
$$
These coordinates are defined on a fixed neighborhood of $p$, so the bounds on the remainder remain fixed as $r_0$ decreases.

\medskip
\noindent\textbf{Step 6. The quantitative estimate.}

For each fixed $r_0$, smooth convergence on the annulus $r_0\leq r\leq2r_0$ and the expansion in {\bf Step 5} imply that
$$
    \limsup_{i\to\infty}M_i(r_0)\leq Cr_0^2.
$$
We substitute this bound into \eqref{eq:fixed-rod-local-cutoff} and use
$$
    t_iR_i^\delta =R_i^{-2+\delta}r_0^2 =t_i^{1-\delta/2}r_0^\delta
$$
to obtain
$$
    \limsup_{i\to\infty}t_i^{-1+\delta/2} \left( \|k_i\|_{C^{2,\alpha}_{2-\delta}(Y_{R_i})} +|c_i| \right) \leq C_\delta r_0^\delta.
$$
Then the conclusion follows.
\end{proof}

\subsection{Vanishing of the boundary integral}

We now pair the Einstein equation with the dilation tensor $\zeta_{h_i}$. The angle contribution is controlled by Lemma~\ref{lem:fixed-rod-angle-cancellation}, while the local estimate bounds the remaining errors.

\begin{proposition}\label{prop:fixed-rod-local-obstruction}
We have
$$
    \mathcal B(H_2)=0.
$$
\end{proposition}

\begin{proof}
We use the metrics $h_i$, tensors $k_i$, angle changes $c_i$, and identifications constructed in the preceding subsection. Thus
$$
    \gamma_i=J_{i,c_i}(h_i+k_i),\qquad \delta_{\gamma_i}b_i=0,\qquad \Ric_{\gamma_i}=\lambda_i t_i\gamma_i.
$$
We fix $r_0>0$ sufficiently small and set $R_i=r_0/\sqrt{t_i}$. Throughout the proof, we first let $i\to\infty$ with $r_0$ fixed, and then let $r_0\to0$.

\medskip
\noindent\textbf{Step 1. The angle contribution.}

We denote the exponent in Lemma~\ref{lem:fixed-rod-angle-cancellation} by $\vartheta$ and choose
$$
    0<\delta<\frac{2\vartheta}{1+\vartheta}.
$$
In particular, $\delta<1$ and $(1-\delta/2)(1+\vartheta)>1$. By Lemma~\ref{lem:fixed-rod-local-estimate}, for each fixed $r_0$ and all sufficiently large $i$, we have
\begin{equation}\label{eq:fixed-rod-weak-rate}
    |c_i|+\|k_i\|_{C^{2,\alpha}_{2-\delta}(Y_{R_i})} \leq Ct_i^{1-\delta/2}r_0^\delta.
\end{equation}
We set $\zeta_i=\zeta_{h_i}$ and $A_{i,j}=\Area_{h_i}(E_j)$. The identities in \eqref{eq:fixed-rod-exact-dilation} hold uniformly for the conical family, so
$$
    \mathscr P_{h_i}\zeta_i=0,\qquad \delta_{h_i}\zeta_i=0,\qquad \tr_{h_i}\zeta_i=0.
$$
We apply Lemma~\ref{lem:pair} at $h_i$, varying one angle at a time. Since $\zeta_i$ represents the negative of the scaling variation modulo a Lie derivative, and the Lie derivative has zero pairing with the compactly supported divergence-free tensor $B_{i,j}$, we obtain
\begin{equation}\label{eq:fixed-rod-angle-dilation}
    \int_Y\langle B_{i,j},\zeta_i\rangle_{h_i} \dvol_{h_i} =8\pi A_{i,j}.
\end{equation}
Thus the pairing with $\zeta_i$ depends on the angle changes through the sum $\sum_jA_{i,j}c_{i,j}$. By Lemma~\ref{lem:fixed-rod-angle-cancellation} and \eqref{eq:fixed-rod-weak-rate}, we have
\begin{equation}\label{eq:fixed-rod-small-angle-pairing}
    \left|\sum_jA_{i,j}c_{i,j}\right| \leq C|c_i|^{1+\vartheta} \leq C\left(t_i^{1-\delta/2}r_0^\delta\right)^{1+\vartheta} =o(t_i).
\end{equation}

\medskip
\noindent\textbf{Step 2. The remaining error terms.}

We expand the Einstein equation in the chosen gauge as
$$
    \mathscr P_{h_i}k_i+\sum_jc_{i,j}B_{i,j} =-\lambda_i t_i(h_i+k_i)+\mathcal Q_i(k_i,c_i),
$$
where $\mathcal Q_i$ denotes the quadratic remainder. On the ALE end, where the angle identifications are the identity, we have
$$
    |\mathcal Q_i(k_i,c_i)|_{h_i} \leq C\left( |k_i|_{h_i}|\nabla^2k_i|_{h_i} +|\nabla k_i|_{h_i}^2 +R^{-2}|k_i|_{h_i}^2 \right).
$$
Here the derivatives are taken with respect to $h_i$, and the terms involving $c_i$ are supported in a fixed compact subset.

We set
$$
    \varepsilon_i=t_i^{1-\delta/2}r_0^\delta.
$$
By \eqref{eq:fixed-rod-weak-rate}, we have
$$
    |\nabla^m k_i|_{h_i} \leq C\varepsilon_i(1+R)^{2-\delta-m}, \qquad 0\leq m\leq2.
$$
Using $|\zeta_i|_{h_i}\leq C(1+R)^{-4}$ and the ALE volume growth, we therefore obtain
$$
    \begin{aligned} &\int_{Y_{R_i}} |\langle \mathcal Q_i(k_i,c_i),\zeta_i\rangle_{h_i}| \dvol_{h_i} \leq C\varepsilon_i^2 \left(1+\int_1^{R_i}R^{1-2\delta} \dd R\right) \leq Ct_ir_0^2+o(t_i),\\
                    &|\lambda_i|t_i\int_{Y_{R_i}} |\langle k_i,\zeta_i\rangle_{h_i}| \dvol_{h_i} \leq Ct_i\varepsilon_i \left(1+\int_1^{R_i}R^{1-\delta} \dd R\right) \leq Ct_ir_0^2+o(t_i). \end{aligned}
$$
Indeed, the contributions from the end satisfy
$$
    \varepsilon_i^2R_i^{2-2\delta}=t_ir_0^2, \qquad t_i\varepsilon_iR_i^{2-\delta}=t_ir_0^2,
$$
while the compact contributions $\varepsilon_i^2$ and $t_i\varepsilon_i$ are $o(t_i)$.

\medskip
\noindent\textbf{Step 3. The boundary integral.}

We set
$$
    \mathcal B_i(R_i) =\int_{Y_{R_i}} \langle\mathscr P_{h_i}k_i,\zeta_i\rangle_{h_i} \dvol_{h_i}.
$$
Since $\mathscr P_{h_i}\zeta_i=0$, integration by parts expresses $\mathcal B_i(R_i)$ as an integral over $\partial Y_{R_i}$. The contributions from small tubes around the rod spheres vanish by regularity in local developments.

We pair the equation in {\bf Step 2} with $\zeta_i$ and use \eqref{eq:fixed-rod-angle-dilation}. Since $\tr_{h_i}\zeta_i=0$, the term $-\lambda_i t_i h_i$ has zero pairing, and we obtain
$$
    \begin{aligned} \mathcal B_i(R_i)+8\pi\sum_jA_{i,j}c_{i,j} =-\lambda_i t_i\int_{Y_{R_i}} \langle k_i,\zeta_i\rangle_{h_i} \dvol_{h_i} +\int_{Y_{R_i}} \langle \mathcal Q_i(k_i,c_i),\zeta_i\rangle_{h_i} \dvol_{h_i}. \end{aligned}
$$
Thus \eqref{eq:fixed-rod-small-angle-pairing} and {\bf Step 2} imply that
\begin{equation}\label{eq:fixed-rod-boundary-small}
    |\mathcal B_i(R_i)| \leq Ct_ir_0^2+o(t_i).
\end{equation}

\medskip
\noindent\textbf{Step 4. Identification of the boundary limit.}

We use the coordinates on the local orbifold cover near $p$ obtained in Lemma~\ref{lem:fixed-rod-local-estimate}, in which
$$
    g=g_{\mathrm{Euc}}+H_2+O_\infty(r^3).
$$
At the scale of $M$, the reference metrics converge to $g_{\mathrm{Euc}}$ on every fixed punctured annulus. The gauge condition therefore passes to the limit as
$$
    \delta_g g_{\mathrm{Euc}}=0.
$$
Taking the term of degree one in this identity, we obtain
$$
    \delta_{g_{\mathrm{Euc}}}H_2=0.
$$
Thus, when we pair $\mathscr P_{g_{\mathrm{Euc}}}H_2$ with $-4H^4$ and integrate by parts, the outer boundary integral over $\mathbb S^3/\Gamma$ equals $\mathcal B(H_2)$.

We now rescale the outer annulus to unit size. In Cartesian components, convergence to $g$ and the uniform expansion of $\zeta_i$ imply that
$$
    \begin{aligned} &k_i(R_i x) =r_0^2H_2(x)+O(r_0^3)+o_i(1),\\
                    &R_i^4\zeta_i(R_i x) =-4H^4(x)+o_i(1), \end{aligned}
$$
with the corresponding differentiated estimates on a fixed annulus. The boundary expression contains one derivative, while the boundary measure has order $R_i^3$. Together with the factor $R_i^{-4}$ in $\zeta_i$, this contributes the factor
$$
    R_i^{-2}=\frac{t_i}{r_0^2}.
$$
Thus
$$
    t_i^{-1}\mathcal B_i(R_i) =\mathcal B(H_2)+O(r_0)+o_i(1).
$$
For each fixed $r_0$, we let $i\to\infty$ and combine the boundary expansion with \eqref{eq:fixed-rod-boundary-small} to obtain
$$
    |\mathcal B(H_2)|\leq Cr_0+Cr_0^2.
$$
Since the coordinates and the constants are independent of the smaller radius $r_0$, we can then let $r_0\to0$ and conclude that $\mathcal B(H_2)=0$.
\end{proof}
\subsection{Completion of the compactness argument}

\begin{proof}[Proof of Theorem~\ref{thm:fixed-rod-smooth-compactness}]
We consider a sequence $[g_i]\in\mathcal M_{\mR}$ and suppose that its curvature is unbounded. Then we obtain a contradiction by Proposition \ref{prop:positivity of pairing} and  Proposition~\ref{prop:fixed-rod-local-obstruction}.
\end{proof}

\bibliographystyle{alpha}
\bibliography{reference}

@article{AndersonSymmetries,
  author       = {M. T. Anderson},
  title        = {{Extension of symmetries on Einstein manifolds with boundary}},
  journal      = {Selecta Mathematica (New Series)},
  volume       = {16},
  year         = {2010},
  pages        = {343--375},
  note         = {Corrected version: arXiv:0704.3373v6}
}

@article{Lafuente-Broder,
  author = {Broder, Kyle and Lafuente, Ramiro A.},
  title = {{$T^2$}-invariant {Einstein} metrics on closed {$4$}-manifolds},
  journal = {arXiv preprint arXiv:2610.18219},
  year = {2026}
}

@article{OzuchIntegrability,
  author  = {Ozuch, Tristan},
  title   = {Integrability of {Einstein} deformations and desingularizations},
  journal = {Communications on Pure and Applied Mathematics},
  volume  = {77},
  number  = {1},
  year    = {2024},
  pages   = {177--220},
  doi     = {10.1002/cpa.22129}
}

@article{ChuStern,
  author = {Chu, Adrian Chun-Pong and Stern, Daniel},
  title = {A variational construction for {Einstein} metrics on connected sums of {$S^2\times S^2$}},
  journal = {arXiv preprint arXiv:2610.29882},
  year = {2026}
}

@article{AA,
  author       = {L. Andersson and B. Araneda},
  title        = {{Infinitesimal rigidity of Hermitian gravitational instantons}},
  journal      = {arXiv preprint},
  year         = {2026},
  eprint       = {2601.17206},
  note         = {arXiv:2601.17206}
}

@article{AL,
  author       = {B. Araneda and J. Lucietti},
  title        = {{All toric Hermitian ALE gravitational instantons}},
  journal      = {arXiv preprint},
  year         = {2025},
  eprint       = {2510.09291},
  note         = {arXiv:2510.09291}
}

@article{AtiyahLeBrun,
  author       = {M. Atiyah and C. LeBrun},
  title        = {{Curvature, cones, and characteristic numbers}},
  journal      = {Mathematical Proceedings of the Cambridge Philosophical Society},
  volume       = {155},
  year         = {2013},
  pages        = {13--37}
}

@article{AtiyahSingerIII,
  author       = {M. F. Atiyah and I. M. Singer},
  title        = {{The index of elliptic operators. III}},
  journal      = {Annals of Mathematics (2)},
  volume       = {87},
  number       = {3},
  year         = {1968},
  pages        = {546--604}
}

@article{Bando,
  author       = {S. Bando},
  title        = {{Bubbling out of Einstein manifolds}},
  journal      = {Tohoku Mathematical Journal},
  volume       = {42},
  year         = {1990},
  pages        = {205--216},
  url          = {https://www.jstage.jst.go.jp/article/tmj1949/42/2/42_2_205/_pdf/-char/en}
}

@article{BKN,
  author       = {S. Bando and A. Kasue and H. Nakajima},
  title        = {{On a construction of coordinates at infinity on manifolds with fast curvature decay and maximal volume growth}},
  journal      = {Inventiones Mathematicae},
  volume       = {97},
  year         = {1989},
  pages        = {313--349},
  doi          = {10.1007/BF01389045}
}

@misc{BiquardDGGA,
  author       = {O. Biquard},
  title        = {{Polycopi\'e on Differential Geometry and Global Analysis}},
  howpublished = {Lecture notes},
  year         = {2008},
  note         = {February 6, 2008},
  url          = {https://webusers.imj-prg.fr/~olivier.biquard/dgga2007.pdf}
}

@article{BiquardI,
  author       = {O. Biquard},
  title        = {{D\'esingularisation de m\'etriques d'Einstein. I}},
  journal      = {Inventiones Mathematicae},
  volume       = {192},
  year         = {2013},
  pages        = {197--252}
}

@article{BG,
  author       = {O. Biquard and P. Gauduchon},
  title        = {{On toric Hermitian ALF gravitational instantons}},
  journal      = {Communications in Mathematical Physics},
  volume       = {399},
  year         = {2023},
  pages        = {389--422}
}

@book{BH,
  author       = {M. R. Bridson and A. Haefliger},
  title        = {{Metric Spaces of Non-Positive Curvature}},
  series       = {Grundlehren der mathematischen Wissenschaften},
  volume       = {319},
  publisher    = {Springer-Verlag},
  address      = {Berlin},
  year         = {1999},
  doi          = {10.1007/978-3-662-12494-9}
}

@incollection{BunchDonaldson,
  author       = {R. S. Bunch and S. K. Donaldson},
  title        = {{Numerical approximations to extremal metrics on toric surfaces}},
  booktitle    = {{Handbook of Geometric Analysis, No. 1}},
  series       = {Advanced Lectures in Mathematics},
  volume       = {7},
  publisher    = {International Press},
  year         = {2008},
  pages        = {1--28}
}

@article{BHHPW,
  author       = {T. Buttsworth and W. Hadden and E. Heyes and D. Platt and T. Wiseman},
  title        = {{The K\"ahler--Einstein metric on the third Del Pezzo surface}},
  journal      = {arXiv preprint},
  year         = {2026},
  eprint       = {2609.32655},
  note         = {arXiv:2609.32655}
}

@article{CFG,
  author       = {J. Cheeger and K. Fukaya and M. Gromov},
  title        = {{Nilpotent structures and invariant metrics on collapsed manifolds}},
  journal      = {Journal of the American Mathematical Society},
  volume       = {5},
  year         = {1992},
  pages        = {327--372}
}

@article{CDS1,
  author       = {X. Chen and S. Donaldson and S. Sun},
  title        = {{K\"ahler--Einstein metrics on Fano manifolds. I: Approximation of metrics with cone singularities}},
  journal      = {Journal of the American Mathematical Society},
  volume       = {28},
  number       = {1},
  year         = {2015},
  pages        = {183--197},
  doi          = {10.1090/S0894-0347-2014-00799-2}
}

@article{CDS2,
  author       = {X. Chen and S. Donaldson and S. Sun},
  title        = {{K\"ahler--Einstein metrics on Fano manifolds. II: Limits with cone angle less than $2\pi$}},
  journal      = {Journal of the American Mathematical Society},
  volume       = {28},
  number       = {1},
  year         = {2015},
  pages        = {199--234},
  doi          = {10.1090/S0894-0347-2014-00800-6}
}

@article{CDSIII,
  author       = {X. Chen and S. Donaldson and S. Sun},
  title        = {{K\"ahler--Einstein metrics on Fano manifolds. III: Limits as cone angle approaches $2\pi$ and completion of the main proof}},
  journal      = {Journal of the American Mathematical Society},
  volume       = {28},
  number       = {1},
  year         = {2015},
  pages        = {235--278},
  doi          = {10.1090/S0894-0347-2014-00801-8}
}

@article{CLW,
  author       = {X. Chen and C. LeBrun and B. Weber},
  title        = {{On conformally K\"ahler, Einstein manifolds}},
  journal      = {Journal of the American Mathematical Society},
  volume       = {21},
  year         = {2008},
  pages        = {1137--1168}
}

@article{ChenLi,
  author       = {X. Chen and Y. Li},
  title        = {{On the geometry of asymptotically flat manifolds}},
  journal      = {Geometry \& Topology},
  volume       = {25},
  number       = {5},
  year         = {2021},
  pages        = {2469--2572},
  doi          = {10.2140/gt.2021.25.2469}
}

@phdthesis{Dammerman,
  author       = {B. Dammerman},
  title        = {{Metrics of special curvature with symmetry}},
  school       = {University of Oxford},
  type         = {D.Phil. thesis},
  year         = {2004}
}

@article{DavisJanuszkiewicz,
  author       = {M. W. Davis and T. Januszkiewicz},
  title        = {{Convex polytopes, Coxeter orbifolds and torus actions}},
  journal      = {Duke Mathematical Journal},
  volume       = {62},
  year         = {1991},
  pages        = {417--451}
}

@article{Delay2012,
  author       = {E. Delay},
  title        = {{Smooth compactly supported solutions of some underdetermined elliptic PDE, with gluing applications}},
  journal      = {Communications in Partial Differential Equations},
  volume       = {37},
  number       = {10},
  year         = {2012},
  pages        = {1689--1716}
}

@article{Derdzinski,
  author       = {A. Derdzi\'nski},
  title        = {{Self-dual K\"ahler manifolds and Einstein manifolds of dimension four}},
  journal      = {Compositio Mathematica},
  volume       = {49},
  number       = {3},
  year         = {1983},
  pages        = {405--433}
}

@incollection{DonaldsonCone,
  author       = {S. K. Donaldson},
  title        = {{K\"ahler metrics with cone singularities along a divisor}},
  booktitle    = {Essays in Mathematics and its Applications},
  publisher    = {Springer},
  address      = {Heidelberg},
  year         = {2012},
  pages        = {49--79},
  doi          = {10.1007/978-3-642-28821-0_4}
}

@article{DonaldsonSun,
  author       = {S. Donaldson and S. Sun},
  title        = {{Gromov--Hausdorff limits of K\"ahler manifolds and algebraic geometry}},
  journal      = {Acta Mathematica},
  volume       = {213},
  number       = {1},
  year         = {2014},
  pages        = {63--106},
  doi          = {10.1007/s11511-014-0116-3}
}

@article{DHHKW,
  author       = {C. Doran and M. Headrick and C. P. Herzog and J. Kantor and T. Wiseman},
  title        = {{Numerical K\"ahler--Einstein metric on the third del Pezzo}},
  journal      = {Communications in Mathematical Physics},
  volume       = {282},
  year         = {2008},
  pages        = {357--393}
}

@article{HallMurphyNumerics,
  author       = {S. J. Hall and T. Murphy},
  title        = {{Numerical approximations to extremal toric K\"ahler metrics with arbitrary K\"ahler class}},
  journal      = {Proceedings of the Edinburgh Mathematical Society (2)},
  volume       = {60},
  year         = {2017},
  pages        = {893--910}
}

@article{Hitchin,
  author       = {N. J. Hitchin},
  title        = {{K\"ahlerian twistor spaces}},
  journal      = {Proceedings of the London Mathematical Society (3)},
  volume       = {43},
  number       = {1},
  year         = {1981},
  pages        = {133--150}
}

@article{KRWY24,
  author       = {M. Khuri and M. Reiris and G. Weinstein and S. Yamada},
  title        = {{Gravitational solitons and complete Ricci flat Riemannian manifolds of infinite topological type}},
  journal      = {Pure and Applied Mathematics Quarterly},
  volume       = {20},
  number       = {4},
  year         = {2024},
  pages        = {1895--1921}
}

@article{LeBrun,
  author       = {C. LeBrun},
  title        = {{On Einstein, Hermitian 4-manifolds}},
  journal      = {Journal of Differential Geometry},
  volume       = {90},
  year         = {2012},
  pages        = {277--302}
}

@article{LiHermitian,
  author       = {M. Li},
  title        = {{Classification results for conformally K\"ahler gravitational instantons}},
  journal      = {arXiv preprint},
  year         = {2023},
  eprint       = {2310.13197},
  note         = {arXiv:2310.13197}
}

@article{LiSun,
  author       = {M. Li and S. Sun},
  title        = {{Gravitational instantons and harmonic maps}},
  journal      = {arXiv preprint},
  year         = {2025},
  eprint       = {2507.15284},
  note         = {arXiv:2507.15284}
}

@article{Liu,
  author       = {T. Liu},
  title        = {{Toric Einstein 4-manifolds with non-negative sectional curvature}},
  journal      = {Annals of Global Analysis and Geometry},
  volume       = {67},
  year         = {2025},
  pages        = {Article 14},
  doi          = {10.1007/s10455-025-09990-3}
}

@article{Lu,
  author       = {Z. Lu},
  title        = {{A note on special K\"ahler manifolds}},
  journal      = {Mathematische Annalen},
  volume       = {313},
  year         = {1999},
  pages        = {711--713},
  url          = {https://lu.math.uci.edu/pdfs/publications/1997-2000/01.pdf}
}

@article{Nilsson,
  author       = {G. Nilsson},
  title        = {{Topology of toric gravitational instantons}},
  journal      = {Differential Geometry and its Applications},
  volume       = {96},
  year         = {2024},
  pages        = {102171},
  doi          = {10.1016/j.difgeo.2024.102171}
}

@article{OrlikRaymond,
  author       = {P. Orlik and F. Raymond},
  title        = {{Actions of the torus on $4$-manifolds. I}},
  journal      = {Transactions of the American Mathematical Society},
  volume       = {152},
  year         = {1970},
  pages        = {531--559},
  doi          = {10.1090/S0002-9947-1970-0268911-3}
}

@article{OzuchI,
  author       = {T. Ozuch},
  title        = {{Noncollapsed degeneration of Einstein 4-manifolds, I}},
  journal      = {Geometry \& Topology},
  volume       = {26},
  year         = {2022},
  pages        = {1483--1528}
}

@article{OzuchII,
  author       = {T. Ozuch},
  title        = {{Noncollapsed degeneration of Einstein 4-manifolds, II}},
  journal      = {Geometry \& Topology},
  volume       = {26},
  year         = {2022},
  pages        = {1529--1634}
}

@article{OzuchStability,
  author       = {T. Ozuch},
  title        = {{Stability of gravitational instantons with a bounded Killing vector field}},
  journal      = {arXiv preprint},
  year         = {2026},
  eprint       = {2608.31008},
  note         = {arXiv:2608.31008}
}

@article{OzuchPositive,
  author       = {T. Ozuch},
  title        = {{Positive Einstein metrics on infinitely many four-manifolds}},
  journal      = {arXiv preprint},
  year         = {2026},
  eprint       = {2610.08052},
  note         = {arXiv:2610.08052}
}

@article{Page,
  author       = {D. N. Page},
  title        = {{A compact rotating gravitational instanton}},
  journal      = {Physics Letters B},
  volume       = {79},
  year         = {1978},
  pages        = {235--238}
}

@article{Siu,
  author       = {Y.-T. Siu},
  title        = {{The existence of K\"ahler--Einstein metrics on manifolds with positive anticanonical line bundle and a suitable finite symmetry group}},
  journal      = {Annals of Mathematics (2)},
  volume       = {127},
  year         = {1988},
  pages        = {585--627}
}

@article{Tian-Yau,
  author       = {G. Tian and S.-T. Yau},
  title        = {{K\"ahler--Einstein metrics on complex surfaces with $c_1>0$}},
  journal      = {Communications in Mathematical Physics},
  volume       = {112},
  year         = {1987},
  pages        = {175--203}
}

@article{Zinger,
  author       = {A. Zinger},
  title        = {{The determinant line bundle for Fredholm operators: construction, properties, and classification}},
  journal      = {Mathematica Scandinavica},
  volume       = {118},
  year         = {2016},
  pages        = {203--268}
}

\end{document}